\documentclass[12pt]{article}
\usepackage{hyperref}
\usepackage{amssymb}
\usepackage{latexsym}
\usepackage{amsmath}
\usepackage{amscd}
\usepackage{theorem}
\usepackage{xy}
\xyoption{all}
\CompileMatrices

\usepackage{graphicx}

\makeatletter

\def\@@seccntfont{\bfseries\slshape}\def\@@secheaderfont{\bfseries\upshape}\def\@@precnt{{\upshape}}\def\@@postcnt{{\upshape}}\def\@@startsection#1#2#3#4#5#6{\if@noskipsec \leavevmode \fi\par
  \@tempskipa #4\relax\@afterindenttrue
  \ifdim \@tempskipa <\z@\@tempskipa -\@tempskipa
  \@afterindentfalse\fi
  \if@nobreak\everypar{}\else\addpenalty\@secpenalty\addvspace\@tempskipa\fi
  \@ifstar{\@dblarg{\@sect{#1}{#2}{#3}{#4}{#5}{#6}}}    {\@dblarg{\@@sect{#1}{#2}{#3}{#4}{#5}{#6}}}}
\def\@@sect#1#2#3#4#5#6[#7]#8{\ifnum
  #2>\c@secnumdepth\let\@svsec\@empty
  \else\refstepcounter{#1}\protected@edef\@svsec{\@seccntformat{#1}\relax}\fi\@tempskipa #5\relax
  \ifdim
  \@tempskipa>\z@\begingroup#6{\@hangfrom{\hskip#3\relax\@svsec}\interlinepenalty 
  \@M #8\@@par}\endgroup\csname #1mark\endcsname{#7}\else
  \def\@svsechd{#6{\hskip #3\relax\@svsec #8}\csname
    #1mark\endcsname{#7}}\fi\@xsect{#5}}
\def\setseccntfmt{\renewcommand{\@seccntformat}[1]{\S
  \csname the##1\endcsname\hspace{1ex}}}
\def\@@seccntfmt{\renewcommand{\@seccntformat}[1]{  {\@@seccntfont\@@precnt\csname
  the##1\endcsname}\@@postcnt\hspace{1ex}}}
\newcommand{\@@secnostar}[1][]{\def\tmpa{}\def\tmpb{#1}  \ifx\tmpa\tmpb\def\tmpa{-1ex}\else\def\tmpa{-\labelsep}\fi\@@seccntfmt\@@startsection
  {\@@name}{\@@level}{0mm}{-\baselineskip}{\tmpa}{\@@secheaderfont}{#1}}
\newcommand{\@@secstar}[2][]{\def\tmpa{}\def\tmpb{#1}  \ifx\tmpa\tmpb\def\tmpa{#2}\else\def\tmpa{#1}\fi\@@seccntfmt\@startsection
  {\@@name}{\@@level}{0mm}{-\baselineskip}{-\labelsep}{\@@secheaderfont}[\tmpa]{#2}}
\def\@@section{\def\@@name{section}\def\@@level{1}\@ifstar{\@@secstar}{\@@secnostar}}
\def\@@subsection{\def\@@name{subsection}\def\@@level{2}\@ifstar{\@@secstar}{\@@secnostar}}
\def\@@subsubsection{\def\@@name{subsubsection}\def\@@level{3}\@ifstar{\@@secstar}{\@@secnostar}}
\def\@@paragraph{\def\@@name{paragraph}\def\@@level{4}\@ifstar{\@@secstar}{\@@secnostar}}
\def\@@subparagraph{\def\@@name{subparagraph}\def\@@level{5}\@ifstar{\@@secstar}{\@@secnostar}}
\let\@@latexsection=\section\let\mysection=\@@section
\let\@@latexsubsection=\subsection\let\mysubsection=\@@subsection
\let\@@latexsubsubsection=\subsubsection\let\mysubsubsection=\@@subsubsection
\let\@@latexparagraph=\paragraph\let\myparagraph=\@@paragraph
\let\@@latexsubparagraph=\subparagraph\let\mysubparagraph=\@@subparagraph
\newif\if@@lastcharstar
\def\@@xxlastcharstar#1{\gdef\@@prevchar{}\@@lastcharstarfalse\@@yylastcharstar#1\end}
\def\@@yylastcharstar#1{\ifx#1\end\def\tmpa{*}\ifx\@@prevchar\tmpa\@@lastcharstartrue\fi
  \let\@@next=\relax\else\def\@@prevchar{#1}\let\@@next=\@@yylastcharstar\fi\@@next}
\gdef\@thm#1#2{\@@xxlastcharstar{#1}\if@@lastcharstar\else\refstepcounter{#1}\fi
  \trivlist\@topsep \theorempreskipamount\@topsepadd
  \theorempostskipamount
  \@ifnextchar [{\@ythm{#1}{#2}}{\@begintheorem{#2}{\csname
  the#1\endcsname}\ignorespaces}}
\gdef\th@nonumplain{\normalfont\itshape
  \def\@begintheorem##1##2{\item[\hskip\labelsep\theorem@headerfont
  ##1]}  \def\@opargbegintheorem##1##2##3{    \item[\hskip\labelsep \theorem@headerfont ##1\ (##3)]}}
\gdef\th@change{  \def\@begintheorem##1##2{\item[\hskip\labelsep{\@@seccntfont
  \@@precnt##2\@@postcnt\hskip 1ex}\theorem@headerfont ##1]}  \def\@opargbegintheorem##1##2##3{\item[\hskip\labelsep{\@@seccntfont
  \@@precnt##2\@@postcnt\hskip 1ex}\theorem@headerfont ##1\ (##3)]}}
\makeatother   
\def\newtheoremset[#1][#2]{{\theoremstyle{change}\newtheorem{-#1}[section]{#2}  \newtheorem{#1}[subsection]{#2}\newtheorem{*#1}[subsubsection]{#2}  \newtheorem{**#1}[paragraph]{#2}\newtheorem{***#1}[subparagraph]{#2}}
  {\theoremstyle{nonumplain}\newtheorem{#1*}{#2}}}
{\theoremheaderfont{\normalfont\bfseries}
{\theorembodyfont{\normalfont\itshape}
\newtheoremset[theorem][Theorem.]\newtheoremset[proposition][Proposition.]
\newtheoremset[lemma][Lemma.]\newtheoremset[corollary][Corollary.]\newtheoremset[properties][Properties.]\newtheoremset[property][Property.]
\newtheoremset[conjecture][Conjecture.]}\newtheoremset[sublemma][Sublemma.]\newtheoremset[assumption][Assumption.]
{\theorembodyfont{\normalfont\rmfamily}
\newtheoremset[definition][Definition.]
\newtheoremset[example][Example.]
\newtheoremset[facts][Facts.]\newtheoremset[fact][Fact.]
\newtheoremset[examples][Examples.]\newtheoremset[claim][Claim.]
\newtheoremset[notation][\itshape\mdseries Notation.]
\newtheoremset[remark][\itshape\mdseries Remark.]
\newtheoremset[remarks][\itshape\mdseries Remarks.]}}
\def\XYmatrix{\xymatrix@M=5pt} \def\ncmd{\newcommand}
\ncmd{\xysubset}[1][r]{\ar@<-2.5pt>@{^(-}[#1]\ar@<2.5pt>@{_(-}[#1]}
\ncmd{\XYmatrixc}[1]{\vcenter{\XYmatrix{#1}}}
\ncmd{\xyto}[1][r]{\ar@{->}[#1]}      \ncmd{\xyinj}[1][r]{\ar@{^(->}[#1]}
\ncmd{\xysurj}[1][r]{\ar@{->>}[#1]}   \ncmd{\xyline}[1][r]{\ar@{-}[#1]}
\ncmd{\xydotsto}[1][r]{\ar@{.>}[#1]}  \ncmd{\xydots}[1][r]{\ar@{.}[#1]}
\ncmd{\xyleadsto}[1][r]{\ar@{~>}[#1]} \ncmd{\xyeq}[1][r]{\ar@{=}[#1]}
\ncmd{\xyequal}[1][r]{\ar@{=}[#1]}    \ncmd{\xyequals}[1][r]{\ar@{=}[#1]}
\ncmd{\xymapsto}[1][r]{\ar@{|->}[#1]}\ncmd{\xyimplies}[1][r]{\ar@{=>}[#1]}
\ncmd{\xytofrom}[1][r]{\ar@{<->}[#1]} 
   
\def\XYTOTO[#1]^#2_#3{\xyto[#1]<0.5ex>^{#2}\xyto[#1]<-0.5ex>_{#3}}
\ncmd{\xytoto}[1][r]{\XYTOTO[#1]}

\newenvironment{proof}
{\noindent {\em Proof.}}
{\hfill $\Box$\medbreak}

\newcommand{\R}{{\mathbb R}}
\newcommand{\Z}{{\mathbb Z}}
\newcommand{\op}{\operatorname}

\newcommand{\sign}{\op{sign}}

\newcommand{\cok}{\op{coker}}

\newcommand{\ind}{\op{ind}}

\newcommand{\qed} {\hfill$\Box$}
\newcommand{\cm}{{\mathcal M}}

\newcommand{\id}{\op{id}}

\newcommand{\p}{\partial}

\newcommand{\lc}{{\mathcal C}}

\begin{document}
\title{Reidemeister Torsion in Floer-Novikov Theory\\ and \\Counting
  Pseudo-holomorphic Tori, II.}
\author{\scshape Yi-Jen Lee\medskip \thanks{supported in part by NSF
    grant DMS 0333163}\\  \small Department of Mathematics, Purdue
  University\\ \small West Lafayette, IN 47907 U.S.A. \\ \small {\tt yjlee@math.purdue.edu}}
\date{{\small This version: April 2005}}

\begin{titlepage}
\maketitle

\abstract{This is the second part of an article in two parts ,
  which builds the foundation of a Floer-theoretic invariant, \(I_F\).
(See \cite{part1} for part I).

Having constructed \(I_F\) and outlined a proof of its invariance
based on bifurcation analysis in part I, in this part we prove a
series of gluing theorems to confirm the bifurcation behavior
predicted in part I. These gluing theorems are different from (and much
harder than) the more conventional versions in that they deal with
broken trajectories or broken orbits connected at degenerate rest
points. The issues of orientation and signs are also settled in the
last section.

This part is strongly {\em dependent} on part I, and is meant 
only for readers familiar with the previous part of
this article.  
}
\newpage

\tableofcontents
\end{titlepage}

\section{Overview.}

This second part forms the main technical core of the present article.

We have not attempted to make this part independent of Part I, and 
shall frequently make use of the definitions, results,
notation, and convention from Part I without repetition. 
Thus, we urge the reader to familiarize him/herself with Part I before
attempting this one, paying particular attention to the convention in I.1.3. 

References in the form of I.* shall refer to section, theorem, or
equation numbers from Part I.

\subsection{A Brief Summary.}
The following summarizes the results contained in this part. Recall
the definitions of (RHFS*), (NEP), and admissible \((J,
X)\)-homotopies from sections 4.3, 4.4, 6.2 of \cite{part1} respectively. 
\begin{theorem*}
Let \(\Lambda=[1,2]\), and \((J^\Lambda, X^\Lambda)\) be an admissible \((J, X)\)-homotopy
connecting two regular pairs, \((J_1, X_1)\), \((J_2, X_2)\).
Then:
\begin{description}\itemsep -1pt
\item[(1) {\sc (Corner structures of parameterized moduli spaces)}]
The properties \(\quad\) (RHFS2c) and (RHFS3c) hold for the CHFS
generated by \((J^\Lambda, X^\Lambda)\);
\item[(2) {\sc (Orientation)}] The parameterized moduli spaces \(\cm_P^{\Lambda, +}\),
\(\bar{\cm}_O^{\Lambda, 1, +}\) may be respectively given coherent and grading
compatible orientations such that (RHFS4) holds;
\item[(3) {\sc (Existence of nonequivariant perturbations)}] (NEP)
  holds for all Type II handleslides in the CFHS
 generated by \((J^\Lambda, X^\Lambda)\).
\end{description}
\end{theorem*}

Combining with Propositions I.4.4.6 and I.6.2.2, this completes the
proof of
the general invariance theorem stated in part I, Theorem I.4.1.1.

Item (1) above follows from the gluing theorems proven in sections
2--6 below. 
Section 7 contains the discussion on orientability of the moduli
spaces, the definitions of coherent and grading-compatible
orientations, and as a
consequence, the proof of item (2) above.
Item (3) is established in sections 6.2--6.3.
There we introduce a class of (possibly nonlocal) perturbations to the
induced flow on the finite-cyclic covers in the statement of (NEP), 
establish the expected regularity and compactness properties of the
moduli spaces of such perturbed flows, and show how the arguments
in the proof of Theorem I.6.2.2 may be adapted to establish the
\(R\)-regularity of parametrized moduli spaces in this context, as
required by (NEP).

Gluing theory is the unifying theme of Part II. Not only is it 
used repeatedly to establish the bifurcation analysis, but it also
appears in the definition of 
coherent orientations in section 7. Linearized versions of the gluing
theorems in sections 2--6, which actually form part of the proofs of
these gluing theorems, play a major role in the verification of signs
for item 2 of Theorem 1.1 above. It is for this reason that we postpone
all discussion of orientations until the gluing results have been
fully treated.

Thus, we begin with 
a quick overview of the general features of gluing theory in next
subsection, then give a more specific outline of the variants
contained in this article in \S1.3.

\subsection{Basics of Gluing Theory.}
This subsection gives only a minimal outline of gluing
theory and its applications in Floer theory. Rather than a general
account, our aim is to set up the basic framework for the proofs of
the gluing theorems contained in this article, and to introduce some basic
notions and terminologies frequently used, some of
which are not conventional.  
The reader may find more details and better-balanced treatments 
in the vast literature on this subject, for example \cite{DK,
  D:floer}, and Floer's original papers. Also, precision will
sometimes be sacrificed here for the overall picture. We shall be
precise in later sections, when we return to the specific context of
this article.

\subsubsection{The four steps of gluing theory.}
Gluing is useful for studying the local structure of a stratified
moduli space, usually coming from compactification. 
Given a space of ``gluing parameters'' \(\Xi({\mathbb S})\) associated with a
codimension \(>0\) stratum \({\mathbb S}\), a typical gluing theorem
constructs a map from \(\Xi(\mathbb{S})\) to a neighborhood of \({\mathbb S}\) in the moduli space
of solutions to a PDE
\[
{\cal F}(w)=0,
\]
which is a local diffeomorphism.

The proof of a typical gluing theorem comprises of the following four
major steps:
\medbreak

\noindent{\sc Step 1.} {\it Constructing the pregluing map and error estimates. }
For each gluing parameter \(\chi\), one constructs an 
approximate solution \(w_\chi\) to the PDE considered, 
which varies smoothly with \(\chi\). 
The {\em pregluing map} \(\chi\mapsto w_\chi\) maps the space of
gluing parameters into a set in the ambient 
configuration space, that is close to the space of solutions.
An explict estimate, referred to as the ``error estimate'' 
is required to show that \({\cal F}(w_\chi)\) is
sufficiently small.
\medbreak

\noindent{\sc Step 2.} {\it Kuranishi structure.}
Let \(\mathfrak{D}_w: E\to F\) 
denote the linearization of \({\cal F}\) at \(w\) (i.e. the
deformation operator). This should be a Fredholm operator, and
ideally, one wants to show that \(\mathfrak{D}_{w_\chi}\) has a right inverse 
bounded uniformly in \(\chi\). 
Namely, there is a \(\chi\)-independent constant \(C_P>0\), and
operators \(P_\chi\) depending continuously on \(\chi\), such that
\[
\mathfrak{D}_{w_\chi}P_\chi=\op{id}, \quad \|P_\chi\|\leq C_P.
\]
For this to hold, judicious choices of normed spaces for \(E\), \(F\)
are often called for. 
\medbreak

\noindent{\sc Step 3.} {\it Obtaining a quadratic bound on the nonlinear part of
\({\cal F}\)}, namely, (\ref{Bq})
below. 
In local coordinates, one may write 
\begin{equation}\label{pde}
{\cal F}(w)={\cal F}(w_\chi)+\mathfrak{D}_{w_\chi}\xi+N_{w_\chi}(\xi)\quad
\text{for \(w=w_\chi+\xi \).}
\end{equation}
Setting \(\xi=P_\chi\eta_\chi\), 
a solution to \({\cal F}(w)=0\) is obtained by solving 
\begin{equation}\label{c-eq}
\eta_\chi=-N_{w_\chi}(P_\chi\eta)-{\cal F}(w_\chi).
\end{equation}
The contraction mapping theorem shows that
\begin{lemma*}
Let \(C_P\) be the upper bound on \(\|P_\chi\|\) as
above, and suppose that there is a \(\chi\)-independent constant \(k\)
such that 
\begin{gather}
\|{\cal F}(w_\chi)\| \leq \frac{1}{10kC_P^2},\nonumber\\
\|N_{w_\chi}(\xi_1)-N_{w_\chi}(\xi_2)\| \leq
k(\|\xi_1\|+\|\xi_2\|)\|\xi_1-\xi_2\|\quad \forall\xi_1, \xi_2.\label{Bq}
\end{gather}
Then there exists a unique $\eta_\chi$ with $\|\eta_\chi\|\leq 1/(5kC_P^2)$
solving (\ref{c-eq}).
Moreover, the solution $\eta_\chi$ varies smoothly with $\chi$, and
$
\|\eta_\chi\| \leq 2 \|{\cal F}(w_\chi)\|.
$
\end{lemma*}
Thus, by assigning to each gluing parameter \(\chi\) the corresponding
\(w_\chi+P_\chi\eta_\chi\), one obtains a smooth map from the space of
gluing parameters to the moduli space. This is the {\em gluing map}.
\medbreak

\noindent{\sc Step 4.} {\em Showing that the gluing map is a local diffeomorphism to
a neighborhood of \(S\)}. 

\subsubsection{Typical pregluing constructions in Floer theory.}
In Floer theory, \(\mathcal{F}=\partial_s+{\cal V}\), and 
\({\mathbb S}\) is a stratum in a moduli space of broken
trajectories or broken orbits. Thus, it is a product of reduced moduli
spaces 
\[\begin{split}
{\mathbb S}=& \hat{\cm}_0\times\hat{\cm}_1\times\cdots\times \hat{\cm}_k, \quad
\text{or}\\
&\hat{\cm}_1\times\hat{\cm}_2\times\cdots\times \hat{\cm}_k/\Z/k\Z.
\end{split}
\]
In the case of a family of Floer theories parameterized by \(\Lambda\),
\({\mathbb S}\) is a fiber product of reduced, parameterized moduli spaces over
\(\Lambda\)
\[
\begin{split}
{\mathbb S}=& \hat{\cm}_0^\Lambda\times_\Lambda\hat{\cm}_1^\Lambda\times_\Lambda\cdots\times_\Lambda \hat{\cm}_k^\Lambda, \quad
\text{or}\\
&\hat{\cm}_1^\Lambda\times_\Lambda\hat{\cm}_2^\Lambda\times_\Lambda\cdots\times_\Lambda \hat{\cm}_k^\Lambda/\Z/k\Z.
\end{split}
\] 
The space of gluing parameters in these cases is \(\Xi({\mathbb S})={\mathbb S}\times
(\Re, \infty)^k\) for certain large \(\Re\), and
the gluing maps map into a reduced moduli space or a reduced,
parameterized moduli space. 

We now describe the typical pregluing construction in these
situations.

Given a (unreduced) flow \(u(s)\) from the critical point \(x\) to
\(y\), we define its truncation
\[u_{[-R_-, R+]}(s): =\begin{cases}
u(s) &\text{when \(-R_-/2\leq s\leq R_+/2\)}\\
\exp(y, \beta(2-2s/R_+)\eta_y(s)) &\text{when \(s\geq R_+/2\)}\\
\exp(x, \beta(2s/R_-+2)\eta_x(s)) &\text{when \(s\leq-R_-/2\),}
\end{cases}\]
where \(\beta\) is a smooth cutoff function with \(\beta(s)=0\) for
\(s\leq0\), and \(\beta(s)=1\) for \(s\geq1\), and \(\eta_y\),
\(\eta_x\) are defined such that
\[u(s)=\begin{cases}\exp (x, \eta_x(s))&\text{for \(s\ll-1\),}\\
\exp (y, \eta_y(s))&\text{for \(s\gg1\).}
\end{cases}
\]
Let \(u_{(-\infty, R+]}\), \(u_{[-R_-,\infty)}\) be similarly
defined, truncated only at the positive/negative end respectively.

Let \(\{\hat{u}_0, \hat{u}_1, \ldots, \hat{u}_k\}\) be a broken trajectory from \(x\) to
\(y\), and \(u_i\) be representatives in the respective
unreduced moduli spaces. Given \((R_1, \ldots, R_k)\in \R_+^k\),
we define the {\em glued trajectory}:
\begin{equation}\label{eq:glue}\begin{split}
u_0\#_{R_1}u_1\#_{R_2}& u_2\cdots\#_{R_k}u_k(s):=\\
&\begin{cases}
u_{0, (-\infty, R_1]}(s) &\text{when \(s\leq R_1\)}\\
\tau_{2R_1}u_{1, [-R_1, R_2]}(s) &\text{when \(s\in [R_1,
  2R_1+R_2]\)}\\
&\vdots \\
\tau_{2\sum_{i=1}^{k-1}R_i}u_{k-1, [-R_{k-1}, R_k]}(s) &\text{when \(s\in [2\sum_{i=1}^{k-2}R_i+R_{k-1}, 2\sum_{i=1}^{k-1}R_i+R_{k}]\)}\\
\tau_{2\sum_{i=1}^kR_i}u_{k, [-R_k, \infty)}(s) &\text{when \(s\in [2\sum_{i=1}^{k-1}R_i+R_k, \infty)\)},\\
\end{cases}\end{split}\end{equation}
where \(\tau_L\) denotes translation by \(L\):
\[\tau_Lw(s):=w(s-L).\]

When \(\{\hat{u}_1, \ldots, \hat{u}_k\}\) is a broken orbit, we may also define the {\em glued orbit}
\begin{equation}\label{eq:cglue}\begin{split}
u_1\#_{R_1}& u_2\#_{R_2}\cdots u_k\#_{R_k}(s):=\\
&\begin{cases}
\tau_{2R_1}u_{1, [-R_1, R_2]}(s) &\text{when \(s\in [R_1,
  2R_1+R_2]\)}\\
&\vdots \\
\tau_{2\sum_{i=1}^{k-1}R_i}u_{k-1, [-R_{k-1}, R_k]}(s) &\text{when \(s\in [2\sum_{i=1}^{k-2}R_i+R_{k-1}, 2\sum_{i=1}^{k-1}R_i+R_{k}]\)}\\
\tau_{2\sum_{i=1}^kR_i}u_{k, [-R_k, R_1]}(s) &\text{when \(s\in [2\sum_{i=1}^{k-1}R_i+R_k, 2\sum_{i=1}^kR_i+R_1]\)},\\
\end{cases}\\
&\quad \qquad \text{for \(s\in \R/(2\sum_{i=1}^kR_i)\Z\).}
\end{split}\end{equation}
We shall sometimes suppress the subscript \(R_i\) from \(\#\) when it
is not important. 

To define the pregluing map, 
in the case of broken trajectories, assign each
\[\chi=\{\hat{u}_0\}\times\cdots \{\hat{u}_k\}\times (R_1, \ldots,
R_k)\in \hat{\cm}_0\times \cdots\times \hat{\cm}_k\times \R^k_+\] 
the \(\R\)-orbit \(\hat{w}_\chi\) of the glued trajectory
\[
w_\chi=u_0\#_{R_1}u_1\#_{R_2} u_2\cdots\#_{R_k}u_k
\]
in the configuration space \({\cal B}_P(x, y)\), taking \(u_i\) to be
{\em centered} representatives of \(\hat{u}_i\).
Similarly for the case of broken orbits or the parameterized case. 
Due to the exponential decay of flows to nondegenerate critical points,
these constructions typically give good approximation to flow lines
when the connecting rest points in the broken trajectory 
are nondegenerate.
In this article, they are used for handleslide bifurcations, and in
the discussion of coherent orientations.
\begin{remark*}
Equivalently, there is an unreduced version of the above construction,
where the gluing map 
maps products of unreduced moduli space to an unreduced
moduli space. 
Namely, take the space of gluing parameters to be an appropriate open
subset 
\[\breve{\Xi}(\mathbb{S})\subset\cm_0\times\cdots \times\cm_k,\]
and let the pregluing be given by the same formulae above, for
{\em fixed}  large \((R_1,\ldots, R_k)\), and {\em not necessarily
  centered} \(u_i\). Notice that there is a free \(\R^{k+1}\) action
on \(\breve{\Xi}(\mathbb{S})\), namely the product of translations on
each factor moduli space \(\cm_i\), and the quotient \(\breve{\Xi}(\mathbb{S})/\R^{k+1}=\mathbb{S}\).

The equivalence is easily seen by observing that,
given \((L_0, \cdots, L_k)\in \R^{k+1}\),
there is a unique \((L, R'_1, \ldots, R'_k)\in \R\times\R_+^k\), so that
\[
\tau_{L_0}u_0\#_{R_1}\tau_{L_1}u_1\#_{R_2}
\tau_{L_2}u_2\cdots\#_{R_k}\tau_{L_k}u_k \text{\small  \, approximates }
\tau_L(u_0\#_{R'_1}u_1\#_{R'_2} u_2\cdots\#_{R'_k}u_k(s)).
\]
(They are equal if \(u_i\) are replaced by their truncations). 
Furthermore, under this identification, a diagonal \(\R\)-translation
\((L_0,\ldots, L_k)\to (L_0+l,\ldots, L_k+l)\) corresponds to an
\(\R\) translation in the first factor \((L, R_1', \ldots, R_k')\to
(L+l, R_1', \ldots, R_k')\). Thus, we have a diffeomorphism
\[
\breve{\Xi}(\mathbb{S})/\R=\Xi(\mathbb{S}), \quad \text{by assigning \((\tau_{L_0}u_0,\ldots,
  \tau_{L_k}u_k) \mod \R \mapsto (\{\hat{u_0},\ldots, \hat{u}_k\}, R_1',\ldots, R_k')\)},
\]
and a commutative diagram
\[\begin{CD}
\breve {\Xi}(\mathbb{S})@>\text{pregluing map}>> {\cal B}\\
@V /\R VV @V /\R VV\\
\Xi(\mathbb{S})@>\text{pregluing map}>> {\cal B}/\R
\end{CD} 
\]
We prefer the reduced perspective in this article, 
because when the connecting
rest points are degenerate, the (reduced) space of gluing paramaters
\(\Xi\) can still be described in a way similar to the above
discussion, while \(\breve{\Xi}\) is no
longer a product of unreduced moduli spaces.
\end{remark*}

\subsubsection{K-models.}\label{sec:K-model}
In general, the deformation operator might not be surjective, and
the gluing theory gives a local description of the moduli space as an analytic variety
in the cokernel of the deformation operator. This is the ``Kuranishi
model''.

For our purpose, it is convenient to introduce a linear variant
of Kuranishi models, which we call ``K-models''.
This notion of K-model will be useful both for Step 2 of the
gluing procedure and in discussing the
orientation issue. 
\begin{definition*}
A {\em K-model} for a Fredholm operator \(\mathfrak{D}: E\to F\), denoted
\([\mathfrak{D}:K\to  C]_B\), or simply \([K\to C]\) when there is no danger of
confusion,
is a triple \(K, C, B\), where \(K, C\) are finite-dimensional subspaces \(K\subset E\),
\(C\subset F\) respectively, and \(B\subset E\) is a closed subspace such
that 
\begin{itemize}\itemsep -1pt
\item  \(\mathfrak{D}|_B: B\to \mathfrak{D}(B)\) is an isomorphism, and
\item there are decompositions
\(E=K\oplus B\),  \(F=C\oplus \mathfrak{D}(B)\) (possibly not orthogonally).
\end{itemize}
An {\em orientation} of a K-model is a choice of orientations for the
spaces \(K\) and \(C\). 
\end{definition*}
\begin{example*}[Standard K-models]
In this article, the ``cokernel'' \(\cok \mathfrak{D}\) refers either to the
quotient space \(F/\op{Image} (\mathfrak{D})\) or an arbitrary subspace of \(F\)
complementary to \(\op{Image}(\mathfrak{D})\).
A trivial example of K-model is \([\mathfrak{D}: \ker \mathfrak{D}\to \cok \mathfrak{D}]_B\), for any
subspace \(B\subset E\) complementary to \(\ker \mathfrak{D}\).
Such will be called a {\em standard K-model} for \(\mathfrak{D}\).
\end{example*}
We shall call \(K\) a {\it ``generalized
kernel''} of \(\mathfrak{D}\), \(C\) a {\it ``generalized cokernel''}, and \(B\)
a {\em ``B-space''}, for lack of better terminology.
The honest kernel and cokernel of \(\mathfrak{D}\) may be described in terms of
\(K\) and \(C\) via the exact sequence:
\begin{equation}\label{K-C-sequence}
0\to \ker \mathfrak{D}\stackrel{\Pi_{K}}{\longrightarrow} K \stackrel{\Pi_C\circ \mathfrak{D}}{\longrightarrow} C\to \cok \mathfrak{D}\to 0, 
\end{equation}
where \(\Pi_K\), \(\Pi_C\) are projections with respect
to the above decompositions of \(E\) and \(F\). 

Here are some other simple examples of K-models frequently encountered
in this article:
\begin{example*}[K-model of a stabilization]
Let \(\hat{\mathfrak{D}}_{\Psi}:\R^k\oplus E\to F\) denote a finite-dimensional extension of the
Fredholm map \(\mathfrak{D}: E\to F\), 
\[\hat{\mathfrak{D}}_\Psi (\vec{r},\xi)=\Psi(\vec{r})+\mathfrak{D}\xi,\] 
where \(\Psi : \R^k\to F\) is a linear map.
We call \(\hat{\mathfrak{D}}_{\Psi}\) a (rank-\(k\)) {\it stabilization} of \(\mathfrak{D}\).

Let \([K\to C]_B\) be a K-model for \(\mathfrak{D}\), and 
\[
\hat{K}:=\R^k\oplus K \subset \R^k\oplus E, \quad \hat{B}:= *\oplus
B\subset \R^k\oplus E,
\]
where \(*\) denotes the trivial vector space. Then 
\([\hat{K}\to C]_{\hat{B}}\) is a K-model for
\(\hat{\mathfrak{D}}_\Psi\), called the {\em stabilization} of \([K\to
C]_B\).
\end{example*}
\begin{example*}[Reductions of K-models]
Let \([\mathfrak{D}: K\to C]_B\) be a K-model, and suppose that there
are subspaces \(Q\subset K\), \(K'\subset K\), \(C'\subset C\) such that \(\Pi_C\circ \mathfrak{D}|_Q\)
is injective, and \(K\), \(C\) decompose as:
\[K=K'\oplus Q;\quad C=C'\oplus \Pi_C(\mathfrak{D}(Q)).\] 
Then \([K'\to C']_{B'}\)
is another K-model for \(\mathfrak{D}\), where \(B'=Q+B\). 
Such K-models will be called {\em reductions} (by \(Q\)) of \([K\to C]\).
\end{example*}

Notice that if two K-models for \(\mathfrak{D}\), 
\([\mathfrak{D}:K_1\to C_1]_{B_1}\), \([\mathfrak{D}:K_2\to C_2]_{B_2}\) have the same
B-space \(B_1=B_2\), then projections of \(K_1\) to \(K_2\)
and \(C_1\) to \(C_2\) (with respect to the decompositions \(E=K_1\oplus
B_1\), \(F=C_1\oplus \mathfrak{D}(B_1)\)) are isomorphisms, and vice versa.
In this case, we say that the two K-models are {\em
  equivalent}. Two oriented K-models are said to be {\em equivalent}
if they are equivalent K-models in the above sense, and the
projections involved are orientation-preserving. 

K-models are particularly useful in family settings.
We adopt the convention of denoting a Banach space bundle over
\(\Lambda\) by \(V^\Lambda\), with the fiber over
\(\lambda\in\Lambda\) denoted as \(V_\lambda\).
Let \(\Lambda\) be a connected manifold, and 
\(E^\Lambda\), \(F^\Lambda\) be Banach space bundles over \(\Lambda\).
Let \(\mathfrak{D}^\Lambda:=\{\mathfrak{\mathfrak{D}}_\lambda|\, \mathfrak{\mathfrak{D}}_\lambda: E_\lambda\to F_\lambda, \, \lambda\in
\Lambda\}\) be a family of uniformly bounded Fredholm
operators, continuous in operator norm. 
A {\em (family) K-model} for \(\mathfrak{D}^\Lambda\), written as
\([\mathfrak{D}^\Lambda: K^\Lambda\to  C^\Lambda]_{B^\Lambda }\), is a triple of Banach
space subbundles \(K^\Lambda\subset E^\Lambda, C^\Lambda\subset
F^\Lambda, B^\Lambda\subset E^\Lambda\),
so that the fibers over each \(\lambda\in\Lambda\), \([K_\lambda\to
C_\lambda]_{B_\lambda}\) form a K-model for \(\mathfrak{D}_\lambda\), and
\(\mathfrak{D}_\lambda|_{B_\lambda}\) has a uniformly bounded left inverse.

If \(\Lambda\) is finite-dimensional and compact, 
such K-models always exist by the Fredholmness of the family 
\(\mathfrak{D}^\Lambda\).
In contrast, \(\bigcup_\lambda \ker \mathfrak{D}_\lambda\), \(\bigcup_\lambda
\cok \mathfrak{D}_\lambda\) may not form bundles as the dimensions of the
kernels and cokernels may jump with \(\lambda\).

Two K-models \([D_1: K_1\to C_1]_{B_1}\), \([D_2: K_2\to C_2]_{B_2}\)
are said to be {\em  correlated} via the family K-model \([D^\Lambda:
K^\Lambda\to C^\Lambda]_{B^\Lambda}\) if they may be identified with
two fibers, \([D_{\lambda_1}: K_{\lambda_1}\to C_{\lambda_1}]_{B_{\lambda
    _1}}\), \([D_{\lambda_2}:K_{\lambda_2}\to C_{\lambda_2}]_{B_{\lambda
    _2}}\) over \(\lambda_1, \lambda_2\in \Lambda\). They are said to
be {\em equivalent} via the family K-model if they are equivalent to
two fibers. Finally, the notions of correlation and equivalence for 
{\em oriented} K-models are obtained by inserting the adjective
``oriented'' before every mention of K-model or family K-model in
the above paragraph.

\begin{example*}
If two Fredholm operators \(\mathfrak{D}, \mathfrak{D}'\) are close in
operator norm, one may always include them in a family 
\[\mathfrak{D}^\Lambda=\{\mathfrak{D}_\lambda
\, |\, \|\mathfrak{D}-\mathfrak{D}_\lambda\|<\varepsilon\} \quad
\text{for an \(\varepsilon\ll1\)}.\] 
Any K-model \([\mathfrak{D}: K\to C]_B\) may be extended into a family
K-model for \(\mathfrak{D}^\Lambda\), with trivial \(B^\Lambda=B\times
\Lambda\). In this case, we shall refer to the equivalence of K-models
for \(\mathfrak{D}\), \(\mathfrak{D}'\) without specifying the family
K-model---a family K-model of the above description will be implied. 
Moreover, if \(\mathfrak{D}\) is surjective, and \(\varepsilon \) is sufficiently small, \([\bigcup_\lambda \ker
\mathfrak{D}_\lambda\to *]\) form a K-model for
\(\mathfrak{D}^\Lambda\). Thus, in this case we shall refer to
correlated orientations of \(\ker \mathfrak{D}\) and \(\ker
\mathfrak{D}'\) without further specifications.
\end{example*}

\subsubsection{Gluing operators and gluing K-models.}
A major motivation to introduce K-models is that, in gluing theory, 
generalized
kernels and cokernels are typically easier to construct and work with
than the honest kernels and cokernels. This subsubsection explains
why.

We summarize the typical properties of the Fredholm operators 
appearring in Floer theories as follows.
A {\em Floer-type operator} is a Fredholm operator of
the form:
\[
\mathfrak{D}=\partial_s+A(s) : E\to F,\quad \text{where:}
\]
\begin{itemize}\itemsep -1pt
\item \(E=W(\R_s\times Y, p_2^*V)\), \(F=L(\R_s\times Y, p_2^*V)\) for
suitable Sobolev norms \(W, L\), 
\item \(V\) is an Euclidean or hermitian bundle over the manifold \(Y\),
\(\R_s\) denotes the real line parameterized by \(s\),
\(p_2: \R_s\times Y\to Y\) denotes the projection.
\item \(A(s): \Gamma(Y; V)\to \Gamma(Y;V)\) is a first order linear
differential operator, which is surjective and \(L^2\)-self-adjoint
when \(|s|\gg1\).
\end{itemize}

A {\em stabilized Floer-type operator} is a stabilization of a
Floer-type operator by multiplication with compactly-supported
functions. 
\medbreak
\noindent{\em Examples.} In Morse theory, \(Y\) is a point. 
In the symplectic Floer
theory considered in this article, \(Y=S^1\), and \(p_2^*V=\R^{2n}\)
(obtained from trivializing some \(u^*K\)). 
\(Y\) is a 3-manifold in Seiberg-Witten or
instanton Floer theories. 
\medbreak

An ordered \(k\)-tuple of Floer-type operators 
\[\mathfrak{D}_1=\partial_s +A_1(s), \ldots , \mathfrak{D}_k=\partial_s
+A_k(s) : E\to F\] 
are said to be {\em glue-able} if 
\begin{itemize}\itemsep -1pt
\item\(A_1(s)\) is constant for large
\(s\), \(A_k(s)\) constant for very negative \(s\), and for \(i=2,
\ldots, k-1\), \(A_i(s)\) is constant in \(s\) for \(|s|\gg1\);
\item\(A_i(\infty)=A_{i+1}(-\infty)\) for \(i=1, \ldots, k-1\).
\end{itemize}
Given a glue-able \(k+1\)-tuple of Floer-type operators \(\mathfrak{D}_0,
\ldots, \mathfrak{D}_k\), and \(k+1\)-tuple of functions \((f_0, \ldots,
f_k)\in E^k\) or \(F^k\), we may define the glued operator
\(\mathfrak{D}_0\#_{R_1}\cdots\#_{R_k}\mathfrak{D}_k\) and glued function
\(f_0\#_{R_1}\cdots \#_{R_k}f_{k}\) via the same formula
(\ref{eq:glue}), replacing \(u_{i, [-R_i, R_{i+1}]}\) there by
\(\mathfrak{D}_i\) and \(f_{i, [-R_i, R_{i+1}]}\) respectively, where
\(f_{i, [-R_i, R_{i+1}]}\) is the truncation
\[
f_{i, [-R_i, R_{i+1}]}=\beta_{[-R_i,
  R_{i+1}]}(s) f_i, \quad \text{where \(\beta_{[-R_i,
  R_{i+1}]}(s)=\beta(2s/R_i+2)\beta(2-2s/R_{i+1})\)},
\] 
with \(R_0\), \(R_{k+1}\) understood as \(-\infty, \infty\)
respectively, and \(\beta_{(-\infty, R]}(s):=\beta(2-2s/R)\),
\(\beta_{[-R, \infty)}(s):=\beta(2s/R+2)\).

Let \(K_i\), \(i=0, \ldots, k\) be subspaces in \(E\) or \(F\). We
denote by \(K_0\#_{R_1}\cdots \#_{R_k}K_k\subset E\) or \(F\) the subspace
\[
K_0\#_{R_1}\cdots \#_{R_k}K_k :=\{f_0\#_{R_1}\cdots \#_{R_k}f_k\, |\, f_i\in K_i, i=0,
\ldots, k\}.
\]
In parallel, let \(\mathfrak{D}_1=\partial_s+A_1(s), \ldots, \mathfrak{D}_k=\partial_s+A_k(s)\) be a
\(k\)-tuple of glueable Floer-type operators, with 
\(A_1(-\infty)=A_k(\infty)\). We call such operators {\em
  cyclically-glueable}. 
In this case, we may define the cyclically glued operator and functions
\(\mathfrak{D}_1\#_ {R_1}\cdots\#_{R_{k-1}}\mathfrak{D}_k\#_{R_k}\),
\(f_1\#_{R_1}\cdots \#_{R_{k-1}}f_k\#_{R_k}\in \Gamma(S^1_{2\sum_iR_i}\times Y;
p_2^*V)\) via the formula (\ref{eq:cglue}), with similar modifications.
The subspace \(K_1\#_{R_1}\cdots \#_{R_{k-1}}K_k\#_{R_k}\subset \Gamma(S^1_{2\sum_iR_i}\times Y;
p_2^*V)\) may also be similarly defined.
Furthermore, gluing and cyclic-gluing extend in an obvious way to
stabilized Floer-type operators. 

We denote by \(\iota_{\#}^jK_j\) the subspace 
\[
\begin{split}
\iota_{\#}^jK_j= & *\#_{R_1}\cdots
*\#_{R_{j-1}}K_j\#_{R_j}*\cdots\#_{R_{k-1}}*\subset
K_1\#_{R_1}\cdots\#_{R_{k-1}}K_k\\
&\text{or} *\#_{R_1}\cdots
*\#_{R_{j-1}}K_j\#_{R_j}*\cdots*\#_{R_{k}}\subset
K_1\#_{R_1}\cdots\#_{R_{k-1}}K_k\#_{R_k}
\end{split}\]
depending on the context. Notice that when \(R_1, \ldots, R_k\) are
sufficiently large and the subspaces \(K_j\) are finite dimensional,
then \(\iota_{\#}^j\) are injective for all \(j\).

Given \(f\in \Gamma(\R\times Y; p_2^*V)\) or
\(\Gamma(S^1\times Y; p_2^*V)\) and \(c\in \R\) or \(S^1\), let
\[\op{res}_c(f):=f|_{\{c\}\times Y}.\] For a subspace \(K\subset
\Gamma(\R\times Y; p_2^*V)\) or \(\Gamma(S^1\times Y; p_2^*V)\), let
\(\op{res}_cK\) denote the subspace
\(\{f|_{\{c\}\times Y}\, |\, f\in K\}\subset \Gamma(Y; V)\).

\begin{lemma*}[Glued K-models]\label{glue-K}
Let \(\mathfrak{D}_1, \ldots, \mathfrak{D}_k\) be a \(k\)-tuple of
glueable Floer-type operators, and \([\mathfrak{D}_i: K_i\to C_i]_{B_i}\) be K-models such that
\(\op{res}_0|_{K_i}: K_i\to \op{res}_0K_i\) is an isomorphism, 
and let \(\op{res}_0B_i\subset\op{res}_0E\) be a complementary
subspace to \(\op{res}_0K_i\). Set \[B_\#:=\Big\{f\, \Big|\,
\op{res}_0(\tau_{-2\sum_{j=1}^{i-1}R_j}f)\in \op{res}_0B_i, i=1, \ldots, k\Big\}.\]
\begin{description}\itemsep -1pt
\item[(1a)] Suppose for \(\Re\gg1\) and \(\vec{R}:=(R_1, \ldots, R_{k-1})\in [\Re,
  \infty)^{k-1}\), 
\begin{equation}\label{ind-sum}
\text{\(\mathfrak{D}_{\#\vec{R}}:=\mathfrak{D}_1\#_{R_1}\cdots\#_{R_{k-1}}\mathfrak{D}_k\)
  is Fredholm of index \(\sum_{i=1}^k\ind \mathfrak{D}_i\), 
.}\end{equation}
and \(K_{\#\vec{R}}:=K_1\#_{R_1}\cdots\#_{R_{k-1}}K_k\),
\(C_{\#\vec{R}}:=C_1\#_{R_1}\cdots\#_{R_{k-1}}C_k\).
Then 
\([\mathfrak{D}_{\#\vec{R}}: K_{\#\vec{R}}\to C_{\#\vec{R}}]_{B_\#}\)
forms a K-model. In fact, these form a
family K-model for the family of operators
\(\{\mathfrak{D}_{\#\vec{R}}\}_{\vec{R}}\). In particular, when \(\mathfrak{D}_i\) are
surjective \(\forall i\), the glued operator has a right inverse
bounded uniformly in \(R_1, \ldots, R_{k-1}\).
\item[(1b)] The same holds for $\vec{R}:=(R_1, \ldots, R_{k})\in [\Re,
  \infty)^{k}$, \(\mathfrak{D}_{\#\vec{R}}=\mathfrak{D}_1\#_{R_1}\cdots
\#_{R_{k-1}}\mathfrak{D}_k\#_{R_k}\),
\(K_{\#\vec{R}}=K_1\#_{R_1}\cdots\#_{R_{k-1}}K_k\#_{R_k}\), \(C_{\#\vec{R}}=
C_1\#_{R_1}\cdots\#_{R_{k-1}}C_k\#_{R_k}\) 
if, in addition, \(\mathfrak{D}_1, \ldots, \mathfrak{D}_k\)
  is cyclically glueable with 
\(
\ind\mathfrak{D}_{\#\vec{R}}=\sum_{i=1}^k\ind \mathfrak{D}_i
.\)
\item[(2)] Furthermore, the projection \(\Pi_{\iota_\#^jC_j}\) (with respect to the
decomposition \(F=\bigoplus_{i=1}^k\iota_\#^iC_i \oplus\mathfrak{D}_{\#\vec{R}}(B_\#)
\)) approximates 
\(\Pi_{C_j}\circ \beta_{[-R_{j-1}, R_j]}\circ
\tau_{-2\sum_{i=1}^{j-1}R_i}\), where the projection \(\Pi_{C_j}\) is
with respect to the decomposition \(F=C_j\oplus \mathfrak{D}_j(B_j)\).
\end{description}
\end{lemma*}
There are many other ways of choosing the B-space \(B_\#\) for the
statement of this Lemma to hold (cf. e.g. \cite{FH} Proposition 9); the one
described above is that which we shall stick to for the gluing
constructions in this article. Notice that
with this choice of \(B_\#\), the projection \(\Pi_{\iota_\#^jK_j}\)
(with respect to the decomposition \(E=
\bigoplus_i\iota_\#^iK_i\oplus B_\#\)) is given by 
\[
\Pi_{\iota_\#^jK_j}=(\op{res}_0|_{K_j})^{-1}\circ \Pi_{\op{res}_0K_j}\circ \op{res}_0 \circ
\tau_{-2\sum_{i=1}^{j-1}R_i},\] where the projection \(\Pi_{\op{res}_0K_j}\) is
with respect to the decomposition \(\op{res}_0E=\op{res}_0K_j\oplus\op{res}_0B_j\).

This gluing procedure also generalizes to family situations
to construct family K-models for glued family of operators from family
K-models of the family of operators to be glued. 
\begin{example*}[K-models of deformation operators at glued
  trajectories/orbits]
Let \(\{\hat{u}_0, \ldots, \hat{u}_k\}\) be a broken trajectory, and
\(u_i\) be centered representatives of \(\hat{u}_i\). Then
\[
E_{u_0\#_{R_1}\cdots \#_{R_k}u_k}=E_{u_{(-\infty, -R_1]}}\#_{R_1}
\cdots\#_{R_k} E_{v_{[R_k, \infty)}}.
\]
When $R_0, \ldots, R_k$ are large enough, \([\ker E_{u_i}, \cok
E_{u_i}]\) is a K-model for
$E_{u_{i, [-R_i, R_{i+1}]}}$. Furthermore, viewing \(\ker
E_{u_i}\) as the solution space of the first order linear differential
equation \(E_{u_i}\xi=0\), we see that \(\op{res}_0|_{\ker E_{u_i}}\)
is an isomorphism. Take \(B_i=\{f| \op{res}_0(f)\in \ker
E_{u_i}^\perp\}\). By the above lemma, 
\[
[\ker E_{u_0}\#_{R_1}\cdots \#_{R_k}\ker E_{u_k}\to \cok E_{u_0}\#_{R_1}\cdots \#_{R_k}\cok E_{u_k}]_{B_\#}
\]
is a K-model for \(E_{u_0\#_{R_1}\cdots \#_{R_k}u_k}\).
Similarly, in the case of broken orbits, we obtain a K-model for the
deformation operator at the glued orbit by cyclically gluing the standard
K-models of the deformation operators at the component trajectories.
\end{example*}

\subsubsection{Proof by contradiction and excision for right-invertibility.}
Though Lemma \ref{glue-K} above is standard, we shall include a proof
here, since it showcases the typical arguments for establishing the
(uniform) right invertibility of \(\mathfrak{D}_{w_\chi}\) required by
Step 2 of gluing: 
In simple situations, one may construct by excision
a right inverse to \(\mathfrak{D}_{w_\chi}\) from right inverses of the deformation
operators \(\mathfrak{D}_{u_i}\) associated to the gluing parameter \(\chi\). (See
e.g. \cite{D:floer, DK, salamon.park}).
In more intricate situations such as those frequently encountered in
this article, it
is often convenient to use an indirect, non-constructive method, which
we refer to as ``proof by contradiction''. This
method starts by choosing a codimension \(\ind \mathfrak{D}_{w_\chi}\) subspace
\(B_\chi\subset E\). By the Fredholmness of \(\mathfrak{D}_{w_\chi}\), if
\(\mathfrak{D}_{w_\chi}|_{B_\chi}\) is injective, then \(\mathfrak{D}_{w_\chi}\) has a
bounded right inverse \(P_\chi: F\to B_\chi\). 
Suppose otherwise, that there is a sequence of unit length \(\xi_\chi\in B_\chi\),
such that \(\mathfrak{D}_{w_\chi}\xi_\chi\to 0\). One then shows that this is
impossible by estimating \(\|\xi_\chi\|\) in terms of
\(\|\mathfrak{D}_{w_\chi}\xi_\chi\|\), showing that the former must go to 0 as the
latter does so. This estimate is usually obtained by breaking \(\xi_\chi\)
into summands \(\xi_i\) supported in different regions, and bounding
the summands using the surjectivity of \(\mathfrak{D}_{u_i}\). 
\medbreak

\noindent {\em Proof of Lemma \ref{glue-K}.} 
{\bf (1):} The proofs of (1a) and (1b) are almost identical; so we
shall focus on (1a). We follow the proof by contradiction framework.
First, fix \(\vec{R}\) and omit it from the subscripts to simplify notation.
Let \(\Psi_j: C_j\to F\) be the inclusion, and let
\(\mathfrak{D}_{\Psi_\#}: \bigoplus _jC_j\oplus B_\#\to F\) be defined
by 
\[
\mathfrak{D}_{\Psi_\#} (v_1, \ldots, v_k, \xi)=\mathfrak{D}_\# \xi+\Psi_1(v_1)\#_{R_1}\cdots\#_{R_{k-1}}\Psi_k(v_k).
\]
Choose \(R_1, \ldots, R_{k-1}\) to be large enough so that \(\dim
K_\#=\sum_i \dim K_i\), and \(\dim C_\#=\sum_i \dim C_i\).
We have the decomposition \(E=K_\#\oplus B_\#\) by construction. To show
that \([\mathfrak{D}_{\#}: K_\#\to C_\#]_{B_\#}\) indeed forms a K-model,
it suffices to show that \(\mathfrak{D}_{\Psi_\#}\) is surjective, in
view of the definition of \(K_\#, C_\#, B_\#\) and the index constraint
(\ref{ind-sum}). Since \(\mathfrak{D}_{\Psi_\#}\) is Fredholm of
index 0, it is equivalent to show that it is injective. 

Suppose the contrary, that there exists \((v_1, \ldots, v_k, \xi)\in
\bigoplus _jC_j\oplus B_\#\) with unit norm such that
\begin{equation}\label{assume-con}
\mathfrak{D}_{\Psi_\#} (v_1, \ldots, v_k, \xi)=0.\end{equation} 
The fact that
\([\mathfrak{D}_j: K_j\to C_j]_{B_j}\) is a K-model implies that the
operator
\[
\mathfrak{D}_{\Psi_j}: C_j\oplus B_j\to F, \quad (v, \eta)\mapsto \mathfrak{D}_j\eta+\Psi_j(v)
\]
has a bounded inverse, and hence
\begin{equation}\label{est:beta-xi}
\begin{split}
\|( &v_j, (\beta_{[-R_{j-1}, R_j]}\circ 
\tau_{-2\sum_{i=1}^{j-1}R_i})\xi\|\leq C\|\mathfrak{D}_j (\beta_{[-R_{j-1}, R_j]}\circ
\tau_{-2\sum_{i=1}^{j-1}R_i})\xi    +\Psi_j(v)\|\\
&\leq C\|\beta_{[-R_{j-1}, R_j]}\circ
\tau_{-2\sum_{i=1}^{j-1}R_i})\mathfrak{D}_{\Psi_\#}(v_1, \ldots, v_k,
\xi)\|\\
&\qquad +C\|\beta_{[-R_{j-1}, R_j]}'\tau_{-2\sum_{i=1}^{j-1}R_i})\xi\|
+C\|(1-\beta_{[-R_{j-1}, R_j]})\Psi_j(v_j)\|\\
&\ll 1,
\end{split}\end{equation}
using the assumption (\ref{assume-con}) for the first term, the fact
that \(\beta'_{[-R_{j-1}, R_j]}<C(\sum_iR_i^{-1})\ll1\) for the second
term, and the fact that \(\|\Psi_j(v_j)\|\) is bounded and \(R_i\) are
large for the last term. 
Meanwhile, observe that \(\beta_j:=\tau_{2\sum_{i=1}^{j-1}R_i}\beta_{[-R_{j-1}, R_j]}
\) are disjointly supported for different
\(j\), and write 
\[
1-\sum_{j=1}^{k}\beta_j=\sum_{l=1}^{k-1}\varphi_l,
\]
where \(\varphi_l\) is a non-negative function supported on
\((2\sum_{i=1}^{l-1}R_i-R_l/2, 2\sum_{i=1}^{l-1}R_i+R_l/2)\). Choose
\(R_1, \ldots, R_k\) to be large enough so that over the support of
\(\varphi_l\),
\(\mathfrak{D}_l=\mathfrak{D}_{l+1}=\partial_s+A_l(\infty)\). Since by
assumption \(\partial_s+A_l(\infty)\) has a bounded inverse, we have:
\begin{equation}\label{est:phi-xi}\begin{split}
\|\varphi_l\xi\|& \leq C'\|(\partial_s+A_l(\infty))(\varphi_l\xi)\|\\
& \leq C'\|\varphi_l\mathfrak{D}_{\Psi_\#}(v_1, \ldots, v_k,
\xi)\|+C'\|\varphi_l'\xi\|
+\Big\|\varphi_l \sum_{j=l}^{l+1}\tau_{2\sum_{i=1}^{j-1}R_i}(\beta_{[-R_{j-1},
  R_j]}\Psi_{j}(v_j))\Big\|\\
&\ll1.
\end{split}
\end{equation}
Summing (\ref{est:beta-xi}) and (\ref{est:phi-xi}) for all \(j\) and
\(l\), we obtain the desired contradiction that
\[
\|(v_1, \ldots, v_k, \xi)\|\leq \sum_j \|( v_j, \beta_j\xi)\|+\sum_l \|\varphi_l\xi\|\ll1.
\]
To see the assertion about family K-models, replace \((v_1, \ldots,
v_k, \xi)\) above by a sequence of unit vectors \(\{(v_1^\nu, \ldots, v_k^\nu,
\xi^\nu)\}_\nu\) with \(\|\mathfrak{D}_{\Psi_\#\vec{R}_\nu}(v_1^\nu, \ldots, v_k^\nu,
\xi^\nu)\|\to 0\). This is impossible by the same estimates,
since the above estimates do not depend on the specific values of
\(R_1, \ldots, R_{k-1}\). 
\medbreak

\noindent{\bf (2):} We now switch to the excision method. Let
\(\chi_j\) be a smooth cutoff function with value 1 on the support 
of \(\beta_j\), and
vanishes outside the support of \(\varphi_{j-1}+\varphi_j\) (with
\(\varphi_0:=0=:\varphi_k\)). Let \(\tilde{\chi}_l\) be a smooth
cutoff function with value 1 on the support of \(\varphi_l\), and
vanishes outside the support of \(\sum_{j=l-1}^l
\beta_j\). We choose these cutoff
functions such that \(|\chi_j'|\), \(|\tilde{\chi}_l'|\) are both
bounded by \(\op{min}_i R_i^{-1}/4\) for all \(j, l\).
Let \(\mathfrak{G}_{\Psi_j}=(g_{\Psi_j}, G_{\Psi_j}): F\to C_j\oplus
B_j\) and \(\tilde{G}_l: F\to E\) be the inverses of
\(\mathfrak{D}_{\Psi_j}\) and \(\partial_s+A_l(\infty)\)
respectively. Let \(\mathfrak{G}_{\Psi_j}^\tau=(g_{\Psi_j}^\tau,
G_{\Psi_j}^\tau ):= (g_{\Psi_j} \tau_{-2\sum_{i=1}^{j-1}R_i},
\tau_{2\sum_{i=1}^{j-1}R_i}G_{\Psi_j}\tau _{-2\sum_{i=1}^{j-1}R_i})\)
and set \(\mathfrak{G}_{\Psi_\#}: F\to \bigoplus _j
C_j\oplus B_\#\) to be
\[
\mathfrak{G}_{\Psi_\#}=\Big(g_{\Psi_1}^\tau\beta_1, \ldots, g_{\Psi_k}^\tau\beta_k, \sum_j\chi_j G_{\Psi_j}^\tau\beta_j+\sum_l \tilde{\chi}_l
\tilde{G}_l \varphi_l\Big).
\]
A straightforward computation shows that
\(\mathfrak{D}_{\Psi_\#}\mathfrak{G}_{\Psi_\#}=1+\Xi\), where \(\Xi\)
is small in operator norm, and so the inverse of
\(\mathfrak{D}_{\Psi_\#}\) is
\(\mathfrak{G}_{\Psi_\#}(1+\Xi)^{-1}\). Now, the projection from \(F\)
to \(\iota_\#^jC_j\) is given by \(\Pi_{C_j}\mathfrak{G}_{\Psi_\#}(1+\Xi)^{-1}\)
while the projection from \(F\) to \(C_j\) is given by
\(\Pi_{C_j}\mathfrak{G}_{\Psi_j}\). Claim (2) of the Lemma follows from
comparing these two.
\hfill\(\Box\)

\subsubsection{Generalizing the gluing map.}\label{K-map}
Suppose the deformation operator
\(\mathfrak{D}_{w_\chi}\) has a K-model \([K\to C]\) with nontrivial \(C\), the
construction of gluing map in \S1.2.1 Step 3 may be generalized as
follows. 

Write in local coordinates near \(w_\chi\) as in \S1.2.1, and
project (\ref{pde}) to the subspaces \(\mathfrak{D}(B)\), \(C\subset
F\) respectively, while decomposing
\[\xi=P_\chi\eta_\chi+\xi_K\quad \text{for \(\xi_K\in K\),
  \(P_\chi\eta_\chi\in B\)},\]
where \(P_\chi: \mathfrak{D}(B)\to B\) being the left inverse of
\(\mathfrak{D}_{w_\chi}|_B\). We have:  
\begin{gather*}
\eta_\chi+\Pi_{\mathfrak{D}(B)}({\cal F}(w_\chi)+\mathfrak{D}_{w_\chi}\xi_K)+\Pi_{\mathfrak{D}(B)}N_{w_\chi}(\xi_K+P_\chi\eta_\chi)=0,\\
\Pi_C({\cal F}(w_\chi)+\mathfrak{D}_{w_\chi}\xi_K)+N_{w_\chi}(\xi_K+P_\chi\eta_\chi))=0.
\end{gather*}
If \(\xi_K\) is sufficiently small, the contraction mapping theorem (Lemma
1.2.1) applies to the first equation above to obtain a solution of
\(\eta_\chi\) depending on \(\xi_K\). Substitute this into the second
equation, we obtain a {\em finite rank} equation in \(\xi_K\), which
is itself in a finite dimensional space. (The function on the LHS of
this equation is the ``Kuranishi map''). Thus, the solution space of
\(\xi\) is now an analytic variety in \(C\). If \([K\to C]\) is a
fiber of a family K-model \([K^\Xi\to C^\Xi]\) for \(\{\mathfrak{D}_{w_\chi}\}_{\chi\in
  \Xi}\), this describes the local structure of moduli space near the
image of the pregluing map as an analytic variety in the finite
dimensional vector bundle \(C^\Xi\). 
\(C^\Xi\) is the so-called ``obstruction bundle'', and 
this is essentially the ``obstruction bundle technique'' pioneered by
Taubes. 

In general, it is difficult to understand the structure of this
analytic variety. An example from this article is the case of gluing
a broken trajectory or orbit involving \(m\)
Type II handleslides, where \(m>1\) (cf. section 6). 
According to Lemma \ref{glue-K}, in this case the glued
K-model has an \(m\)-dimensional generalized cokernel.
Our inability to describe the local structure of
\(\hat{\cm}^{\Lambda, 1, +}_P\) near the stratum \(T_{P, hs-m}\) or
that of \(\hat{\cm}^{\Lambda, 1, +}_O\) near the stratum \(T_{O, hs-m}\) 
is precisely due to the lack of understanding on the relevant analytic variety
in this bundle of generalized cokernels.

\subsubsection{Typical arguments for Step 4 in Floer theory.}\label{step4}
Typically, it follows directly from the discussion on Kuranishi
structure in Step 2 that the gluing map is a local diffeomorphism.
For example, let
\(\chi=\{\hat{u}_0, \ldots, \hat{u}_k\}\times(R_1, \ldots, R_{k-1})\),
and \(\check{\chi}=\{u_0\}\times\cdots\times \{u_k\}\) for corresponding
representatives \(u_i\) of \(\hat{u}_i\) given in Remark 1.2.2. 
When \(\hat{u}_i\) are all nondegenerate, Lemma \ref{glue-K}
asserts that \(\ker \mathfrak{D}_{w_\chi}\) is isomorphic to \[\ker \mathfrak{D}_{u_0}\#
\cdots \# \ker
\mathfrak{D}_{u_k}\simeq T_{\check{\chi}}\check{\Xi}(\mathbb{S})\simeq
T_\chi\Xi(\mathbb{S})\times
\R w_\chi',\]
where the first isomorphism in the above expression 
is the differential of the pregluing map,
and the second isomorphism is due to Remark 1.2.2 and the fact 
that \(\mathfrak{D}_L(\tau_Lw_\chi)=w_\chi'\). 
On the other hand, the pregluing \(w_\chi\) is close to the
corresponding image of the gluing map, \(w\). 
Thus, \(\ker \mathfrak{D}_{w_\chi}\simeq\ker \mathfrak{D}_w=T_w
\cm_P\). These together imply that the differential of the gluing map
is an isomorphism from 
\(T_\chi\Xi\) to \(T_w\hat{\cm}_P\).

To show that the gluing map is actually surjective to a neighborhood
of \({\mathbb S}\) in \(\hat{\cm}_P^+\), one starts with
the following simple consequence of the implicit function theorem:
\begin{lemma*}
In the above situation, let \({\cal T}_\chi\subset T_{w_\chi}{\cal
  B}_P=E\) be the image of the differential of the pregluing map
at \(\chi\). Suppose the following hold for all \(\chi\in \Xi\):
\begin{itemize}\itemsep -1pt
\item \({\cal T}_\chi\), and \(w'_\chi\)
vary smoothly with \(\chi\);
\item    \(\exists \) subspaces \(B_\chi\subset E\) forming fibers of
  a bundle \(B^\Xi\to \Xi\), such that \(E\) decomposes as \(E=B_{\chi}\oplus{\cal T}_{\chi}\oplus \R
  w_{\chi}',\)
and the projections to the summands are bounded {\em uniformly} in
\(\chi\).
\end{itemize}
Let \(\exp (w_\chi, b_\chi)\in {\cal B}_P\) denote the element of
coordinates \(b_\chi\) in the local chart centered at \(w_\chi\). 
Then there is a diffeomorphism from a small tubular neighborhood of
\(\{((\chi,0), 0)\} \subset B^{\Xi}\times \R\) to a small tubular
neighborhood, \(U_\epsilon=\{\exp (\tau_Lw_{\chi}, \xi) |\,
\|\xi\|_E<\epsilon\}\subset {\cal B}_P\) defined by 
\[
((\chi, b_\chi), \tau) \mapsto \tau_L(\exp (w_\chi,b_\chi)). 
\]
\end{lemma*}
In other words, \(B^{\Xi}\) gives a
good coordinate system of a slice of the \(\R\)-action in
\(U_\epsilon\).
In our context, \({\cal B}_\chi\) is the B-space defined in Step 2, and the projection
\(\Pi_{B_\chi}=P_\chi \mathfrak{D}_{w_\chi}\). 
Proofs of analogous statements in the harder 
gauge-theoretic context, where the
\(\R\) action is replaced by the action of an infinite dimensional
gauge group, may be found in \cite{DK} 7.3, and \cite{D:floer} pp
97--99.

Together with the contraction mapping theorem stated in \S1.2.1, this
lemma implies that the gluing map surjects to a tubular neighborhood
({\em in \(E\)-norm}) in the moduli space.
However, the moduli space of broken trajectories is endowed
with the coarser chain topology instead.
Thus, a major task in Step 4 is to show that any
flowline in a chain topology neighborhood of \({\mathbb S}\) in fact
lies in \(U_\epsilon\). This requires a decay
estimate of the flow lines near the connecting rest points.

In the case where the connecting rest points are nondegenerate,
the relevant exponential decay estimate is akin to the decay estimate
for flows ending at \(y\), which has been used to derive (global)
compactness of \(\cm_P\) from Gromov (local) compactness. 
Proposition 4.4 of \cite{D:floer} is recommended for a 
well-written account of this estimate (in the gauge-theoretic
context).

\subsection{Gluing Flowlines Ending in Degenerate Critical Points.}

To verify the prediction of (RHFS2c, 3c) on the corner structure of
\(\hat{\cm}_P^{\Lambda, 1,+}\) or \(\hat{\cm}_O^{\Lambda, 1,+}\) near
\(\mathbb{T}_{P,db}\), \(\mathbb{J}_P\) or \(\mathbb{T}_{O,db}\),
one needs to glue flow lines ending at a degenerate critical point. 

Let \((J^\Lambda, X^\Lambda)\) be an admissible \((J, X)\)-homotopy,
and let \((0,y)\in {\cal P}^{\Lambda, deg}(J^\Lambda, X^\Lambda)\).
In sections 2--4, we set
\({\mathbb S}=\mathbb{T}_{P,db}\) or \(\mathbb{T}_{O,db}\), 
which consists of broken trajectories or orbits with 
all the connecting rest points being \(y\). In section 5,
\(\mathbb{S}\) is the subset in \(\mathbb{J}_P\) consisting of
connecting flow lines starting or ending in \(y\).
The space of gluing parameters in both cases will be
\(
\Xi({\mathbb S})={\mathbb S}\times S, \)
where \(S\) is an open interval in \(\Lambda\) with left or right end
\(0\).

The gluing theory in these cases differ
from the ``standard'' case outlined in \S1.2 in many aspects. 
Much of the additional complication arises from the fact that, instead
of the usual configuration space modeled locally on Sobolev spaces or
exponentially-weighted Sobolev spaces,
the moduli spaces of flows to \(y\) now embed in configuration spaces
modeled on the polynomially-weighted \(W_u\)-norm,
and the deformation operator is between the \(W_u\) and \(L_u\) spaces
introduced in \S I.5.
The main difference between working with these polynomially-weighted
spaces and the more commonly seen exponentially weighted ones is that,
the range space \(L_u\) now has larger weights in the longitudinal
direction than the domain space \(W_u\). This often implies that all
the estimates in the gluing theory need to be particularly precise in
the longitudinal direction, especially near \(y\), where the weight is
large. Below is a quick outline of the strategies adopted in sections
2--5. 

\subsubsection{Constructing pregluing.} 
Let \(\chi=(\{\hat{u}_1, \ldots, \hat{u}_k\}, \lambda)\in \Xi(\mathbb{S})\).
Due to the aforementioned
problem with large weights in the longitudinal direction, 
one needs a more delicate
pregluing construction instead of the typical one explained in \S1.2.1.

Let \(u_{\lambda, i}\) be a centered representative of \(\hat{u}_i\)
or a suitable cut-off version of it (to be specified later).
Noticing that a variation in parametrization (by \(s\)) 
of an element in \({\cal B}_P\) or \({\cal B}_O\) gives rise to a
variation of the element in the longitudinal direction, 
a natural solution to the above problem 
is to find (\(\lambda\)-dependent) diffeomorphisms
\(\gamma_{u_i}: I_i\to \R\), such that: 
\begin{itemize}\itemsep -1pt
\item setting the pregluing
  \(w_\chi(s,t)=u_{\lambda,i}(\gamma_{u_i}(s),t)\) over \(I_i\times
  S^1\), the error \(\mathcal{F}(w_\chi)\) projects trivially to the
  longitudinal direction (i.e. the direction of \(w_\chi'\)) where
  \(w_\chi(s, \cdot)\) is close to \(y\), and \(\gamma_{u_i}'=1\) elsewhere.
\item\(s_i<s_j\) if \(s_i\in I_i\),
\(s_j\in I_j\), and \(i<j\), and the closures of \(I_i\times S^1\)
cover the domain of \(w_\chi\), \(\Theta\).
\end{itemize}

The above condition gives an ODE which determines
\(\gamma_{u_i}\). Furthermore, from the ODE one may derive various 
behaviors of \(\gamma_{u_i}\), which will be important for the
estimates throughout the proof. For instance, the length of \(I_i\) is
of order \(|\lambda|^{-1/2}\) if \(u_i\) is not the first or last component of
a broken trajectory.

\subsubsection{\(\lambda\)-dependent \(W\)-norms and partitioning of \(\Theta\).}
In these settings, the gluing map to be constructed takes values in 
parameterized moduli spaces endowed with the ordinary
\(L^p_1\)-topology.
However, instead of the ordinary \(L^p_1\)-norms, we shall work with
certain weighted norms \(W_\chi\), \(L_\chi\), because the right
inverses of the deformation operator at \(w_\chi\) is not 
bounded {\em uniformly} in the ordinary Sobolev norms.  
These weighted norms are defined similarly to the \(W_u\) and \(L_u\)-norms in
\S I.5.2, and are in some sense
a combination of the \(W_{u_i}\)- or \(L_{u_i}\)-norms of the
components \(u_i\); thus, when \(u_i\) are all nondegenerate, the
right inverse of the deformation operator at \(w_\chi\) is expected to
have a uniform bound in terms of the norms of the right inverses of
the deformation operators at \(u_i\).
They are all commensurate with the
usual Sobolev norms, though dependent on the gluing parameter \(\chi\). 

When performing estimates, we typically partition \(\Theta\) into several
regions depending on whether \(\gamma_{u_i}'\) is close to 1, and
estimate over each region separately. Over the region \(\Theta_{u_i}\), the
values of \(\gamma_{u_i}'\) is close to 1, and hence \(w_\chi'\) approximates \(\partial_\gamma
u(\gamma)\), the \(W_\chi\)-norm approximates the
\(W_{u_i}\)-norm, and the deformation operator at \(w_\chi\) may be
approximated by that at \(u_i\). The length of these regions are
typically of order \(|\lambda|^{-1/2}\) or infinite, and the estimates
over these regions are similar to those in \S I.5.

In the case considered in sections 2--4, the other
regions are \(\Theta_{yj}\). They have lengths of order
\(|\lambda|^{-1/2}\), and estimates over these regions often use the facts that
 on \(\Theta_{yj}\), \(w_\chi(s, \cdot)\) is close to \(y\) (of
distance \(\leq C|\lambda|^{1/2}\) for some positive constant \(C\)),
and that \(\gamma_{u_i}(s)\) grows polynomially as positive multiples of \((|\lambda |(\mathfrak{l}-s))^{-1}\).

In the case considered in section 5, the other regions are
\(\Theta_{y\pm}\). These are of infinite length, but \(\gamma_u'\) and hence
also \(w_\chi'\) decay exponentially in the form \(C_\pm\exp(\mp
\mu_\pm|\lambda|^{1/2}s)\), \(C_\pm, \mu_\pm\) being positive constants of \(O(1)\).
In addition to this, we also often use the fact that over this region,
 \(w_\chi(s, \cdot)\) is close to the new critical points \(y_{\lambda\pm}\) (of
distance \(\leq C|\lambda|^{1/2}\) for some positive constant \(C\)),
and the estimates about \(y_{\lambda\pm}\) in \S I.5.3.

\subsubsection{K-models} 

\noindent{\sc (a) Choosing the triple \(K, C, B\).}
The deformation operators for parameterized
moduli spaces are stabilizations of those for \(\cm_P\),
\(\cm_O\), \(\mathfrak{D}_u=E_u\) or \(\tilde{D}_{u}\) respectively.
Thus, it suffices to construct K-models for the latter. 
Similar to the case in \S\ref{glue-K}, we shall always take the B-space to be
\(W_\chi'\), the subspace of \(W_\chi\) consisting of those \(\xi\)
such that \(\op{res}_{\gamma_{u_i}^{-1}(0)}\xi\) is \(L^2_t\)-orthogonal to
\(\op{res}_0\ker E_{u_i}\) \(\forall i\). The generalized
kernel will be the sum of the subspaces
\(\gamma_{u_i}^*\ker E_{u_i}=\{\gamma_{u_i}^*f\, |\, f\in \ker
E_{u_i}\}\). The generalized cokernel is
trivial in the case of section 5, but it is nontrivial in the
case of sections 2-4. In fact, by additivity of indices, its dimension
is precisely the number of connecting rest points of the broken
trajectory/orbit \(\{\hat{u}_1, \ldots, \hat{u}_k\}\). 

In this case, we choose the generalized cokernel to be spanned by 
\(\{\mathfrak{f}_j\}\), where \(\mathfrak{f}_j\) is a positive
multiple of the product of the characteristic function of \(\Theta_{yj}\)
with a unit vector in the longitudinal direction. If one requires
\(\mathfrak{f}_j\) to be of unit \(L_\chi\)-norm, the
\(L^\infty\)-norm of \(\mathfrak{f}_j\) would be of order
\(|\lambda|^{1+1/(2p)}\). Heuristically, this choice is natural in the following
sense:
\begin{enumerate}\itemsep -1pt
\item In this case, \(\mathfrak{D}_{LL}\) 
is modeled on the operator \(d/ds: L^p_1([\gamma_{u_{1}}^{-1}(0), \gamma_{u_{k}}^{-1}(0)])\to L^p([\gamma_{u_{1}}^{-1}(0), \gamma_{u_{k}}^{-1}(0)])\) while
\(B=W'_\chi\) models on the subspace of functions vanishing at the
points \(\gamma_{u_i}^{-1}(0)\). Thus, \(\mathfrak{D}(B)\) models on
the space of functions integrating to 0 on all the intervals
\([\gamma_{u_{i}}^{-1}(0), \gamma_{u_{i+1}}^{-1}(0)]\). A natural
choice for the complementary space $C$ is the subspace 
spanned by characteristic functions over these intervals. 
\item Let \(\lambda, \lambda'\) be respectively in a small death/birth
  neighborhood of \(0\), \(\chi=(\{\hat{u}_1, \ldots, \hat{u}_k\},
  \lambda)\), and \(\tilde{u}_i \in\cm_{P, \lambda'}\) be the flow
  line close to \(u_i\), \(\tilde{y}\in \cm_{P, \lambda'}\) be the
  short flow line from \(y_{\lambda' +} \) to \(y_{\lambda' -}\) close
  to the constant flow line \(\bar{y}(s)=y\) \(\forall s\). Let 
\(\tilde{y}^{inv}(s):=\tilde{y}(-s)\).
There is a glued trajectory or orbit
\(w_\#=\tilde{u}_{1}\#\tilde{y}^{inv}\#\tilde{u}_2\#\tilde{y}^{inv}\#\cdots\)
that approximates \(w_\chi\). Note that \(\ker
E_{\tilde{y}^{inv}}\simeq\cok E_{\tilde{y}}\), \(\cok
E_{\tilde{y}^{inv}}\simeq\ker E_{\tilde{y}}\); the former being
trivial, while the latter approximates the 1-dimensional space of
constant functions in the longitudinal direction (cf. \S5.3.1
below). Thus, the glued K-model for \(E_{w_\#}\)
constructed in Example \ref{glue-K} also form a K-model for
\(E_{w_\chi}\), in which the general cokernel is spanned by
\(\{*\#\cdots*\#\ker E_{\tilde{y}}\#*\cdots\}\), which approximates \(\{\mathfrak{f}_j\}\).
\end{enumerate} 
\bigbreak
\noindent{\sc (b) Proving the isomorphism.} To verify that the above choices do
give rise to a desired K-model, we need to show that the following
operators are isomorphisms with uniformly bounded inverses:
\begin{itemize}\itemsep -1pt
\item in the case of section 5, \(\mathfrak{D}_{w_\chi}|_{W'_\chi}:
  W'_\chi\to L_\chi\), 
\item in the case of sections 2-4, the stabilization \(\tilde{\mathfrak{D}}_{w_\chi}:
  \R^m\oplus W'_\chi\to L_\chi\), \(\tilde{\mathfrak{D}}_{w_\chi}(
  \iota_1, \ldots, \iota_m, \xi):=\mathfrak{D}_{w_\chi}\xi+\sum_j\iota_j \mathfrak{f}_j\).
\end{itemize}
The general outline of the proofs follows the ``proof by
contradiction'' framework sketched in \S1.2.5, estimating over different
regions in \(\Theta\) separately according to the partition outlined
in \S1.3.2, and incoporating several extra ingredients including:
\begin{itemize}\itemsep -1pt
\item  variants of Floer's lemma (cf. e.g. Lemma \ref{claimF}), which
  gives a \(L^\infty\)-bound on \(|\lambda|^{-1/2}\xi\) over
  \(\Theta_{yj}\). This is useful for ensuring that, in spite of the
  potential problem with large weights, the extra term
  \(\beta' \xi_T\) introduced by the cutoff function (as in
  (\ref{est:beta-xi})) when estimating the transversal component
  \(\xi_T\) is still sufficiently small. (A different method is needed
  for the longitudinal component, where the problem with large weights
  is worse). This estimate is also useful for bounding the \(W_\chi\)
  norm of \(\xi_T\) over \(\Theta_{yj}\).
\item estimates for \(\iota_j\) and \(\xi_L\) over \(\Theta_{yj}\). In
  contrast to estimates over \(\Theta_{u_i}\), the estimates over
  \(\Theta_{yj}\) differ substantially from the stereotype exemplified by
  the proof of Lemma \ref{glue-K}, especially for the longitudinal
  direction, since \(\partial_s+A_y\) is not surjective, or even
  Fredholm. Since \(\mathfrak{D}_{LL}\) in this region is modeled on
  \(\partial_s\), a basic tool of these estimates is a simple Lemma
  (Lemma \ref{fact}) bounding the \(L^p\) norm of a real-valued
  function \(f\) over an interval \(I\) in terms of the
  \(\|f'\|_{L^p(I)}\), the value of \(f\) at an end point of \(I\),
  and the length of \(I\). The latter are in turn bounded via 
\(\|\tilde{\mathfrak{D}}_{w_\chi}(\iota_1, \ldots, \iota_m,
\xi)\|_{W_\chi}\), the vanishing of \(\xi_L\) at the points
\(\gamma_{u_i}^{-1}(0)\), and the length estimate of \(\Theta_{yj}\).
\end{itemize}
\bigbreak
\noindent{\sc (c) Understanding the Kuranishi map.} As explained above, in the
case of sections 2-4, the Kuranishi model is more interesting, as the
Kuranishi map is nontrivial. To understand the Kuranishi
map, one needs a better description of the projection \(\Pi_C\). In
general this is not easy to compute when the decomposition
\(C\oplus \mathfrak{D}(B)\) is not orthogonal. Fortunately, due to
the special property of our \(\mathfrak{D}\) and our choice of \(C\),
there is a relatively simple way of computing \(\Pi_C\): very roughly speaking,
modulo certain typically ignorable terms and multiplication by
positive scalars, \(\Pi_{\mathfrak{f}_j}\) is
given by integrating the longitudinal component  over the interval
\([\gamma_{u_{j-1}}^{-1}(0), \gamma_{u_{j}}^{-1}(0)]\). (See Lemma
\ref{prop:P-j} for the precise statement). Notice that this conforms with the
heuristic picture sketched in item 1 of part
{\sc (a)} above.

\subsubsection{Surjectivity of gluing map.}
As explained in \S1.2.7, the main task of this step is a decay
estimate for the flow line near \(y\), which has to be particularly
precise when \(y\) is degenerate, due to the polynomially weighted
norms adopted. This will be done via various
refinements of the decay estimate in section I.5.
Given \(w=\exp (w_\chi, \xi)\in \cm_P\) in a chain-topology
neighborhood of the pregluing \(w_\chi\), we estimate the transversal and
longitudinal components of \(\xi \) separately. 
First, reparameterize \(w_\chi\) to get \(\tilde{w}\), such that the
difference between \(w\) and \(\tilde{w}\) is transversal near
\(y\). This difference satisfies a differential equation which is used
to obtain its pointwise estimate. On the other hand, comparing this
parametrization with \(\gamma_{u_i}\), which was used in the
definition of \(w_\chi\), one may estimate the difference between the
two parametrizations via an ODE, which in turn gives a pointwise
bound on the difference between \(w_\chi\) and \(\tilde{w}\) (note
that this is longitudinal). The desired bound on \(\|\xi\|_{W_\chi}\)
is obtained using 
the transversal and longitudal pointwise estimates above.

\section{Gluing at Deaths I: Pregluing and Estimates.}

The following three sections give a detailed proof of 
Proposition 2.1 below, following the outline in \S1. 

This section contains the pregluing construction, the definitions of the
Banach spaces as the domain and range of the relevant deformation
operator, the error estimates, and estimates for the nonlinear term. 
Namely, Steps 1 and 3 of the gluing construction sketched in \S1. 

\subsection{Statement of the Gluing Theorem.}
The following Proposition describes the appearance of new trajectories
and closed orbits near a death-birth bifurcation, by gluing broken
trajectories and broken orbits at a death-birth. 
These trajectories all appear for \(\lambda\) in a death-neighborhood;
for this reason, we call this a ``gluing theorem at deaths'', in
contrast to the gluing theorems in section 5, where the images of the
gluing maps project via \(\Pi_\Lambda\) to birth-neighborhoods.

\begin{proposition*}\label{8.1}
Let \((J^\Lambda, X^\Lambda)\) be an admissible \((J, X)\)-homotopy
connecting two regular pairs, and \({\bf x}, {\bf z}\) be two path
components of \({\cal P}^\Lambda\backslash {\cal P}^{\Lambda, deg}\).
Then:
\begin{description}\itemsep -1pt
\item[(a)] a chain-topology neighborhood of \({\mathbb T}_{P, db} ({\bf x}, {\bf
  z}; \Re)\) in \(\hat{\cm}_P^{\Lambda, 1, +}({\bf x}, {\bf z}; \op{wt}_{-\langle {\cal Y}\rangle,
  e_{\cal P}}\leq \Re)\) is l.m.b. along \({\mathbb T}_{P, db} ({\bf x}, {\bf
  z}; \Re)\);
\item[(b)] a chain-topology neighborhood of \({\mathbb T}_{O, db}
  (\Re)\) in \(\hat{\cm}_O^{\Lambda, 1, +}(\op{wt}_{-\langle {\cal Y}\rangle,
  e_{\cal P}}\leq \Re)\) is l.m.b. along \({\mathbb T}_{O, db} (\Re)\).
\end{description}
Furthermore, \(\Pi_\Lambda\) maps these neighborhoods to death-neighborhoods. 
\end{proposition*}

We shall focus on the proof of part (a), since the proof of part (b)
is very similar: in fact, only the discussion in section 4 on gluing
maps needs slight modification. The necessary
modification for part (b) will be briefly indicated in \S4.3. 

Recall that the admissibility of \((J^\Lambda, X^\Lambda)\) implies
that elements in \({\cal P}^{\Lambda, deg}\) satisfy (RHFS1i), and lie
in standard d-b neighborhoods, namely, satisfy the conditions
described in Definition I.5.3.1. Thus, 
by possibly restricting to a sub-homotopy and/or reversing the orientation
of \(\Lambda\), we may assume without loss of
generality that \({\cal P}^{\Lambda, deg}\) contains exactly one
point, \(y\), which is a death. Namely, the constant \[C'_y>0\]
in Definition I.5.3.1 (2b). 
We may also assume without loss of generality that 
\[
\Pi_\Lambda y=0.
\]

We now begin the construction of a gluing map from 
\(\Xi({\mathbb S})\) to \(\hat{\cm}_P^{\Lambda, 1}({\bf x}, {\bf z}; \op{wt}_{-\langle {\cal Y}\rangle,
  e_{\cal P}}\leq \Re)\), where in this case
\[
{\mathbb S}={\mathbb T}_{P, db} ({\bf x}, {\bf z}; \Re); \quad
\Xi({\mathbb S})={\mathbb S}\times (0, \lambda_0) \quad \text{for a
  small \(\lambda_0>0\).}
\]

As \({\mathbb T}_{P, db} ({\bf x}, {\bf z}; \Re)\) consists of
finitely many isolated points, we may focus on a broken trajectory 
\(\{\hat{u}_0, \ldots, \hat{u}_{k+1}\}\) in \({\mathbb T}_{P, db} ({\bf x}, {\bf z}; \Re)\).
As usual, \(u_i\) will denote the centered representative of \(\hat{u}_i\).

\subsection{The Pregluing.}
Let
$\chi:=(\{\hat{u}_0, \ldots,\hat{u}_{k+1}\}, \lambda)\in{\mathbb T}_{P,
  db} ({\bf x}, {\bf z}; \Re)\times (0,\lambda_0)$. 
Choose the representatives \(u_i\), \(i=0, \ldots, k+1\) such that
\(u_i(0)\) lies away from the neighborhood of \(y\) mentioned in
Definition I.5.3.1 (2a) and (2d).
Let \[\delta_\lambda {\cal V}:={\cal V}_{X_\lambda}-{\cal V}_{X_0}.\]
By Definition I.5.3.1 (2a), this is given by
\(\check{\theta}_{X_\lambda}-\check{\theta}_{X_0}\) when
\(\lambda<\lambda_0\) is sufficiently small. That is, when 
\(\lambda_0\) is so small such that \(J_\lambda\) is constant in
\(\lambda\) for \(\lambda\in (-\lambda_0, \lambda_0)\). We choose
\(\lambda_0\) so that this is the case, and shall
{\it simply write \(J_\lambda=J\) for such \(\lambda\).}

\begin{*lemma}\label{lemma:gamma}
Let \({\mathfrak l}_{0}=-\infty\), \(\mathfrak{l}_1=0\), and
\(\mathfrak{l}_{k+2}=\infty\). Then there exist \(\mathfrak{l}_i\in \R\),
\(i=2, \ldots, k+1\), and homeomorphisms 
\[
\gamma_{u_i}: (\frak{l}_i, \mathfrak{l}_{i+1})\to \R\quad \forall
i\in\{0, \ldots, k+1\},
\]
so that the configuration \(\underline{w}_{\chi}\in {\cal
  B}_P(x_0, z_0)\) defined by
\begin{equation}\label{equation:pregluing}
\underline{w}_{\chi}(s):=
\begin{cases}
u_i(\gamma_{u_i} (s)) & \text{for $s\in({\mathfrak l}_{i},{\mathfrak l}_{i+1})$, \(i=0, \ldots, k+1\);}\\
y &\text{for $s={\mathfrak l}_j, j=1, \ldots, k+1$}
\end{cases}
\end{equation}
satisfies 
\begin{equation}
\begin{cases}
\langle \underline{w}_{\chi}'(s),
\bar{\partial}_{J X_{\lambda}}\underline{w}_{\chi}(s)\rangle_{2,t}=0;&
\text{on
  \(\bigcup_{i=0}^{k+1}\, [\gamma_{u_i}^{-1}(0), \gamma_{u_{i+1}}^{-1}(0)]\)}\\
\gamma_{u_i}'=1 &\text{otherwise.}
\end{cases}
\end{equation}
Furthermore, 
\begin{gather*}
C_0\lambda^{-1/2}\leq -\gamma_{u_0}^{-1}(0) \leq C'_0\lambda^{-1/2}\\
C_i\lambda^{-1/2}\leq {\mathfrak l}_{i+1}-{\mathfrak l}_{i}\leq
C'_i\lambda^{-1/2} \quad \text{for \(i=1, 2, \ldots, k\).}\\
C_{k+1}\lambda^{-1/2}\leq \gamma_{u_{k+1}}^{-1}(0)-{\mathfrak l}_{k+1} \leq C'_{k+1}\lambda^{-1/2}.
\end{gather*}
\end{*lemma}
\begin{notation*}
To avoid confusion, we write
\(u_{\gamma}=\partial_{\gamma}u\) and reserve \(u'=u_s\) for
\(\partial_s u\).

\(\gamma_{u_i}\) will also be used to denote \(\gamma_{u_i}\times\id:
(\mathfrak{l}_i, \mathfrak{l}_{i+1})\times S^1\to \R\times S^1\).

\end{notation*}
\begin{proof}
We'll focus on the case of \(i=0\), since the cases with other \(i\)'s
are similar.  
From the definition, $\gamma_{u_0}(s)$ satisfies 
\[
\frac{d}{ds}\gamma_{u_0}(s)=h_{u_0}(\gamma_{u_0}(s)),
\]
where $h_{u_0}: \R\to \R$ is defined as:
\begin{equation}\label{h-}
h_{u_0}(\gamma):=
\begin{cases}
-\langle (u_0)_{\gamma}(\gamma), (\bar{\partial}_{J_\lambda X_{\lambda}}u_0)(\gamma)\rangle_{2,t}
\|(u_0)_{\gamma}(\gamma)\|_{2,t}^{-2}+1 &\text{when \(\gamma\geq0\);}\\
1 &\text{when \(\gamma<0\).}
\end{cases}
\end{equation}
Our choice of the representatives \(u_i\) ensures that \(h_{u_i}\) is
continuous. 
From the decay estimates in Proposition I.5.1.3 and the fact that
\(y\) is in standard d-b neighborhood, for large
\(\gamma\) we have
$$
(\bar{\partial}_{JX_{\lambda}})u_0(\gamma)=\delta_\lambda{\cal V}(u_0(\gamma))
=T_{y,u_0(\gamma)}\Big(\lambda C'_y {\bf e}_y\Big)+O(\lambda \gamma^{-1})+O(\lambda^2).
$$

On the other hand, $(u_0)_{\gamma}(\gamma)$ approaches the
direction $-{\bf e}_y$ for large $\gamma$; therefore there are
\(\lambda\)-independent positive constants \(A, A'\), such that
\begin{equation}\label{h-est}
A'\lambda\gamma^2\geq h_{u_0}(\gamma)\geq A\lambda \gamma^2 \qquad
\mbox{for \(\gamma\gg1\).} 
\end{equation}

We see that the inverse function of $\gamma_{u_0}(s)$, given by integration
\begin{equation}\label{gamma-int}
\int_{\gamma_{u_0}}^{\infty} \frac{d\gamma}{h_{u_0}(\gamma)}=\int_{s(\gamma_{u_0})}^0 ds'
\end{equation}
is well defined where $\gamma_{u_0}$ is large. On the other hand, $h_{u_0}$ is always
positive and goes to $1$ when $\gamma_{u_0}$ becomes negative; we see that $\gamma_{u_0}(s)$
defines a homeomorphism from $\R_-$ to $\R$.
\end{proof}
\begin{*definition}\label{def:pregluing}
The {\em pregluing} associated with the gluing data
\(\chi\) above is \((\lambda,w_{\chi})\in {\cal B}^\Lambda_P({\bf x},
{\bf z})\), where 
\[
w_{\chi}:=e_{R_-, R_+}(0,\underline{w}_{\chi}; \lambda,0),
\]
\(e_{R_-, R_+}\) are defined in I.(62), 
and 
\[R_-=\gamma_{u_0}^{-1}(-C\lambda^{-1/2});\quad
R_+=\gamma_{u_{k+1}}^{-1} (C'\lambda^{-1/2})\quad \text{for fixed
  positive constants \(C,C'\).}\]
\end{*definition}
\begin{remark*}
In general, more complicated pregluing constructions are needed if
Definition I.5.3.1 (2b) is not assumed.
\end{remark*}

The following estimates for \(R_\pm\) in terms of \(\lambda\) will be
very useful.
\begin{*lemma}\label{est:R-pm}
$C'_{\pm}\lambda^{-1/2}\leq  R_{\pm}\leq C_{\pm}\lambda^{-1/2}$ for
some \(\lambda\)-independent positive constants \(C_\pm\), \(C'_\pm\).
\end{*lemma}
\begin{proof}
We shall only demonstrate the inequalities about $R_-$, since those
for \(R_+\) are similar. 

Choose a large enough $\gamma_0$ such that
when $s\geq \gamma_0$, the decay estimate in Proposition I.5.1.3 for $u_0(s)$ and
$u_0'(s)$ holds, and 
$$\|\delta_\lambda{\cal V}(u_0(s))-T_{y, u_0(s)}(\lambda {\bf e}_y)\|_{2,t}\leq
C_{\gamma_0}\lambda\|u_0(s)\|_{2,t}.$$

This implies that when $\gamma\geq \gamma_0\lambda^{-1/2}$, 
$A'\lambda\gamma^2\geq h_{u_0}(\gamma)\geq A\lambda\gamma^2\gg1$ for some $\lambda$-independent constants
$A, A'$. Thus 
\begin{equation}\label{gamma-0}
C_2\lambda^{-1/2}\leq -\gamma_{u_0}^{-1}(\gamma_0\lambda^{-1/2})
\leq \int_{\gamma_0\lambda^{-1/2}}^{\infty} \frac{d\gamma}{A\lambda \gamma^2}
=C_1\lambda^{-1/2}.
\end{equation}
On the other hand, $\frac{d\gamma_{u_0}(s)}{ds}=h_{u_0}(\gamma_{u_0}(s))\geq 1$ always, so
\[
\gamma_{u_0}^{-1}(\gamma_0\lambda^{-1/2})-\gamma_{u_0}^{-1}(-C\lambda^{-1/2})
\leq C_2\lambda^{-1/2}.
\]
Combining the above two inequalities we get the claimed inequality for
$R_-$.
\end{proof}

\subsection{The Weighted Norms.}

Define the {\em weight function} $\sigma_{\chi}: \R\to \R^+$ by
\begin{equation}
\sigma_{\chi}(s):=
\begin{cases}
\|w_{\chi}'(s)\|_{2,t}^{-1} 
&\text{when $\gamma_{u_0}^{-1}(0)\leq s\leq \gamma_{u_{k+1}}^{-1}(0)$;}\\
\|w_{\chi}'(\gamma_{u_0}^{-1}(0))\|_{2,t}^{-1} 
&\text{when $s\leq\gamma^{-1}_{u_0}(0)$.}\\
\|w_{\chi}'(\gamma_{u_{k+1}}^{-1}(0))\|_{2,t}^{-1} 
&\text{when $ s\geq \gamma^{-1}_{u_{k+1}}(0)$.}
\end{cases}
\end{equation}
Let $\xi\in \Gamma (w_{\chi}^*K)$, define its `{\it longitudinal'} component as
\[
\xi_L(s):=\beta (s-\gamma_{u_0}^{-1}(0))\beta (\gamma_{u_{k+1}}^{-1}(0)-s)
\sigma_{\chi}(s)^2\langle w'_{\chi}(s), \xi(s)\rangle_{2,t} w'_{\chi}(s),
\]
where $\beta: \R\to [0,1]$ is the smooth cutoff function supported on
\(\R^+\) such that \(\beta(s)=1\) \(\forall s\geq1\) (cf. I.3.2.3).

The norms for the domain and range of $E_{w_{\chi}}$ are
defined as follows.
\begin{*definition}\label{w-norms}
For \(\xi\in \Gamma (w_\chi^*K)\), 
\begin{gather*}
\|\xi\|_{L_{\chi}}:=\|\sigma_{\chi}^{1/2}\xi\|_p+\|\sigma_{\chi}\xi_L\|_p;\\
\|\xi\|_{W_{\chi}}:=\|\sigma_{\chi}^{1/2}\xi\|_{p,1}+\|\sigma_{\chi}\xi'_L\|_p.
\end{gather*}
\end{*definition}
As usual, we also use \(W_{\chi}\), \(L_{\chi}\) to denote the Banach
spaces which are \(C^{\infty}\)-completion with respect to these norms.

We shall extend the norm $W_{\chi}$ to a norm on $T_{(\lambda,
  w_{\chi})}B_P^{\Lambda}({\bf x},{\bf z})=\R\oplus W_{\chi}$ in a
way such that $\hat{E}_{(\lambda, w_\chi)}$ is uniformly bounded. 
(Cf. Lemma \ref{Ebdd}).
\begin{*definition}
Define the following norm on $\hat{W}_\chi:=\R\oplus W_{\chi}$ (denoted by the same notation):
\[
\|(\alpha, \xi)\|_{\hat{W}_\chi}:=\|\xi\|_{W_{\chi}}+\lambda^{-1/(2p)-1}|\alpha |.
\]
\end{*definition}

\subsection{The Error Estimate.}
The main goal of this subsection is to obtain the following estimate:
\begin{*proposition}\label{err-est}
In the notation of \S2.2, 2.3,
\[
\|\bar{\partial}_{J_\lambda X_\lambda}w_{\chi}\|_{L_{\chi}}
\leq C\lambda^{1/2-1/(2p)}.
\]
\end{*proposition}
\begin{proof}
By direct computation, \(\bar{\partial}_{J_\lambda X_\lambda}w_{\chi}\) is
supported on
\((-R_--1, R_++1)\times S^1\), on which it is given by 
\begin{equation}\label{err}
T_{\underline{w}_{\chi},
  w_{\chi}}\Big(\tilde{\Pi}_{\underline{w}_{\chi}'}^{\perp}
\delta_\lambda{\cal V}(\underline{w}_{\chi})\Big)+r_{\lambda}(x, z), \quad \text{where}
\end{equation}
\begin{itemize}
\item \(T_{\underline{w}_{\chi}, w_{\chi}}\) is as in Notation
I.5.2.6;
\item $r_{\lambda}(x, z)$ is a `remainder term' supported on
\((-R_--1, R_++1)\backslash(-R_-, R_+)\times S^1\) which consists of terms involving
$\beta(-R_--s)\bar{x}_{\lambda}^{0,\underline{w}_{\chi}}$,
$\beta(s-R_+)\bar{z}_{\lambda}^{0,\underline{w}_{\chi}}$
and their derivatives (cf. I.(62) for notation);
\item letting \(\Pi_{\underline{w}_{\chi}'(s)}^{\perp}\) denote the \(L^2_{t}\)-orthogonal projection to the orthogonal complement of
  \(\R w_{\chi}'(s)\), \[
\tilde{\Pi}_{\underline{w}_{\chi}'(s)}^{\perp}=\begin{cases}
\Pi_{\underline{w}_{\chi}'(s)}^{\perp}&\text{  for \(s\in [\gamma_{u_0}^{-1}(0),
  \gamma_{u_{k+1}}^{-1}(0)]\);}\\
 \op{Id} &\text{otherwise.}\end{cases}\]
\end{itemize}

To estimate the terms in (\ref{err}), note:
\begin{*lemma}\label{lem:err-ptws}
When $-R_-\leq s\leq R_+$, there is a constant $C$ independent of
$\lambda$ and $s$, such that
\[\|\sigma_{\chi}(s)\tilde{\Pi}_{w_{\chi}'}^{\perp}
\delta_\lambda{\cal V}(w_{\chi}(s))\|_{\infty,t}
\leq C\lambda^{1/2} \quad \text{\(\forall \) sufficiently small $\lambda$.}\]
\end{*lemma}

Combining this Lemma with Lemma \ref{est:R-pm}, one may bound
the contribution to \(\|\bar{\p}_{J_\lambda X_{\lambda}}w_{\chi}\|_{L_{\chi}}\) from the first term
in (\ref{err}) by \(C_1\lambda^{1/2-1/(2p)}\).

The contribution from the second term can be bounded by \(C\lambda\),
using the $C^2$ bound on $J$ and $X$, and the following estimates:
\begin{eqnarray*}
\lefteqn{\sum_{k=0}^1\sup_{s\in [-R_--1,
  -R_-]}\|\partial_s^k\bar{x}_\lambda^{0,\underline{w}_{\chi}}(s)\|_{2,1,t}}\\
&\leq C_2\zeta_\lambda^x
&+C_3\sum_{k=0}^1\sup_{s\in [-R_--1,
  -R_-]}\|\mu(s)\|_{2,1,t}\\
&\leq C'_2\lambda&+ 
C_3'e^{-C_4\lambda^{1/2}}\\
&\leq C_x \lambda,&
\end{eqnarray*}
where \(\mu(s), \zeta_\lambda^x\) are defined by
\(\exp(\underline{w}_{\chi}(s), \mu(s))=x_0\), \(\exp(x_0, x_\lambda)=\zeta_\lambda^x\),
and the second inequality follows from 
the exponential decay of \(\underline{w}_{\chi}\) to
\(x_0\), and the estimate for \(R_-\) in Lemma
\ref{est:R-pm}.
Similarly, \(\sum_{k=0}^1\sup_{s\in [R_+,
  R_++1]}\|\partial_s^k\bar{z}_\lambda^{0,\underline{w}_{\chi}}(s)\|_{2,1,t}\leq C_z \lambda\).

These together implies the Proposition.
\end{proof}

\noindent {\em Proof of Lemma \ref{lem:err-ptws}.}
By Sobolev embedding it suffices to estimate the $L^2_{1,t}$ norm.
Again we will estimate only the $s\leq 0$ part, since the other parts are
entirely similar.
Let $\gamma_0$ be as in Lemma \ref{est:R-pm}. Consider
the following two cases separately. 
Case 1: $-R_-\leq s\leq \gamma_{u_0}^{-1}(\gamma_0)$; 
Case 2: $\gamma_{u_0}^{-1}(\gamma_0)\leq s \leq 0$.

\noindent \underline{\sc Case 1}: In this region 
$\|\tilde{\Pi}_{w_{\chi}'}^{\perp} (\delta_\lambda {\cal V} (w_{\chi}(s))\|_{2,1,t}\leq C\lambda$. 
On the other hand on this region $\sigma_{\chi}\leq C$; in sum we have
$\|\sigma_\chi\tilde{\Pi}_{w_{\chi}'}^{\perp} (\delta_\lambda{\cal V}(w_{\chi}(s))\|_{2,1,t}\leq C_3 \lambda$.

\noindent\underline{\sc Case 2}: In this region, the fact that \(y\)
is in a standard d-b neighborhood plus the decay estimates in
Proposition I.5.1.3 imply that for small enough $\lambda$,
\[
\begin{split}
\|\Pi_{(u_0)_{\gamma}}^{\perp}\Big(\delta_\lambda{\cal V}(u_0(\gamma_{u_0})\Big)\|_{2,t}
&\leq \lambda C_0(\|\Pi_{(u_0)_{\gamma}}^{\perp}T_{y, u_0(\gamma_{u_0})}{\bf e}_y\|_{2,t}+\|\mu(\gamma_{u_0})\|_{2,t})\\
&\leq \lambda ((1-(1+C_1^2\|b(\gamma_{u_0})\|_{2,t}^2)^{-1})^{1/2}+C''\gamma_{u_0}^{-1}\\
&\leq \lambda (C_1\|b(\gamma_{u_0})\|_{2,t}+C''\gamma_{u_0}^{-1})\\
&\leq C_4\lambda \gamma_{u_0}^{-1},
\end{split}
\] 
where \(\mu, b\) are defined by \(\exp(y, \mu(\gamma))=u_0(\gamma)\),
\(b(\gamma)=\Pi_{{\bf e}_y}\mu(\gamma)\), as in \S I.5. Meanwhile,  
\[
\begin{split}
&\Big\|\p_t[\Pi_{(u_0)_{\gamma}}^{\perp}\delta_\lambda{\cal V}(u_0(\gamma_{u_0}))]\Big\|_{2,t}\\
&\qquad\leq \|\p_t[\delta_\lambda{\cal V}(u_0(\gamma_{u_0}))]\|_{2,t}
+\Big\|\p_t(\Pi_{u_{\gamma_{u_0}}})(\delta_\lambda{\cal V}(u_0(\gamma_{u_0})))\Big\|_{2,t}\\
&\qquad\leq \lambda C_5\Big(\|\p_t(u_0(\gamma_{u_0}))\|_{2,t}+\|\sigma_{u}\p_t(u_0)_{\gamma}(\gamma_{u_0})\|_{2,t}\Big)\\
&\qquad\leq C_6\lambda \gamma_{u_0}^{-1}.
\end{split}
\]

On the other hand, we have from direct computation:
\begin{equation}\label{w'}\begin{split}
\sigma_{\chi}^{-1}(s)&=\|w'_{\chi}(s)\|_{2,t}\\
&=\Pi_{(u_0)_{\gamma}(\gamma_{u_0}(s))}\Big(\delta_\lambda{\cal
  V}(u_0(\gamma_{u_0}(s)))\Big)+\|(u_0)_{\gamma}(\gamma_{u_0}(s))\|_{2,t}\quad \text{when $\gamma_{u_0}^{-1}(0)\leq s\leq 0$.}\end{split}
\end{equation}
In particular, 
\begin{equation}\label{sigma}
C'(\lambda+\gamma_{u_0}(s)^{-2})\geq \|w'_{\chi}(s)\|_{2,t}\geq
C_7(\lambda+\gamma_{u_0}(s)^{-2})\quad \text{ when $0\geq s\geq \gamma_{u_0}^{-1}(\gamma_0)$.}
\end{equation}
In sum, in case 2
\[
\sigma_{\chi}(s)\Big\|\Pi_{w_{\chi}'}^{\perp}\delta_\lambda{\cal V}(w_{\chi}(s))\Big\|_{2,1,t}
\leq \frac{C_8\lambda\gamma_{u_0}(s)^{-1}}{\lambda+\gamma_{u_0}(s)^{-2}}
\leq C_9\lambda^{1/2}.
\]
The last step above is obtained by a simple estimate of the critical value of 
the rational function.
Combining the two cases, we have proved the lemma.
\hfill\(\Box\)

\subsection{Bounding Linear and Nonlinear Terms.} 
In the previous subsection, we obtained the estimate for the 0-th
order term of the expansion (\ref{pde}). We estimate the linear and
nonlinear terms in this subsection. In our context, this means bounding
\(\hat{E}_{(\lambda, w_{\chi})}\) and \(\hat{n}_{(\lambda,
  w_{\chi})}\). These are done respectively in Lemmas \ref{Ebdd} and
\ref{nonlinear:bound} below.
\begin{*lemma}\label{Ebdd}
With respect to the norms \(\hat{W}_\chi, L_\chi\) of \S2.3, the
deformation operator 
$\hat{E}_{(\lambda, w_{\chi})}$ is bounded uniformly in $\lambda$.  
\end{*lemma}
\begin{proof}
The uniform boundedness of $E_{w_{\chi}}$ follows from simple
adaptation of I.5.2.3. We therefore just have to estimate the \(L_{\chi}\) norm of 
\[\hat{E}_{(\lambda, w_{\chi})}(1,0)=Y_{(\lambda, w_{\chi})},\]
where \(Y_{(\lambda, w_{\chi})}\) is as in I.(63). By the properties
of \(Y_{(\lambda, w_{\chi})}\) listed following I.(63),  
\(Y_{(\lambda, w_{\chi})}\) is supported on \((-R_--1, R_++1)\times
S^1\), over which it has a \(\lambda\)-independent \(C^{\infty}_\epsilon\)
bound. Also, from (\ref{w'}) we have $\sigma_{\chi}\leq C\lambda^{-1}$.
These, together with Lemma \ref{est:R-pm}, imply
\[
\|\hat{E}_{(\lambda, w_{\chi})}(1,0)\|_{L_{\chi}}\leq C'\lambda^{-1-1/(2p)}=C'\|(1,0)\|_{\hat{W}_\chi}.
\]
for a \(\lambda\)-independent positive constant $C'$.
\end{proof}

Given $(\alpha, \xi)\in T_{(\lambda,w_{\chi})}{\cal B}_P^\Lambda({\bf x},{\bf z})$, let
\(\hat{n}_{(\lambda,w_{\chi})}(\alpha, \xi)\) be:
\[
\hat{n}_{(\lambda,w_{\chi})}(\alpha, \xi):=
T^{-1}_{w_\chi, e (\lambda,
w_{\chi};\alpha,\xi)}\bar{\p}_{J^\Lambda, X^\Lambda}(\lambda+\alpha, e (\lambda,
w_{\chi};\alpha,\xi))-\bar{\p}_{J_\lambda X_{\lambda}}w_{\chi}-\hat{E}_{(\lambda,
  w_{\chi})}(\alpha, \xi).
\]

\begin{*lemma}\label{nonlinear:bound}
There is a \(\lambda\)-independent constant \(C_n\) such that
for any \(\hat{\xi}=(\alpha, \xi), \hat{\eta}=(\alpha', \eta)\in\hat{W}_{\chi}\),
$$
\|\hat{n}_{(\lambda, w_{\chi})}(\hat{\xi})-\hat{n}_{(\lambda,
  w_{\chi})}(\hat{\eta})\|_{L_{\chi}}\leq C_n(\|\hat{\xi}\|_{\hat{W}_{\chi}}+\|\hat{\eta}\|_{\hat{W}_{\chi}})\|\hat{\xi}-\hat{\eta}\|_{\hat{W}_{\chi}}.
$$ 
\end{*lemma}
\begin{proof}
These follow from direct computations, via the
\(C^{\infty}_{\epsilon}\)-bounds of \(J, X\). First, observe the
pointwise estimate
\[
\begin{split}
&|\hat{n}_{(\lambda, w_{\chi})}(\hat{\xi})-\hat{n}_{(\lambda,
  w_{\chi})}(\hat{\eta})|\\
&\quad \leq
C_1(|\xi|+|\eta|)(|\xi-\eta|+|\nabla(\xi-\eta)|)+C_2(|\nabla\xi|+|\nabla\eta|)|\xi-\eta|\\
&\qquad +(|\alpha|+|\alpha'|)(|\alpha-\alpha'|)|Z_{\lambda\lambda}|
+\Big((|\alpha|+|\alpha'|)|\xi-\eta|+|\alpha-\alpha'|(|\xi|+|\eta|)\Big)|Z_{\lambda
  w}|,
\end{split}
\]
where \(Z_{\lambda\lambda}, Z_{\lambda w}\) are both supported on
\((-R_--1, R_++1)\times S^1\), over which they are 
\(\partial_{\lambda}^2\check{\theta}_{X_\lambda}(w_{\chi})\),
\(\partial_{\lambda}\nabla\check{\theta}_{X_\lambda}(w_{\chi})/2\) respectively, up to ignorable
terms. Estimating similarly to the proof of Lemma \ref{Ebdd}, we have
\[
\|\sigma_{\chi}Z_{\lambda\lambda}\|_{p,1}+\|\sigma_{\chi}Z_{\lambda w}\|_{p,1}\leq C'\lambda^{-1-1/(2p)}.
\]
Thus 
\[\begin{split}
&\|\hat{n}_{(\lambda, w_{\chi})}(\hat{\xi})-\hat{n}_{(\lambda,
  w_{\chi})}(\hat{\eta})\|_{L_{\chi}}\leq\Big\|\sigma_{\chi}\Big(\hat{n}_{(\lambda, w_{\chi})}(\hat{\xi})-\hat{n}_{(\lambda,
  w_{\chi})}(\hat{\eta})\Big)\Big\|_{p,1}\\
&\quad\leq
C_1(\|\sigma_{\chi}^{1/2}\xi\|_{\infty}+\|\sigma_{\chi}^{1/2}\eta\|_{\infty})\|\xi-\eta\|_{W_{\chi}}
+C_2(\|\xi\|_{W_{\chi}}+\|\eta\|_{W_{\chi}})\|\sigma_{\chi}^{1/2}(\xi-\eta)\|_{\infty}\\
&\qquad+C'\lambda^{-1-1/(2p)}(|\alpha|+|\alpha'|)(|\alpha-\alpha'|)\\
&\qquad+C'\lambda^{-1-1/(2p)}\Big((|\alpha|+|\alpha'|)\|\xi-\eta\|_{\infty}+|\alpha-\alpha'|(\|\xi\|_{\infty}+\|\eta\|_{\infty})\Big)\\
&\quad \leq
C_n(\|\hat{\xi}\|_{\hat{W}_{\chi}}+\|\hat{\eta}\|_{\hat{W}_{\chi}})\|\hat{\xi}-\hat{\eta}\|_{\hat{W}_{\chi}},
\end{split}
\]
using a Sobolev inequality to bound the \(L^{\infty}\)-norm by \(L^p_1\)-norm.
\end{proof}

\section{Gluing at Deaths II: the Kuranishi Structure.}

The purpose of this section is to introduce a K-model for the
operator \(E_{w_\chi}\). By stabilization, this also yields a K-model
for \(\hat{E}_{(\lambda,w_{\chi})}\).
The main result is summarized in Proposition
\ref{UI} below. 

\subsection{The Generalized Kernel and Generalized Cokernel.}
Given \(y\in {\cal P}\), we denote by \(\bar{y}\) the constant flow
\(\bar{y}(s)=y\) \(\forall s\).

We first partition the domain \(\Theta=\R\times S^1\) into several
regions, over which \(w_\chi\) approximates either one of \(u_i\), or
\(\bar{y}\).

For a subdomain \(\Theta'\subset\Theta\) and some norm \(L\), we
denote by \(\|\xi\|_{L(\Theta')}:=\|\xi\Big|_{\Theta'}\|_L\).

\begin{*definition}[Partioning \(\Theta\)]\label{def:partition}
Fix a small positive number \(\epsilon>\lambda\). For \(i=0, \ldots, k+1\),
let \[\mathfrak{r}_i:=(2C_{u_i}\epsilon)^{1/2}(\lambda C_y')^{-1/2},\]
where $C_{u_i}$ is the constant in the bound $\|u'_i(s)\|_2\geq
C_{u_i}/s^2$ (cf. I.5.1.3). For \(j=1, \ldots, k+1\), define 
\[\begin{split}
\Theta_{yj} :=&[\mathfrak{s}_{j-}, \mathfrak{s}_{j+}]\times S^1,\quad \text{where}\\
& \mathfrak{s}_{j-}:=\gamma_{u_{j-1}}^{-1}(\mathfrak{r}_{j-1}), \quad 
\mathfrak{s}_{j+}=
\gamma_{u_j}^{-1}(-\mathfrak{r}_j).
\end{split}\]
Let \(\Theta_{u_i}\) denote the \((i+1)\)-th component of
\(\Theta\backslash \bigcup_j\Theta_{yj}\), and let
\(\Theta'_{yj}=(\gamma_{u_{j-1}}^{-1}(\mathfrak{r}_{j-1}-1), 
\gamma_{u_j}^{-1}(-\mathfrak{r}_j+1))\times S^1\supset \Theta_{yj}\).
\end{*definition}
Notice that the `length' of the region \(\Theta_{yj}\),
\(\mathfrak{s}_{j+}-\mathfrak{s}_{j-}\), is bounded as
\begin{equation}\label{Theta0:length}
C_1(\epsilon\lambda)^{-1/2}\leq \mathfrak{s}_{j+}-\mathfrak{s}_{j-}\leq C_2(\epsilon\lambda)^{-1/2}.
\end{equation}
These inequalities follow from the arguments leading to
(\ref{gamma-0}), using respectively inequalities of the type of 
the left and the right inequalities in (\ref{h-est}).
The length of \(\Theta_{yj}'\) satisfies similar bounds, with the
constants \(C_1, C_2\) above replaced by different constants \(C_1',
C_2'\).

\begin{*definition}[Bases for generalized kernel/cokernel]
For \(i=1, \ldots, k+1\), let 
\[\mathfrak{e}_{u_i}:=\gamma_{u_i}^*u_i' \in
W_\chi.
\]
For \(j=1, \ldots, k+1\), define the following elements in \(L_\chi\):
\[
\mathfrak{f}_{j}:=C_j|\lambda|^{1+1/(2p)}\vartheta_{\Theta_{yj}}w_\chi'\|w_\chi\|_{2,t}^{-1},
\]
where \(\vartheta_{\Theta_{yj}}\) is a characteristic function
supported on $\Theta_{yj}$, and $C_j$ are constants chosen such that
$\|{\mathfrak f}_j\|_{L_\chi}=1$.
\end{*definition}
Let  
\[\begin{split}
K_\chi  & :=\op{Span}\{\mathfrak{e}_{u_i}\}_{i\in\{0, \ldots, k+1\}}\subset W_\chi;\\
\hat{K}_\chi &:=\op{Span}\{(1, 0), (0, \mathfrak{e}_{u_i})\}_{i\in\{0, \ldots, k+1\}}\subset \hat{W}_\chi;\\
C_\chi &:=\op{Span}\{\mathfrak{f}_{j}\},\quad \text{and}\\
W_\chi' &:=\Big\{\xi \, | \langle
(\gamma_{u_i}^{-1})^*\xi(0), \eta(0)\rangle_{2,t}=0 \,
\forall \eta\in \ker E_{u_i}\, \forall i\in \{0,
1, \ldots, k+1\}\Big\}\subset W_\chi.\\
\end{split}\]
(Note that linearly independent elements in $\ker E_{u_i}$ restrict
linearly independently to the circle $s=0$, since they satisfy 
a homogeneous first order differential equation.)

A quick computation shows that the \(W_\chi\)-norm on \(K_\chi\) and
the \(L_\chi\)-norm on \(C_\chi\) are commensurate with the standard
norm on Euclidean spaces with respect to the bases given above.

These are respectively fibers of Banach spaces bundles over the space
of gluing parameters \(\Xi(\mathbb{S})\), \(K^\Xi, \hat{K}^\Xi, C^\Xi, W^{'\Xi}\), and
\(\hat{W}^{'\Xi}\).

Obviously, \(W_\chi=K_\chi\oplus W_\chi' \) and
\(\hat{W}_\chi=\hat{K}_\chi\oplus W_\chi' \).
Let \(\tilde{W}_\chi:=\R^{k+1}\oplus W_\chi'\), with the standard metric on
\(\R^{k+1}\). (As usual, we denote the norm on it by
the same notation). Let
\[
\tilde{E}_\chi: \tilde{W}_\chi \to L_\chi, \quad
\tilde{E}_\chi( \iota_1, \ldots,
\iota_{k+1}, \xi):=E_{w_\chi}\xi+\sum_{j=1}^{k+1}\iota_j{\mathfrak f}_j.
\]
A quick computation using (\ref{Theta0:length}) shows 
that this is a bounded operator.
The rest of this section is devoted to proving the following
\begin{*proposition}\label{UI}
For sufficiently small \(\lambda_0\), the triples \(K^\Xi, C^\Xi, W^{'\Xi}\), and \(\hat{K}^\Xi, C^\Xi\),
\(W^{'\Xi}\) are respectively
K-models for the families of operators \(\{E_{w_\chi}\}_{\chi\in
  \Xi(\mathbb{S})}\) and \(\{\hat{E}_{(\lambda, w_\chi)}\}_{\chi\in
  \Xi(\mathbb{S})}\).

In particular, there is an inverse $\tilde{G}_{\chi}:
L_{\chi}\to\tilde{W}_{\chi}$ of $\tilde{E}_{\chi}$, which is bounded
uniformly in \(\lambda\).
\end{*proposition}
We shall concentrate on proving the existence of a uniformly bounded
\(\tilde{G}_\chi\), since the rest of the assertions follow in a
straightforward manner. The proof follows the ``proof by
contradiction'' framework outlined in \S1.2.3: since
\(\tilde{E}_\chi\) is Fredholm with \(\ind \tilde{E}_\chi=0\),
it suffices to show that there exists a \(\lambda\)-independent
constant \(C\), such that \(\|\tilde{\xi}\|_{L_\chi}\leq C\|\tilde{E}_\chi
\tilde{\xi}\|_{\tilde{W}_\chi}\) \(\forall \tilde{\xi}\in \tilde{W}_\chi\).

Suppose the contrary: that there exists a sequence
\(\{\tilde{\xi}_\lambda=(\iota_{1,\lambda}, \ldots,
\iota_{k+1, \lambda}, \xi_{\lambda})\in \tilde{W}_{\chi}\}\), with
\begin{gather}
\|\tilde{\xi}_{\lambda}\|_{\hat{W}_{\chi}}=1; \nonumber\\
\|\hat{E}_{(\lambda, w_\chi)}(\tilde{\xi}_{\lambda})\|_{L_{\chi}}=:\varepsilon_E(\lambda)\to 0\quad
\text{where $\lambda\to 0$.}\label{As:xi}
\end{gather}
We shall estimate \(\tilde{\xi}_\lambda\) in
terms of \(\varepsilon_E(\lambda)\) over the various domains
introduced in \S3.1.1, to obtain a contradiction.

\subsection{Estimates over $\Theta_{u_i}$.}
Given a diffeomorphism \(\gamma: I\to \R\), let 
\[{\cal T}^{\gamma}_{w,\underline{w}}:=(\gamma^{-1})^*T_{w,\underline{w}}: \Gamma(\Theta,
w^*K)\to \Gamma(\gamma(\Theta), (\gamma^{-1})^*\underline{w}^*K).\]
By construction, for \(\xi\in \Gamma(w_{\chi}^*K)\), 
\({\cal T}_{w_{\chi}, \underline{w}_{\chi}}^{\gamma_{u_i}}\xi \in\Gamma(u_i^*K)\).

Since the discussion in this subsection holds for all \(i\), we shall
often drop the index \(i\). For instance, \(u=u_i\) for some \(i\).

In this subsection, the `transversal' or `longitudinal components' shall 
refer to the respective components of elements in \(\Gamma(u^*K)\).

\mysubsubsection[Comparing \(W_\chi\), \(L_\chi\)-norms and \(W_u\), \(L_u\)-norms.]
According to computation in 
the proof of Lemma \ref{lemma:gamma} and the definition
of \(\Theta_{u}\), in this region $|h_{u}(\gamma_{u})-1|\leq \epsilon$.
Thus, $\sigma_u$ and $\sigma_{\chi}$ are close in this region, and by
direct computation we have:
\begin{lemma*}\label{epsilon} 
Suppose $\xi\in \Gamma (w_{\chi}^*K)$ is supported on $\Theta_u$, and
let \(\epsilon, \lambda\) be as in Definition \ref{def:partition}. Then
\begin{gather*}
(1-2\epsilon)\|\xi\|_{W_{\chi}}\leq \|{\cal T}_{w_{\chi}, \underline{w}_{\chi}}^{\gamma_u}\xi\|_{W_{u}}
\leq (1+2\epsilon)\|\xi\|_{W_{\chi}};\\
(1-2\epsilon)\|\xi\|_{L_{\chi}}\leq \|{\cal T}_{w_{\chi}, \underline{w}_{\chi}}^{\gamma_u}\xi\|_{L_{u}}
\leq (1+2\epsilon)\|\xi\|_{L_{\chi}},
\end{gather*}
and for some constant $C$,
\[
\|\hat{E}_{(\lambda, w_\chi)}(\alpha, \xi)\|_{L_{\chi}}
\geq (1+2\epsilon)^{-1}\|E_u({\cal T}_{w_{\chi},
  \underline{w}_{\chi}}^{\gamma_u}\xi)+\alpha {\cal T}_{w_{\chi},
  \underline{w}_{\chi}}^{\gamma_u}Y_{(\lambda, w_\chi)}\|_{L_u}
-(C\lambda+\epsilon)\|\xi\|_{W_{\chi}}.
\]
\end{lemma*}
\begin{remark*} The fact that \(y\) is in a standard d-b neighborhood
(more precisely, the condition Definition I.5.3.1(2c)) is used here.
In general, 
the last term on the RHS of the above inequality would be larger.
\end{remark*}

\mysubsubsection[From \(\xi_\lambda\in W'_\chi\) to \(\bar{\xi}_{u_i,
  \lambda}\in W_{u_i}'\).]
For \(i=1,\ldots, k\), and \(\xi_\lambda\in W'_\chi\), let 
\begin{equation}\label{xi-bar}
\bar{\xi}_{u_i, \lambda}:= \begin{cases}\begin{split}
(\gamma_{u_i}^{-1})^*\xi_\lambda& -\beta _i((\gamma_{u_i}g^{-1})^*\xi_\lambda)_T -\theta_{i+}
((\gamma_{u_i}^{-1})^*\xi_\lambda)_L-c_{i+}u'_i)\\
& \quad -\theta_{i-}((\gamma_{u_i}^{-1})^*\xi_\lambda)_L-c_{i-}u'_i)\end{split} 
&\text{on $(-\mathfrak{r}_i,
  \mathfrak{r}_i)\times S^1$}\\
0 &\text{outside}.
\end{cases} \end{equation}
where: 
\begin{itemize}
\item$\beta_i$ is a smooth cutoff function in \(s\) supported away from 
$(-{\mathfrak r}_{i }+1,{\mathfrak r}_{i}-1)$, being 
1 outside \((-{\mathfrak r}_{i}, {\mathfrak r}_{i})\). 
\item$\theta_{i\pm}$ are characteristic functions of \((-\infty
  -\mathfrak{s}_i)\) and \((\mathfrak{s}_i, \infty)\) respectively. 
\item $c_{i\pm}$ are constants defined by
\begin{equation}\label{c-minus}
((\gamma_{u_i}^{-1})^*\xi_{\lambda })_L(\pm\mathfrak{r}_i)=c_{i\pm}
u'(\pm\mathfrak{r}_i).
\end{equation}
\end{itemize}
For \(i=0\) or \(k+1\) and similarly defined constants \(c_0,
c_{k+1}\), let
\begin{gather*}
\bar{\xi}_{u_0, \lambda}:=\begin{cases}
\begin{split}
 {\cal T}^{\gamma_{u_0}}_{w_{\chi}, \underline{w}_{\chi}}\xi_{\lambda}&
 -\beta(s-\mathfrak{r}_0+1)({\cal T}^{\gamma_{u_0}}_{w_{\chi}, \underline{w}_{\chi}}\xi_{\lambda})_T\\
& -\theta(s-\mathfrak{r}_0)(({\cal T}^{\gamma_{u_0}}_{w_{\chi}, \underline{w}_{\chi}}\xi_{\lambda})_L-c_0 u'_0)
\end{split}&\text{on $(-\infty, \mathfrak{r}_0)\times S^1$}\\
0 &\text{outside};
\end{cases}\\
\bar{\xi}_{u_{k+1}, \lambda}:=\begin{cases}
\begin{split}
 {\cal T}^{\gamma_{u_{k+1}}}_{w_{\chi}, \underline{w}_{\chi}}\xi_{\lambda}&
 -\beta(\mathfrak{r}_{k+1}-1-s)({\cal T}^{\gamma_{u_{k+1}}}_{w_{\chi}, \underline{w}_{\chi}}\xi_{\lambda})_T\\
& -\theta(\mathfrak{r}_{k+1}-s)(({\cal T}^{\gamma_{u_{k+1}}}_{w_{\chi}, \underline{w}_{\chi}}\xi_{\lambda})_L-c_{k+1} u'_{k+1})
\end{split}&\text{on $(\mathfrak{r}_{k+1}, \infty)\times S^1$}\\
0 &\text{outside},
\end{cases}\\
\end{gather*}
where \(\beta\) is the smooth cutoff function as in Part I and
\S1.2.2, \(\theta\) is the characteristic function of \(\R^+\).

\begin{remark*}
The point of the above definition is 
to introduce cutoff on $\xi_\lambda$, while keeping
the extra terms (arising from the cutoff)
in $E_u(\bar{\xi}_{\lambda})$ ignorable. The usual smooth cutoff works
for the transversal direction, but not for the longitudinal
component, over which the weight function is greater. Instead, we replace the
longitudinal component over the cutoff region by a suitable multiple
of $u'$ determined by the matching condition (\ref{c-minus}), 
and make use of the fact that $E_u (u')=0$. 
\end{remark*}

\mysubsubsection[Estimating \(E_u\bar{\xi}_{u, \lambda}\).] The
estimate for all \(i\) are similar. Taking \(i=0\) for example, a
straightforward computation yields:
\[
\begin{split}
E_{u_0}(\bar{\xi}_{u_0, \lambda})=&(1-\beta(s-\mathfrak{r}_0+1))E_{u_0}({\cal T}^{\gamma_{u_0}}_{w_{\chi}, \underline{w}_{\chi}}\xi_{\lambda
  T})\\
& \quad +(1-\theta(s-\mathfrak{r}_0))E_{u_0}
({\cal T}^{\gamma_{u_0}}_{w_{\chi},
  \underline{w}_{\chi}}\xi_{\lambda L}(\mathfrak{r}_0))\\
& \quad -\beta'(s-\mathfrak{r}_0+1) {\cal T}^{\gamma_{u_0}}_{w_{\chi},
  \underline{w}_{\chi}}\xi_{\lambda T}\\
& \quad -\delta(s-\mathfrak{r}_0)\Big({\cal T}^{\gamma_{u_0}}_{w_{\chi},
  \underline{w}_{\chi}}\xi_{\lambda L}
                                -c_0 u'_0)\Big).
\end{split}
\]
The last term above has vanishing $L_u$-norm because of the condition
(\ref{c-minus}); by Lemma \ref{claimF} below, the $L_u$-norm of the
penultimate term can be
bounded by $C\varepsilon _0(\lambda)$, which goes to 0 as $\lambda\to 0$.
Thus, by the previous lemma, we have for small
$\lambda$ that 
\begin{equation}\label{Theta-bound}
\begin{split}
\|E_u\bar{\xi}_{u, \lambda}\|_{L_u(\gamma_u(\Theta_u))}& \leq \|E_u({\cal T}^{\gamma_u}_{w_{\chi}, \underline{w}_{\chi}}\xi_{\lambda})
\|_{L_u(\gamma_u(\Theta_u))}+C\varepsilon_0(\lambda)\\
&\leq (1+2\epsilon)\|E_{w_\chi}\xi_\lambda\|_{L_{\chi}(\Theta_u)}
+(C'\lambda+\epsilon)\|\xi_{\lambda}\|_{W_{\chi}(\Theta_u)}+C\varepsilon_0(\lambda).\\
\end{split}
\end{equation}
In the last expression, the first term goes to zero because of
(\ref{As:xi}) and the fact that over \(\Theta_u\),
\(\tilde{E}_\chi\tilde{\xi}_\lambda=E_{w_\chi}\xi_\lambda\). The second term
is small since $\|\xi_{\lambda}\|_{W_{\chi}}\leq 1$.

\mysubsubsection[Estimating \(\bar{\xi}_{u_i, \lambda}\).]
Since \(\bar{\xi}_{u,\lambda}\in W_u'\), where 
\[
W'_u:=\Big\{\xi|\, \xi\in W_u, \, \langle u'(0), \xi(0)\rangle_{2,t}=0\Big\},
\]
by the right-invertibility of $E_u$, 
$\|\bar{\xi}_{u, \lambda}\|_{W_u}\leq C\|E_u\bar{\xi}_{u, \lambda}\|_{L_u}
\leq \varepsilon$. In particular 
\begin{equation}\label{xi-minus}
\|{\cal T}^{\gamma_u}_{w_{\chi}, \underline{w}_{\chi}}\xi_{\lambda}
\|_{W_u(\gamma_u(\Theta_u))}
\leq \varepsilon_u\quad \text{when $\lambda\leq \lambda_0$ is
  sufficiently small,}
\end{equation} 
where \(\varepsilon_u\) is of the form
\[
\varepsilon_u=C\varepsilon_0(\lambda)+C_2(\lambda+\epsilon)+2\varepsilon_E(\lambda),
\] 
which can be made arbitrarily small by choosing the
small constants $\epsilon,\lambda$ appropriately.

\subsection{Estimates over $\Theta_{yj}$.}
The estimates over \(\Theta_{yj}\) for different \(j\) are similar; so
we shall drop the subscript \(j\) in the discussion below.

First, note that from the computation of (\ref{sigma}) that there exist
\(\lambda\)-independent constants $C_M, C_m$ such that
\begin{equation}\label{est:sigma-w}
C_m\lambda^{-1}\leq \sigma_{\chi}(s)\leq C_M\lambda^{-1}\quad \text{on
\(\Theta_y\).}
\end{equation}
We may therefore replace (modulo multiplication by 
a constant) the weights in the
$W_{\chi}$ and $L_{\chi}$ norms by $\lambda^{-1}$.

\mysubsubsection[Estimating the transversal component.]
In the transversal direction, the estimates are again similar to the
standard case:
By looking at the limit of \((\alpha_{\lambda}, \xi_\lambda)\), one has:
\begin{lemma*}[Floer] \label{claimF}
Let $(\alpha_{\lambda},\xi_{\lambda})$ be as in (\ref{As:xi}). Then 
for all sufficiently small
$\lambda$, 
\[
\|\xi_\lambda\|_{L^\infty(\Theta'_y)}\leq \varepsilon_0(\lambda)\lambda^{1/2}
\] 
where \(\varepsilon_0(\lambda)\) is a small positive number,
\(\lim_{\lambda\to 0}\varepsilon_0(\lambda)=0\).
\end{lemma*}
\begin{proof}
Let $(s_{\lambda}, t_{\lambda})$ be a maximum of
\(|\xi_{\lambda }|\) in \(\Theta'_y\). Consider a slight enlargement of \(\Theta'_y\), \(\Theta''_y\supset
\Theta'_y\), and let \(C>0\) be such that
\([-C^{-1}(\epsilon\lambda)^{-1/2},C^{-1}(\epsilon\lambda)^{-1/2}]\times
S^1\subset \Theta''_y\). Define
\[
\varsigma_{\lambda}(s,t):=\lambda^{-1/2}\tilde{\beta}(C(\epsilon\lambda)^{1/2}s)\,
T_{w_\chi\bar{y}}\xi_{\lambda
}(s+s_{\lambda},t)\quad \text{on \(\Theta''_y\)},
\]
where \(\tilde{\beta}\) is a smooth cutoff function supported on \((-1,1)\)
which equals \(1\) on \((-1/2,1/2)\). By (\ref{As:xi}),
\(\|\varsigma_{\lambda}\|_{p,1}\) is uniformly bounded and
thus by Sobolev embedding \(\varsigma_{\lambda}\) converges in \(C_0\)
(taking a subsequence if necessary) to a \(\varsigma_{0}\), which satisfies
\(E_{\bar{y}}\varsigma_{0}=0\). (Note that the term involving \(\iota\)
dropped out because of the assumption \(|\iota|\leq
\lambda^{1+1/(2p)}\).) Such a $\varsigma_{0}$ must be
  identically zero (cf. \cite{floer.jdg} pp.542--543); so
\[
\|\varsigma_{\lambda}\|_{L^\infty([-1,1]\times S^1)}
\to 0 \quad \text{as  \(\lambda\to 0\),}
\]
and thus
\(\|\xi_{\lambda}\|_{L^\infty(\Theta'_y)}<\varepsilon_0(\lambda)\lambda^{1/2}\).
\end{proof}
\begin{remark*}
In fact, one may be more precise about the longitudinal component: by
part (a) of Lemma \ref{fact},
\(
\|\xi_{\lambda L}\|_{L^\infty(\Theta_y')}\leq C\lambda^{1/2+1/(2p)}\epsilon^{1/(2p)-1/2}.
\)
\end{remark*}

Lemma \ref{claimF} tells us that 
$\|E_{w_\chi}(\beta_y\xi_{\lambda })_T\|_{L_\chi}\to 0$,
where $\beta_y$ is a smooth cutoff function supported on 
$\Theta''_y
$ with value $1$ on $\Theta_y'$, since contribution from the extra term due to the cutoff function
goes to zero. Thus since \(E_{w_\chi}\) is right-invertible (being close to a
conjugation of \(E_{\bar{y}}\)) on the transversal subspace,
\begin{equation}\label{xi-0T}
\|(\xi_{\lambda })_T \|_{W_\chi(\Theta_y')}\leq C\varepsilon_E(\lambda)+C'\varepsilon_0(\lambda)\to 0
\quad\text{as \(\lambda\to 0\)}.
\end{equation}

\mysubsubsection[A useful normalizing function.]
The estimates for the longitudinal components hinge on the observation
that, after certain normalization, \(E_{w_\chi}\) behaves like the simple
operator \(d/ds\) over the longitudinal components.
\begin{definition*}
Let \(\ell(s)\) be the positive real function such that 
\begin{gather}\label{equation:h-0}
\langle w'_\chi, E_{w_\chi}(\ell {\bf e}_w)\rangle_{2,t}=0 \quad
\text{and $\ell(\gamma_{u_0}^{-1}(0))=1$, where}\\
{\bf e}_w(s):=\|w'_\chi\|^{-1}_{2,t}(s) w'_\chi(s).\nonumber
\end{gather}
\end{definition*}
\begin{lemma*}\label{h0:ext}
The function $\ell$ is always positive, and
there are positive $\lambda$-independent constants $C_1, C_2$ such that 
\begin{equation}\label{est:h-0}
0<C_1\leq \ell(s)\|w'_\chi\|_{2,t}(s)^{-1}\leq C_2 \quad \forall s\in
[\gamma_{u_0}^{-1}(0), \gamma_{u_{k+1}}^{-1}(0)].
\end{equation}
Furthermore, over \([{\mathfrak s}_{j-}, {\mathfrak s}_{j+}]\)
\(\forall j\),
\begin{equation}\label{h-0'}
0<C_1\lambda\leq \ell\leq C_2 \lambda ; \quad |\ell'|\leq C|\lambda|^{1/2} |\ell|.
\end{equation}
\end{lemma*}
\begin{proof}
$\ell$ satisfies a first order linear differential equation, so its
existence and uniqueness is obvious. It also follows that $\ell$ has no
zeros, because otherwise it would be identically zero. The condition
that $\ell(\gamma_{u_0}^{-1}(0))=1$ therefore implies that $\ell$ is always
positive. 
\medbreak

(\ref{est:h-0}) follows from the next
Claim by observing that $\ell\|w'_\chi\|_{2,t}^{-1}(\gamma_{u_0}^{-1}(0))$ is
$\lambda$-independent.
\medbreak
\noindent {\bf Claim.} 
Let \({\mathfrak
  s}_{0+}:=\gamma_{u_0}^{-1}(0)\); \({\mathfrak s}_{k+2\,
  -}:=\gamma_{u_{k+1}}^{-1}(0)\). Then 
\[\Big|\ln (\ell\|w'_\chi\|_{2,t}(r_1))-\ln(\ell\|w'_\chi\|_{2,t}(r_2))\Big|\leq
C\] 
for a constant $C$ independent of $r_1,
r_2$ and $\lambda$ when $r_1, r_2$ are both in: 
\begin{description}\itemsep -1pt
\item[(a)]
$[{\mathfrak s}_{i+},{\mathfrak s}_{i+1\,-}]$ for some \(i\in \{0,1,
\ldots, k+1\}\) or
\item[(b)] \([{\mathfrak s}_{j-}, {\mathfrak s}_{j+}]\) for some \(j\in \{1,
2, \ldots, k+1\}\).
\end{description}
\medbreak

\noindent{\it Proof of the Claim:}
In case (a), set \(u=u_i\) and drop the index \(i\). In this case, 
\(|\gamma_u|\leq C_i\lambda^{-1/2}\), $\gamma'_{u}$ is
close to 1, and $\|w'_\chi\|_{2,t}$ can be
approximated by $\|u_\gamma(\gamma_{u})\|_{2,t}$. We will therefore estimate
$\ell\|u_{\gamma}\|_{2,t}^{-1}$ instead. In this region, rewrite (\ref{equation:h-0})
as:
$$
\frac{d}{ds}(\ln (\ell\|u_{\gamma}\|_{2,t}^{-1}))=(\gamma'_u-1)\frac{d}{d\gamma}(\ln \|u_{\gamma}\|_{2,t}^{-1})
$$
and integrate over \(s\). Using the estimates in section 2, it is easy to see that in this region, the $L^{\infty}$ norm of the right hand side of
the above equation can be bounded by $C\lambda|\gamma|\leq C'\lambda^{1/2}$. On
the other hand, the distance between \(r_1\) and
$r_2$ can be bounded by a multiple of $\lambda^{-1/2}$, so
the claim is verified in this case.

In case (b), set \(u=u_{j-1}\) or \(u_j\) depending on whether \(s\)
is smaller or larger than \(\mathfrak{l}_j\), and again drop the index
\(j-1\) or \(j\). In this case, \(|\gamma_u|\geq C_i\lambda^{-1/2}\), and $\lambda\|w'_\chi\|_{2,t}^{-1}$ can be bounded above and below
independently of $\lambda$ (cf. (\ref{est:sigma-w})) in this region,
so it suffices to estimate the variation in $\ell$. We write (\ref{equation:h-0}) in the form:
\[
\frac{d}{ds} (\ln
\ell)=-\|u_{\gamma}\|_{2,t}^{-1}(u_{\gamma\gamma})_{L}
\]
in this case and again integrate over $s$. 
The Claim then follows from the bound
\begin{equation}\label{equation:h-0-y}
\left|\frac{d}{ds} (\ln \ell)\right|\leq C_3|\gamma_u|^{-1}\leq
C_3'\lambda^{1/2}
\end{equation}
and bound of the distance between $r_1$ and $r_2$ can be bounded by 
$C_4\lambda^{-1/2}$.\hfill\qed
\medbreak

Continuing the proof of the Lemma, 
the first inequality in (\ref{h-0'}) is the consequence of (\ref{est:h-0})
and (\ref{est:sigma-w}). The second inequality of (\ref{h-0'}) follows
directly from (\ref{equation:h-0-y}) and the first inequality.
\end{proof}

\mysubsubsection[Estimating the longitudinal direction.]
It is convenient to introduce:
\begin{definition*}
For \(j=1, \ldots, k+1\), let the $\R$-valued function $f_j(s)$ be the unique
solution of \begin{gather}
E_{w_\chi}(f_j\ell{\bf e}_w)=f'_j\ell{\bf e}_w=\mathfrak{f}_j,\nonumber\\
\ell f_j(\mathfrak{s}_{j-})=0. \label{f0:initial}
\end{gather}
Also, let  $\phi_{\lambda}(s)$, $\psi_{\lambda  , i}(\gamma)$ be the $\R$-valued functions defined  
respectively by 
\[\begin{split}
\ell\phi_{\lambda }(s)&=\langle \xi_{\lambda }(s), {\bf e}_w(s)\rangle_{2,t},\\ 
\psi_{\lambda , i}(\gamma)&=\langle{\cal T}^{\gamma_{u_i}}_{w_{\chi},
  \underline{w}_{\chi}}\xi_{\lambda}(\gamma), {\bf
  e}_{u_i}(\gamma)\rangle_{2,t},\quad
\text{where \({\bf
  e}_{u_i}=\|u_i'\|_{2,t}^{-1}u'_i\).}\\ 
\end{split}\]
\end{definition*}

The estimates for the longitudinal components will be based on the 
following elementary lemma.
\begin{lemma*}\label{fact}
If $q\in L^p_1([0,l])$, then\newline
(a) $\|q\|_{\infty}\leq C_1l^{1-1/p}\|q'\|_p+C_2l^{-1/p}\|q\|_p.$

If furthermore \(q(0)=0\), then in addition:\newline
(b) $\|q\|_{\infty}\leq Cl^{1-1/p}\|q'\|_p;$\newline
(c) $\|q\|_p\leq C'l\|q'\|_p$. 

The positive constants $C, C', C_1, C_2$  are independent of 
$q$ and $l$. 
\end{lemma*}

Let \(\bar{\phi}_{\lambda,j}:=\phi_{\lambda }+\iota_{\lambda , j}f_j\). Then
by (\ref{h-0'}), (\ref{est:sigma-w}), (\ref{Theta0:length}), and part
(c) of the above Lemma,
\begin{equation}\label{xi-0L-}\begin{split}
&\|(\xi_\lambda)_L+\iota_{\lambda ,j}\ell f_j{\bf
  e}_w\|_{W_\chi(\Theta_{yj})}\\
&\quad\leq
C_1\left(\lambda^{1/2}\|\bar{\phi}_{\lambda, j}\|_{L^p(\Theta_{y
    j})}+\|\bar{\phi}_{\lambda, j}'\|_{L^p(\Theta_{y, j})}\right)\\
&\quad\leq C_2 \left(\epsilon^{-1/2}\|\bar{\phi}_{\lambda, j}'\|_{L^p(\Theta_{y,
    j})}+\lambda^{-1/2-1/(2p)}\epsilon^{-1/(2p)}|\psi_{\lambda, j-1}|(\mathfrak{r}_{j-1})\right)\\
&\quad\leq C_3\left(\epsilon^{-1/2}\|(\tilde{E}_\chi\tilde{\xi}_{\lambda
  L})_L\|_{L_\chi(\Theta_{y j})}+\lambda^{-1/2-1/(2p)}\epsilon^{-1/(2p)}|\psi_{\lambda, j-1}|(\mathfrak{r}_{j-1})\right)\\
&\quad\leq C_4\left(\epsilon^{-1/2}\varepsilon_E+\epsilon^{-1/2-1/(2p)}\lambda^{1/2-1/(2p)}\varepsilon_0)+\lambda^{-1/2-1/(2p)}\epsilon^{-1/(2p)}|\psi_{\lambda, j-1}|(\mathfrak{r}_{j-1})\right).
\end{split}
\end{equation}
In the last line above, \(\varepsilon_0\) comes from Lemma
\ref{claimF} and an estimate for
\((E_{w_\chi}\xi_{\lambda T})_L\) via a computation similar to
I.(47).
To estimate the last term above, note that from Lemma \ref{epsilon}, 
we have
\[
\begin{split}
&\|E_{u_i}(\psi_{\lambda , i}\mathbf{e}_{u_i})\|_{L_{u_i}([0,\mathfrak{r}_i]\times S^1)}\\
&\quad \leq
(1+2\epsilon)\|\tilde{E}_{\chi}(\tilde{\xi}_{\lambda})\|_{L_{\chi}}+(C\lambda+\epsilon)\|\xi_{\lambda
  L}\|_{W_{\chi}}+C'\|{\cal
  T}^{\gamma_{u_i}}_{w_{\chi},
  \underline{w}_{\chi}}\xi_{\lambda T}\|_{W_{u_i}([0,\mathfrak{r}_i]\times S^1)}\\
&\quad \leq C\varepsilon'_{u_i} \quad \to 0\quad \text{as $\lambda\to 0$.}
\end{split}
\]
In the above, we used 
(\ref{xi-minus}) to estimate $\|{\cal
  T}^{\gamma_{u_i}}_{w_{\chi},
  \underline{w}_{\chi}}\xi_T\|_{W_{u_i}([0,\mathfrak{r}_i]\times S^1)}$.
On the other hand, 
$$
\|E_{u_i}(\psi_{\lambda , i}\mathbf{e}_{u_i})\|_{L_{u_i}([0,\mathfrak{r}_i]\times S^1)}
\geq C\|(\sigma_{u_i}(\psi_{\lambda, i})_\gamma\|_{L^p([0,\mathfrak{r}_i])}. 
$$
Using Lemma \ref{fact} (b) and the fact that
\(\psi_{\lambda, i}(0)=0\) (because \(\xi_\lambda\in W_\chi'\)), the
previous two inequalities imply: 
\begin{equation}\label{sum:begin}
\begin{split}
|\psi_{\lambda , i} (\mathfrak{r}_i)|
&\leq C\lambda\mathfrak{r}_i^{1-1/p}\varepsilon'_{u_i}\\
& \leq C'\lambda^{1/2+1/(2p)}\epsilon^{1/2-1/(2p)}\varepsilon'_{u_i}.
\end{split}
\end{equation}
Combining this with (\ref{xi-0L-}), we have
\begin{equation}\label{xi-0L}
\|(\xi_\lambda)_L+\iota_{\lambda ,j}\ell f_j{\bf
  e}_w\|_{W_\chi(\Theta_{yj})}\leq 
C_5\left(\epsilon^{-1/2}\varepsilon_E+\epsilon^{-1/2-1/(2p)}\lambda^{1/2-1/(2p)}\varepsilon_0+\epsilon^{1/2-1/p}\varepsilon'_{u_{j-1}}\right).
\end{equation}

\subsection{Estimating $\iota_j$ and $f_j$.}

This subsection fills in the last ingredients for the proof of
Proposition \ref{UI}: estimates for \(\iota_j\) and \(f_j\) (lemmas
\ref{lemma:alpha} and \ref{lemma:f} respectively). Combining these
estimates with the estimates obtained in previous subsections, we
finish the proof of Proposition \ref{UI} in \S\ref{UI:end}. 
\begin{*lemma}\label{lemma:alpha}
Let $\lambda,\, \epsilon$ be small positive numbers as before.
Then 
\[
|\iota_{\lambda, j}|\leq \varepsilon_{\iota, j}(\lambda, \epsilon),
\]
where $\varepsilon _{\iota, j}>0$ can be made arbitrarily small as
\(\lambda\to 0\) by choosing $\epsilon$ appropriately.
\end{*lemma}

\begin{proof}
This lemma follows from a lower bound on 
\(\ell f_{j}({\mathfrak s}_{j+})\), and a upper bound on 
\(\iota_{\lambda , j}\ell f_{j}({\mathfrak
  s}_{j+})\), given respectively in (\ref{alpha1}), (\ref{UpBdd}) below.

Note from the defining equation for $f_j$ and (\ref{h-0'})
that $|f_j'|\geq C\lambda^{1/(2p)}$ for some $\lambda$-independent constant
$C$. Therefore by (\ref{Theta0:length}), (\ref{h-0'}),
and the initial condition of \(f_j\) (\ref{f0:initial}),
\[\begin{split}
\ell f_j(\mathfrak{s}_{j+})&=\ell(\mathfrak{s}_{j+})\Big(f_j(\mathfrak{s}_{j+})-f_j(\mathfrak{s}_{j-})\Big)\\
&\geq C_1\lambda^{1+1/(2p)}(\mathfrak{s}_{j+}-\mathfrak{s}_{j-})\\
&\geq C_j\epsilon^{-1/2}\lambda^{1/2+1/(2p)}.
\end{split}\]
A similar calculation establishes an analogous upper bound, and we have:
\begin{equation}\label{alpha1}
C_j'\epsilon^{-1/2}\lambda^{1/2+1/(2p)}\geq \ell f_j({\mathfrak s}_{j+})\geq
C_j \epsilon^{-1/2}\lambda^{1/2+1/(2p)}.
\end{equation}

On the other hand, 
\[
 |\iota_{\lambda, j}|\ell f_j(\mathfrak{s}_{j+})
\leq |\psi_{\lambda,
  j-1}(\mathfrak{r}_{j-1})|+|-\psi_{\lambda, j}(\mathfrak{r}_j)|
+\Big|\ell\bar{\phi}_{\lambda }(\mathfrak{s}_{j+})-\ell\bar{\phi}_{\lambda }(\mathfrak{s}_{j-})\Big|.
\]
The first two terms on the RHS are already estimated in (\ref{sum:begin});
the third term can be bounded by
\[
\Big|\ell(\mathfrak{s}_{j+})\Big(\bar{\phi}_{\lambda }(\mathfrak{s}_{j+})-\bar{\phi}_{\lambda }(\mathfrak{s}_{j-})\Big)\Big|
+\Big|\Big(\ell(\mathfrak{s}_{j+})-\ell(\mathfrak{s}_{j-})\Big)\bar{\phi}_\lambda(\mathfrak{s}_{j-})\Big|,
\]
in which the first term may be bounded via (\ref{h-0'}), Lemma \ref{fact}
(b) by
\[
C\lambda(\epsilon\lambda)^{-1/2+1/(2p)}\|\bar{\phi}_\lambda'\|_{L^p([\mathfrak{s}_{j-},
  \mathfrak{s}_{j+}])}\leq C'\lambda^{1/2+1/(2p)}(\epsilon^{-1/2+1/(2p)}\varepsilon_E+\epsilon^{-1/2}\lambda^{1/2-1/(2p)}\varepsilon_0),
\]
according to the computation in (\ref{xi-0L-}), line 3-5.

The second term, via (\ref{sum:begin}), the initial condition
(\ref{f0:initial}), and (\ref{h-0'}), may be bounded by
\[
C''\lambda^{1/2+1/(2p)}\epsilon^{1/2-1/(2p)}\varepsilon'_{u_{j-1}}.
\]
Summing up, we have:
\begin{equation}\label{UpBdd}
 |\iota_{\lambda, j}|\ell f_j(\mathfrak{s}_{j+})\leq
 C_0\epsilon^{-1/2}\lambda^{1/2+1/(2p)}
(\epsilon^{1/(2p)}\varepsilon_E
+(\lambda+\epsilon)\epsilon^{1-1/(2p)}+\varepsilon_0 (\lambda^{1/2-1/(2p)}+\epsilon^{1-1/(2p)})).
\end{equation}
Comparing with (\ref{alpha1}), we have an estimate for
\(\iota_{\lambda, j}\) as asserted in the Lemma, with 
\[
\varepsilon(\epsilon,\lambda)=C_1(\epsilon^{1/(2p)}\varepsilon_E
+(\lambda+\epsilon)\epsilon^{1-1/(2p)}+\varepsilon_0 (\lambda^{1/2-1/(2p)}+\epsilon^{1-1/(2p)})),
\] 
which can be made arbitrarily small 
by choosing \(\lambda, \epsilon \) appropriately. \end{proof}

\begin{*lemma}\label{lemma:f}
For \(j=1, \ldots, k+1\), 
\begin{gather*}
\lambda^{-1/2}\|\ell f_j{\bf e}_w\|_{L^p(\Theta_{yj})}+\lambda^{-1}\|(\ell f_j)'{\bf e}_w\|_{L^p(\Theta_{yj})}\leq C_j''\epsilon^{-1/2-1/(2p)}.
\end{gather*}
\end{*lemma}
\begin{proof}
Note that since \(f_j'=\ell^{-1}\langle {\bf
  e}_w, \mathfrak{f}_j\rangle_{2,t}>0\), \(f_j\) is increasing, and thus the
estimates leading to (\ref{alpha1}) imply that
\[\|\ell f_j\|_{L^\infty(\Theta_{yj})}\leq
C_1\epsilon^{-1/2}\lambda^{1/2+1/(2p)}.\]
On the other hand, this and (\ref{h-0'}) yield
\[\begin{split}
\|(\ell f_j)'\|_{L^\infty(\Theta_{yj})}&\leq \|\ell f_j'\|_{L^\infty(\Theta_{yj})}+\|\ell' f_{yj}\|_{L^\infty(\Theta_{yj})}\\
&\leq C_2(\lambda^{1+1/(2p)}+\lambda^{1/2}\|\ell
f_j\|_{L^\infty(\Theta_{yj})})\\
&\leq C_3\epsilon^{-1/2}\lambda^{1+1/(2p)}.
\end{split}\]
These two \(L^\infty\) bounds together with the length estimate for
\(\Theta_{yj}\) imply the lemma.
\end{proof}

\subsubsection{Concluding the proof of Proposition \ref{UI}.}\label{UI:end}
Now we have all the ingredients to finish the proof of the
Proposition.

By Lemmas \ref{lemma:alpha}, \ref{lemma:f}, the
\(\iota\)-terms in (\ref{xi-0L})
are ignorable as \(\lambda\to 0\): They are bounded by expressions of
the form
\[
C_1(\varepsilon_E\epsilon^{-1/2}+(\lambda+\epsilon)\epsilon^{1/2-1/p}+\varepsilon_0(\lambda^{1/2-1/(2p)}\epsilon^{-1/2-1/(2p)}+\epsilon^{1/2-1/p})),
\]
which can be made arbitrarily small by requiring, e.g.,
\begin{equation}\label{choice:epsilon}\begin{split}
&\text{\(\lambda=\lambda(\epsilon)\) is small enough such that}\\
&\lambda<\epsilon^5,
\varepsilon_E(\lambda)+\varepsilon_0(\lambda)<\epsilon^3, \quad \text{and}\\
&\epsilon\to 0.\end{split}
\end{equation}

Applying the comparison lemma
\ref{epsilon} to (\ref{xi-minus}) for all \(i\), and adding them to 
(\ref{xi-0T}) and (\ref{xi-0L}) for all \(j\),
we see that
\[
\|\xi_{\lambda}\|_{W_{\chi}}\leq C_4(\epsilon^{-1/2}\varepsilon_E +\lambda+\epsilon+\varepsilon_0(1+\epsilon^{-1/2-1/(2p)}\lambda^{1/2-1/(2p)}))
\ll1
\]
by the same choice (\ref{choice:epsilon}). Combining this with 
Lemma \ref{lemma:alpha}, we arrive at a contradiction to
(\ref{As:xi}). \hfill\qed

\section{Gluing at Deaths III: the Gluing Map.}
Sections 4.1 and 4.2 finish the proof of Proposition \ref{8.1} (a): 
by further analyzing the Kuranishi map associated to the K-model of
last section, we obtain a smooth gluing map, which is a local
diffeomorphism. We show that this gluing map surjects to a
neighborhood of the stratum \(\mathbb{S}\) in \S4.2.

In \S4.3, we discuss the minor modification needed to obtain part (b)
of Proposition \ref{8.1}, which glues
broken orbits at \(\lambda=0\) to creat new closed orbits for small
\(\lambda>0\). 

\subsection{Understanding the Kuranishi Map.}
In the previous section, we constructed the K-model for the family of
deformation operators, \([K^\Xi\to C^\Xi]\). According to the
discussion in \S\ref{K-map}, this yields a local description of the
moduli space as the zero locus of an analytic map. In this subsection,
we analyze this analytic map in more details; this analysis enables us
to show that the moduli space is in fact (Zariski) {\em smooth}.

Recall that Proposition \ref{UI} shows that we have decompositions:
\begin{equation}\label{decomposition:L-chi}
W_\chi=\bigoplus _i \R\mathfrak{e}_{u_i}\oplus W_\chi', \quad
L_\chi=\bigoplus_j\R {\mathfrak f}_j\oplus E_{w_\chi}(W_\chi'), 
\end{equation}
and the (non-orthogonal) projection of \(L_\chi\) to the
${\mathfrak f}_j$ direction is given by \[P_{{\mathfrak f}_j}=\Pi_j \tilde{G}_\chi,\]
where \(\Pi_j\) is the projection to the \(j\)-th \(\R\)-component of
\(\tilde{W}_\chi\).
By Proposition \ref{UI}, \(P_{{\mathfrak f}_j}\) has
uniformly bounded operator norm.
However, we need a finer estimate for
\(P_{{\mathfrak f}_j}\) to understand the Kuranishi map. 
The next Lemma is a useful tool for this purpose.
\begin{*lemma}[Projection via integration]\label{prop:P-j}
Let \(\eta\in L_\chi\), and as usual denote
\[
\underline{\eta}_L(s):=\|w'_\chi(s)\|_{2,t}^{-1}\langle w'_\chi(s),
\eta(s)\rangle_{2,t} \quad \text{for \(s\in [\gamma_{u_0}^{-1}(0), \gamma_{u_{k+1}}^{-1}(0)]\)}.
\]
Then the projection \(P_{{\mathfrak f}_j}\eta\) is bounded above and
below by expressions of the form 
\begin{equation}\label{est:P-j}
C_{1\pm}\lambda^{1/2-1/(2p)}\int_{\gamma_{u_{j-1}}^{-1}(0)}^{\gamma_{u_j}^{-1}(0)}\ell ^{-1}\underline{\eta}_L\, ds
-C_{2\pm}\lambda^{1/2-1/(2p)}\|\eta\|_{L_\chi}
\end{equation}
for \(\lambda\)-independent constants \(C_{1\pm}, C_{2\pm}\). 
\end{*lemma}
In our later applications of this lemma, the
second term in the above expression is typically dominated by
the first term, and hence ignorable.
\medbreak

\begin{proof}
In accordance with the decomposition (\ref{decomposition:L-chi}), write
\begin{equation}\label{eta-decomposition}
\eta=E_{w_\chi}\xi+\sum_{j=1}^{k+1}\iota_j{\mathfrak f}_j
\end{equation}
for \((\iota_1, \ldots, \iota_{k+1}, \xi)\in
\tilde{W}_\chi\). Thus \begin{equation}\label{compute:p-j}
\begin{split}
\int _{{\mathfrak s}_{j-}}^{{\mathfrak s}_{j+}}\ell ^{-1}ds\,
\iota_j&=C_j^{-1}\lambda^{-1-1/(2p)}\Big(\int
_{\gamma_{u_{j-1}}^{-1}(0)}^{\gamma_{u_j}^{-1}(0)}\ell ^{-1}\underline{\eta}_Lds\,
-\int_{\gamma_{u_{j-1}}^{-1}(0)}^{\gamma_{u_j}^{-1}(0)}\ell ^{-1}\underline{(E_{w_\chi}\xi)}_L\,ds\Big)\\
&=C_j^{-1}\lambda^{-1-1/(2p)}\Big(\int
_{\gamma_{u_{j-1}}^{-1}(0)}^{\gamma_{u_j}^{-1}(0)}\ell ^{-1}\underline{\eta}_Lds\,
-\int
_{\gamma_{u_{j-1}}^{-1}(0)}^{\gamma_{u_j}^{-1}(0)}\ell ^{-1}\underline{(E_{w_\chi}\xi_T)}_L\,ds\Big).
\end{split}\end{equation}
The second identity above is due to the fact that
\(\xi_L(\gamma_{u_{j-1}}^{-1}(0))=0=\xi_L(\gamma_{u_j}^{-1}(0))\):
writing \(\xi_L(s)=\ell \phi w_\chi'\|w_\chi'\|_{2,t}^{-1}(s)\),
we see from the definition of \(\ell \) that
\(\ell ^{-1}\underline{(E_{w_\chi}\xi_L)}_L=\phi'\). Now integrate by
parts, using the fact that since \(\ell (s)\neq0\, \forall s\), 
\(\phi(\gamma_{u_{j-1}}^{-1}(0))=0=\phi(\gamma_{u_j}^{-1}(0))\).

By (\ref{h-0'}) and (\ref{Theta0:length}), we have
\begin{equation}\label{h-0:integral}
C_h\lambda^{-3/2}\leq \int
_{{\mathfrak s}_{j-}}^{{\mathfrak s}_{j+}}\ell ^{-1}\,ds\leq
C_h'\lambda^{-3/2};
\end{equation}
on the other hand, a computation similar to that leading to I.(47) yields
\[
\Big|\ell ^{-1}\underline{(E_{w_\chi}\xi_T)}_L\Big|(s)\leq
C_1\Big(\sigma_\chi\|u_\gamma\|_{2,t}^{-1}(\gamma_u'+1)\|(u_{\gamma\gamma})_T\|_{2,t}+\lambda\Big)\|\xi_T\|_{2,t}(s)\leq
C_2\|\xi_T\|_{2,t}(s),
\]
where \(u=u_{j-1}\) or \(u_j\) depending on whether \(s<{\mathfrak
  l}_j\) or \(>{\mathfrak l}_j\). So by the estimates for
\(\|w_\chi'\|_{2,t}\) and
\(\gamma_{u_i}^{-1}(0)-\gamma_{u_{i-1}}^{-1}(0)\), and the uniform
boundedness of \(\tilde{G}_\chi\), 
\begin{equation}\label{est:E-xi}
\Big|\int
_{\gamma_{u_{j-1}}^{-1}(0)}^{\gamma_{u_j}^{-1}(0)}\ell ^{-1}\underline{(E_{w_\chi}\xi_T)}_L\,ds\Big|\leq
C_3\|\sigma_\chi^{-1/2}\|_{L^q((\gamma_{u_{j-1}}^{-1}(0),
\gamma_{u_j}^{-1}(0))\times S^1)}\|\xi_T\|_{W_\chi}\leq C_4\|\eta\|_{L_\chi},
\end{equation}
where \(q^{-1}:=1-p^{-1}\). Putting (\ref{compute:p-j}),
(\ref{h-0:integral}), (\ref{est:E-xi}) together, we arrive at
(\ref{est:P-j}).
\end{proof}

Next, applying the recipe of \S1.2.5 to the K-models given by
Proposition \ref{UI}, we look for solutions \((\alpha, \varphi_0,
\ldots, \varphi_{k+1}, \xi )\in W'_\chi\oplus \R\oplus  \R^{k+2}\) of:
\begin{gather}
P^c\Big(\bar{\p}_{J X_\lambda}w_\chi+\hat{E}_{(\lambda, w_\chi)}(\alpha,
\xi)+\sum_{i=0}^{k+1}\varphi_iE_{w_\chi}{\mathfrak
  e}_{u_i}+\hat{n}_{(\lambda, w_\chi)}(\alpha, \xi+\sum_{i=0}^{k+1}\varphi_i{\mathfrak e}_{u_i})\Big)=0;\label{equation:image}\\
P_{{\mathfrak f}_j}\Big(\bar{\p}_{J
  X_\lambda}w_\chi+\hat{E}_{(\lambda,
  w_\chi)}(\alpha, \sum_{i=0}^{k+1}\varphi_i{\mathfrak
e}_{u_i})+\hat{n}_{(\lambda, w_\chi)}(\alpha, \xi+\sum_{i=0}^{k+1}\varphi_i{\mathfrak e}_{u_i})\Big)=0,\label{equation:cok}
\end{gather}
where \(P^c:=1-\sum_jP_{{\mathfrak f}_j}\).
\begin{*lemma}[Solving the infinite dimensional equation]\label{est:im}
Given \[\hat{\varphi}:=(\alpha, \varphi_0, \ldots, \varphi_{k+1})\in \R\oplus\R^{k+2},
\quad \text{with \(|\hat{\varphi}|^2:=|\lambda^{-3/2}\alpha|^2+\sum_i|\varphi_i|^2\ll1\), }\] 
(\ref{equation:image}) has a unique solution \(\xi(\hat{\varphi})\), with 
\begin{equation}\label{xi-West}
\|\xi(\hat{\varphi})\|_{W_\chi}\leq C_1\lambda^{1/2-1/(2p)}+C_2\lambda^{1/2-1/(2p)}|\hat{\varphi}|+C_3\lambda^{1-1/(2p)}|\hat{\varphi}|^2.
\end{equation}
Furthermore, the solution \(\xi(\hat{\varphi}')\) corresponding to
another \(\hat{\varphi}'=(\alpha', \varphi_0', \ldots, \varphi_{k+1}')\) satisfies
\begin{equation}\label{delta-xi-West}
\|\xi(\hat{\varphi})-\xi(\hat{\varphi}')\|_{W_\chi}\leq C_4\lambda^{1/2-1/(2p)}|\hat{\varphi}-\hat{\varphi}'|.
\end{equation}
\end{*lemma}
The significance of the factor \(|\lambda|^{2/3}\) associated to \(\alpha\) in the
definition of \(|\hat{\varphi}|\) will become clear in (\ref{alpha-small}).
\medbreak
\begin{proof}
To apply the usual contraction mapping argument (Lemma 1.2.1)
to (\ref{equation:image}), we need to estimate the ``error''
\(\mathcal{F}\) and the nonlinear term \(N\), and to show that the
linearization has a uniformly bounded right inverse.

In this context, the error \(\mathcal{F}\) consists of 
\[P^c\bar{\p}_{J
  X_\lambda}w_\chi+P^c\sum_{i=0}^{k+1}\varphi_iE_{w_\chi}{\mathfrak
  e}_{u_i}+\alpha P^c Y_{(\lambda, w_\chi)}+P^c \hat{n}_{(\lambda, w_\chi)}(\alpha, \sum_{i=0}^{k+1}\varphi_i{\mathfrak e}_{u_i}),\]
which we estimate term by term below.
We shall drop all $P^c$ from the terms, since by Proposition
\ref{UI}, it has uniformly bounded operator norm, and thus 
only affect the estimate by a \(\lambda\)-independent factor.

For the first term, note that \(\|\bar{\p}_{J
  X_\lambda}w_\chi\|_{L_\chi}\) is readily estimated by Proposition \ref{err-est}.

For the second term, we claim:
\begin{equation}\label{Eest:xi-u}
\|E_{w_\chi}\gamma_{u_i}^*({u_i})_\gamma\|_{L_{\chi}}\leq C\lambda^{1/2-1/(2p)}.
\end{equation}
We shall again suppress the subscript \(i\) below. Note that 
\begin{equation}\label{E:xi-u}
E_{w_\chi}\gamma_{u}^*u_\gamma=(\gamma_{u}'-1)u_{\gamma\gamma}(\gamma_u)+Z(u)u_\gamma(\gamma_u),
\end{equation}
where \(Z\) arises from the difference between \(X_\lambda\) and
\(X_0\), and hence \(\|Z(u)\|_{\infty}\leq C\lambda\). Thus, by
I.5.3.1 (2c) and routine estimates, the
\(L_\chi\)-norm of the second term above is also bounded by \(C'\lambda\).

For the first term, note that the length of $[\gamma_u^{-1}(-\gamma_0), \gamma_u^{-1}(\gamma_0)]$ is
bounded independently of $\lambda$; therefore the $L_{\chi}$ norm of
it in this region is bounded by $C\lambda$. On the other hand, the length
of the intervals where \(|\gamma_u|\geq \gamma_0\) is bounded by
$C'\lambda^{-1/2}$. When $s$ is in these intervals, by the
computations in \S\ref{lem:err-ptws} case (2),
\(\|\sigma_\chi(\gamma_{u}'-1)u_{\gamma\gamma}(\gamma_u)\|_{\infty}\leq
C_3\lambda^{1/2}\). Together with the length estimate above, we see
that the \(L_\chi\) norm in this region is bounded by
$C''\lambda^{1/2-1/(2p)}$. (\ref{Eest:xi-u}) is verified. 

For the third term, recall the following estimate obtained 
in the proof of Lemma \ref{Ebdd}:
\begin{equation}\label{a-bdd}
\|\alpha Y_{(\lambda, w_\chi)}\|_{L_\chi}\leq C\|(\alpha,
0)\|_{\hat{W}_\chi}\leq C\lambda^{1/2-1/(2p)}|\hat{\varphi}|.
\end{equation}

For the last term in the error, we have:
\[
\begin{split}
&\|\hat{n}_{(\lambda, w_\chi)}(\alpha,
  \sum_{i=0}^{k+1}\varphi_i{\mathfrak e}_{u_i})\|_{L_\chi};\\
&\quad \leq
  C_2\sum_{i=0}^{k+1}\Big(\|n_{w_\chi}(\varphi_i\mathfrak{e}_{u_i})\|_{L_\chi}+\|\alpha\varphi_i\nabla_{{\mathfrak e}_{u_i}}Y_{(\lambda, w_\chi)}\|_{L_\chi}\Big)+C_2'\alpha^2\|\p_\lambda Y_{(\lambda, w_\chi)}\|_{L_\chi}+\text{\small higher order terms}\\
&\quad \leq
C_3\sum_{i=0}^{k+1}\Big(\lambda|\varphi_i|^2+|\alpha|\lambda^{-1/(2p)}|\varphi_i|\Big)+C_3'|\alpha|^2\lambda^{-1-1/(2p)}\\
&\quad \leq C_4\lambda|\hat{\varphi}|^2.\end{split}
\label{est:hatN}\]
For the first inequality above, we used the fact 
that for different \(i, j\), \({\mathfrak e}_{u_i}, {\mathfrak
  e}_{u_j}\) have disjoint supports.
For the second inequality, we used the invariance of the
flow equation under translation, which implies 
\begin{equation}\label{n-vanishing}
n_{u_i}^{J X_0}(\varphi _i(u_{i})_{\gamma})=0\quad \forall  \varphi_i\in\R.
\end{equation}

Next, 
the linear term in (\ref{equation:image}) is of the form \(E'_\chi\xi\), where \(E'_\chi\) is
\(E_{w_\chi}\)  perturbed by a term coming from 
$\hat{n}_{(\lambda, w_\chi)}(\alpha, \sum_{i=0}^{k+1}\varphi_i{\mathfrak
  e}_{u_i}+\xi)$, which has operator norm bounded by
\[C(\sum_{i=0}^{k+1}|\varphi_i|+\lambda^{-1/2}|\alpha|)\leq C|\hat{\varphi}|.\]
So with the assumption that
\(|\hat{\varphi}|\ll1\), \(E'_\chi\) is uniformly right invertible as
$E_{w_\chi}$ is. So by
contraction mapping theorem and the error estimates above we have an
$\xi(\hat{\varphi})$ satisfying (\ref{xi-West}).

The estimate for the nonlinear term is not very different from 
that in Lemma \ref{nonlinear:bound}, which we shall omit. 

Finally, to estimate \(\xi-\xi'\), where \(\xi:=\xi(\hat{\varphi}); \xi':=\xi
(\hat{\varphi}')\), notice that it satisfies 
\[\begin{split}
E_{w_\chi}(\xi-\xi')=-& P^c\Big(\sum_{i=0}^{k+1}(\varphi_i-\varphi'_i)E_{w_\chi}{\mathfrak
  e}_{u_i})+(\alpha-\alpha')Y_{(\lambda, w_\chi)}\\
&\quad +\hat{n}_{(\lambda, w_\chi)}(\alpha,
\xi+\sum_{i=0}^{k+1}\varphi_i{\mathfrak e}_{u_i})-\hat{n}_{(\lambda,
  w_\chi)}(\alpha', \xi'+\sum_{i=0}^{k+1}\varphi'_i{\mathfrak e}_{u_i})\Big).
\end{split}\]
Thus, by Proposition \ref{UI}, 
\[\begin{split}
\|\xi-\xi'\|_{W_\chi} &\leq
C'\Big(\sum_{i=0}^{k+1}|\varphi_i-\varphi'_i|\|E_{w_\chi}{\mathfrak
  e}_{u_i}\|_{L_\chi}+|\alpha-\alpha'|\|Y_{(\lambda,
  w_\chi)}\|_{L_\chi}\\
&\quad +\Big\|\hat{n}_{(\lambda, w_\chi)}(\alpha,
\xi+\sum_{i=0}^{k+1}\varphi_i{\mathfrak e}_{u_i})-\hat{n}_{(\lambda,
  w_\chi)}(\alpha, \xi'+\sum_{i=0}^{k+1}\varphi'_i{\mathfrak e}_{u_i})\Big\|_{L_\chi}\Big)
\end{split}\]
for a \(\lambda\)-independent constant \(C'\). The first two terms inside
the parenthesis may be bounded by
\(C_1\lambda^{1/2-1/(2p)}|\hat{\varphi}-\hat{\varphi}'|\) according to
(\ref{Eest:xi-u}) and (\ref{a-bdd}). The third term, by direct computation and
(\ref{n-vanishing}) again, may be bounded by
\begin{equation}\label{nonlinear}
\begin{split}
&C_2\Big(\|\xi\|_{W_\chi}+\|\xi'\|_{W_\chi}+\sum_{i=0}^{k+1}(|\varphi_i|+|\varphi'_i|)\|{\mathfrak
  e}_{u_i}\|_{W_\chi}+ (|\alpha|+|\alpha'|)\lambda^{-1/2}\Big)\|\xi-\xi'\|_{W_\chi}\\
 &\quad+\sum_{i=0}^{k+1}\Big(\|\xi\|_{W_\chi}+\|\xi'\|_{W_\chi}+\lambda(|\varphi_i|+|\varphi'_i|)\|{\mathfrak
  e}_{u_i}\|_{W_\chi}+ (|\alpha|+|\alpha'|)\lambda^{-1/2}\Big)\\
&\qquad \qquad \qquad\cdot C_3 \|{\mathfrak
  e}_{u_i}\|_{W_\chi}|\varphi_i-\varphi'_i|\\
&\quad +  \Big(\lambda^{-1/2}(\|\xi
\|_{W_\chi}+\|\xi'\|_{W_\chi})
+\lambda ^{-1-1/(2p)}(|\alpha|+|\alpha'|)+\sum_{i=0}^{k+1}(|\varphi_i|+|\varphi'_i|)\lambda^{-1/2}\|{\mathfrak
  e}_{u_i}\|_{W_\chi}\Big)\\
&\qquad\qquad \qquad \cdot C_4 |\alpha-\alpha'|\\
  &\leq C_2'\varepsilon\|\xi-\xi'\|_{W_\chi}+C_3'\lambda^{1/2-1/(2p)}|\hat{\varphi}-\hat{\varphi}'|,
\end{split}\end{equation}
where \(0<\varepsilon\ll1\), and we have used (\ref{xi-West}) and the
fact that \(|\hat{\varphi}|\ll1\)
above. Now, the first term in the last expression above can be got rid
of by a rearrangement argument, and we arrive at (\ref{delta-xi-West}).
\end{proof}

Next, substitute \(\xi(\hat{\varphi})\) back in (\ref{equation:cok})
to solve for \(\hat{\varphi}\). To understand the behavior of the
solutions, we estimate each term in the Kuranishi map in turn.
\begin{*lemma}[Terms in the Kuranishi map]\label{est:cok}
Let \(q^{-1}:=1-p^{-1}\). Then:
\begin{description}
\item[(a)] \(|P_{{\mathfrak f}_j}\bar{\p}_{J
    X_\lambda}w_\chi|\leq C'\lambda^{1/q}\) for any \(j\in \{1,
  \ldots, k+1\}\);

\item[(b)] For any \(i\in \{0, \ldots, k+1\}\), \(j\in \{1, \ldots,
  k+1\}\), 
\begin{gather*}
|P_{{\mathfrak f}_j}E_{w_\chi}{\mathfrak e}_{u_i}|\leq C\lambda^{1/q}\quad
\text{if \(j\neq i\)
    or \(i+1\)}; \\
-C'_-\lambda^{1/(2q)}\geq P_{{\mathfrak f}_{i}}(E_{w_\chi}{\mathfrak e}_{u_i})\geq
-C_-\lambda^{1/(2q)};\\
C'_+\lambda^{1/(2q)}\geq P_{{\mathfrak f}_{i+1}} (E_{w_\chi}{\mathfrak e}_{u_i})\geq
  C_+\lambda^{1/(2q)}.
\end{gather*}

\item[(c)] Let \(\hat{\varphi}, \hat{\varphi}'\),
  \(\xi:=\xi(\hat{\varphi}), \xi':=\xi(\hat{\varphi}')\) be as in the
  previous lemma. Then \(\forall j\), 
\[ \Big|P_{{\mathfrak f}_j}\Big(\hat{n}_{(\lambda, w_\chi)}(\alpha, \xi+\sum_{i=0}^{k+1}\varphi_i{\mathfrak
  e}_{u_i})-\hat{n}_{(\lambda, w_\chi)}(\alpha', \xi'+\sum_{i=0}^{k+1}\varphi_i'{\mathfrak
  e}_{u_i})\Big)\Big|\leq C_n\lambda^{1/q}|\hat{\varphi}-\hat{\varphi}'|.\]
\end{description}
\end{*lemma}
\begin{proof}
\underline{(a):} Apply (\ref{est:P-j}) with \(\eta=\bar{\p}_{JX_\lambda}w_\chi\).
The integrals in the first terms 
vanish because by our definition of pregluing, \(\bar{\p}_{J X_\lambda}w_\chi\)
has no longitudinal component in this region. On the other hand,
\(\|\eta\|_{L_\chi}\leq C\lambda^{1/(2q)}\) by Proposition \ref{err-est}.

\underline{(b):} Let \(\eta=E_{w_\chi}{\mathfrak e}_{u_i}\) in (\ref{est:P-j}). By
(\ref{Eest:xi-u}), the second term of (\ref{est:P-j}) (multiple of
\(\lambda^{1/(2q)}\|\eta\|_{L_\chi}\)) contributes a multiple of
\(\lambda^{1/q}\). 
If \(j\neq i, i+1\), the first term of (\ref{est:P-j}) (multiples of
integrals) vanishes, because
\(\eta\) is supported away from the interval of integration. 
These together imply the first line of (b).

For the other cases (\(u=u_j\) or \(u_{j-1}\)), we shall show that 
\[
\pm C'\leq  \int_{\gamma_{u_{j-1}}^{-1}(0)}^{\gamma_{u_{j}}^{-1}(0)}
\ell^{-1}\underline{E_{w_\chi}\mathfrak{e}_u}_L\leq \pm C\quad
\text{(\(-\) when \(u=u_j\), \(+\) when \(u=u_{j-1}\)).}
\]
This would imply that the first terms of (\ref{est:P-j}) is bounded
below and above by positive multiples of \(\pm \lambda^{1/(2q)}\),
dominating the second term. The other two cases of (b) would then follow.

To see this, recall the computation of \(E_{w_\chi}\mathfrak{e}_u\)
from (\ref{E:xi-u}), and note that the longitudinal component of the
second term vanishes because of I.5.3.1 (2c). Thus, 
\[
\int_{\gamma_{u_{j-1}}^{-1}(0)}^{\gamma_{u_{j}}^{-1}(0)}
\ell^{-1}\underline{E_{w_\chi}\mathfrak{e}_u}_L=\int_{\gamma_{u_{j-1}}^{-1}(0)}^{\gamma_{u_{j}}^{-1}(0)}
\ell^{-1}(\gamma_{u}'-1)\underline{u_{\gamma\gamma}(\gamma_u)}_L
\]
To estimate the integral on the RHS, note that on the interval of
integration, \(\gamma_{u_j}, \gamma_{u_{j-1}}\) are negative/positive
respectively. Also, 
\[
\begin{cases}
|\ell ^{-1}(\gamma_u'-1)\underline{u_{\gamma\gamma}}_L|\leq
C_1\lambda&\text{when \(|\gamma_u|\leq \gamma_0\)};\\
C_2'\lambda|\gamma_u|\leq|\ell ^{-1}(\gamma_u'-1)\underline{u_{\gamma\gamma}}_L|\leq
C_2\lambda|\gamma_u|&\text{when \(\gamma_0\leq|\gamma_u|\leq
  \epsilon^{1/2}\lambda^{-1/2}\)};\\
C_3'|\gamma_u|^{-1}\leq|\ell ^{-1}(\gamma_u'-1)\underline{u_{\gamma\gamma}}_L|\leq
C_3|\gamma_u|^{-1}&\text{when \(|\gamma_u|\geq \epsilon^{1/2}\lambda^{-1/2}\)};
\end{cases}\]
Furthermore, when \(|\gamma_u|\geq\gamma_0\), the sign of
\(\ell ^{-1}(\gamma_u'-1)\underline{u_{\gamma\gamma}}_L\) is the sign of
\(\gamma_u\) . Thus, the contribution to the first integral from
the two regions where \(|\gamma_u|\geq\gamma_0\) is bounded above and below
by expressions of the form 
\[
\sign (\gamma_u) \, C_2''\Big( \int_{\gamma_0}^{\epsilon^{1/2}\lambda^{-1/2}}\lambda\gamma(\gamma')^{-1}\,
d\gamma\,
+\int_{\epsilon^{1/2}\lambda^{-1/2}}^\infty\gamma^{-1}(\gamma')^{-1}\,
d\gamma\Big).
\]
By our estimate for \(\gamma'\) in section 2,
this is in turn bounded above and below by \(\sign (\gamma_u) \,  C\), where \(C>0\)
is a \(\lambda\)-independent constant, and the sign is \(-\) when \(u=u_j\); \(+\) when
\(u=u_{j-1}\). Meanwhile, the contribution from the region where
\(|\gamma_u|\leq\gamma_0\) is bounded by \(C'\lambda\); therefore
ignorable. In sum, we have the claimed estimate, and hence 
the assertion (b).

\underline{(c):} Let \(\eta=\hat{n}_{(\lambda, w_\chi)}(\alpha, \xi+\sum_{i=0}^{k+1}\varphi_i{\mathfrak
  e}_{u_i})-\hat{n}_{(\lambda, w_\chi)}(\alpha', \xi'+\sum_{i=0}^{k+1}\varphi_i'{\mathfrak
  e}_{u_i})\) in (\ref{est:P-j}). The second term in it, 
by (\ref{nonlinear}) and (\ref{delta-xi-West}), are bounded
by \(C\lambda^{1/(2q)}\lambda^{1/(2q)}|\hat{\varphi}-\hat{\varphi}'|\). On
the other hand, by H\"{o}lder inequality and the same direct computation
that appeared in the end of the proof of last lemma, 
the first term of (\ref{est:P-j}) can be bounded in absolute value by
\[\begin{split}
&C\lambda^{1/(2q)}\Big( \sum_i|\varphi_i-\varphi'_i|\Big(\lambda(|\varphi_i|+|\varphi'_i|)\int_{\gamma_{u_{j-1}}^{-1}(0)}^{\gamma_{u_{j}}^{-1}(0)}\ell
^{-1}\|(u_i)_\gamma\|_{2,t}^2ds\,\\
&\qquad\qquad\qquad\qquad
+(|\alpha|+|\alpha'|)\int_{\gamma_{u_{j-1}}^{-1}(0)}^{\gamma_{u_{j}}^{-1}(0)}\ell^{-1}
\underline{(u_i)_\gamma}_Lds\Big)\\
&\quad+(\|\xi\|_{W_\chi}+\|\xi'\|_{W_\chi})\Big(\|\xi-\xi'\|_{W_\chi}\lambda^{-1/2+1/p}+\sum_i|\varphi_i-\varphi'_i|\|\ell^{-1/2}(u_i)_\gamma\|_{L^q((\gamma_{u_{j-1}}^{-1}(0),
  \gamma_{u_{j}}^{-1}(0))\times S^1)}\\
&\qquad\qquad\qquad\qquad + |\alpha-\alpha'|\,
\|\ell^{-1/2}\|_{L^q((\gamma_{u_{j-1}}^{-1}(0),
  \gamma_{u_{j}}^{-1}(0)))}
 \Big)\\
&\quad +|\alpha-\alpha'|\Big((|\alpha|+|\alpha'|)\int_{\gamma_{u_{j-1}}^{-1}(0)}^{\gamma_{u_{j}}^{-1}(0)}\ell ^{-1}\, ds\, +\sum_i(|\varphi_i|+|\varphi'_i|)\int_{\gamma_{u_{j-1}}^{-1}(0)}^{\gamma_{u_{j}}^{-1}(0)}\ell^{-1}
\underline{(u_i)_\gamma}_Lds\Big)\\
&\quad
+\|\xi-\xi'\|_{W_\chi}\Big(\sum_i(|\varphi_i|+|\varphi'_i|)\|\ell ^{-1/2}(u_i)_\gamma\|_{L^q((\gamma_{u_{j-1}}^{-1}(0),
  \gamma_{u_{j}}^{-1}(0))\times S^1)}\\
&\qquad\qquad\qquad\qquad +(|\alpha|+|\alpha'|)\|\ell^{-1/2}\|_{L^q((
    \gamma_{u_{j-1}}^{-1}(0),
    \gamma_{u_{j}}^{-1}(0))}\Big) \Big)\\
&\leq
C'\lambda^{1/(2q)}\Big(\lambda^{1/(2q)}|\hat{\varphi}-\hat{\varphi}'|+(\varepsilon_1+C''\lambda^{1/(2p)})\|\xi-\xi'\|_{W_\chi}\Big)\\
&\leq C_3\lambda^{1/q}|\hat{\varphi}-\hat{\varphi}'|.
\end{split}
\]
In the above we again used (\ref{delta-xi-West}), the estimate for
\(\ell\) in Lemma \ref{h0:ext}, and the estimates for
\(\gamma_{u_{i}}\) and \(\sigma_\chi\) in section 2.
Summing up, this gives us assertion (c).
\end{proof}

\subsubsection{Constructing the gluing map.}\label{gluing-map:2.1}
It follows immediately from the previous Lemma that the linearization
of the Kuranishi map is surjective, and hence the moduli space is
(Zariski) smooth. More precisely, choose 
\[Q_\chi:=\op{Span}\{{\frak e}_{u_1},
\ldots, {\frak e}_{u_{k+1}}\}\subset K_\chi \quad Q^\Xi:=\bigcup_\chi Q_\chi.\]
The reductions of the K-models 
\([K_\chi\to C_\chi]_{W'_\chi}\), \([\hat{K}_\chi\to
C_\chi]_{W'_\chi}\) by \(Q_\chi\) give respectively the 
standard K-models for \(E_{w_\chi}\) and \(\hat{E}_{(\lambda,
  w_\chi)}\): 
\[
[\ker E_{w_\chi} \to *]_{Q_\chi\oplus W'_\chi},\quad  [\ker \hat{E}_{(\lambda,
  w_\chi)}\to *]_{Q_\chi\oplus W'_\chi}.
\]

Indeed, from Lemma \ref{est:cok} (b) we see that the \((k+1)\times(k+1)\)-matrix \(E=(E_{ji})\),
\[
E_{ji}:=\lambda^{-1/(2q)}P_{{\mathfrak f}_{j}}(E_{w_\chi}{\mathfrak
  e}_{u_i}),\quad i, j \in \{1,\ldots, k+1\} 
\]
is, up to ignorable terms, of the form
\[
\left(\begin{array}{ccccc}
 - &0&\cdots& \cdots&0\\
+ &-& 0&\cdots &0\\
 0 &+ &\ddots&\cdots & 0 \\
\vdots &0 &\ddots & \ddots&\vdots \\
 0 &\cdots &\cdots& + &-
\end{array}\right)\quad \text{\small (\(+/-\) denote positive/negative
numbers of \(O(1)\).)} 
\] 
Thus, it 
has a uniformly bounded inverse, denoted \((G_{ij})\). Restricted to
\(Q_\chi\) (i.e. setting \(\alpha=\varphi_0=0\)), (\ref{equation:cok}) can be rewritten in the form 
\begin{equation}\label{f.d.contraction}
\vec{\varphi}=\Psi(\vec{\varphi}),\quad\text{where \(\vec{\varphi}:=(\varphi_1, \ldots, \varphi_{k+1})\),}
\end{equation}
and \(\Psi: \R^{k+1}\to \R^{k+1}\) is the map given by
\[
(\Psi(\vec{\varphi}))_i=-\sum_{j=1}^{k+1}G_{ij}\lambda^{-1/(2q)}P_{{\mathfrak f}_{j}}\Big(\bar{\p}_{JX_\lambda}w_\chi+n_{w_\chi}\Big(\xi(\vec{\varphi})+\sum_{l=1}^{k+1}\varphi_l{\mathfrak
  e}_{u_l}\Big)\Big).
\]
Note from the uniform boundedness of \((G_{ij})\) and Lemmas
\ref{est:im}, \ref{est:cok} (a) that 
\begin{equation}\label{Psi-0}
\begin{split}
|\Psi(\vec{0})|&\leq
\sum_jC\lambda^{-1/(2q)}\Big(|P_{{\mathfrak
    f}_j}\bar{\p}_{JX_\lambda}w_\chi|+\|n_{w_\chi}(\xi(\vec{0}))\|_{L_\chi}\Big)\\
&\leq
C_2\lambda^{-1/(2q)}(\lambda^{1/q}+\|\xi(\vec{0})\|_{W_\chi}^2)\\
&\leq
C_2'\lambda^{1/(2q)}\ll1.\end{split}\end{equation}
On the other hand, Lemmas \ref{est:cok} (c), \ref{est:im} and the uniform boundedness
of \((G_{ij})\) again imply:
\[
|\Psi(\vec{\varphi})-\Psi(\vec{\varphi'})|\leq
K|\vec{\varphi}-\vec{\varphi'}| \quad \text{for a positive constant
  \(K\leq C\lambda^{1/(2q)}\ll1\).}
\]
Thus by the contraction mapping theorem, we have a unique solution of
(\ref{f.d.contraction}) among all small enough \(\vec{\varphi}\). 

To summarize, for sufficiently small \(\lambda_0>0\), there is a
universal positive constant \(C_w\), such that for
all \(\chi\in \Xi(\mathbb{S})\), there is a unique 
\[(\vec{\varphi}_\chi, \xi_\chi)\in Q_\chi\oplus W_\chi'\quad \text
{with \(\|\xi_\chi\|_{W_\chi}^2+|\vec{\varphi}_\chi|^2\leq
C_w\),}\]
 which solves (\ref{equation:image}), (\ref{equation:cok}). In
fact, the solution satisfies 
\[\|\xi_\chi+\sum_i\varphi_{\chi, i}\mathfrak{e}_{u_i}\|_{W_\chi}^2
\leq \|\xi_\chi\|_{W_\chi}^2+C'|\vec{\varphi}_\chi|^2\leq C\lambda^{1/q},\]
because of (\ref{Psi-0}) and (\ref{xi-West}).

We define the gluing map to be the map from \(\Xi(\mathbb{S})\) to
\(\hat{\cm}^\Lambda_P\) sending 
\[
\chi\mapsto \exp \Big(w_\chi, \xi_\chi+\sum_{i=1}^{k+1} \varphi_{\chi, i}\mathfrak{e}_{u_i}\Big).
\]

\subsection{Surjectivity of the Gluing Map.}
With the gluing map constructed above, 
the standard arguments outlined in \S\ref{step4} shows that it is a
local diffeomorphism onto \(\hat{\cm}_P^\Lambda(\mathbf{x},
\mathbf{z})\cap U_\varepsilon\), the intersection of the moduli space with 
a tubular neighborhood \(U_\varepsilon \) of the image of the pregluing map, in {\em
  \({\cal B}\)-topology}. 

As our goal is instead to show that the gluing map is a local
diffeomorphism onto (the interior of) a neighborhood of 
\(\mathbb{S}\subset \hat{\cm}_P^{\Lambda, +}(\mathbf{x},\mathbf{z})\)
in the coarser {\em chain topology}, 
we need to show that the latter neighborhood in fact lies in
\(U_\varepsilon\). This is done via the following variant of decay
estimates for flows near \(y\).

\mysubsubsection[\(w\) in terms of \(\tilde{w}\) and \(\tilde{\xi}\).]
Let \((\lambda, \hat{w})\in \hat{\cm}_P^{(-\lambda_0,\lambda_0),
  1}(\mathbf{x},\mathbf{z}; \op{wt}_{-\langle {\cal Y}\rangle, e_{\cal P}}\leq\Re)\) be
in a chain topology neighborhood of 
\(\mathbb{S}\subset \hat{\cm}_P^{\Lambda, +}(\mathbf{x},\mathbf{z})\).
Namely, \(|\lambda|\ll1\), and there is a broken trajectory \(\hat{\chi}:=\{u_0,
u_1, \ldots, u_{k+1}\}\in \mathbb{S}\), which is close to \(\hat{w}\) in chain
topology. Let \(\chi:=(\hat{\chi}, \lambda)\in \Xi(\mathbb{S})\). We
may find a representative \(w\) of \(\hat{w}\), such that:
\[
w(s)=\exp (\tilde{w}(s), \tilde{\xi}(s)),
\] 
where \(w\), \(\tilde{w}\) are chosen such that:
\begin{itemize}
\item\(\tilde{w}(s)=y\) at \(s=\tilde{\mathfrak l}_1, \ldots,
\tilde{\mathfrak l}_{k+1}\); \(\tilde{\mathfrak l}_1<\tilde{\mathfrak l}_2<\cdots<
\tilde{\mathfrak l}_{k+1}\) subdivide \(\R\) into \(k+2\) open
intervals \(I_i, \, i=0, 1, \ldots, k+1\);
\item \(\tilde{w}(s)=u_i(\tilde{\gamma}_{u_i}(s))\) over \(I_i\),
  where \(\tilde{\gamma}_{u_i}: I_i\to\R\) are homeomorphisms
  determined by:
\begin{equation}\label{calibration}
\begin{cases}
\Pi_{{\bf e}_y}\tilde{\zeta}(s)=\Pi_{{\bf e}_y}\zeta(s) &\text{
for \(s\in \bigcup_i[\tilde{\gamma}_{u_{i-1}}^{-1}(\gamma_0),
\tilde{\gamma}_{u_i}^{-1}(-\gamma_0)]\), }\\
\tilde{\gamma}'_{u_i}(s)=1 &\text{for \(s\in \R\backslash  \bigcup_i[\tilde{\gamma}_{u_{i-1}}^{-1}(\gamma_0),
\tilde{\gamma}_{u_i}^{-1}(-\gamma_0)]\),}
\end{cases}\end{equation}
with \(\zeta\), \(\tilde{\zeta}\) given by \(\tilde{w}(s)=\exp(y,
\tilde{\zeta}(s))\), \(w(s)=\exp(y, \zeta(s))\),
and \(\gamma_0\) the large positive constant in \S2.4.
\item \begin{equation}\label{after-transl}
\gamma_{u_0}^{-1}(0)=\tilde{\gamma}_{u_0}^{-1}(0); \quad
\langle(u_0)_\gamma(0),
\tilde{\gamma}_{u_0}^*\tilde{\xi}(0)\rangle_{2,t}=0.
\end{equation}
\end{itemize}
Because of elliptic regularity and the fact that
\((\lambda, w)\) is close to \(\mathbb{S}\) in chain topology, we may
assume without loss of generality that 
\[\begin{split}
\|\tilde{\xi}(s)\|_{\infty,1,t}& +\|\tilde{\xi}'(s)\|_{\infty,1,t}<
  \varepsilon\,\,  \forall s \quad \text{for \(|\lambda|^{1/2}<\varepsilon\ll1\)};\\
\tilde{\gamma}_{u_j}^{-1}(0)
&-\tilde{\gamma}_{u_{j-1}}^{-1}(0)>\varepsilon^{-1}\quad \forall j\in
\{1, \ldots, k+1\}.
\end{split}\]

\mysubsubsection[Estimating \(\tilde{\xi}\).] Because of the large
weights near \(y\), we need the following more refined pointwise
estimate for \(\tilde{\xi}\) near \(y\).
\begin{lemma*}
Let \((\lambda, w)\in \cm_P^\Lambda\) be close to \(\mathbb{S}\) in chain topology, and
let \(y\) be a death as before. Then \(\lambda>0\).
Furthermore, in the notation of \S 4.2.1, there is a small positive
constant \(\varepsilon_0=\varepsilon_0(\gamma_0^{-1}, \varepsilon)\)
independent of \(w\), such that
\begin{equation}\label{calibrated}
\|\tilde{\xi}(s)\|_{2,2,t}+\|\tilde{\xi}'(s)\|_{2,1,t}\leq
\varepsilon_0(\|\Pi_{{\bf
    e}_y}\tilde{\zeta}(s)\|_{2,1,t}^4+|\lambda|)\quad \forall s\in
\bigcup_j\, [\tilde{\gamma}_{u_{j-1}}^{-1}(\gamma_0),
\tilde{\gamma}_{u_j}^{-1}(-\gamma_0)].
\end{equation}
\end{lemma*}
\begin{proof}
Let \(s\in \bigcup_j[\tilde{\gamma}_{u_{j-1}}^{-1}(\gamma_0),
\tilde{\gamma}_{u_j}^{-1}(-\gamma_0)]\) throughout this proof. In
fact, it suffices to consider only one \(j\).

To estimate \(\xi\), it is equivalent to estimate \(\zeta-\tilde{\zeta}\), which we denote
by \(c\). The assumption (\ref{calibration}) implies that \(c\in \ker
A_y^\perp\).

Write \(b:=\Pi_{\ker A_y}\tilde{\zeta}\), and let \(Z:\ker
A_y\to \ker A_y^{\perp}\) be such that
\(\tilde{\zeta}=(1+Z)b\). Similar to I.(38), I.(39), the flow equation
can be re-written as:
\begin{gather}
\label{B:flow}
-\frac{d\tilde{\zeta}}{ds}=(1+\nabla_{b}Z)(\lambda C'_y{\bf e}_y+\Pi_{\ker A_y} \hat{n}_{(0,y)}(\lambda,\tilde{\zeta}+c));\\
\label{C:flow}
-\frac{dc}{ds}=A_yc+(1-\Pi_{\ker A_y}-\nabla_{b}Z\Pi_{\ker A_y})(\hat{n}_{(0,y)}(\lambda,\tilde{\zeta}+c)-n_y(\tilde{\zeta}))-\nabla_{b}Z(\lambda C'_y{\bf e}_y).
\end{gather}
Taking the \(L^2_t\)-inner product of (\ref{C:flow}) with \(c\) and
rearranging like the proof of Sublemma I.5.1.7, we get (adopting
the notation of I.5.1)
\begin{gather*}
\frac{d\|c_+\|_{2,t}}{ds}\geq\mu_+\|c_+\|_{2,t}-\epsilon_+\|c\|_{2,t}-C_+|\lambda|\|\tilde{\zeta}\|_{2,1,t},\\
\frac{d\|c_-\|_{2,t}}{ds}\leq-\mu_-\|c_-\|_{2,t}+\epsilon_-\|c\|_{2,t}+C_-|\lambda|\|\tilde{\zeta}\|_{2,1,t}.
\end{gather*}
Subtracting a suitable multiple of the first inequality from the
second, one obtains:
\[
(\|c_-\|_{2,t}-\epsilon'_-\|c_+\|_{2,t})'\leq -\mu_-'\|c_-\|_{2,t}+C'_-\gamma_0^{-1}|\lambda|;
\]
Taking convolution product with the integral kernel of \(d/ds
+\mu_-'\) on both sides, one gets
\[
\|c_-\|_{2,t}-\epsilon'_-\|c_+\|_{2,t}\leq C_-\varepsilon e^{-\mu'(s-\tilde{\gamma}^{-1}_{u_{j-1}}(\gamma_0))}+C''_-\gamma_0^{-1}|\lambda|,
\] 
and similarly,
\[
\|c_+\|_{2,t}-\epsilon'_+\|c_-\|_{2,t}\leq C_+\varepsilon e^{-\mu'(\tilde{\gamma}^{-1}_{u_j}(-\gamma_0)-s)}+C''_+\gamma_0^{-1}|\lambda|.
\]
Adding the above two inequalities, we get
\[
\|c\|_{2,t}\leq C\varepsilon\Big(e^{-\mu'(s-\tilde{\gamma}^{-1}_{u_{j-1}}(\gamma_0))}+e^{-\mu'(\tilde{\gamma}^{-1}_{u_j}(-\gamma_0)-s)}\Big)+C''\gamma_0^{-1}|\lambda|.
\]
This may be improved to give a similar estimate for \(\|c\|_{2,1,t}\) using
(\ref{C:flow}) by the same elliptic bootstrapping and Sobolev embedding
argument as in I.5.1.7.

On the other hand, write \(b(s)=\underline{b}(s) {\bf e}_y\) as usual,
and notice that by taking \(\Pi_{\ker A_y}\) of (\ref{B:flow}), \(b(s)\) satisfies:
\begin{equation}\label{flow:b}
-b'(s)=\lambda C'_y{\bf e}_y+\Pi_{\ker A_y} \hat{n}_{(0,y)}(\lambda,\tilde{\zeta}+c).
\end{equation}
Integrating this equation, it is easy to see that \(\lambda<0\) would
contradict the fact that, due to the proximity of \(w\) and
\(\tilde{w}\), \[
\underline{b}(\tilde{\gamma}_{u_{j-1}}^{-1}(\gamma_0))>0, \, \, 
\underline{b}(\tilde{\gamma}_{u_{j}}^{-1}(-\gamma_0))<0, \quad 
\tilde{\gamma}_{u_{j-1}}^{-1}(\gamma_0)< \tilde{\gamma}_{u_{j}}^{-1}(-\gamma_0).\]
Thus, \(\lambda\) must be positive. 

On the other hand, as \(\lambda>0\), (\ref{flow:b}) implies that 
\(\underline{b}(s)\) decreases monotonically with \(s\). We now claim that
\[e^{-\mu'(s-\tilde{\gamma}^{-1}_{u_{j-1}}(\gamma_0))}+e^{-\mu'(\tilde{\gamma}^{-1}_{u_j}(-\gamma_0)-s)}\leq
C_e(\underline{b}(s)^4+|\lambda|).\] Combined with the above estimates for \(c\), this would then imply the
second assertion of the Lemma. 

To prove the claim, note that by symmetry and the decay/growth
behavior of the two terms on the LHS, it suffices to show that 
\begin{gather*}
\underline{b}^4(s)\geq
C_1e^{-\mu'(\tilde{\gamma}^{-1}_{u_j}(-\gamma_0)-s)}\quad \text{when \(-1\ll\underline{b}(s)\leq-|\lambda|^{1/2}\);}\\
\underline{b}^4(s)\geq
C_2e^{-\mu'(s-\tilde{\gamma}^{-1}_{u_{j-1}}(\gamma_0))}\quad \text{when \(1\gg\underline{b}(s)\geq|\lambda|^{1/2}\)}
\end{gather*} 
for \(s\)-independent constants \(C_1, C_2\). We shall only
demonstrate the second inequality since the first is similar. 
When \(\underline{b}(s)\geq |\lambda|^{1/2}\), (\ref{flow:b}) together
with the above estimate for \(\|c\|_{2,1,t}\) imply that 
\[
(\underline{b}^4)'\geq -\mu'\underline{b}^4-C_b\varepsilon^4e^{-4\mu'(s-\tilde{\gamma}^{-1}_{u_{j-1}}(\gamma_0))}.
\]
Taking convolution product with the integral kernel of \(d/ds+\mu'\), we get in this region 
\[
\underline{b}^4(s)\geq C_6e^{-\mu'(s-\tilde{\gamma}^{-1}_{u_{j-1}}(\gamma_0))}-C_7\varepsilon^4e^{-4\mu'(s-\tilde{\gamma}^{-1}_{u_{j-1}}(\gamma_0))},
\]
and hence the claim.
\end{proof}

\mysubsubsection[From \(\tilde{w}\) to \(w_\chi\).] Next, notice that 
\(\tilde{w}\) differs from the pregluing \(w_\chi\) by a
reparametrization. We shall estimate the difference between
\(\tilde{w}\) and \(w_\chi\) by estimating the difference between
\(\tilde{\gamma}_{u_i}^{-1}\) and \(\gamma_{u_i}^{-1}\).

Similar to \(\gamma_{u_i}\) (see \S2.2.1), \(\tilde{\gamma}_{u_i}\) satisfies:
\[
\Big\langle u_{\gamma}(\tilde{\gamma}_{u_i}),
-(\tilde{\gamma}_{u_i}'-1)u_{\gamma}(\tilde{\gamma}_{u_i})+\lambda
Y_{(0, \tilde{w})}+E_{\tilde{w}}\tilde{\xi}+\tilde{n}(\lambda,
\tilde{\xi})\Big\rangle_{2,t}=0,
\]
where \(\tilde{n}\) is some nonlinear term in \(\lambda, \tilde{\xi}\).
By (\ref{calibration}), when \(\tilde{\gamma}_{u_i}(s)\geq |\gamma_0|\),
\[
\|u_{\gamma}(\tilde{\gamma}_{u_i})\|_{2,t}^{-1} |\langle
u_{\gamma}(\tilde{\gamma}_{u_i}), E_{\tilde{w}}\tilde{\xi}+\tilde{n}(\lambda,
\tilde{\xi})\rangle_{2,t}|
\leq
C_8(\|\tilde{\xi}\|_{2,1,t}+\|\tilde{\xi}'\|_{2,t})/\tilde{\gamma}_{u_i}+C'_8(\|\tilde{\xi}\|_{2,1,t}+\lambda)^2.
\]

Write \(\gamma_{\lambda, u_i}=\gamma_{u_i}\) in (\ref{equation:pregluing}) to
emphasize the parameter \(\lambda\) used in the definition, and let
\[\Delta_{s, i}(\gamma):=\gamma^{-1}_{\lambda, u_i}(\gamma)-\tilde{\gamma}^{-1}_{u_i}(\gamma).\]
Comparing the defining equations for \(\gamma_{\lambda , u_i}\) and
\(\tilde{\gamma}_{\lambda, u_i}\), and using (\ref{calibrated}) and the above estimate for \(\|u_{\gamma}(\tilde{\gamma}_{u_i})\|_{2,t}^{-1} |\langle
u_{\gamma}(\tilde{\gamma}_{u_i}), E_{\tilde{w}}\tilde{\xi}+\tilde{n}(\lambda,
\tilde{\xi})\rangle_{2,t}|\),
we find that for all \(i\in \{0, \ldots, k\}\), 
\begin{equation}\label{Delta-s}\begin{split}
|\Delta_{s, i} &(\gamma_+) -\Delta_{s, i}(\gamma_-)|\leq\\
&\begin{cases}
\Big|\int_0^{\gamma_0}(O(\varepsilon)+O(\lambda))\, d\gamma\, \Big|\leq
C_1\varepsilon &\text{when \(0\leq\gamma_+, \gamma_-\leq\gamma_0\)};\\
\Big|\int_{\gamma_0}^{\mathfrak{r}_i}C_7(\varepsilon\gamma (\lambda+\gamma^{-4})+(\lambda
+\gamma^{-4})^2\gamma^2)\,d\gamma\Big|\leq C_6\varepsilon &\text{when \(\gamma_0<\gamma_+, \gamma_-\leq\mathfrak{r}_i\)};\\
\Big|\int_{\mathfrak{r}_i}^{\infty}\frac{C_5(\varepsilon\gamma(\lambda+\gamma^{-4})+(\lambda
+\gamma^{-4})^2\gamma^2)\,
  d\gamma}{(1+\lambda\gamma^2)^2}\Big|\leq C_4\varepsilon
&\text{when \(\gamma_+, \gamma_-\geq \mathfrak{r}_i\)}.
\end{cases}\end{split}
\end{equation}
Similar estimates hold for negative \(\gamma_+, \gamma_-\) when \(i\in
\{1, \ldots, k+1\}\). For \(i=0\) and any two negative
\(\gamma_+\), \(\gamma_-\), or for \(i=k+1\) and any two positive
\(\gamma_+\), \(\gamma_-\), the estimate in the first case above
holds.

Combining this with the initial conditions from (\ref{after-transl}):
\[
\Delta_{s, 0}(0)=0; \quad \Delta_{s, i-1}(\infty)=\Delta_{s,
  i}(-\infty) \, \, \forall  i\in \{1, \ldots, k+1\},\]
we have
\[
|\Delta_{s, i} (\gamma)|\leq C\varepsilon \, \, \forall \gamma, i\quad \text{for a
  \(\lambda\)-independent constant \(C>0\).}
\]
Applying the the mean value theorem and
the estimate for $w'_{\chi}$ in (\ref{sigma}), and recalling
the assumption (\ref{after-transl}), we see that
\(\tilde{w}=\exp (w_\chi, \tilde{\xi}_\chi)\) for an
\(\tilde{\xi}_\chi\) with:
\begin{equation}\label{t-xi-chi}
\begin{split}
& \langle  (u_0)_\gamma, (0), (\gamma_{u_0}^{-1})^*\tilde{\xi}_\chi(0)\rangle_{2,t}=0;\\
\|\tilde{\xi}_\chi\|_{2,1, t}\leq &
C'(\lambda+\varepsilon\gamma^{-2}_{u_i}) \quad \text{on \(I_i\cap \bigcup_j[\gamma_{u_{j-1}}^{-1}(\gamma_0),
\gamma_{u_j}^{-1}(-\gamma_0)]\,
  \forall i\)}
\end{split}\end{equation}

\mysubsubsection[From pointwise estimates to \(\hat{W}_\chi\) estimates.]
Recall that our goal is to show that given \((\lambda, \hat{w})\in
\hat{\cm}^\Lambda_P\) in a chain topology neighborhood of
\(\mathbb{S}\),  as prescribed in \S4.2.1, we may write
\begin{equation}\label{w-w-chi}
(\lambda,w)=e(\lambda',w_{\chi'};\alpha_{\chi'},\xi_{\chi'})\quad
\text{for some
  \(\chi'=(\underline{\chi}, \lambda')\in \Xi(\mathbb{S})\),
  \(\underline{\chi}:=\{\hat{u}_0, \ldots, \hat{u}_{k+1}\}\),}
\end{equation}
with \((\alpha_{\chi'},\xi_{\chi'})\) satisfying
\begin{equation}\label{surj-est}
(1) \, \|(\alpha_{\chi'},\xi_{\chi'})\|_{\hat{W}_{\chi'}}\leq
C\varepsilon;\quad
(2)\, (\alpha_{\chi'},\xi_{\chi'})\in B_{\chi'},
\end{equation}
where \(C\) is a \(\lambda\)-independent constant, and \(B_{\chi'}\)
is the B-space chosen so that \([\ker \hat{E}_{(\lambda',
  w_{\chi'})}\to *]_{B_{\chi'}}\) forms a K-model for \(\hat{E}_{(\lambda',
  w_{\chi'})}\).
\begin{lemma*}
Suppose the \((\alpha_{\chi'}, \xi_{\chi'})\) given in (\ref{w-w-chi}) satisfies
\begin{equation}\label{pt-wise}
\begin{split}
|\alpha_{\chi'}|& \leq C'\lambda^{3/2}\varepsilon;\\
\|\xi_{\chi'}\|_{2,1, t} &\leq
C'_\xi(\lambda+\varepsilon\gamma^{-2}_{u_i}) \quad \text{on \(I_i\cap
  \bigcup_j\, [\gamma_{u_{j-1}}^{-1}(\gamma_0),
\gamma_{u_j}^{-1}(-\gamma_0)]\,\,
  \forall i\)}
\end{split}\end{equation}
for \(\lambda\)-independent constants \(C', C'_\xi\). Then
(\ref{surj-est}.1) holds. 
\end{lemma*}
\begin{proof}
First, notice that the assumption on \(\alpha_{\chi'}\) implies 
\(\|(\alpha_{\chi'}, 0)\|_{\hat{W}_{\chi'}}\leq C_1\lambda^{1/2-1/(2p)}\varepsilon\).
On the other hand, the assumption that \((\lambda, \hat{w})\) is close
to \(\chi'\) in chain topology implies that 
over \(\Theta^c:=\Theta\backslash\bigcup_j[\gamma_{u_{j-1}}^{-1}(\gamma_0),
\gamma_{u_j}^{-1}(-\gamma_0)]\times S^1\),
\[
\|\xi_{\chi'}\|_{W_{\chi'}(\Theta^c)}\leq
C_2\|\xi_{\chi'}\|_{L^p_1(\Theta^c)}\leq C_3\varepsilon \quad \text{for \(\lambda\)-independent constants \(C_2, C_3\).
}\] 
Thus, it remains to estimate \(\|\xi_{\chi'}\|_{W_{\chi'}([\gamma_{u_{j-1}}^{-1}(\gamma_0),
\gamma_{u_j}^{-1}(-\gamma_0)]\times S^1)}\). We shall focus on
estimates on the region \([\gamma_{u_{i}}^{-1}(\gamma_0),
\gamma_{u_{i}}^{-1}(\infty))\times S^1\) for an \(i\in \{0, \ldots,
k\}\), since estimates on the rest are similar.

By the definition of $\gamma_{u_i}$, on this region the flow
equation has the form: 
\begin{equation}\label{pcr-uni}
(\bar{\p}_{J X_{\lambda'}} w_{\chi'})_{T_w} +\hat{E}_{(\lambda', w_{\chi'})} (\alpha_{\chi'},\xi_{\chi'})
  +\hat{n}_{(\lambda', w_{\chi'})}(\alpha_{\chi'},\xi_{\chi'})=0.
\end{equation}
($T_w$ here means the transverse component with respect to
$w_{\chi'}'$, in contrast to \(T_y\) below).

Subdivide the region again into \([\gamma_{u_i}^{-1}(\mathfrak{r}_i),
\gamma_{u_i}^{-1}(\infty))\times S^1\) and the rest.

{\em Over the first region}, namely when
\(\gamma_{u_i}\geq\mathfrak{r}_i\), in place of \(\xi_{\chi'}\) and its \(W_{\chi'}\)-norm,
it is equivalent to estimate
\begin{gather*}
\xi_0:=T_{w_{\chi'}, \bar{y}}\xi_{\chi'}\in \Gamma(\bar{y}^*K)\quad
\text{in the norm}\\
\|\xi_0\|_{W_0^{\lambda'}}:=(\lambda')^{-1/2}\|\xi_0\|_{p,1}+(\lambda')^{-1}\|(\xi_{0})'_{L_y}\|_p,
\end{gather*}
where \((\xi_{0})_{L_y}(s)=\langle{\bf e}_y, \xi_0(s)\rangle_{2,t}{\bf
  e}_y\) is the `longitudinal direction with respect to \(y\)'. Notice
that \({\bf e}_y\) differs from the original
longitudinal direction
\(T_{w_{\chi'}(s),y}w'_{\chi'}(s)\|w'_{\chi'}(s)\|_{2,t}^{-1}\) by an
ignorable factor of \(C(\lambda'/\epsilon)^{1/2}\). Let the
transversal direction \(T_y\) and the \(L_0^{\lambda'}\)-norm 
be similarly defined. 

Rewriting the flow equation in terms of the above
transverse and longitudinal directions, we have:
\begin{equation}\label{split:pcr}
\begin{split}
E_{\bar{y}}(\xi_{0T_y})&=-\alpha _{\chi'}(T_{w_{\chi'},
  \bar{y}}Y_{(\lambda, w_{\chi'})})_{T_y}-(T_{w_{\chi'}, \bar{y}}\hat{n}_{(\lambda',w_{\chi'})}(\alpha_{\chi'},\xi_{\chi'}))_{T_y}+Z_{T_y}+\Upsilon_{T_y}; \\
E_{\bar{y}}(\xi_{0L_y})&=-\alpha _{\chi'}(T_{w_{\chi},
  \bar{y}}Y_{(\lambda, w_{\chi'})})_{L_y}-(T_{w_{\chi'},
  \bar{y}}\hat{n}_{(\lambda',w_{\chi'})}(\alpha_{\chi'},\xi_{\chi'}))_{L_y}+Z_{L_y}+\Upsilon_{L_y}, 
\end{split}\end{equation}
where:
\begin{itemize}
\item $Z_{T_y}, Z_{L_y}$ come from the difference between
$E_{\bar{y}}$ and $E_{w_{\chi'}}$. Thus,
their $L^2_{1,t}$-norms are bounded by 
$C\|\xi_{\chi'}\|_{\infty,1,t} \gamma^{-1}_{u_i}$; 
\item $\Upsilon_T$, $\Upsilon_L$ are
terms coming from $(\bar{\p}_{J X_{\lambda'}} w_{\chi'})_{T_w}$. 
The computation in \S\ref{lem:err-ptws} shows that 
$\|\Upsilon_{T_y}\|_{2,1,t}$ is bounded by $C_1\lambda' \gamma^{-1}_{u_i}$,
while $\|\Upsilon_{L_y}\|_{2,1,t}$ is bounded by
$C_2\lambda'\gamma_{u_i}^{-2}$. 
\end{itemize}

Now, the length estimate for \([\gamma_{u_i}^{-1}(\mathfrak{r}_i),
\gamma_{u_i}^{-1}(\infty))\times S^1\) (cf. (\ref{Theta0:length})) and
the assumption (\ref{pt-wise}) yield
\[
(\lambda')^{-1/2} (\|\xi_{0}\|_{L^p([\gamma_{u_i}^{-1}(\mathfrak{r}_i),
\gamma_{u_i}^{-1}(\infty))\times S^1)}+\|\dot{\xi}_{0}\|_{L^p([\gamma_{u_i}^{-1}(\mathfrak{r}_i),
\gamma_{u_i}^{-1}(\infty))\times S^1)})\leq C'(\lambda')^{1/2-1/(2p)}\ll\varepsilon. 
\]
In addition, the second line of (\ref{split:pcr}) and the above
estimates for terms therein, combined with (\ref{pt-wise}) 
and the length estimate for this region yield \[(\lambda')^{-1}\|\xi_{0L_y}'\|_{L^p([\gamma_{u_i}^{-1}(\mathfrak{r}_i),
\gamma_{u_i}^{-1}(\infty))\times S^1)}\leq
C_L'(\lambda')^{1/2-1/(2p)}\ll \varepsilon.\] 
In sum, we have \[
\|\xi_{\chi'}\|_{W_{\chi'}([\gamma_{u_i}^{-1}(\mathfrak{r}_i),
\gamma_{u_i}^{-1}(\infty))\times S^1)}\leq C_1\|\xi_0\|_{W^{\lambda'}_0([\gamma_{u_i}^{-1}(\mathfrak{r}_i),
\gamma_{u_i}^{-1}(\infty))\times S^1)}\ll\varepsilon.\]
{\em To estimate on the second region}, namely on
$[\gamma_{u_i}^{-1}(\gamma_0),
\gamma_{u_i}^{-1}(\mathfrak{r}_i)]\times S^1$,
let \(\beta^+_i(s)\) be a smooth cutoff function with:
\begin{itemize} 
\item support on \([\gamma^{-1}_{u_i}(\gamma_0)-1,
(1+\epsilon)\gamma_{u_i}^{-1}(\mathfrak{r}_i)]=:\Theta_{\beta_i^+}\);
\item value 1 over \([\gamma_{u_i}^{-1}(\gamma_0),
\gamma_{u_i}^{-1}(\mathfrak{r}_i)]\), and 
\item \(
|(\beta_i^+)'|<C'\lambda^{1/2}\quad \text{on \([\gamma^{-1}_{u_1}(\mathfrak{r}_i),
(1+\epsilon)\gamma_{u_i}^{-1}(\mathfrak{r}_i)]\).}\)
\end{itemize}
Notice that  $(\gamma_{u_i}^{-1})^*(\beta_i^+\xi_{\chi'})\in
W_{u_i}'$, and since 
$E_{u_i}|_{W_{u_i}'}$ an isomorphism,
we have from Lemma \ref{epsilon} and (\ref{pcr-uni}) that
\[\begin{split}
\|\beta_i^+\xi_{\chi'}\|_{W_{\chi'}} &\leq 
\|E_{w_{\chi}} (\beta_i^+\xi_{\chi'})\|_{L_{\chi'}}\\
 \leq |\alpha   _{\chi'}| \|\beta_i^+ &Y_{(\lambda', w_{\chi'})}\|_{L_{\chi'}}+\|\beta_i^+(\bar{\p}_{J
  X_{\lambda'}} w_{\chi'})_{T_w}\|_{L_{\chi'}} +\|\beta_i^+
\hat{n}_{(\lambda',w_{\chi'})}(\alpha_{\chi'},\xi_{\chi'})\|_{L_{\chi'}}
+\|(\beta_i^+)'\xi_{\chi'}\|_{L_{\chi'}}.
\end{split}
\]
By the assumption (\ref{pt-wise}) and the length estimate
\(
|\gamma_{u_i}^{-1}(\mathfrak{r}_i)-\gamma_{u_i}^{-1}(\gamma_0)|\leq
C'(\lambda')^{-1/2}, 
\)
we may bound each terms on the RHS as follows:\
\begin{itemize}
\item The first term may be bounded by \(C_1|\lambda'|^{1/2-1/(2p)}\).
\item The second term is 
already estimated to be small in Proposition \ref{err-est}. 
\item The computation in the proof of Lemma \ref{nonlinear:bound}
  shows that
\[\begin{split}
&\|\beta_i^+\hat{n}_{(\lambda',w_{\chi'})}(\alpha_{\chi'},\xi_{\chi'})\|_{L_{\chi'}}\\
&\quad\leq
C_n(\|(\alpha_{\chi'},0)\|_{\hat{W}_{\chi'}}^2+\|(\alpha_{\chi'},0)\|_{\hat{W}_{\chi'}}\|\beta_i^+\xi_{\chi'}\|_{W_{\chi'}}+\|\sigma_{\chi'}^{1/2}\xi_{\chi'}\|_{L^\infty(\Theta_{\beta_i^+})}\|\beta_i^+\xi_{\chi'}\|_{W_{\chi'}})\\
&\quad \leq
C'_n\Big(C_1(|\lambda'|^{1/2-1/(2p)})^2+C_2(|\lambda'|^{1/2-1/(2p)}+\varepsilon\gamma_0^{-1}+\lambda^{1/2})\|\beta_i^+\xi_{\chi'}\|_{W_{\chi'}}\Big).
\end{split}\]
\item By the defining properties of
\(\beta_i^+\),
\[
\|(\beta^+_i)'\xi_{\chi'}\|_{L_{\chi'}}\leq C''|\lambda'|^{1/2-1/(2p)}.
\]
\end{itemize}
Collecting all the above and rearranging, we obtain
\[
\|\xi_{\chi'}\|_{W_{\chi'}([\gamma^{-1}_{u_i}(\gamma_0),
\gamma_{u_i}^{-1}(\mathfrak{r}_i)]\times S^1)} \leq C_i\varepsilon.
\]
Now that we have the estimates for the \(W_{\chi'}\)-norm over all the
various regions, we conclude \(\|\xi_{\chi'}\|_{W_{\chi'}}\leq
C\varepsilon\), and hence the claim of the Lemma.
\end{proof}

\mysubsubsection[Concluding the proof of Proposition \ref{8.1} (a).]
Recall from \S4.1.4 the K-model for \(\hat{E}_{(\lambda',
  w_{\chi'})}\): \([\ker \hat{E}_{(\lambda',
  w_{\chi'})}\to *]_{B_{\chi'}}\), where \(B_{\chi'}\) was chosen to
be the following subspace of \(\hat{W}_{\chi'}\):
\[
B_{\chi'}=\Big\{(0, \xi_{\chi'}) \, \Big|\, \langle  (u_0)_\gamma(0), (\gamma_{u_0}^{-1})^*\xi_{\chi'}(0)\rangle_{2,t}=0\Big\}.
\]
Thus, setting \(\lambda'=\lambda\) and \(\chi'=\chi\),
\((\alpha_{\chi'}, \xi_{\chi'})=(0, \xi_\chi)\), and \(\xi_\chi\) can
be expressed in terms of \(\tilde{\xi}_\chi\) (cf. \S4.2.3) and \(\tilde{\xi}\)
(cf. \S4.2.2). In particular, by (\ref{after-transl}) and the first
line of (\ref{t-xi-chi}), (\ref{surj-est}.2) holds.
On the other hand, combining Lemma 4.2.2 and the second line of
(\ref{t-xi-chi}), we see that the assumption (\ref{pt-wise}) holds, 
and therefore Lemma 4.2.4 implies the validity of (\ref{surj-est}.1).
The arguments in \S\ref{step4} then complete the last step of the
proof of Proposition \ref{8.1} (a).\qed

\subsection{Gluing Broken Orbits.}
We now discuss the modification needed for the proof of Proposition
\ref{8.1} (b).

Given a broken orbit \(\{\hat{u}_1, \hat{u}_2, \ldots, \hat{u}_k\}\)
connected at \(y\), and an \(\lambda\in (0, \lambda_0)\), the
pregluing \(w_\chi\) associated to \(\chi=(\{\hat{u}_1, \hat{u}_2,
\ldots, \hat{u}_k\}, \lambda)\in \Xi(\mathbb{S})\) is given by 
\[
w_\chi=\underline{w}_\chi,
\]
where \(\underline{w}_\chi\) is given by the same formula
(\ref{equation:pregluing}), except that now \(i\in \{1, \ldots, k\}\)
only, and instead of taking values in \(\R\), \(s\) now takes value in \(\R /
T_\chi\Z\), where \[T_\chi:=\mathfrak{l}_{k+1}.\]

With this explained, the material in sections 2 and 3 transfers
directly to the case of broken orbits, but the discussion in \S4.1,
4.2 above requires the following modification.

\mysubsubsection[Constructing the gluing map.]\label{gluing-map:2.1b}
At a closed orbit \((\lambda, (T, w))\in {\cal
  B}_O^\Lambda=\Lambda\times {\cal B}_O\), the deformation operator is
\(\hat{D}_{(\lambda, (T, w))}: \R_\alpha\oplus \R_\varrho\oplus
L^p_1(w^*K)\to L^p(w^*K)\), 
\[
\hat{D}_{(\lambda, (T, w))}(\alpha, \varrho, \xi)=\alpha
\partial_\lambda \check{\theta}_{X_\lambda}+\tilde{D}_{(T, w)} (\varrho, \xi),
\]
namely, it is a rank 2 stabilization of \(D_w\) (cf. \S3.3.1).
(\(\R_\alpha\), \(\R_\varrho\) above respectively parametrize
variation in \(\lambda\) and in the period \(T\)).
In our context, this operator is a map between the weighted spaces
\[
\hat{W}_\chi:=\R_\alpha\oplus \R_\varrho\oplus W_\chi \quad \text{and}
\quad L_\chi.
\]
Let \(K_\chi=\op{Span}\{\mathfrak{e}_{u_i}\}_{i=1}^k\), \(C_\chi=\op{Span}\{\mathfrak{f}_{j}\}_{j=1}^k\).
The analog of Proposition \ref{UI} shows that
\([K_\chi\to C_\chi]\) forms a K-model for \(D_{w_\chi}\),
which induces a K-model for \(\hat{D}_{(\lambda, (T_\chi, w_\chi))}\) by stabilization.

However, since \(s\) is now periodic instead of real, the matrix
\(
\lambda^{-1/2+1/(2p)}(P_{\mathfrak{f}_j} D_{w_\chi} \mathfrak{e}_{u_i})
\)
is no longer (approximately) triangular, and hence not clearly 
uniformly invertible. Consequently, it is no longer clear that,
with the choice of \(Q_\chi\) in \S4.1, reduction by
\(Q_\chi \) gives another K-model for the deformation operator.
Instead, use the following subspace \(Q_{O, \chi}\subset \hat{W}_\chi\):
\[
Q_{O, \chi}=\R_\alpha\oplus *\oplus \op{Span}\{\mathfrak{e}_{u_i}\}_{i=1}^{k-1}.
\]
Note that from (\ref{a-bdd}) and the uniform boundedness of \(P_{{\mathfrak
    f}_j}\) that:
\[
C_{\alpha-}\lambda^{1/2-1/(2p)}\leq\lambda^{3/2}P_{{\mathfrak
    f}_j}\p_\lambda \check{\theta}_{X_\lambda}(w_\chi)\leq C_{\alpha
  +}\lambda^{1/2-1/(2p)} \quad \text{\small for \(\lambda\)-independent
  constants \(C_{\alpha\pm}>0\)}.
\]
Supplementing Lemma
\ref{est:cok} (b) with this additional estimate, we see that the matrix
representation of the operator 
\[\lambda^{-1/2+1/(2p)}\Pi_{C_\chi}\hat{D}_{(\lambda, (T_\chi,
  w_\chi))}|_{Q_{O,\chi}}\] 
with respect to the bases
\[\Big\{(\lambda^{3/2},0,0), (0, 0, {\mathfrak e}_1), \ldots, (0, 0, {\mathfrak e}_{k-1})
\Big\} ,\quad 
\{\mathfrak{f}_1, \ldots, \mathfrak{f}_k\}\] 
is, modulo ignorable terms, of the form
\[
\left(\begin{array}{ccccc}
 +&- &0&\cdots& 0\\
+&+ &-& 0&0\\
 + &0 &+ &\ddots&\vdots\\
\vdots &\vdots &\ddots & \ddots& - \\
 +&0 &\cdots &0& +
\end{array}\right)\quad \text{\small (\(+/-\) denote positive/negative
numbers of \(O(1)\)),} 
\] 
which is easily seen to have a uniformly bounded right inverse. 
Thus, the rest of \S4.1 may be repeated with \(Q_\chi\) replaced by
\(Q_{O, \chi}\) to define a gluing map in this case, which is also a
local diffeomorphism.

\mysubsubsection[Surjectivity of the gluing map.]
As the choice of \(Q_{\chi'}\) is changed, the definition of the
\(B\)-space \(B_{\chi'}\) changes accordingly. In this situation, 
\[
B_{\chi'}=\Big\{(\alpha', 0, \xi_{\chi'}) \, \Big|\, \langle
(u_k)_\gamma(0),
(\gamma_{u_k}^{-1})^*\xi_{\chi'}(0)\rangle_{2,t}=0\Big\}\subset \hat{W}_{\chi'}.
\]
(Note that in the case of broken orbits, \(i\in \Z/k\Z\), thus
\(u_0=u_k\)). The work in \S4.2 needs corresponding modification.

Given a \((\lambda, (T, w))\in \hat{\cm}_O^\Lambda\) close to the
broken orbit \(\{\hat{u}_1, \ldots, \hat{u}_k\}\), one may define \(\tilde{w}\) and
\(\tilde{\gamma}_{u_i}\) in essentially the same way as \S4.2.1, and
the estimates in \S4.2.2 still hold. However, for \(\chi=(\{\hat{u}_1,
\ldots, \hat{u}_k\}, \lambda)\), the period \(T_\chi\) of the 
pregluing \(w_\chi\) differs from those of \(\tilde{w}\) or
\(w\). Thus, instead of comparing with \(w_\chi\), we compare \(w\) or
\(\tilde{w}\) with \(w_{\chi'}\), where \(\chi'=(\{\hat{u}_1, \ldots,
\hat{u}_k\}, \lambda')\), and \(\lambda'=\lambda+\alpha'\) is chosen so that the period
of \(w_{\chi'}\) agrees with the period of \(w\) (which is also the
period of \(\tilde{w}\)). 
With this choice of \(\chi'\), the assumption (\ref{after-transl}),
together with the definition of \(B_{\chi'}\) above, imply (\ref{surj-est}.2). 

Moreover, the length estimates in \S2.2.1 show that 
\[
C_-\lambda^{-1/2}\leq T_\chi\leq C_+\lambda^{-1/2};
\]
combining with the estimate for the difference in periods of \(w\) and
\(w_\chi\) given by (\ref{Delta-s}), we have:
\begin{equation}\label{alpha-small}
|\alpha'|\leq C \lambda^{3/2}\varepsilon.\end{equation}
For such \(\lambda'=\lambda+\alpha'\), the difference \(\gamma_{\lambda',
  u_i}^{-1}-\tilde{\gamma}^{-1}_{u_i}\) satisfies estimates similar to
(\ref{Delta-s}). Thus, (\ref{pt-wise}) holds for this choice of
\(\chi'\), which in turn implies (\ref{surj-est}.1), via Lemma 4.2.4.

\section{Gluing at Births.}

The purpose of this section is to prove Proposition \ref{thm:birth}
below. The proof is in many ways similar to the proof of Proposition
\ref{8.1}, but simpler in Step 2, since here we glue only a single flow
line, and the generalized cokernel is
this case is trivial.

\subsection{Statement of the Gluing Theorem.}
The next Proposition verifies part of 
(RHFS2c, 3c) for admissible \((J,X)\)-homotopies.

\begin{proposition*}\label{thm:birth}
Let \((J^\Lambda, X^\Lambda)\) be an admissible \((J, X)\)-homotopy
connecting two regular pairs, and \({\bf x}, {\bf z}\) be two path
components of \({\cal P}^\Lambda\backslash {\cal P}^{\Lambda, deg}\).
Then a chain-topology neighborhood of \({\mathbb J}_{P} (\Lambda, {\bf x}, {\bf
  z}; \Re)\) in \(\hat{\cm}_P^{\Lambda, 1, +}({\bf x}, {\bf z}; \op{wt}_{-\langle{\cal Y}\rangle,
  e_{\cal P}}\leq \Re)\) is l.m.b. along \(\mathbb{J}_{P} (\Lambda , {\bf x}, {\bf
  z}; \Re)\).

Furthermore, \(\Pi_\Lambda\) maps these neighborhoods 
to birth-neighborhoods. 
\end{proposition*}

We shall restrict our attention to 
\(\hat{u}\in \hat{\cm}_{P,\lambda}^{0}((\underline{\bf x, [w]}), 
(\underline{\bf z, [v]}))
\cap {\mathbb J}_{P} (\Lambda, {\bf x}, {\bf z}; \Re) \), where one of
\(x_\lambda\) and \(z_\lambda\) is a death-birth, and
\(\op{gr}_+((x_\lambda, [w_\lambda]), (z_\lambda, [v_\lambda]))=\op{gr} (({\bf
  x}, [{\bf w}]), ({\bf z}, [{\bf v}]))\). The other cases follow either from standard
gluing theory or structure theory of parameterized moduli spaces,
since the flow lines decay exponentially to the critical points in these cases.

Without loss of generality, assume as in sections 2--4 that
\(\lambda=0\), that \(z_\lambda=y\)
is in a standard death-birth neighborhood, and that the \((J, X)\)-homotopy
is oriented such that \(C'_y>0\). Under these assumptions, a birth
neighborhood is \((-\lambda_0, 0)\subset \Lambda\), for a small \(\lambda_0>0\).

Our goal is thus to construct a gluing map from 
\(\Xi({\mathbb S})\) to \(\hat{\cm}_P^{\Lambda, 1}({\bf x}, {\bf z}; \op{wt}_{-\langle {\cal Y}\rangle,
  e_{\cal P}}\leq \Re)\), where
\[
{\mathbb S}=\hat{\cm}_{P,0}^{0}(\underline{\bf x},
\underline{\bf y};\op{wt}_{-\langle {\cal Y}\rangle, e_{\cal P}}\leq\Re); \quad
\Xi({\mathbb S})={\mathbb S}\times (-\lambda_0,0) \quad \text{for a
  small \(\lambda_0>0\).}
\]

We shall again focus on a single \(\hat{u}\in
\mathbb{S}\), since in this case \({\mathbb S}\) also consists of
finitely many isolated points.  Notice that when \(x_0=y\), \(\hat{u}\) can be the
constant flow at \(y\), \(\bar{y}\). The argument required for this
case is somewhat different from the other cases. We discuss this case
in \S5.3, and the other cases in \S5.2.

\subsection{When $u\neq\bar{y}$.}
Assume without loss of generality that 
\(x_0\neq y\) is nondegenerate, so that we may concentrate on
the region where \(s>0\).

\mysubsubsection[Pregluing.]
Let \(\chi:=(\lambda,\hat{u})\in\Xi(\mathbb{S})\) as above, and let
\(u\) be a centered representative.
Write \[u(s)=\exp(y, \mu(s))\quad \text{for large \(s\),}\] 
and as in I.5.3.2, let \[y_{\lambda -}=\exp (y, \eta_{\lambda-})\in {\cal
  P}_\lambda\] be the critical point near \(y\) of index \(\ind_-(y)\). 
Note that $\langle{\mathbf{e}_y}, \mu\rangle_{2,t}(s)>0$ is a decreasing
function for large \(s\), sending \((s_0, \infty)\) to \((C, 0)\) for
some positive numbers \(s_0, C\). Since $\langle{\mathbf{e}_y},
\eta_{\lambda -}\rangle_{2,t}$ is a small positive number, it equals 
$\langle\mathbf{e}_y, \mu(s)\rangle_{2,t}$ for certain large
$s=\breve{\gamma}_{\lambda}$. From the
estimates in Lemma I.5.3.2 and
Proposition I.5.1.3, we have
\[
C|\lambda|^{-1/2}\leq \breve{\gamma}_\chi\leq C'|\lambda|^{-1/2}.
\]

Let $R<\breve{\gamma}_\chi-1$ be a $\lambda$-independent large positive number
such that $u(s)$ is close to $y$ for $s\geq R$, and set
$R_{\pm}=\pm C_0|\lambda|^{-1/2}$ for some \(\lambda\)-independent constant \(C_0>0\). 
Define $u_{\lambda}\in \Gamma((-\infty, \breve{\gamma}_\chi)\times
S^1, p_2^*T_f)$ by
\begin{equation}\label{u-lambda}
u_{\lambda}(s):=\begin{cases}
e_{R_-, R_+}(0, u;\lambda, 0) &\text{when \(s\leq R/2\),}\\
\exp \Big(y, \mu(s)+\beta(s-R)\Pi_{\ker A_y}^{\perp}\eta_{\lambda -}\Big)&\text{when \(s\geq R/2\)},
\end{cases}
\end{equation}

\begin{lemma*}\label{lemma:gamma2}
There is a function $\gamma_{\chi}(s)$ defining a homeomorphism from
$\R$ to $(-\infty, \breve{\gamma}_\chi)$, such that 
\begin{equation}
\begin{cases}
\langle w_{\chi}'(s),
\bar{\partial}_{J,X_{\lambda}}w_{\chi}(s)\rangle_{2,t}=0;&
\text{when 
  \(s\in [\gamma_\chi^{-1}(0), \infty)\)};\\
\gamma_{\chi}'=1 &\text{otherwise.}
\end{cases}
\end{equation}
\end{lemma*}
\begin{proof}
Write
\(
\frac{d\gamma_{\chi}}{ds}=h(\gamma_{\chi})\),
where 
\[
h(\gamma):=\begin{cases}
1-\langle (u_{\lambda})_\gamma, \bar{\partial}_{JX_{\lambda}}(u_{\lambda})\rangle_{2,t}\|(u_{\lambda})_\gamma\|_{2,t}^{-2}&
\text{when 
  \(s\in [\gamma_\chi^{-1}(0), \infty)\)};\\
1 &\text{otherwise.}
\end{cases}
\]
We now examine the behavior of 
$\bar{\partial}_{JX_{\lambda}}u_{\lambda}(\gamma)$ near $\gamma=\breve{\gamma}_\chi$.
Here since $u_{\lambda}$ is close to $y_{\lambda-}$, expanding 
$\bar{\partial}_{JX_{\lambda}}$ about $y_{\lambda-}$ and writing
$\exp (y_{\lambda -}, \mu_\lambda(\gamma))=u_{\lambda}(\gamma)$, we have:
\[
T_{u_\lambda, y_{\lambda -}}\bar{\partial}_{JX_{\lambda}}u_{\lambda}(\gamma)=(\mu_\lambda)_\gamma
+A_{y_{\lambda -}}\mu_\lambda
+n_{y_{\lambda -}}(\mu_\lambda).
\]
By definition, $\mu_\lambda(\breve{\gamma}_\chi)=0$; hence
$h(\breve{\gamma}_\chi)=0$. Thus, 
\begin{equation}\label{h-lambda}\begin{split}
h(\gamma)&=
\Big\langle (u_{\lambda})_\gamma(\gamma), T_{y_{\lambda -}, u_\lambda}A_{y_{\lambda -}}T_{y_{\lambda -}, u_\lambda}^{-1}((\breve{\gamma}_\chi-\gamma)((u_{\lambda})_\gamma(\gamma))\Big\rangle_{2,t}\|(u_{\lambda})_\gamma (\gamma)\|_{2,t}^{-2}\nonumber\\
&\qquad +O\Big(|\breve{\gamma}_\chi-\gamma|^2\|(u_{\lambda})_\gamma(\gamma)\|_{2,1,t}\Big).\end{split}
\end{equation}
By the estimate for minimal eigenvalue of \(A_{y_{\lambda-}}\) in I.5.3.2, this 
is bounded above and below by multiples of $|\lambda|^{1/2}(\breve{\gamma}_\chi-\gamma)$. Integrating like (\ref{gamma-int}), we see that
for large \(s\)
\begin{equation}\label{delta-gamma:est}
C'_5e^{-c_6'|\lambda|^{1/2}s}\leq \breve{\gamma}_\chi-\gamma_\chi\leq C_5e^{-c_6|\lambda|^{1/2}s},
\end{equation}
while on the other end 
$\gamma_{\chi}(s)= s+c_{\lambda}$ for some constant $c_{\lambda}$. 
We define $\gamma_{\chi}$ such
that $\gamma_{\chi}(s)=s$ for \(s<0\).
\end{proof}
\begin{definition*}
The {\em pregluing} \(w_\chi\) corresponding to gluing data
$\chi=(\lambda, u)$ is
\[
w_{\chi}(s):=u_{\lambda}(\gamma_{\chi}(s)).
\]
\end{definition*}

\mysubsubsection[The weighted norms.]
The norms $W_{\chi}$, $L_{\chi}$ here are defined by the same formulae in Definition
\ref{w-norms}, with the {\em weight function} $\sigma_{\chi}$ replaced by
\[
\sigma_{\chi}(s):=
\begin{cases}
\|w'_{\chi}(\gamma_{\chi}^{-1}(0))\|_{2,t}^{-1} &\text{when $s\leq \gamma_{\chi}^{-1}(0)$,}\\
\|w'_{\chi}(s)\|_{2,t}^{-1} &\text{when $\gamma_{\chi}^{-1}(0)\leq s\leq \gamma_{\chi}^{-1}(\mathfrak{r}_\chi)$,}\\
\|w'_{\chi}(\gamma_{\chi}^{-1}(\mathfrak{r}_\chi))\|_{2,t}^{-1} &\text{when $s\geq\gamma_{\chi}^{-1}(\mathfrak{r}_\chi)$,}
\end{cases}
\]
where $\mathfrak{r}_\chi=\mathfrak{r}_\chi(\lambda, \epsilon)=C_{\tau}(\lambda/\epsilon)^{-1/2}<\breve{\gamma}_\chi$ is chosen such that
$1-h(s)\leq \epsilon$ where $s\leq \mathfrak{r}_\chi$ for a small
positive number $\epsilon$.

\medbreak

We shall frequently call on the following useful 
\begin{facts*}\begin{description}
\item{(a)} In this case $\gamma'_{\chi}\leq 1$.  
\item{(b)} \(\|\Pi_{\ker A_y}^{\perp}\eta_{\lambda -}\|_{2,2,t}\leq
C|\lambda|\) by I.(55), Lemma I.5.3.2, and the decay estimates in Proposition
I.5.1.3. 
\item{(c)} $\sigma_{\chi}\leq C|\lambda|^{-1}$. \end{description}\end{facts*}
In particular, Fact (b) often implies that in addition to estimates
analogous to those in the proof of Proposition 2.1, the extra
terms introduced by the cutoff function $\beta$ in the definition of 
$u_{\lambda}$ is usually ignorable.

\mysubsubsection[Error estimate.]\label{err-est2} Proceeding to Step 1
of the gluing theory, we have:
\begin{lemma*}
$\|\bar{\partial}_{JX_{\lambda}}w_{\chi}\|_{L_{\chi}}\leq C\lambda^{1/2-1/(2p)}$.
\end{lemma*}
\begin{proof}
Consider the two regions 
(a) $\gamma_{\chi}^{-1}(-R_-)\leq s\leq \gamma_{\chi}^{-1}(R)$, 
(b) $\gamma_{\chi}^{-1}(R)\leq s\leq
\infty$ separately. The point is to expand
$\bar{\partial}_{JX_{\lambda}}w_{\chi}(s)=\tilde{\Pi}_{u_{\lambda}'}^{\perp}(\bar{\partial}_{JH_{\lambda}}u_{\lambda}(\gamma_{\chi}(s)))$ differently in the two regions:
expand $u_{\lambda}$ about $u$ in region (a), and about $y_{\lambda}$ in
region (b).

{\em In region (a)}, modulo terms coming from $\beta(\gamma_\chi-R)\Pi_{\ker A_y}^{\perp}\eta_{\lambda -}$, 
the estimate of the norm is entirely parallel to that in Proposition \ref{err-est}:
The time $w_{\chi}$ spends in this region is
$$\gamma_{\chi}^{-1}(R)-\gamma_\chi^{-1}(R_-)\leq C|\lambda|^{-1/2};$$ 
the $L^{\infty}$ norm of
$\sigma_{\chi}\bar{\partial}_{JX_{\lambda}}w_{\chi}$
can be estimated as in Case 1 of Proposition \ref{err-est}, with 
$\gamma_0$ replaced by $R$.

On the other hand, since $\sigma_{\chi}$ has a
$\lambda$-independent uniform bound in this region, and the
$L^p_t$-norm of the contribution to 
$\bar{\partial}_{JX_{\lambda}}w_{\chi}$ from the extra terms
introduced by $\beta \Pi_{\ker A_y}^{\perp}\eta_{\lambda -}$ can be bounded by
$C\|\Pi_{\ker A_y}^{\perp}\eta_{\lambda -}\|_{2,1,t}\leq C'|\lambda|$, the contribution from these
terms to $\|\sigma_{\chi} \bar{\partial}_{JX_{\lambda}}w_{\chi}\|_p$
is thus bounded by $C|\lambda|$.

{\em For region (b)}, $w_{\chi}$ spends infinite amount of time here; however
\[\begin{split}
\Pi_{(u_{\lambda})_\gamma}^{\perp}\bar{\partial}_{JX_{\lambda}}u_{\lambda}(\gamma)&=\Pi_{(u_{\lambda})_\gamma}^{\perp}\Big((\delta\gamma)T_{y_{\lambda-},
  u_\lambda}A_{y_{\lambda}}(T_{y_{\lambda-},
  u_\lambda})^{-1}(u_\lambda)_\gamma(\bar{\gamma}))\\
&\qquad +T_{y_{\lambda-},
  u_\lambda}n_{y_{\lambda}}((\delta\gamma )(T_{y_{\lambda-},u_\lambda})^{-1}(u_{\lambda})_\gamma(\bar{\gamma}))\Big), 
\end{split}\]
where $\delta \gamma:=\breve{\gamma}_\chi-\gamma$; $\gamma\leq
\bar{\gamma}\leq \breve{\gamma}_\chi$. On the other hand in this region
$\sigma_{\chi}(s)\leq C|\lambda|^{-1}$. 
Thus by Lemma I.5.3.2, on this region $\|\sigma_{\chi}
\bar{\partial}_{JX_{\lambda}}w_{\chi}\|_p$ is bounded by 
$C\|(\delta\gamma)\|_p|\lambda|^{1/2}\leq C'|\lambda|^{1/2-1/(2p)}$,
since $\delta\gamma\leq C_5e^{-C_6|\lambda|^{1/2}s}$ by (\ref{delta-gamma:est}).
\end{proof}

\mysubsubsection[Existence and uniform boundedness of the right inverse $G_{\chi}:L_{\chi}
\to W_{\chi}$ of $E_{w_\chi}$.] We now proceed to Step 2 of the proof. 
In this case $W_{\chi}'\subset W_{\chi}$ is 
\[
W_{\chi}':=\Big\{
\xi\in W_{\chi}\, |\, \langle\nu(0),\xi(\gamma_\chi^{-1}(0))\rangle_{2,t}=0
\;\text{for all $\nu\in \ker E_u$}\Big\},
\]
and we aim to show that there is a uniformly bounded isomorphism
\(G_\chi:L_\chi\to W_\chi'\) which is a right inverse of \(E_{w_\chi}\).
Assume the opposite, that there is a sequence 
\(\{\xi_{\lambda}\in W'_{\chi}\}_\lambda\) satisfying
\begin{gather}
\|\xi_{\lambda}\|_{W_{\chi}}=1; \nonumber\\
\|E_{w_\chi}\xi_{\lambda}\|_{L_{\chi}}=:\varepsilon_E(\lambda)\to 0\quad
\text{when $\lambda\to 0$.}\label{As:xi2}
\end{gather} 

Divide $\Theta=\R\times S^1$ into two parts $\Theta_u, \Theta_{y-}$,
separated by the line $s=\gamma_{\chi}^{-1}(\mathfrak{r}_\chi)$. 
Let $\Theta_u':=(-\infty, \gamma_{\chi}^{-1}(\mathfrak{r}_\chi)+1)\times S^1\supset\Theta_u$;  
$\Theta_{y-}':=(\gamma_{\chi}^{-1}(\mathfrak{r}_\chi-1), \infty)\times S^1\supset
\Theta_{y-}$.

\underline{On $\Theta_u'$}, we define $\xi_{\lambda , u}\in \Gamma(u^*K)$ by 
$$
T_{u, u_\lambda}\xi_{\lambda, u}(\gamma_{\chi}(s))=\xi_\lambda(s).
$$
Let $(\xi_{\lambda, u})_L$ be the projection of $\xi_{\lambda , u}$ to the direction of
$u'$ and let
\((\xi_{\lambda, u})_T=\xi_{\lambda, u}-(\xi_{\lambda, u})_L\). 
Let \(\beta_u\) be a smooth cutoff function supported on \(\gamma_\chi(\Theta_u')\)
with value 1 on \(\gamma_\chi(\Theta_u)\).
Arguing as in the proof of Proposition \ref{UI}, one obtains:
\begin{eqnarray}
\lefteqn{\|\xi_\lambda\|_{W_{\chi}(\Theta_u)}
\leq C\|\xi_{\lambda , u}\|_{W_u(\gamma_{\chi}(\Theta_u))}  }\nonumber\\
&&\leq
C'\|\beta_uE_u\xi_{\lambda , u}\|_{L_u(\gamma_{\chi}(\Theta_u'))}+C\|\beta'_u(\xi_{\lambda, u})_T\|_{L_u(\gamma_{\chi}(\Theta_u'))}\nonumber\\
&&\leq
C'(1+2\epsilon)\|E_{w_\chi}(\xi_\lambda)\|_{L_{\chi}(\Theta')}+C''(\epsilon+|\lambda|^{1/2})\|\xi_\lambda\|_{W_{\chi}}+C\|\beta'_u(\xi_{\lambda, u})_T\|_{L_u(\gamma_{\chi}(\Theta_u'))}\nonumber\\
&&\leq 2C'\varepsilon_E+C''(\epsilon+|\lambda|^{1/2})+C\varepsilon_0.\label{Theta-Bdd}
\end{eqnarray}
(Note in comparison with (\ref{Theta-bound}), the 2nd term in the 3rd line
above has a worse factor of $|\lambda|^{1/2}$ instead of
\(|\lambda|\); this arises from the difference between \(u\) and \(u_\lambda\).)

On the other hand, \underline{on $\Theta_{y-}'$} we consider
$\xi_{\lambda , y-}\in \Gamma(\bar{y}_{\lambda -}^*K)$ defined by
$$
T_{\bar{y}_{\lambda -}, w_\chi}\xi_{\lambda , y-}=\xi_\lambda.
$$
Let \({\bf e}_{y_{\lambda -}}\) be the unit eigenvector associated with
the minimal eigenvalue of \(A_{y_{\lambda-}}\), which goes
to \({\bf e}_y\) as \(\lambda \to 0\). By I.(57), \({\bf e}_{y_{\lambda -}}\) differs from $T_{w_\chi,
  y_{\lambda-}}w_\chi'/\|w'_\chi\|_{2,t}(s)$ by
$O(|\lambda|^{1/2})$ for \(s\in \Theta_{y-}'\). 
Let \[(\xi_{\lambda , y-})_L:=\Pi_{{\bf e}_{y_{\lambda -}}}\xi_{\lambda , y-}=\underline{(\xi_{\lambda , y-})}_L{\bf e}_{y_{\lambda -}},\]
and \((\xi_{\lambda , y-})_T=\xi_{\lambda , y-}-(\xi_{\lambda , y-})_L\).
The above observation about \({\bf e}_{y_{\lambda -}}\), together with
the fact that in this region $\sigma_{\chi}$ is bounded
above and below by multiples of $|\lambda|^{-1}$ imply that to
estimate \(\|\xi\|_{W_\chi(\Theta_{y-})}\) or \(\|\xi\|_{L_\chi(\Theta_{y-})}\), it is equivalent to
estimate \(\|\xi_{y-}\|_{W_{y-}(\Theta_{y-})}\) or \(\|\xi_{y-}\|_{L_{y-}(\Theta_{y-})}\), where
\[
\|\xi_{y-}\|_{W_{y-}}:=|\lambda|^{-1/2}\|\xi_{y-}\|_{p,1}+|\lambda|^{-1}\|(\xi_{y-})'_L\|_p,
\quad \|\xi_{y-}\|_{L_{y-}}:=|\lambda|^{-1/2}\|\xi_{y-}\|_{p,1}+|\lambda|^{-1}\|(\xi_{y-})_L\|_p.
\]
We have a refined version of Lemma \ref{claimF} in this case:
\begin{lemma*}[Refining Floer's lemma]\label{claimF2} 
Let $\xi_{\lambda}$ be as in (\ref{As:xi2}). Then 
for all sufficiently small
$\lambda$, 
\[
|\lambda|^{1/(2p)}\|\xi_{\lambda , y-}\|_{L^\infty(\Theta'_{y-})}+
\|(\xi_{\lambda , y-})_L\|_{L^\infty(\Theta'_{y-})}\leq \varepsilon_0(\lambda)|\lambda|^{1/2+1/(2p)},
\] 
where \(\varepsilon_0(\lambda)\) is a small positive number, with
\(\lim_{\lambda\to 0}\varepsilon_0(\lambda)=0\).
\end{lemma*}
\begin{proof}
The estimate for \(\|\xi_{\lambda , y-}\|_{L^\infty(\Theta'_{y-}})\) follows
easily from the argument for Lemma \ref{claimF}.
The longitudinal component has a more refined bound because it has a
better bound on the Sobolev norm. Let 
\[
\tilde{\varsigma}_{\lambda}(\tau):=\lambda^{-1/2-1/(2p)}\underline{(\xi_{\lambda , y-})}_L(\lambda^{-1/2}(\tau+s_{\lambda})) \quad \text{over \([1,\infty)\),}
\]
where \(s_\lambda\) are constants chosen such that \(\lambda^{-1/2}(1+s_{\lambda})=\gamma_{\chi}^{-1}(\mathfrak{r}_\chi-1)\).
Then by (\ref{As:xi2}), \(\|\tilde{\varsigma}_\lambda\|_{L^p_1([1, \infty))}\)  is
bounded (note the rescaling contributes a factor of \((\lambda)^{-1/(2p)}\) to the \(L^p_1\) norm).
Thus (again after possibly taking a subsequence)
\(\tilde{\varsigma}_\lambda\) converges in \(C_0\) to \(\tilde{\varsigma}_0\),
and \(s_\lambda\to s_0\). \(\tilde{\varsigma}_0\) satisfies an equation of the
form
\begin{equation}\label{equation:tilde-zeta}
\frac{d\tilde{\varsigma}_0}{d\tau}+\chi\tilde{\varsigma}_0=0, \quad \text{where \(\chi\sim C_1+C_2e^{-\nu'\tau}\).}
\end{equation}
 (The assumption of being in a standard d-b neighborhod
is used to simplify the differential equation above. Notice also 
that $(\xi_\lambda)_T$ does not
appear in this equation, because by the \(L^\infty_t\) estimate for \(\xi_\lambda\), its
contribution vanishes as \(\lambda\to0\).)
Thus,
\(\|\tilde{\varsigma}_0\|_{\infty}\leq|\tilde{\varsigma}_0(1)|\). Meanwhile, 
\(\tilde{\varsigma}_0(1)=0\) since by the argument for (\ref{sum:begin}),
\(\|(\xi_{\lambda})_L(\gamma_{\chi}^{-1}(\mathfrak{r}_\chi-1))\|_{\infty,t}\leq
C\lambda^{1/2+1/(2p)}\varepsilon_u\).
\end{proof}

Let $\beta_{y-}$ be a cutoff function on $\R$ which vanishes in $(-\infty,
\mathfrak{r}_\chi-1]$ and is 1 on $[\mathfrak{r}_\chi, \infty)$. 
We may estimate the longitudinal component as:
\begin{eqnarray}\label{longi}
\lefteqn{|\lambda|^{-1}\|(\xi_{\lambda , y-})_L'\|_{L^p(\Theta_{y-})}
+|\lambda|^{-1/2}\|(\xi_{\lambda , y-})_L\|_{L^p(\Theta_{y-})}}\nonumber\\&&\leq C|\lambda|^{-1} \|E_{y_{\lambda}}(\beta_{y-}(\gamma_{\chi})(\xi_{\lambda , y-})_L)\|_{L^p(\Theta_{y-}')}\nonumber\\
&&\leq C|\lambda|^{-1} \|\beta_{y-}(\gamma_{\chi})(E_{y_{\lambda-}}(\xi_{\lambda , y-})_L)_L\|_{L^p(\Theta_{y-}')}
+C'|\lambda|^{-1}\|(\beta_{y-})_\gamma(\gamma_\chi)\gamma_{\chi}'(\xi_{\lambda , y-})_L\|_{L^p(\Theta_{y-}')}\nonumber\\
&& \leq C_1\|\beta_{y-}\circ\gamma_{\chi}E_{w_{\chi}}\xi_\lambda\|_{L_{\chi}(\Theta_{y-}')}
+C_2|\lambda|^{-1/2}\|\beta_{y-}(\gamma_{\chi})e^{-C_6|\lambda|^{-1/2}(s-\gamma_\chi^{-1}(\mathfrak{r}_\chi))}(\xi_{\lambda , y-})_L\|_{L^p(\Theta_{y-}')}\nonumber\\
&&\qquad
+C'|\lambda|^{-1}\|(\beta_{y-})_\gamma(\gamma_\chi)\gamma_{\chi}'(\xi_{\lambda , y-})_L\|_{L^p(\Theta_{y-}')}+C_3\|\beta_{y-}(\gamma_{\chi})e^{-C_6|\lambda|^{-1/2}(s-\gamma_\chi^{-1}(\mathfrak{r}_\chi))}(\xi_{\lambda , y-})_T\|_{L^p(\Theta_{y-}')}\nonumber\\
&&\leq C_1'\varepsilon_E+C_2'\varepsilon_0.
\end{eqnarray}
The first inequality above follows from the eigenvalue estimate for
\(A_{y_{\lambda -}}\) in I.5.3.2.
The second term in the penultimate expression above comes from the
difference between \(E_{w_\chi}\) and (a conjugate of)
\(E_{y_{\lambda-}}\), while the last term arises from \((T_{\bar{y}_{\lambda -}, w_\chi}^{-1}E_{w_\chi}T_{\bar{y}_{\lambda -}, w_\chi}(\xi_{\lambda , y-})_T)_L\).
(Note that this term would have an extra factor of
\(|\lambda|^{-1/2}\) if I.5.3.1 (1c) is not assumed). We have also
used Lemma \ref{claimF2} and the estimates that in this region,
\[\begin{split}
|(\beta_{y-})_\gamma\gamma_{\chi}'|\leq C|\lambda|^{1/2}\exp
(-C_6|\lambda|^{1/2}(s-\gamma_\chi^{-1}(\mathfrak{r}_\chi)))\quad \text{and that }\\
\|\mu_\lambda (\gamma_\chi(s))\|_{2,2,t}\leq C'|\lambda|^{1/2}\exp
(-C_6|\lambda|^{1/2}(s-\gamma_\chi^{-1}(\mathfrak{r}_\chi)), 
\end{split}\]
which in turn
follows from the computation in the proof of Lemma \ref{lemma:gamma2}.

Similarly, the transversal direction can be estimated by:
\begin{equation}\label{transv}
\begin{split}
&|\lambda|^{-1/2}(\|(\xi_{\lambda , y-})_T'\|_{L^p(\Theta_{y-})}+\|(\xi_{\lambda , y-})_T\|_{L^p(\Theta_{y-})})\\
&\qquad\leq C\|\beta_{y-}(\gamma_{\chi})E_{w_\chi}\xi_\lambda\|_{L_{\chi}(\Theta'_{y-})}
+C'|\lambda|^{1/2-1/(2p)}\varepsilon_0(\lambda)\leq C''\varepsilon_E+C'|\lambda|^{1/2-1/(2p)}\varepsilon_0.
\end{split}
\end{equation}
Combining (\ref{longi}), (\ref{transv}) and (\ref{Theta-Bdd}), we obtain
$\|\xi_\lambda\|_{W_{\chi}}\ll1$ for all large enough $\lambda$, and
hence the desired contradiction. 

\mysubsubsection[Surjectivity of the gluing map.] Estimates for the
nonlinear terms required for Step 3 in this case are not very
different from those discussed in \S 2.5, and hence will be omitted. 
The argument in \S1.2.1 then defines a gluing map, which is a local
diffeomorphism onto a \({\cal B}\)-topology neighborhood of the image of
pregluing map. Again, we need to show that the latter neighborhood
contains a chain-topology neighborhood of \(\mathbb{S}\).

To adapt the proof in \S4.2, given 
\((\lambda, \hat{w})\in \hat{\cm}^{1,\Lambda}({\bf x}, {\bf
  y}_{\lambda -})\) close to \(\hat{u}\in \mathbb{S}\) in the chain
topology neighborhood, 
we may again choose a representative \(w\) and \(\tilde{w}\) as in
\S4.2.1, satisfying conditions similar to (\ref{calibration}) and
(\ref{after-transl}):
\begin{itemize}
\item\(\tilde{w}(s):=u_{\lambda}(\tilde{\gamma}_\chi(s))\), where 
\(\tilde{\gamma}_\chi: \R\to(-\infty, \breve{\gamma}_\chi)\) is a
homeomorphism determined by 
\begin{equation}\label{new-calib}
\Pi_{{\bf e}_y}\tilde{\zeta}_\lambda(s)=\Pi_{{\bf e}_y}\zeta(s) \quad\forall s\in
[\tilde{\gamma}_\chi^{-1}(R+1), \infty),
\end{equation}
and \(\zeta\), \(\tilde{\zeta}_\lambda\) are defined by \(w(s)=\exp (y, \zeta(s))\), \(\tilde{w}(s)=\exp(y,
\tilde{\zeta}_\lambda(s))\) as in \S4.2.1. 
\item\(\gamma_{\chi}^{-1}(0)=\tilde{\gamma}_{\chi}^{-1}(0); \quad
\langle u_\gamma(0),\tilde{\gamma}_{\chi}^*\tilde{\xi}(0)\rangle_{2,t}=0\).
\end{itemize}
(\ref{calibrated}) is
in this case replaced by:
\begin{lemma*}\label{lemma:calibrate}
\(\forall s\in[\tilde{\gamma}_\chi^{-1}(R+1), \infty)\),
\begin{equation}\label{calibrated2}
\|\tilde{\xi}(s)\|_{2,2,t}+\|\tilde{\xi}'(s)\|_{2,1,t}\leq C(|\lambda|+\|\Pi_{{\bf
    e}_y}(\tilde{\zeta}_\lambda(s)-\tilde{\zeta}_\lambda(\infty))\|_{2,t}^3)\|\Pi_{{\bf
    e}_y}(\tilde{\zeta}_\lambda(s)-\tilde{\zeta}_\lambda(\infty))\|_{2,t}.
\end{equation}
\end{lemma*}
\begin{proof}
Write \(u(\tilde{\gamma}_\chi(s))=\exp(y, \tilde{\zeta}(s))\), and let
\(b(s):=\Pi_{{\bf e}_y}\tilde{\zeta}(s)\),
\(c(s):=\zeta(s)-\tilde{\zeta}(s)\). 
Note that \(\Pi_{{\bf e}_y}\tilde{\zeta}=\Pi_{{\bf
    e}_y}\tilde{\zeta}_\lambda\), and on this region
\(\tilde{\zeta}_\lambda-\tilde{\zeta}=\eta_{\lambda -}\) \(\forall s\).
The functions \(b(s), c(s)\) still satisfy (\ref{flow:b}),
(\ref{C:flow}).
However, we want to estimate instead
\[
c_d(s):=\zeta(s)-\tilde{\zeta}_\lambda(s)=c(s)-\Pi_{\ker A_y}^{\perp}\eta_{\lambda -}:
\]
From the definitions, estimates for \(c_d\) would imply similar
estimates for \(\tilde{\xi}\).

Let \(\underline{b}_d(s):=\langle {\bf e}_y,
\tilde{\zeta}_\lambda(s)-\tilde{\zeta}_\lambda(\infty)\rangle_{2,t}\);
\(b_d(s)=\underline{b}_d(s) {\bf e}_y\). Noting that 
\[\begin{split}
&-A_y\Pi_{\ker A_y}^{\perp}\eta_{\lambda -}\\
&\quad =(1-\Pi_{\ker A_y}-\nabla_{b(\infty)}Z\Pi_{\ker A_y})\Big(\hat{n}_{(0,y)}(\lambda,\tilde{\zeta}(\infty)+c(\infty))-n_y(\tilde{\zeta}(\infty))\Big)-\nabla_{b(\infty)}Z(\lambda C'_y{\bf e}_y),
\end{split}\]
we see that (\ref{C:flow}) may be rewritten in terms of \(c_d\) as:
\begin{equation}\label{eqn:c-d}
\begin{split}
-c_d'=&A_yc_d\\
&\quad+(1-\Pi_{\ker A_y}-\nabla_{b}Z\Pi_{\ker
  A_y})\Big(\hat{n}_{(0,y)}(\lambda,\tilde{\zeta}_\lambda+c_d)-n_y(\tilde{\zeta})-\hat{n}_{(0,y)}(\lambda,\tilde{\zeta}_\lambda(\infty))+n_y(\tilde{\zeta}(\infty))\Big)\\
&\quad -\nabla_{b_d}Z\Big(\lambda C'_y{\bf e}_y+\Pi_{\ker A_y}(\hat{n}_{(0,y)}(\lambda,\tilde{\zeta}_\lambda(\infty))-n_y(\tilde{\zeta}(\infty)))
\Big).
\end{split}
\end{equation}
By the nature of \(u\), \(u_\lambda\), and \(\hat{n}_{(0,y)}\), this leads to the
familiar estimates:
\begin{gather}
\|c_{d+}\|_{2,t}'\geq\nu_+\|c_{d+}\|_{2,t}-\epsilon_+\|c_d\|_{2,t}-C_+|\lambda\underline{b}_d|;\nonumber\\
\|c_{d-}\|_{2,t}'\leq -\nu_-\|c_{d-}\|_{2,t}+\epsilon_-\|c_d\|_{2,t}+C_-|\lambda\underline{b}_d|.\label{c-d-}
\end{gather}
Subtracting the two inequalities, we get
\[
(\|c_{d+}\|_{2,t}-\|c_{d-}\|_{2,t})'\geq\nu'(\|c_{d+}\|_{2,t}-\|c_{d-}\|_{2,t})-C'|\lambda\underline{b}_d|;
\]
Taking convolution product with the integral kernel of \(d/ds-\nu'\), we find that for \(s\geq s_0\)
\[
\|c_{d+}\|_{2,t}(s)-\|c_{d-}\|_{2,t}(s)\geq
(\|c_{d+}\|_{2,t}(s_0)-\|c_{d-}\|_{2,t}(s_0))e^{\nu'
  (s-s_0)}-C'\int_{s_0}^s|\lambda\underline{b}_d(\underline{s})|e^{\nu'(s-\underline{s})}\, d\underline{s},
\]
and since \(\underline{b}_d(s)>0\) decreases with \(s\), this implies
that for all large enough \(s\), 
\begin{equation}\label{c-d-plus}
\|c_{d+}\|_{2,t}(s)\leq \|c_{d-}\|_{2,t}(s)+C''|\lambda\underline{b}_d(s)|,
\end{equation}
otherwise \(\|c_{d+}\|_{2,t}(s)-\|c_{d-}\|_{2,t}(s)\) would be growing
exponentially as \(s\to\infty\), contradicting the fact that by
construction, \(\lim_{s\to\infty}\|c_d(s)\|_{2,t}=0\).

Plugging in this back to (\ref{c-d-}), we get
\begin{equation}\label{c-d-0}
\|c_{d-}\|_{2,t}'\leq -\nu_-'\|c_{d-}\|_{2,t}+C'_-|\lambda\underline{b}_d|,\end{equation}
where \(\nu_-'\) is a positive numbers close to \(\nu_-\). 
Taking convolution product with the integral kernel of \(d/ds+\nu'\),
\begin{equation}\label{c-d-1}
\|c_{d-}(s)\|_{2,t}(s)\leq
C_0e^{-\nu_-'s}+\int_{s_0}^s|\lambda|\underline{b}_d(\underline{s})e^{\nu_-'(\underline{s}-s)}\,
d\underline{s}.
\end{equation}
We claim that there is a positive constant \(\nu''_-\) slightly smaller than \(\nu'_-\)
such that
\begin{equation}\label{claim:b}
\underline{b}_d(\underline{s})\leq
2\underline{b}_d(s)e^{\nu_-''(s-\underline{s})/4}\quad \text{for
  \(s_0\leq \underline{s}\leq s\).}
\end{equation}
Using this in the integrand in (\ref{c-d-1}), we arrive at 
\[\begin{split}
\|c_{d-}(s)\|_{2,t}(s)&\leq
C_0e^{-\nu_-'s}+C_0'|\lambda|\underline{b}_d(s)\\
&\leq
C_d\underline{b}_d(s)(\underline{b}^3_d(s)+|\lambda|).
\end{split}
\]
(In the second step above we used (\ref{claim:b}) again to bound
\[e^{-\nu_-'s}=(e^{-\nu_-''s/4})^4\leq
C_8\underline{b}_d(s)^4. \quad \text{for large \(s\).})\] 
Combining with (\ref{c-d-plus}), we obtain a similar estimate for
\(\|c_d(s)\|_{2,t}\): 
\begin{eqnarray}
\|c_{d}(s)\|_{2,t}(s)&\leq
C_1e^{-\nu_-'s}+C_1'|\lambda|\underline{b}_d(s)\label{c-d-first}\\
&\leq
C_d'\underline{b}_d(s)(\underline{b}^3_d(s)+|\lambda|).
\end{eqnarray}
We now returns to verify the claim (\ref{claim:b}). To see this, note
that projecting the flow equation to \(\ker A_y\), we have
\[
-b_d'=\Pi_{\ker A_y}\Big(\hat{n}_{(0,y)}(\lambda,\tilde{\zeta}_\lambda+c_d)-\hat{n}_{(0,y)}(\lambda,\tilde{\zeta}_\lambda(\infty))\Big)
\]
Then the properties of \(\hat{n}_{(0,y)}\), \(\tilde{\zeta}_\lambda(\infty)\) and \(u\) again give the
estimate:
\begin{equation}\label{b-d0}
\underline{b}_d'\geq-\varepsilon''(\underline{b}_d+\|c_d\|_{2,t})
\end{equation}
for a small positive constant \(\varepsilon''\).
Subtracting a small multiple of this from (\ref{c-d-0}) and using (\ref{c-d-plus}), we have
\begin{equation*}\label{eqn:b-d}
(\|c_{d-}\|_{2,t}-\varepsilon_1\underline{b}_d)'\leq -\nu''(\|c_{d-}\|_{2,t}-\varepsilon_1\underline{b}_d).
\end{equation*}
Taking convolution product with the integral kernel of \(d/ds+\nu''\),
we have
\[
\|c_{d-}\|_{2,t}\leq \varepsilon_1\underline{b}_d+C_1e^{-\nu''s}.
\]
Plug this back in (\ref{b-d0}) and (\ref{c-d-plus}), we get
\[
\underline{b}_d'\geq -\frac{\nu''}{4}\underline{b}_d-\varepsilon_2e^{-\nu''s}.
\]
Now taking convolution product with the integral kernel of
\(d/ds+\nu''/4\), we have for \(\underline{s}<s\):
\[
\underline{b}_d(s)\geq\underline{b}_d(\underline{s})e^{\nu''(\underline{s}-s)/4}-\frac{4\varepsilon_2}{3\nu''}(e^{-3\nu''\underline{s}/4}-e^{-3\nu''s/4})e^{-\nu''s/4}\geq\frac{1}{2}\underline{b}_d(\underline{s})e^{\nu''(\underline{s}-s)/4}.
\]
To obtain the second inequality above, first use the first inequality
and the fact that \(s>\underline{s}\geq s_0\gg1\) to obtain
\(\underline{b}_d(s)\geq Ce^{-\nu''s}\); then use this (with \(s\)
replaced by \(\underline{s}\)) to estimate 
\[
\frac{4\varepsilon_2}{3\nu''}(e^{-3\nu''\underline{s}/4}-e^{-3\nu''s/4})\leq
\frac{\underline{b}_d(\underline{s})}{2}e^{\nu''\underline{s}/4}.
\]
Claim verified.

Next, to get estimates for higher derivatives of \(c_d\) from
(\ref{c-d-first}) we need to elliptic bootstrap
using (\ref{eqn:c-d}) and apply Sobolev embedding as in the proof of
Lemma I.5.1.7.To obtain the estimates claimed in the Lemma, we
need to bound the average of \(\underline{b}_d\) in an interval about
\(s\) in terms of \(\underline{b}_d(s)\). This is obtained using
(\ref{claim:b}) and the fact that \(\underline{b}_d\) is decreasing.
\end{proof}

Next we compare \(\tilde{w}(s)\) with the pregluing
\(w_\chi(s)\) to get a pointwise estimate of \(\xi(s)\), as in \S4.2.3. In this case,
the first two formulas of (\ref{Delta-s}) are still valid (with
\(\mathfrak{r}_i\) there replaced by \(\mathfrak{r}_\chi\)) by arguments similar to
those in \S4.2, but the third needs to be modified. 
In this region (where \(\gamma\geq \mathfrak{r}_\chi\)), we need to 
expand about \(y_{\lambda}\) instead of \(u_{\lambda}\) as in the
proof of Lemma \ref{lemma:gamma2}, keeping in mind that \(\mu_\lambda\) is of order \(\lambda^{1/2}\) while
\((\mu_\lambda)_\gamma\) is of order \(\lambda\). Recall that
\(\gamma_\chi\) satisfies the equation
\(\gamma_\chi'=h(\gamma_\chi)\), with \(h\) given by
(\ref{h-lambda}). The function \(\tilde{\gamma}_\chi\)
satisfies a similar equation:
\[\tilde{\gamma}_\chi'=h(\tilde{\gamma}_\chi)+\|(\mu_{\lambda})_{\gamma}(\tilde{\gamma}_\chi)\|_{2,t}^{-2}\Big\langle(\mu_{\lambda})_\gamma(\tilde{\gamma}_\chi),
E_{y_\lambda}T_{w_\chi
  ,y_\lambda}\tilde{\xi}(s)+o(\|\tilde{\xi}(s)\|_{2,1,t})\Big\rangle_{2,t}.\] 
By (\ref{calibrated2}) and (\ref{new-calib}), the absolute value of
this can be bounded by
\[\begin{split}
&C_1|\lambda|^{-1/2}(\|\mu_\lambda(\tilde{\gamma}_\chi)\|_{2,1,t}^3+|\lambda|)\|\mu_\lambda(\tilde{\gamma}_\chi)\|_{2,1,t}\\
&\quad \leq
C_2|\lambda|^{3/2}(\breve{\gamma}_\chi-\tilde{\gamma}_\chi)(1+\lambda^2(\breve{\gamma}_\chi-\tilde{\gamma}_\chi)^3).\end{split}\]
Recall also the estimate for \(h\) from \S\ref{lemma:gamma2}; we then obtain
\[\begin{split}
|\Delta_s(\gamma)-\Delta_s(\mathfrak{r}_\chi)|&\leq C_3\int_{\mathfrak{r}_\chi}^\gamma
\Big|\frac{|\lambda|^{3/2}(\breve{\gamma}_\chi-\gamma)(1+\lambda^2(\breve{\gamma}_\chi-\gamma)^3)}{|\lambda|(\breve{\gamma}_\chi-\gamma)^2}\,
d\gamma\Big|\\
&\leq
C_4|\lambda|^{1/2}\Big(\lambda^2(\breve{\gamma}_\chi-\gamma)^2+|\ln
(\breve{\gamma}_\chi-\gamma)|\Big)\Big|_{\mathfrak{r}_\chi}^\gamma.
\end{split}\]
Using this and the facts that in this region
\begin{equation}\label{delta-gamma:est2}
\begin{split}
\|w_\chi'\|_{2,2,t}\leq
C_5|\lambda|e^{-C_6|\lambda|^{1/2}(s-\gamma_\chi^{-1}(\mathfrak{r}_\chi))}\quad
\text{and}\\
\breve{\gamma}_\chi-\gamma_\chi(s)\leq
C_5'|\lambda|^{-1/2}e^{-C_6|\lambda|^{1/2}(s-\gamma_\chi^{-1}(\mathfrak{r}_\chi))},
\end{split}
\end{equation}
we can bound \[\|\tilde{\xi}_\chi\|_{2,2,t}\leq
\varepsilon_7|\lambda|e^{-C'_6|\lambda|^{1/2}(s-\gamma_\chi^{-1}(\mathfrak{r}_\chi))},\]
where \(\tilde{\xi}_\chi\) is defined by \(\tilde{w}(s)=\exp (w_\chi(s),
\tilde{\xi}_\chi(s))\) as in \S4.2. Recall also that \(w(s)=\exp
(w_\chi(s), \xi_\chi(s))\).
Combining the above estimate with (\ref{calibrated2}) and the other two lines of
(\ref{Delta-s}), we have
\begin{equation}\label{xi-small2}
\|\xi_\chi(s)\|_{2,2,t}\leq \begin{cases}
\varepsilon_7|\lambda|e^{-C'_6|\lambda|^{1/2}(s-\gamma_\chi^{-1}(\mathfrak{r}_\chi))} &\text{when \(s\in
  [\gamma_\chi^{-1}(\mathfrak{r}_\chi), \infty)\)},\\
\varepsilon_8|\lambda| &\text{when \(s\) is between
  \(\gamma_{\chi}^{-1}(\mathfrak{r}_\chi)\) and \(\gamma_{\chi}^{-1}(
  R)\),}\\
\varepsilon(\lambda) &\text{otherwise.}
\end{cases}
\end{equation}
So over \((-\infty,\gamma_\chi^{-1}(\mathfrak{r}_\chi)]\times S^1\subset
\Theta \), we can estimate \(\|\xi_\chi\|_{W_\chi}\) by the same argument
as in \S4.2. The estimate over \(\Theta_{y-}:=[\gamma_\chi^{-1}(\mathfrak{r}_\chi),
\infty)\times S^1\) is replaced by the following:
By (\ref{xi-small2}), 
\begin{equation}\label{W-bound:xi-T}
|\lambda|^{-1/2} \|\xi_{\chi}\|_{L^p(\Theta_{y-})}+|\lambda|^{-1/2} \|\dot{\xi}_{\chi}\|_{L^p(\Theta_{y-})}\leq\varepsilon_9|\lambda|^{1/2-1/(2p)}\ll1. 
\end{equation}
Next, to estimate \(\xi_\chi'\), it is equivalent to estimate
\(\xi_{\chi, y-}'\), which is obtained by expanding the flow equation
about \(y_{\lambda-}\): Here we have an equation similar to (\ref{split:pcr}),
with \(\bar{y}\) replaced by \(\bar{y}_{\lambda-}\), and
\(\alpha=0\). Using the error estimate in this region in the proof of
Lemma \ref{err-est2}, (\ref{delta-gamma:est2}), Lemma I.5.3.2,
we find
\[
|\lambda|^{-1/2} \|(\xi_{\chi, y-})'_{T}\|_{L^p(\Theta_{y-})}
|\lambda|^{-1} \|(\xi_{\chi, y-})'_{L}\|_{L^p(\Theta_{y-})}\leq C_9|\lambda|^{1/2-1/(2p)}\ll1. 
\]
This together with (\ref{W-bound:xi-T}) shows that
\(\|\xi_\chi\|_{W_\chi(\Theta_{y-})}\ll1\). Now one may follow the argument in
\S1.2.7 to complete Step 4 of the proof of the gluing theorem.

\subsection{When $u=\bar{y}$.}

We now assume that \(x_0=z_0=y\), and \(u=\bar{y}\), the constant flow
at \(y\).

\mysubsubsection[The pregluing.]
Let $\underline{b}^0_{\lambda}(s)$ be the solution of 
\begin{equation}\label{equation:b-0}
-(\underline{b}^{0}_{\lambda})'=C_y'\lambda +C_y
(\underline{b}^0_{\lambda})^2, \quad \text{
with $\underline{b}^0_{\lambda}(0)=0$,}\end{equation}
where \(C_y, C_y'\) are as defined in I.5.3.1. In other words, there
are positive \(\lambda\)-independent constants \(C_0, C'\), such that
\begin{equation}\label{b-tanh}
\underline{b}^0_{\lambda}(s)=C_0|\lambda|^{1/2}\tanh (C'|\lambda|^{1/2}s).
\end{equation}
Let \(b^0_\lambda(s):=\underline{b}^0_{\lambda}(s){\bf e}_y\).
Denote by $\underline{b}^{0\pm}_{\lambda}=\lim_{s\to\pm\infty}\underline{b}^0_{\lambda}(s)=\pm
C_0|\lambda|^{1/2}$. Let 
\[
\tilde{\zeta}_{\lambda}:=b^0_{\lambda}+\beta_+(\eta_{\lambda+}-b^{0+}_{\lambda})+\beta_-(\eta_{\lambda-}-b^{0-}_{\lambda}),
\]
where $\eta_{\lambda\pm}$ are defined by $\exp
(y,\eta_{\lambda\pm})=y_{\lambda\pm}$, and 
$\beta_-(s):=\beta (|\lambda|^{-1}+s)$; $\beta_+(s):=\beta(|\lambda|^{-1}-s)$. 

We define the {\em pregluing} \(w_\lambda\) in this case by
\[
w_{\lambda}:=\exp(y, \tilde{\zeta}_{\lambda}),
\]

\mysubsubsection[The weighted norms.] 
Recall the definition of \(W_{y-}\), \(L_{y-}\) from
\S\ref{claimF2}. Let \(W_{y+}\), \(L_{y+}\) and \(W_{y}\), \(L_{y}\)
be similarly defined for elements in \(\Gamma(\bar{y}_{\lambda +}^*K)\) and
\(\Gamma(\bar{y}^*K)\) respectively, with longitudinal
directions given by \(\mathbf{e}_{y_{\lambda+}}\) and
\(\mathbf{e}_y\). 

Via the map \(T_{y, w_\lambda}: \Gamma(\bar{y}^*K)\to
\Gamma((w_\lambda)^*K)\), the norms \(W_{y}\), \(L_{y}\) on
\(\Gamma(\bar{y}^*K)\) induce norms on \(\Gamma((w_\lambda)^*K)\),
which we denote by $W_{\lambda}$, $L_{\lambda}$. The associated spaces shall be the
domain and range for \(E_{w_\lambda}\).

By the estimates for \(\eta_{\lambda \pm}\), it is easy to see that
the induced norms on \(\Gamma((w_\lambda)^*K)\) via
\(T_{y_{\lambda\pm}, w_\lambda}\) from \(W_{y\pm}\), \(L_{y\pm}\) are
commensurate with \(W_\lambda, L_\lambda\).

\mysubsubsection[Error estimates.]
Divide $\Theta$ into three regions: $\Theta_a$, $\Theta_b$, $\Theta_c$
corresponding to
$s<-|\lambda|$, $|s|<|\lambda|$ and $s>|\lambda|$ respectively. We
will expand $\bar{\p}_{J X_{\lambda}}w_{\lambda}$ around $y_{\lambda+}, y, y_{\lambda-}$
respectively in the three regions. 

{\it Over $\Theta_b$}, using (\ref{equation:b-0}) and the fact of
\(y\) being in a standard d-b neighborhood, we have
\begin{eqnarray}
\lefteqn{(T_{y, w_\lambda})^{-1}\bar{\p}_{J X_{\lambda}}w_{\lambda}}\nonumber\\
&&=E_y\tilde{\zeta}_{\lambda}+n_y(\tilde{\zeta}_{\lambda})
+(T_{y, w_\lambda})^{-1}\delta_\lambda {\cal V}(w_{\lambda})\nonumber\\
&&=-\beta'_{+}(\eta_{\lambda+}-b^{0+}_{\lambda})+\beta'_{-}(\eta_{\lambda -}-b^{0-}_{\lambda})+\beta_+A_y(\eta_{\lambda +}-b^{0+}_{\lambda})
+\beta_-A_y(\eta_{\lambda -}-b^{0-}_{\lambda})\nonumber\\
&&\qquad +C_y\langle {\bf e}_y, 2b_{\lambda}^0
+\delta\rangle_{2,t}\delta+\Pi_{{\bf e}_y}^{\perp}O((\|\tilde{\zeta}_\lambda\|_{2,1,t}+|\lambda|)^2),\label{Theta-b:est}
\end{eqnarray}
where
$\delta:=\beta_+(\eta_{\lambda +}-b^{0+}_{\lambda})+\beta_-(\eta_{\lambda -}-b^{0-}_{\lambda})$.

From the estimates for \(\eta_{\lambda\pm}\) in the proof of Lemma
I.5.3.2, one sees that                \begin{equation}\label{delta-est}
|\lambda|^{-1}\|\delta_L\|_{2,1,t}+|\lambda|^{-1/2}\|\delta_T\|_{2,1,t}\leq
C|\lambda|^{1/2}\quad \forall s,
\end{equation} 
and therefore from (\ref{Theta-b:est})
\[
\|\bar{\p}_{J X_{\lambda}}(w_{\lambda})\|_{L_{\lambda}(\Theta_b)}
\leq C|\lambda|^{1/2-1/p}.
\]

The estimates on $\Theta_a$ and $\Theta_c$ are similar; so we shall
focus on $\Theta_c$. In this region, writing
$w_{\lambda}(s)=\exp(y_{\lambda+},\mu_{\lambda+}(s))$, we have
\begin{equation}\label{Theta-c:est}
(T_{y, w_\lambda})^{-1}\bar{\p}_{J X_{\lambda}}w_{\lambda}
=E_{y_{\lambda+}}\mu_{\lambda+}+n_{y_{\lambda+}}(\mu_{\lambda+}).
\end{equation}
From the definition of $\mu_{\lambda+}$,
$\|\mu_{\lambda+}(s)\|_{p,1,t}\leq
C\|b^{0+}_{\lambda}-b^0_{\lambda}(s)\|_{p,1,t}$.
So by (\ref{Theta-c:est}), (\ref{b-tanh}), 
\[
\|\bar{\p}_{J X_{\lambda}}(w_{\lambda})\|_{L_{\lambda}(\Theta_c)}
\leq C|\lambda|^{-1/2-1/p}e^{-C'|\lambda|^{-1}}\to 0 \quad \text{as $\lambda\to 0$.}
\]

\mysubsubsection[Existence and uniform boundedness of right inverse to
\(E_{w_\lambda}: W_\lambda\to L_\lambda\).]
In this case, let
\[
W'_\lambda:=\{\xi\in W_\lambda |\, \xi_L(0)=0\}.
\]
Again assume the existence of a sequence \(\{\xi_\lambda\in
W_\lambda'\}_\lambda\) satisfying (\ref{As:xi2}), with the obvious modification.

Divide $\Theta$ into three regions $\Theta_y$, $\Theta_{y\pm}$
corresponding respectively to the three possibilities: 
$|s|\leq|\lambda|^{-1/2}$, \(\pm s\geq |\lambda|^{-1/2}\). Let $\Theta'_y\supset
\Theta_y$ be the region in which
$|s|<(1+\varepsilon)|\lambda|^{-1/2}$; let
$\Theta'_{y\pm}\supset\Theta_{y\pm}$ be the region in which $\pm
s>(1-\varepsilon)|\lambda|^{-1/2}$, where \(1\gg\varepsilon>0\). 
Instead of estimating \(\|\xi_\lambda\|_{W_\lambda}\), we shall
estimate
\[\xi_{\lambda, y}:=T_{w_\lambda, y}\xi_\lambda, \quad \xi_{\lambda, y\pm}:=T_{w_\lambda, y\pm}\xi_\lambda
\]
in \(W_y\) or \(W_{y\pm}\) norm over \(\Theta_y\) and
\(\Theta_{y\pm}\) respectively.
First, observe the following  analog of Lemma \ref{claimF2} over \(\Theta_y\) and
\(\Theta_{y\pm}\).
\begin{lemma*}[Analog of Floer's Lemma]\label{claimF3}
Let $\xi_{\lambda}$ be as in (\ref{As:xi2}) with \(W_\chi, W_\chi',
L_\chi\) replaced by \(W_\lambda, W_\lambda', L_\lambda\) respectively. Then 
for all sufficiently small
$\lambda$, there is a small positive number, \(\varepsilon_0(\lambda)\), 
\(\lim_{\lambda\to 0}\varepsilon_0(\lambda)=0\), such that
\[\begin{split}
(a)\, & |\lambda|^{1/(2p)}\|\xi_{\lambda , y}\|_{L^\infty( \Theta'_y)}+
\|(\xi_{\lambda , y})_L\|_{L^\infty(\Theta'_y)}\leq \varepsilon_0(\lambda)|\lambda|^{1/2+1/(2p)};\\
(b) \,& |\lambda|^{1/(2p)}\|\xi_{\lambda, y\pm}\|_{L^\infty(\Theta'_{y\pm})}+
\|(\xi_{\lambda , y\pm})_L\|_{L^\infty(\Theta'_{y\pm})}\leq \varepsilon_0(\lambda)|\lambda|^{1/2+1/(2p)}.
\end{split}\] 
\end{lemma*}
\begin{proof}
The \(L^{\infty}_t\)-estimate for \(\xi_\lambda\) (and hence
\(\xi_{\lambda, y}\) and \(\xi_{\lambda, y\pm}\)) is now routine. The estimates for the longitudinal components
follow the rescaling argument in the proof of Lemma \ref{claimF2},
with the following modifications: 

On \(\Theta_y'\), one may similarly define a sequence
\(\tilde{\varsigma}_\lambda\) of \(L^p_1\)-bounded functions on
\([-1,1]\), which converges to \(\tilde{\varsigma}_0\) which satisfy also
an equation of the form (\ref{equation:tilde-zeta}), but now
\(\chi\sim C\tanh (C'\tau)\). Because \((\xi_\lambda)_L(0)=0\),
\(\tilde{\varsigma}_0(0)=0\), and thus \(\tilde{\varsigma}_0=0\). This proves
part (a) above.

On \(\Theta_{y\pm}'\), we have another version of
\(\tilde{\varsigma}_\lambda\) and \(\tilde{\varsigma}_0\) (which are now
functions on \([1, \infty)\), \((-\infty, -1]\) respectively), and the argument in
the proof of Lemma \ref{claimF2} again gives a bound on
\(\|\tilde{\varsigma}_0\|_{\infty}\) by \(|\tilde{\varsigma}_0|(\pm1)\), which
vanishes by the estimate for \(\|(\xi_\lambda)_L\|_{L^\infty(\Theta_y')}\)
obtained in part (a). This proves part (b).
\end{proof}

We now return to estimate \(\|\xi_{\lambda, y}\|_{W_y}(\Theta_y)\) and
\(\|\xi_{\lambda, y\pm}\|_{W_{y\pm}}(\Theta_{y\pm})\).
\medbreak

{\it On $\Theta_y'$}: let $\chi_y$ be a smooth cutoff function which
vanishes outside $(-1,1)$ and let
$\beta_y(s):=\chi_y(s/|\lambda|^{-1/2})$.
Since \((\xi_{\lambda, y})_L(0)=0\), applying Lemma \ref{fact} (c) to the
longitudinal component, one may estimate:
\begin{eqnarray}\label{xi-0:est3}
\lefteqn{\|\beta_y\xi_{\lambda, y}\|_{W_{y}(\Theta_y)}\leq 
C\|E_{\bar{y}}(\beta_y \xi_{\lambda, y})\|_{L_{\lambda}(\Theta_y')}}\nonumber\\
&&\leq
C_1\|E_{w_\lambda}\xi_\lambda\|_{L_{\lambda}(\Theta_y')}+C_2(|\lambda|^{1/2}\|\beta_y\xi_{\lambda}\|_{W_{\lambda}(\Theta'_y)}+\|\beta_y(\xi_{\lambda,
y})_L\|_{W_{y}(\Theta'_y)})
+C_3\|\beta_y'\xi_\lambda\|_{L_{\lambda}(\Theta'_y)};\nonumber\\
&&\leq C_1\varepsilon_E+C_4(|\lambda|^{1/2}+\varepsilon_0).
\end{eqnarray}
In the above, the second term in the penultimate line came from the difference between $E_{w_\chi}$
and the conjugate of $E_y$, using the fact that
$\|\tilde{\zeta}_{\lambda}\|_{\infty ,1,t}\leq C|\lambda|^{1/2}$.
The last line used Lemma \ref{claimF3} (a) and the equation for
\((\xi_{\lambda, y})_L\).
\medbreak

{\it On $\Theta_{y\pm}'$} one may estimate similarly. Let
\(\beta_{y\pm}\) be smooth cutoff functions supported on
\(\Theta_{y\pm}'\) with value 1 over \(\Theta_{y\pm}\) and
\(|\beta_{y\pm}'|\leq C|\lambda|^{1/2}\). 
By the eigenvalue estimate for \(A_{y_{\lambda\pm}}\) in I.5.3.2,
\[
\|\beta_{y\pm}\xi_{\lambda, y\pm}\|_{W_{y\pm}}\leq 
C\|E_{y_{\lambda\pm}}(\beta_{y\pm}\xi_{\lambda, y\pm})\|_{L_{\pm}}.
\]
The RHS can be estimated like (\ref{longi}, \ref{transv}) using Lemma
\ref{claimF3} (b) below.

Finally, from the estimates for $\xi_{\lambda, y}$ and $\xi_{\lambda,
  y\pm}$ above we obtain the
desired contradiction that $\|\xi_\lambda\|_{W_{\lambda}}\ll1$.

\mysubsubsection[Surjectivity of the gluing map.] We have the routine
estimate for the nonlinear term to define the gluing map. The main
isssue is again to show that the gluing map surjects to a neighborhood
of \(\mathbb{S}\) in the parameterized moduli space of broken
trajectories.

Let \((\lambda, \hat{w})\in \hat{\cm}_P^{\Lambda, 1}({\bf y}_{+}, {\bf
  y}_{-})\) be in a chain topology neighborhood of \(\bar{y}\in
\hat{\cm}_P^{\Lambda, 1, +}\). Choose a representative \(w\) of
\(\hat{w}\) such that, writing 
\[
w(s)=\exp(y, \zeta(s)), \quad w_\lambda(s)=\exp(y,
\tilde{\zeta}(s)),
\] 
the difference \(\eta:=\zeta-\tilde{\zeta}\) satisfies
\(\eta_L(0)=0\). Writing \(w(s)=\exp(w_\lambda(s), \xi(s))\),
we want to show that \(\|\xi\|_{W_\lambda}\ll1\);
equivalently, it suffices to estimate \(\eta\).

Let \(\epsilon'\gg|\lambda|^{1/2}\) be a small positive number Consider the three regions \(\Theta_a^{\epsilon'}=(-\infty,
-\epsilon'|\lambda|^{-1}]\times S^1\),
\(\Theta_b^{\epsilon'}=[-\epsilon'|\lambda|^{-1}, \epsilon'|\lambda|^{-1}]\times
S^1\), \(\Theta_c^{\epsilon'}=[\epsilon'|\lambda|^{-1}, \infty)\times S^1\) separately.

From the flow equation and the definition of \(w_\lambda\), we find
that \(\eta_L, \eta_T\) satisfy respectively:
\begin{gather}
\underline{\eta}_L'+C_y(\underline{\tilde{\zeta}}_L+\underline{\zeta}_L)\underline{\eta}_L=O\Big((|\lambda|+\|\eta_T(s)\|_{2,1,t})^2+\|\zeta_L(s)\|_{2,1,t}^4\Big);\label{equation:eta-L}\\
\eta_T'(s)+A_y\eta_T(s)=O((|\lambda|+\|\zeta\|_{2,1,t})^2).\label{equation:eta-T}
\end{gather}
(In the usual notation, \(\eta_L=:\underline{\eta}_L{\bf e}_y\);
\(\zeta_L=:\underline{\zeta}_L{\bf e}_y\).)

Equation (\ref{equation:eta-T}) and the fact that
\(\eta(\infty)=\eta(-\infty)=0\) imply:
\begin{equation}\label{est:eta-T}
\|\eta_T\|_{2,2,t}\leq
C_1(|\lambda|+\|\zeta_L\|_{\infty})^2\quad\forall s.
\end{equation}

The argument to get this estimate should be by-now familiar to the
reader (cf. e.g. the proof of Lemma
\ref{lemma:calibrate}): Take \(L^2_t\)-inner product of
(\ref{equation:eta-T}) with \(\eta_{T-}\), \(\eta_{T+}\) respectively,
and integrate over \(s\), one obtains
\[
\|\eta_T\|_{2}\leq C_2(|\lambda|+\sup_s\|\zeta_L\|_{2,t}(s))^2.
\]
Then apply the usual elliptic bootstrapping and Sobolev embedding to
get estimates on higher derivatives. Finally, observe that on the
1-dimensional subspace of longitudal direction, the various norms are
all commensurate. 
 
On the other hand, \(\zeta_L\) satisfies 
\[
\underline{\zeta}_L'+C_y\underline{\zeta}_L^2+C'_y\lambda=O\Big((|\lambda|+\|\eta_T\|_{2,1,t}+\|\delta_T\|_{2,1,t}+\|\zeta\|_{2,1,t}^2)^2\Big),
\]
with \(\zeta_L'(\infty)=0=\zeta'_L(-\infty)\) (so when
\(|\underline{\zeta}_L|(s)\) reaches maximum,
\(\underline{\zeta}'_L=0\)). Combining this with (\ref{est:eta-T})
and (\ref{delta-est}), we have \[|\zeta_L|_{\infty}\leq C_L|\lambda|^{1/2}.\]
Plugging this back in (\ref{est:eta-T}), we get
\begin{equation}\label{est:eta-T2}
\|\eta_T\|_{2,2,t}\leq C_1|\lambda|.
\end{equation}
Using these \(L^\infty\) estimates for \(\eta_T, \zeta_L\) and multiplying
(\ref{equation:eta-L}) with \(\underline{\eta}_L\), we
obtain
\[
-C'|\lambda|^2\leq|\underline{\eta}_L|'+C_y(\underline{\tilde{\zeta}}_L+\underline{\zeta}_L)|\underline{\eta}_L|\leq
C|\lambda|^2. \]
Now since \(\underline{\tilde{\zeta}}_L\), \(\underline{\zeta}_L\) are both
\(\geq 0\) when \(s\geq0\), and are both \(\leq0\) when \(s\leq0\), we
see that
\[
\begin{split}
|\underline{\eta}_L|'\leq C|\lambda|^2 &\quad \text{when \(s\geq0\);}\\
-|\underline{\eta}_L|'\leq C|\lambda|^2 &\quad \text{when \(s\leq0\).}
\end{split}
\]
Integrating using the initial condition that \(\eta_L(0)=0\), we see
that
\begin{equation}\label{est:eta-L}
\|\eta_L\|_{L^\infty(\Theta_b^{\epsilon'})}\leq C_1\epsilon'|\lambda|.
\end{equation}
Combining this with (\ref{est:eta-T2}) we get
\[
\|\xi\|_{W_\lambda(\Theta_b^{\epsilon'})}\leq
C_2\epsilon'^{1+1/(2p)}|\lambda|^{1/2-1/p}\to 0\quad \text{as
  \(\lambda\to 0\)}.
\]

We now turn to estimating \(\xi\) over \(\Theta_a^{\epsilon'},
\Theta_c^{\epsilon'}\). We shall only consider \(\Theta_c^{\epsilon'}\)
since the other works by analogy. On this region, writing
\(w_\lambda(s)=\exp(y_{\lambda -}, \tilde{\zeta}_{\lambda -}(s))\), we
have from the definition of \(w_\lambda\) that
\[
\|\tilde{\zeta}_{\lambda -}\|_{2,2,t}(s)\leq
C_3|\lambda|^{3/2}e^{-C_6|\lambda|^{1/2}(s-\epsilon'|\lambda|^{-1})}+\mbox{terms
  involving \(\delta\)}.
\]
We may ignore the terms involving \(\delta\) since by
(\ref{delta-est}), its contribution to the \(W\)-norm is at most of
order \(|\lambda|^{1/2}\). On the other hand, by standard estimates
\(w(s)\) has the same exponential decay behavior in this region,
and so combining the estimate for
\(\zeta_{\lambda-}(\epsilon'|\lambda|^{-1})\) above and
(\ref{est:eta-T2}), (\ref{est:eta-L}), we have on this region\[
\|\zeta_{\lambda-}\|_{2,2,t}(s)\leq
C_4\epsilon'|\lambda|e^{-C_6'|\lambda|^{1/2}(s-\epsilon'|\lambda|^{-1})},
\]
where \(\zeta_{\lambda-}\) is defined by \(\exp(y_{\lambda-},
\zeta_{\lambda -})=w(s)\).
Together with the previous expression, we obtain a pointwise estimate
for \(\xi(s)\) on this region, which, when combined with
(\ref{equation:eta-L}), (\ref{equation:eta-T}), yields
\[
\|\xi\|_{W_\lambda(\Theta_c^{\epsilon'})}\leq C_5\epsilon'|\lambda|^{1/2-1/(2p)}\ll1.
\]

\section{The Handleslide Bifurcation.}
The purpose of this section is to verify the bifurcation behavior at 
handle-slides predicted in I.4.3, namely, Propositions
\ref{handleslide}, \ref{prop:nep} below.

\subsection{Summary of Results.}

Combined with the previous gluing theorems: Propositions \ref{8.1}, \ref{thm:birth}, the following proposition completes the verification of (RHFS2c),
(RHFS3c) for admissible \((J,X)\)-homotopies.
\begin{*proposition}\label{handleslide}
Let \((J^\Lambda, X^\Lambda)\) be an admissible \((J, X)\)-homotopy
connecting two regular pairs, and \({\bf x}, {\bf z}\) be two path
components of \({\cal P}^\Lambda\backslash {\cal P}^{\Lambda, deg}\).
Then:
\begin{description}\itemsep -1pt
\item[(a)] a chain-topology neighborhood of \({\mathbb T}_{P, hs-s} ({\bf x}, {\bf
  z}; \Re)\) in \(\hat{\cm}_P^{\Lambda, 1, +}({\bf x}, {\bf z}; \op{wt}_{-{\cal Y},
  e_{\cal P}}\leq \Re)\) is lmb along \({\mathbb T}_{P, hs-s} ({\bf x}, {\bf
  z}; \Re)\);
\item[(b)] a chain-topology neighborhood of \({\mathbb T}_{O, hs-s}
  (\Re)\) in \(\hat{\cm}_O^{\Lambda, 1, +}(\op{wt}_{-{\cal Y},
  e_{\cal P}}\leq \Re)\) is lmb along \({\mathbb T}_{O, hs-s} (\Re)\).
\end{description}
\end{*proposition}
The proof follows the standard gluing construction outlined in \S1.2,
and shall be omitted. A description of the relevant K-models will be
given in \S7.3.2 and 7.3.3. A result analogous to part (a) above is also
given by \cite{floer.jdg} Proposition 4.2.

The rest of this section will be devoted to the proof of:
\begin{*proposition}\label{prop:nep}
Let \((J^\Lambda, X^\Lambda)\) be an admissible \((J, X)\)-homotopy,
and \(u\in \hat{\cm}_P^{\Lambda, 0}({\bf x}, {\bf x})\). 
Then (NEP) holds for \(u\).
\end{*proposition}
Without loss of generality, we restrict our attention to a
\(J|X\)-homotopy without death-birth bifurcations 
throughout this section.

\subsection{Nonequivariant Perturbations on Finite-cyclic Covers.}

This subsection contains the main body of the proof of Proposition
\ref{prop:nep}.
We first discuss a simpler situation in which the
non-equivariant perturbation may be obtained from a vector field on a
finite cyclic covering of \(M\). In general, we need to resort to
non-local perturbations.

\subsubsection{A special case: Local perturbations from
finite-cyclic covers of \(M\).}

If a finite-cyclic cover \(\lc^{\nu, m}\to \lc\) is (a path component
of) the pull-back
bundle of a finite-cyclic cover
\(\hat{M}\to M\) via \(e_f: \lc\to M\) (cf. I.3.1.1), 
then a non-equivariant function or
vector field on \(\hat{M}\) may induce a 
non-equivariant function or
vector field on \(\lc^{\nu_\lc, m}\).

\begin{example*}
Assume the conditions of Corollary I.2.2.5 (namely, \(M\) is monotone,
\(f\) is symplectic isotopic to \(\op{id}\), and \(\gamma_0\) is the
trace of a point under the symplectic isotopy). We claim that in this
case, for any \(m\in \Z^+\) not dividing \(\op{div}([u])\),
there exists a \(u\)-breaking \(m\)-cyclic cover of \(\lc\)
via the above pull-back construction. Thus, in this case Proposition
\ref{prop:nep} may be proven by
simply repeating the argument for Proposition I.6.2.2 for
non-equivariant Hamiltonian perturbations over finite-cyclic covers of
\(M\). (In fact, only Lemma I.6.2.5 needs to be redone).

To see the claim, recall that in this case, 
\[{\frak H}=H_1(\lc;\Z)=\pi_2(M)\oplus
H_1(M;\Z), \quad \text{and} \quad e_{f*}=0\oplus \op{id}\] 
with respect to this decomposition. 
Notice that \(e_{f*}([u])\) is a non-torsion element in
\(H_1(M;\Z)\). Otherwise, by the commutative diagram from I.(12), 
\[k[u]=b\in \ker c_1\Big|_{\pi_2(M)}\quad \text{for some \(k\in
  \Z^+\).}\]  
But then by monotonicity of \(M\),
\[
[{\cal Y}_X](k[u])=\omega(b)-e_f^*\theta_X(k[u])=0,
\]
contradicting the fact that \(u\) has positive
energy.

Thus, for any \(m\in \Z^+\) not dividing \(\op{div}([u])\), one
may simply set \(\nu_M\in H^1(M;\Z)\) to be a primitive class with
\(\nu_M(e_{f*}([u]))=\op{div}([u])\), and take
\(\lc^{\nu, m}=e_f^*M^{\nu_M, m}\). Furthermore, such a finite-cyclic
cover is always \({\frak H}\)-adapted, and \(u\)-breaking if
\(m\) does not divide \(\op{div}([u])\).
\end{example*}

However, this simple construction does not give all the
\(u\)-breaking finite covers we need.

\subsubsection{The general case: Non-local perturbations.}

Let \(\lc^{\nu, m}\) be a \(u\)-breaking, \({\frak
  H}\)-adapted \(m\)-cyclic cover of
\({\mathcal C}\) introduced in I.4.4.5.

We shall often make use of the following convenient description of
\(\lc^{\nu, m}\):
\[
\lc^{\nu, m}=\Big\{(z, [\mu])\, |\,  z\in \lc, \mu: [0,1]\to \lc; \mu(0)=\gamma_0, \mu(1)=z\Big\}/\sim,
\]
where $(z_1,[\mu_1])\sim(z_2,[\mu_2])$ iff \(z_1=z_2\)  and $\nu( [\mu_1-\mu_2])=0 \mod
m$. Such an equivalence class shall be denoted by a pair \((z, [\mu]_m)\).

Recall that \(\nu\in \op{Hom}({\frak H}, \Z)\). The fact that
\(\lc^{\nu, m}\) is \(u\)-breaking implies that \(\nu\) is
non-torsion.
Thus, one may find a class \(\nu_2\in H^2(T_f;\R)\) extending \(\nu\)
by linearity, that is, satisfying \(\nu_2((i_{\frak H}\ker\nu)\otimes
\R)=0\), and \(\nu_2(i_{\frak H}[u])=\op{div}(u)\) where \(i_{\frak
  H}: {\frak H}\hookrightarrow H_2(T_f;\Z)\) is the inclusion.
Let \(\omega_\nu\) be a smooth closed 2-form on \(T_f\) in the
cohomology class \(\nu_2\).

The 2-form \(\omega_\nu\) defines an \(\R\)-valued function
\(\Omega_\nu\) on \(\tilde{\lc}\), by setting 
\[
\Omega_\nu(z, [\mu]):=\int_{[0,1]\times S^1_1}\mu^*\omega_\nu.
\]
This induces an \(\R/m\Z\)-valued
function on \(\lc^{\nu,m}\), which we shall denote by the same notation.

\begin{definition*}[A class of nonlocal perturbations]
Let \(\chi:\R/m\Z\to \R\) be a smooth function, 
and let \(P\in {\cal H}\). We define the formal vector field \(\wp_{\chi P}\) on \(\lc^{\nu, m}\) by
\begin{equation}\label{add-pcr}
\wp_{\chi P}(z, [\mu]_m):=\chi(\Omega_\nu(z, [\mu]_m))\nabla P(z).
\end{equation}
For a path \((u(s), [\mu(s)]_m)\) in \(\lc^{\nu, m}\), let \[\bar{\p}_{JX}^{\chi
  P}(u,[\mu]_m):=\bar{\p}_{JX}u+\wp_{\chi P}(u, [\mu]_m).\]
Similarly, for a smooth function \(\chi^\Lambda: \Lambda\times
\R/m\Z\to \R\) and \(P^\Lambda\in {\cal H}^\Lambda\), one may define a
path of formal vector fields \(\{\wp_{\chi_\lambda
  P_\lambda}\}_{\lambda\in \Lambda}\) and the section
\(\bar{\p}_{J^\Lambda X^\Lambda}^{\chi^\Lambda
  P^\Lambda}\) on \({\cal B}_P^\Lambda\) or \({\cal B}_O^\Lambda\).
\end{definition*}
For the rest of this section, a ``{\em \(\chi P\)-perturbed flow}'' or
simply a ``{\em perturbed flow}'' will refer to a solution of \(\bar{\p}_{JX}^{\chi
  P}(u,[\mu]_m)=0\). One may define the moduli spaces of such flows,
\(\cm_{P;\nu,m}(J, X;\chi, P), \cm_{O;\nu, m}(J, X;\chi , P)\) etc.,
and their parameterized versions, in
the usual manner (cf. I.2.1.2, I.4.3.1). 
Notice that if one chooses \(P\in V^k_\delta(J, X)\) and
\(P^\Lambda\in V^{\Lambda; k, \kappa}_\delta(J^\Lambda, X^\Lambda)\), then 
\[{\cal P}(X; \chi,P)={\cal P}(X); \quad {\cal P}^\Lambda(X^\Lambda;
\chi^\Lambda,P^\Lambda)={\cal P}^\Lambda(X^\Lambda),\]
and in both equalities, the former is nondegenerate iff the latter is.
We shall show in the next subsection that in this case, when
\(\chi\), \(\chi^\Lambda\) are sufficiently small, and if \((J, X)\)
is regular and \((J^\Lambda, X^\Lambda)\) admissible, then the moduli
spaces of \(\chi P\)-perturbed flows and their parameterized versions 
satisfy all the usually expected regularity and
compactness properties, as described by (FS2), (FS3) and (RHFS2*), (RHFS3*).

\medbreak
\noindent{\it Proof of Proposition \ref{prop:nep}.} Let \(\Re\in \R^+\)
and \(\lc^{\nu, m}\) be fixed as in the statement of (NEP).
Without loss of generality, assume \(\Pi_\Lambda u=0\).

The admissible \((J,X)\)-homotopy \((J^\Lambda, X^\Lambda)\) induces a
homotopy of formal flows on \(\lc^{\nu, m}\), which satisfies all the
Properties listed in I.6.2.3 for admissibility, except for Property
(8) (injectivity of \(\Pi_\Lambda|_{\hat{\cm}^{\Lambda, 0}_P}\)): at
\(\lambda=0\), there are \(m\) distinct elements in
\(\hat{\cm}_P^{0}(J_\lambda, X_\lambda)\), which are precisely the
\(m\) different lifts of \(u\).

We write this induced homotopy of vector fields as \(\{{\cal V}^{\nu,
  m}(J_\lambda, X_\lambda)\}_{\lambda\in \Lambda}\).

To achieve Property (8), we shall consider homotopy of vector fields
on \(\lc^{\nu, m}\) of the form \[\{{\cal V}^{\nu,
  m}(J_\lambda, X_\lambda;\chi_\lambda, P_\lambda)\}_{\lambda\in
  \Lambda}:=\{{\cal V}^{\nu,
  m}(J_\lambda, X_\lambda)+\wp_{\chi _\lambda P_\lambda}\}_{\lambda\in
  \Lambda},\]
where \begin{equation}\label{cond:P1}P^\Lambda\in V^{\Lambda; k, \kappa}_\delta(J^\Lambda, X^\Lambda).\end{equation}
In fact, since \(\hat{\cm}_P^{\Lambda, 0}(J^\Lambda, X^\Lambda;
\op{wt}_{-\langle{\cal Y}\rangle, e_{\cal P}}\leq \Re)\) consists of finitely many
points, each projecting under \(\Pi_\Lambda\) to distinct values, we may
assume that \begin{equation}\label{cond:P2}P_\lambda=0\quad \text{for \(\lambda\in
  \Lambda\backslash S\)},\end{equation}
where \(S\) is a small interval about \(\Pi_\Lambda(u)=0\), so
that 
\[
S\cap (\Lambda_{db}\cup \Pi_\Lambda(\hat{\cm}_P^{\Lambda, 0}(J^\Lambda, X^\Lambda;
\op{wt}_{-\langle{\cal Y}\rangle, e_{\cal P}}\leq \Re)\backslash\{u\})=\emptyset.
\]
Such perturbed homotopy of formal flows might no longer be
co-directional, however, Properties (1)--(6) of I.6.2.3 are
preserved. Moreover, we shall see in the next subsection
that as long as \(\chi\) is sufficiently small in \(C_\epsilon\)-norm,
the parameterized moduli spaces remain \(\Re\)-regular
(i.e. \(\Re\)-truncated version of I.6.2.3 (7) holds).

We now describe an explicit choice of \(\chi^\Lambda\), \(P^\Lambda\)
among all those satisfying both (\ref{cond:P1}), (\ref{cond:P2}), so
that I.6.2.3 (8) may be achieved. For this purpose, the argument in
the proof of Lemma I.6.2.5 is revised as follows. 

Replace \(u_n\) there by \(u\), let \(B\) be a small
neighborhood in \(Q_1\cap Q_2\subset \R\times S^1\). Let
\(P_0\in V^k_\delta(J_0, X_0)\) be supported in a
small neighborhood \({\frak B}\subset T_f\), such that
\(u^{-1}({\frak B})\subset B\), similar to the definition of
\(\underline{H}_{\lambda_n}\) in I.6.2.5. Let \(P^\Lambda\) be
an extension of \(P_0\) satisfying (\ref{cond:P1}) and
(\ref{cond:P2}), which is in turn the analog of
\(\underline{H}^\Lambda\) in I.6.2.5. 

Let \(\tilde{u}_1, \ldots, \tilde{u}_m\) be the \(m\) distinct lifts
of \(u\) in \(\lc^{\nu, m}\). With the above choice of
\(P_0\), the perturbation \(\wp_{\chi _0 P_0}(\tilde{u}_i(s))\) is nontrivial only
when \(s\) is in the small interval \[I_B:=\op{pr}_1(B),\] 
where \(\op{pr}_1: \R\times S^1\to \R\) denotes the projection.
By construction, the values of
\(\Omega_\nu\) at different lifts of a point in \(\lc\) differ
by multiples of \(m\). Thus if \(I_B\) is sufficiently
small, the image \(\Omega_\nu(\tilde{u}_i(I_B))\) form disjoint intervals
in \(\R\) for different \(\tilde{u}_i\). We denote the interval
corresponding to \(\tilde{u}_i\) by \({\cal I}_i\), and 
choose \(\chi_0\) such that 
\[\chi_0(\phi)=C_i \quad \text{when \(\phi\in {\cal I}_i\), \(i=1, 2, \ldots, m\),}\] 
where \(C_i\) are distinct
  constants, and \(\chi_0\) is very small in \(C_\epsilon\)-norm. 
With this choice, \[\wp_{\chi_0 P_0}(\tilde{u}_i(s))=C_i\nabla P_0(\tilde{u}_i(s)), \]
and the analog of I.(71) now reads
\[E_{\tilde{u}_i}\xi+\alpha_i Y_{\tilde{u}_i}+C_i\nabla P_0(\tilde{u}_i)=0,\]
this, together with the contraction mapping theorem, shows that 
\(\tilde{u}_i\) perturbs into a $\tilde{v}_i$, so that 
\(\Pi_\Lambda\tilde{v}_i-\Pi_\Lambda u\)
are, up to higher order correction terms, proportional to
\(C_i\). Hence the perturbed flows have distinct values under
$\Pi_\Lambda$. 

As remarked before, the regularity 
of \(\Re\)-truncated moduli spaces is unaffected by
this perturbation, and thus \(\{{\cal V}^{\nu,
  m}(J_\lambda, X_\lambda; \chi _\lambda ,P_\lambda)\}_{\lambda\in
  \Lambda}\) satisfies all the \(\Re\)-truncated versions of I.6.2.3
(1)--(8). In particular, it has all the \(\Re\)-truncated versions of
the properties (RHFS*), except for (RHFS2c), (RHFS3c) and (RHFS4).

To see that the remaining properties also hold, we need to verify
that the gluing theorems proven in previous sections still
hold. (The arguments for (RHFS4) to appear in section 7 below depend
on the perturbation only through the existence of Fredholm theory, and
the linear gluing theorem \ref{glue-K}.)

By construction, \(S\cap \Lambda_{db}=\emptyset\).
Thus, no gluing for births or deaths (as in sections 2-5)
is necessary.

The proofs of standard gluing theorems such as Proposition
\ref{handleslide}, Lemma \ref{glue-K} do require updates. 
However, because of our choice of \(P^\Lambda\), \(\wp_{\chi_\lambda
  P_\lambda}\) vanishes near the critical point \(x_\lambda\). Thus,
we have the usual exponential decay of flows to critical points, and
the same error estimates. Only two facts need to be verified for the
standard arguments sketched in \S1.2 to go through:
\medbreak
(1) Fredholm theory and surjectivity of the perturbed version of
deformation operators \(E_{\tilde{v}}^{J_\lambda, X_\lambda;
  \chi_\lambda, P_\lambda}\), \(\hat{E}_{\tilde{u}}^{J_\lambda,
  X_\lambda; \chi_\lambda, P_\lambda}\), where \(\tilde{u}\),
\(\tilde{v}\) are the perturbed flows to be glued;

(2) The usual quadratic bound on the nonlinear term \(N_{w_\chi}\),
namely, (\ref{Bq}).
\medbreak

We shall verify these in the next subsection.\qed

\subsection{Properties of $\Re$-truncated Moduli Spaces of Perturbed Flows.}

The structure theory of the moduli spaces of such perturbed flows is
not covered in the literature, or in the discussion of Part I. We need
to start from scratch and check the foundation of this more general
theory. The major components of the expected structure theory are
examined in turn below.

We have already mentioned the following basic fact: 
\begin{*fact}[Exponential decay]\label{pertDecay}
A perturbed flow decays exponentially to a nondegenerate critical
point.
\end{*fact}
This is due to our choice that \(P^\Lambda\in V^{\Lambda;
  k,\kappa}_\delta (J^\Lambda, X^\Lambda)\), in particular, \(P_\lambda\)
vanishes up to \(k\)-th order at the critical points.
In fact, this also shows that
a perturbed flow decays polynomially to a good minimally degenerate
critical point as described in \S I.5. However, we do not
need this fact. 

\mysubsubsection[Fredholm theory.] 
Consider the linearization of \(\bar{\p}_{JX}^{\chi
  P}(u,[\mu]_m)\). We denote it by \(E^{J,X; \chi ,P}_{(u,[\mu]_m)}\) or
\(D^{J,X; \chi ,P}_{(u,[\mu]_m)}\), depending on whether
\((u,[\mu]_m)\) is an connecting flow line or an orbit. In addition to
the well-understood \(E^{J,X}_{(u,[\mu]_m)}(\xi)\) or
\(D^{J,X}_{(u,[\mu]_m)}(\xi)\), it has the following extra terms due to
\(\wp_{\chi P}\):
\begin{equation}\label{add:deformOp}
\begin{split}
&\chi(\Omega_\nu(u, [\mu]_m))\nabla_\xi\nabla P(u)\\
&\quad
+\chi'(\Omega_\nu(u, [\mu]_m))\int_{S^1}\omega_\nu(\xi,
\p_t u)\, dt\, \nabla P(u).
\end{split}\end{equation}
Observing that the first term is a 0-th order multiplicative operator,
and the second term is a mixture which is infinitely smoothing in
\(t\) and 0-th order in \(s\), this implies that \(D^{J,X; \chi ,P}_{(u, [\mu]_m)}\) is still
Fredholm. To see that \(E^{J,X; \chi ,P}_{(u,
  [\mu]_m)}\) is Fredholm, we use in addition Fact \ref{pertDecay}
above, and 
the usual excision argument for Fredholmness in this situation.

The deformation operators for parameterized moduli spaces are
finite-rank stabilizations of the above operators, and their
Fredholmness is thus evident from the above discussion.

\mysubsubsection[Estimating the nonlinear term.]
The contribution of the perturbation to the nonlinear term 
\(N^{J,X;\chi, P}_{(u,[\mu]_m)}(\xi)\) is
\[
\begin{split}
&\chi(\Omega_\nu(u, [\mu]_m))\nabla_\xi\nabla_\xi\nabla P(u)\\
&\quad
+2\chi'(\Omega_\nu(u, [\mu]_m))\int_{S^1}\omega_\nu(\xi,
\p_t u)\, dt\, \nabla_\xi\nabla P(u)\\
&\quad
+\chi''(\Omega_\nu(u, [\mu]_m))\Big(\int_{S^1}\omega_\nu(\xi,
\p_t u)\, dt\Big)^2\, \nabla P(u)\\
&\quad
+\chi'(\Omega_\nu(u, [\mu]_m))\int_{S^1}\omega_\nu(\xi,
\p_t \xi)\, dt\, \nabla P(u).
\end{split}
\]
It is straightforward to check that each term above may be bound by
\[
C\|\xi\|_{C^0}\|\xi\|_{L^p_1}\leq C' \|\xi\|_{L^p_1}^2.
\]
We omit the straightforward estimate for the parameterized version.

\mysubsubsection[Compactness.]  
Let's go over the main ingredients in the usual proof one by one:
\begin{itemize}
\item {\it elliptic regularity.} By the above estimate on the
  nonlinear term, and the form of  (\ref{add:deformOp}), the elliptic bootstrapping
argument still hold, provided a \(C^0\) bound can be established. The
latter relies on the Gromov compactness.
\item {\it Gromov compactness.} Going through the rescaling argument, we
  note that the extra term \(\wp_{\chi P}\) disappears in the limit,
  and therefore again (local) compactness is lost only through
  bubbling off honest holomorphic spheres. This possibility is eliminated via
  transversality as in section 3 of Part I.
\item {\it energy bound.} With this definition of nonlocal perturbations, 
there might not be a good action functional for the perturbed
flows\footnote{We may easily modify the definition of \(\wp_{\chi P}\)
  so that there is; however we would run into difficulty with Gromov
  compactness.}. However, we still have the requisite energy bound for
perturbed flows with weight \(\leq \Re\). Let \((u, [\mu]_m)\) be such a
\(\chi P\)-perturbed flow. Then
\begin{equation}\label{pertEnergy}
\begin{split}{\cal E}(u, [\mu]_m)&=\|\p_s u\|_2^2\\
&=\alpha  \op{wt}_{-\langle{\cal Y}\rangle, e_{\cal P}}(u, [\mu]_m)
+\int\Big\langle\p_s u,
\chi(\Omega_\nu(u, [\mu]))\nabla P(u)\Big\rangle_{2,t}(s)\, ds.
\end{split}\end{equation}
On the other hand, 
\begin{lemma*}
Let \((u, [\mu]_m)\) be a \(\chi P\)-perturbed flow (either a
connecting flow line or
an orbit). Then there
is a constant \(C\) independent of \(s\) or \((u,[\mu]_m)\), such that
\[
\|\nabla P(u)\|_{2,t}(s)\leq C\|\p_s u\|_{2,t}(s) \, \, \forall  s.
\]
\end{lemma*}
\begin{proof}
This follows from the \(L^{\infty}_{1}\)-boundedness of \(P\), the
fact that \(P\) vanishes to up order \(k>2\) at the critical points, and
the following 
\paragraph{Palais-Smale condition.}
{\it There exists an \(\varepsilon'>0\) such that for any \((z,
  [\mu]_m)\in\lc^{\nu, m}\) with \(\|J(z)\partial_tz+\check{\theta}_{X}(z)+\wp_{\chi P}(z,
  [\mu]_m)\|_{2,t}\leq \varepsilon'\), there is a critical point \(z_0\)
  such that \(z=\exp (z_0, \xi)\) for a small \(\xi\), and 
\[\begin{split}
\| &J(z) \partial_tz +\check{\theta}_{X}(z) +\wp_{\chi P}(z,
  [\mu]_m)\|_{2,t}\geq \\
& \begin{cases}
C_1\|\xi\|_{2,t} &\text{when \(z_0\) is nondegenerate}\\
C_2\|\xi\|_{2,t}^2 &\text{when \(z_0\) is minimally degenerate 
in a standard d-b neighborhood.} 
\end{cases}
\end{split}\]}
This in turn follows from the Ascoli-Arzela argument as in the
unperturbed case, since by our condition on \(P\), \(\wp_{\chi P}\)
can be ignored near critical points. 
\end{proof}
Thus, if \(\|\chi\|_{C_{\epsilon}}\leq \varepsilon\), the absolute
value of the last term in (\ref{pertEnergy}) can be bounded by
\(C\varepsilon\|\p_s u\|_2^2\), and by rearranging,
\[
{\cal E}(u, [\mu]_m)\leq (1-C\varepsilon)^{-1}\alpha \op{wt}_{-\langle{\cal
    Y}\rangle, e_{\cal P}}(u, [\mu]_m)\leq (1-C\varepsilon)^{-1}\alpha  \Re.
\]

\item {\it global compactness} (for \(\hat{\cm}_P, \hat{\cm}_P^\Lambda\)). As in the
  unperturbed case, to go from local compactness to global
  compactness, we just need in addition Fact \ref{pertDecay}.
\end{itemize}

\mysubsubsection[Transversality.]
The transversality arguments in Part I uses a unique continuation
theorem extensively, however, Aronszajn's theorem
or the Carleman similarity principle used in \cite{trans} is no longer
applicable as the nonlocal term is introduced. 
While it might be possible to prove a unique continuation result
in this situation, we choose not to develop a general theory here.
Instead, for the purpose of proving Proposition \ref{prop:nep}, we
only need the following
\begin{claim*}
{\it Let \((J, X)\) be regular, \(\chi\) be sufficiently small (in
\(C_\epsilon\) norm), \(P\in V^k_\delta(J, X)\) and \(i\leq2\). Then
for \(\cm=\cm_P\) or \(\cm_O\), \(\cm^i_{\nu, m}(J, X; \chi, P;
\op{wt}_{-\langle {\cal P}\rangle, e_{\cal P}}\leq \Re)\) is (Zariski)
smooth. Similarly for the parameterized versions.}
\end{claim*}

Take \(\cm_P\) for example; the arguments for \(\cm_O\)
or the parameterized versions are similar. 
Due to Lemma \ref{glue-K}, for regular \((J, X)\) and
\(u\) in a neighhorhood of a lower dimensional strata
of \(\hat{\cm}_P^{i-1, +}(J, X; \op{wt}_{-\langle {\cal P}\rangle,
  e_{\cal P}}\leq \Re)\), the deformation operator
\(E_u\) has a uniformly bounded right inverse. Combining this with 
the compactness of 
\(\hat{\cm}_P^{i-1, +}(J, X; \op{wt}_{-\langle {\cal P}\rangle,
  e_{\cal P}}\leq \Re)\), there is a small number
\(\delta>0\)  such that any element in \(\{\mathfrak{D}\, |\, \mathfrak{D}-E_{v}^{JX}\|<\delta,
v\in\hat{\cm}_P^{i-1}(J,X; \op{wt}_{-\langle {\cal P}\rangle,
  e_{\cal P}}\leq \Re)\}\) is surjective. 
In particular, there is a
\(\delta'=\delta'(\delta)\), such that for any element \(w\) in 
\[\Big\{\exp (v,\xi)\, |\,  \|\xi\|_{\infty,1}<\delta', \, \, v\in
\hat{\cm}_P^{i-1}(J, X; \op{wt}_{-\langle {\cal P}\rangle , e_{\cal
    P}}\leq \Re )\Big\},\]
\(E_{(w,[\mu_w]_m)}^{JX; \chi P}\) is surjective for any lift \((w,
[\mu_w]_m)\) of \(w\) in \(\lc^{\nu, m}\), \(\forall \chi\) with \(\|\chi\|_{C_\epsilon}<\delta'\).
Thus, the claim follows from the following
\begin{lemma*}
Fix \(P\) \(i\), and \(\Re\) as above. 
Then there is an \(\varepsilon'>0\) such that for all \(\chi\) with
\(\|\chi\|_{C_{\epsilon}}<\varepsilon'\), any element \((u,
[\mu])\in\hat{\cm}_{P;\nu, m}^{i-1}(J, X;\chi, P, \op{wt}_{-\langle {\cal P}\rangle , e_{\cal
    P}}\leq \Re)\) is close to
\(\hat{\cm}_P^{i-1}(J, X, \op{wt}_{-\langle {\cal P}\rangle , e_{\cal
    P}}\leq \Re)\) in the sense that 
\[(*) \quad u=\exp (v, \xi), \, \, \|\xi\|_{\infty,1}<\delta'
\quad\text{for some \(v\in\hat{\cm}_P^{i-1}(J, X, \op{wt}_{-\langle {\cal P}\rangle , e_{\cal
    P}}\leq \Re)\).}\]
\end{lemma*}
\begin{proof}
Suppose the contrary. Then there exists a sequence \(\{\chi_n\}\),
\(\lim_{n\to \infty}\|\chi_n\|_{C_\epsilon}=0\), and a sequence
\(\{(u_n, [\mu_n]_m) \, |\,  (u_n, [\mu_n]_m)\in \hat{\cm}_{P,m}^{k}(J,
X;\chi, P; \op{wt}_{-\langle {\cal P}\rangle , e_{\cal
    P}}\leq \Re)\}\)
such that none of them satisfies (*).
By Gromov compactness, such a sequence \((u_n,
[\mu_n]_m)\) must weakly converge to an element \(v\) in
\(\hat{\cm}_P^{k}(J, X; \op{wt}_{-\langle {\cal P}\rangle , e_{\cal
    P}}\leq \Re)\) 
together with some bubbles. Since by the regularity of \((J,X)\),
there is no such bubble,
\((u_n, [\mu_n]_m)\) are close to \(v\), contradicting our assumption.
\end{proof}

This also shows that when \(\chi\) is sufficiently small and \((J,
X)\) regular, these \(\chi P\)-perturbed
flows avoid pseudo-holomorphic spheres, as in the case before perturbation.
\medbreak

\section{Orientation and Signs.}

In this section, we tie up the last loose end of this article by
addressing all orientation issues so far ignored: we verify (FS4) for
the Floer theory described in \S I.3, and show that an admissible \((J,
X)\)-homotopy satisfies (RHFS4). 

In \S7.2, we show that the various moduli spaces \(\cm_P(x,
y)\), \(\cm_O^1\) and their parameterized variants 
are orientable; furthermore, we introduce the notions of coherent
orientations for \(\cm_P, \cm_P^\Lambda\) and grading-compatible
orientations for \(\cm_O^1, \cm_O^{\Lambda, 2}\), and show that these
moduli spaces may be endowed with such orientations. This completes
the verification that the formal flow associated to a regular pair
\((J, X)\) forms a Floer system. The coherence of orientation is
determined by linearized versions of the gluing theorems proven in the
previous sections; this is in fact why we have postponed this
discussion. 
Compared with the full gluing theory, major simplifications for these
linearized versions arise from
the fact that we may substitute the complicated polynomially-weighted 
Sobolev spaces used in sections 2--5 
by larger, exponentially-weighted versions, due to the removal of 
constraints from non-linear aspects of general gluing theory.
Furthermore, deformation operators between these exponentially
weighted spaces are conjugate to deformation operators between 
the usual \(L^p_k\) spaces with perturbation by
asympotically constant 0-th order terms, making it possible to work
only with the ordinary Sobolev norms throughout this section. 

In \S7.3, we verify the signs in the expressions for \(\mathbb{T}_{P,
  db}, \ldots, \mathbb{T}_{O, hs-s}\) given in (RHFS4) (cf. \S
I.4.3.7). This is obtained by examining the orientations of the K-models
used in the proofs of the gluing theorems in previous sections. 
With this done, the verification
that admissible \((J,X)\)-homotopies satisfy the assumptions of
Proposition I.4.6.3 is complete, which in turn implies the general invariance
theorem, Theorem I.4.1.1.

\subsection{Basic Notions and Conventions.}

We first review some basic material to fix terminology and conventions.

\mysubsubsection[Orientation for direct sums and and exact sequences.]
Given a direct sum of an ordered
\(k\)-tuple of oriented vector spaces, \(E=E_1\oplus\cdots \oplus
E_k\), we orient it by
\(e_1\wedge\cdots \wedge e_k\in \det E\), where
\(e_i\in \det E_i\) orients \(E_i\).

An exact sequence of finite dimensional vector spaces
\[0\to E_1\stackrel{i_1}{\to} F_1\stackrel{j_1}{\to} E_2\stackrel{i_2}{\to} F_2\cdots\stackrel{j_{n-1}}{\to} E_n\stackrel{j_n}{\to} F_n\to 0\] 
determines an isomorphism
\(
\bigotimes_k\det E_k\simeq\bigotimes_k\det F_k,
\)
by writing 
\[
E_k= B_k^E\oplus j_{k-1}B^F_{k-1},\quad
F_k= i_kB_k^E\oplus B_k^F
\]
for appropriate oriented 
subspaces $B_k^E, B_k^F$, over which $i_k$, $j_k$
restrict to isomorphisms.

\mysubsubsection[Orientation for determinant lines and K-models.]
Given a Fredholm operator \(\mathfrak{D}: E\to F\), the determinant line
\[
\det \mathfrak{D}:=\det\ker \mathfrak{D}\otimes\det(\cok \mathfrak{D})^*.
\]
It is well known that for a continuous family (in operator norm) of
Fredholm operators, the determinant lines above form a real line bundle
over the parameter space, namely the determinant line bundle.
We use \(\op{\frak{or}}(\mathfrak{D})\) to denote the space of
possible orientations for \(\det \mathfrak{D}\) when it is orientable, and similarly, \(\op{\frak{or}}(\mathfrak{D}^\Lambda)\) denotes the space of possible
trivializations of the determinant line bundle for the family
\(\mathfrak{D}^\Lambda\) when it is orientable. These are affine spaces
under \(\Z/2\Z\). If \(\mathfrak{D}_{\lambda_1},
\mathfrak{D}_{\lambda_2}\) are elements of the family of operators
\(\mathfrak{D}^\Lambda\), we say that \(\mathfrak{o}_1\in
\op{\frak{or}}(\mathfrak{D}_{\lambda_1})\) and \(\mathfrak{o}_2\in
\op{\frak{or}}(\mathfrak{D}_{\lambda_2})\) are {\em correlated} via \(\mathfrak{D}^\Lambda\) if
they are restrictions of the same trivialization \(\mathfrak{O}\in
\op{\frak{or}}(\mathfrak{D}^\Lambda)\). They are said to be of {\em
  relative sign} \(\rho\in \{\pm1\}\) (with respect to \(\mathfrak{D}^\Lambda\)), denoted
\([\mathfrak{o}_1/\mathfrak{o}_2]\), if \(\mathfrak{o}_1\) and \(\rho\mathfrak{o}_2\)
are correlated. 

It is convenient to describe the orientation of \(\det \mathfrak{D}\) in
terms of K-models. 
Recall the definition of oriented K-models and 
the exact sequence (\ref{K-C-sequence}) from \S\ref{sec:K-model}.
This exact sequence induces the isomorphism:
\begin{equation*}\label{app-det}
\det \mathfrak{D}\simeq \det K\otimes \det C^*.
\end{equation*}
Thus, an orientation of a K-model for $\mathfrak{D}$ decides an orientation for
\(\det \mathfrak{D}\). Given an orientation of \(\det \mathfrak{D}\), an oriented K-model
\([K; C]\) of \(\mathfrak{D}\) is
said to be {\em compatible} with this orientation, if the orientation
of \([K; C]\) induces the orientation of \(\det \mathfrak{D}\).

Two K-models of $\mathfrak{D}$ are said to be {\em co-oriented} if they give rise to the
same orientation of $\det \mathfrak{D}$. 
Let \([\mathfrak{D}_{\lambda_1}: K_{\lambda_1}\to  C_{\lambda _1}]_{B_{\lambda_1}}\),
\([\mathfrak{D}_{\lambda_2}:K_{\lambda_2}\to
C_{\lambda_2}]_{B_{\lambda_2}}\) be fibers of a family K-model for \(\mathfrak{D}^\Lambda\).
They are said to be {\em mutually co-oriented} via the
family \(\mathfrak{D}^\Lambda\) if they are with respectively
compatible with orientations of \(\det\mathfrak{D}_{\lambda_1}\) and
\(\det\mathfrak{D}_{\lambda_2}\) correlated by \(\mathfrak{D}^\Lambda\).

\mysubsubsection[Induced orientation of a stabilization.]
Let \(\hat{\mathfrak{D}}_{\Psi}:\R^k\oplus E\to F\) be a stabilization
of \(\mathfrak{D}: E\to F\); recall the definition of stabilized
K-models from \S\ref{sec:K-model}.

Given an orientation \({\frak o}\in \op{\frak{or}}(\mathfrak{D})\), we define the
{\em induced orientation} \(\hat{\frak o}\in \op{\frak{or}}(\hat{\mathfrak{D}}_\Psi)\) from \({\frak o}\) as follows. Given an oriented 
K-model \([\mathfrak{D}: K\to C]_B\) compatible with \({\frak o}\), let
\(\hat{\frak o}\) be the orientation given by the stabilization
\([\hat{\mathfrak D}_\Psi:\hat{K}\to C]_{\hat{B}}\), where \(\hat{K}\)
is oriented as
\[
\hat{K}=(-1)^{k\ind \mathfrak{D}}\, \R^k\oplus K.
\]

\mysubsubsection[Reduction of oriented K-models.]
Let the K-model \([K'\to C']\)
be a reduction of another K-model, \([K\to C]\), by \(Q\)
(cf. \S\ref{sec:K-model}). Then the orientation of one K-model induces an
orientation of the other via writing 
\[K=K'\oplus Q;\quad C=C'\oplus \Pi_C\mathfrak{D}(Q)\]
as oriented spaces.  Note that changing the orientation of \(Q\)
results in a co-oriented K-model.

\mysubsubsection[Orientation for glued K-models.] 
Recall the definitions and notations in \S\ref{glue-K}. 
Given an ordered \(k\)-tuple of finite dimensional subspaces \(K_1, \ldots,
K_k\) in \(E\) or \(F\), and sufficiently large \(R_1, \ldots, R_k\),
we orient the glued space \(K_1\#_{R_1}\cdots \#_{R_{k-1}}K_k\) or
\(K_1\#_{R_1}\cdots \#_{R_{k-1}}K_k\#_{R_k}\) by its natural
isomorphism with \(K_1\oplus K_2\oplus \cdots\oplus K_k\). 

Let \(\mathfrak{D}_1\), \(\mathfrak{D}_2\) be an ordered pair of glue-able Floer type operators,
and let \(\mathfrak{D}\) be a cyclically glueable Floer-type operator. Given \({\frak o}_1\in
\op{\frak{or}}(\mathfrak{D}_1)\), \({\frak o}_2\in \op{\frak{or}}(\mathfrak{D}_2)\), \({\frak o}\in
\op{\frak{or}}(\mathfrak{D})\), we define \({\frak o}_1\#_R{\frak o}_2\in \op{\frak{or}}(\mathfrak{D}_1\#_R\mathfrak{D}_2)\), \(\frak{o}_{\#_R}\in \op{\frak{or}}(\mathfrak{D}\#_{R})\) as follows.
Let \([K_1\to C_1]\), \([K_2\to
C_2]\), \([K\to C]\) be oriented \(K\)-models compatible with \({\frak o}_1\),
\({\frak o}_2\), and \({\frak o}\)
respectively. Then the {\it induced orientation}, \({\frak
  o}_1\#_R{\frak o}_2\) and \({\frak o}_{\#R}\),
are respectively the orientation given by the oriented K-models
\begin{equation}\label{lm:glue}
\Big[(-1)^{\dim (C_1) \cdot \ind \mathfrak{D}_2}K_1\#_R K_2\to 
C_1\#_RC_2\Big], \quad
[K_{\#R} \to C_{\#R}].
\end{equation}

The orientations for the generalized kernels and generalized cokernels
of stabilized, reduced, or glued K-models given above are chosen such 
that co-oriented K-models give rise to co-oriented stabilized,
reduced, or glued K-models. Thus, we have well-defined homomorphisms
of affine spaces under \(\Z/2\Z\):
\[
s_\Psi: \op{\frak{or}}(\mathfrak{D}) \to
\op{\frak{or}}(\hat{\mathfrak{D}}_\Psi), \quad
{\frak g}_R: \op{\frak{or}}(\mathfrak{D}_1)\times \op{\frak{or}}(\mathfrak{D}_2)\to
\op{\frak{or}}(\mathfrak{D}_1\#_R\mathfrak{D}_2), \,\, 
\frak{sg}_R: \op{\frak{or}}(\mathfrak{D})
\to \op{\frak{or}}(\mathfrak{D}_{\#R}),
\]
sending \({\frak o}\) to \(\hat{\frak o}\), \(\mathfrak{o}_1\times
\mathfrak{o}_2\) to \(\mathfrak{o}_1\#_R\mathfrak{o}_2\), and
\(\mathfrak{o}\) to \(\mathfrak{o}_{\#R}\) respectively. We call
\(s_\Psi\) the {\em stabilization isomorphism}, and
\(\mathfrak{g}_R\), \(\mathfrak{sg}_R\) the {\em gluing
  homomorphisms}. 
As a consequence of the independence of K-models, the above
constructions also work in the family setting to define induced orientations
for the determinant line bundles of stabilized or glued operators.
In addition, the gluing homomorphisms above may be extended to 
be defined for arbitrary \(k\)-tuple of glueable or cyclically
glueable Floer type operators, by observing that any glued operator or
cyclically glued operator can be obtained by a combination of
translation and the two gluing operations discussed above. 
Morever, with the above definition, it is straightforward to check
that the oriented K-model for the same glued operator obtained from
different combinations are actually the same.  

\begin{remark*}
(1) Alternatively, one may define induced orientation for
stabilization by the oriented K-model \([\R^k\oplus K\to C]\)
instead. We have so chosen our definition because in our context, 
\(\det \mathfrak{D}_\Psi\) gives the orientation of a fiber bundle,
where \(\R^k\) corresponds to the tangent space of the base. We prefer
the ``fiber-first'' convention for orienting a fiber bundle. 
With our definitions, the gluing homomorphism commutes with the
(rank k) stabilization isomorphism on \(\mathfrak{D}_2\), but commutes with
stabilization on \(\mathfrak{D}_1\) modulo the sign \((-1)^{k\ind\mathfrak{D}_2}\).

\noindent (2) The definitions of the orientation for a stabilization
and glued operators in \cite{FH} differ from ours. 
Their definitions have the following disadvantage: 
Given an orientation of a determinant line bundle for a family
\(\{\mathfrak{D}_\lambda\}_{\lambda\in \Lambda}\), the stabilization isomorphism
of \cite{FH} 
gives a {\em possibly discontinuous}, nowhere vanishing section of
the determinant line bundle of the stabilized family
\(\{\hat{\mathfrak{D}}_{\lambda, \Psi}\}_{\lambda\in \Lambda}\).
Furthermore, the gluing morphisms in \cite{FH} commute with stabilization
{\it only up to a sign} depending on the dimension of
\(\R^k\). 
\end{remark*}

\subsection{Orienting Moduli Spaces.}
This subsection addresses the orientability issues required by (FS4)
and (RHFS4).

By an {\em orientation} of a moduli space \(\cm=\cm_P\) or \(\cm_O\), we mean the
following. Notice that the
configuration spaces \({\cal B}_P^k(x, y)\), \({\cal B}_O^k\)
parameterize families of deformation operators, \(\{E_u \,
|\, u\in {\cal B}_P^k(x, y)\}\), \(\{\tilde{D}_{(T, u)}\, |\, (T,
u)\in {\cal B}_O^k\}\). Thus, they carry determinant line bundles, which
we denote by \(L{\cal B}_P^k(x, y)\), \(L{\cal B}_O\).
The moduli space \(\cm\subset {\cal B}={\cal B}_P^k(x, y)\) or \({\cal
  B}_O^k\) parameterizes a subfamily of deformation operators,
and thus carries a determinant line bundle \(L{\cal M}\),
which is the pull-back of \(L{\cal B}\).
An orientation of \(\cm\) will mean a
trivialization of \(L\cm\).
In this article, this will always be the
pull-back of a trivialization of \(L{\cal B}\), and we shall therefore
focus on orienting \(L{\cal B}\) for various configuration spaces
\({\cal B}\). Similarly, parameterized moduli spaces
\(\cm^\Lambda\) will be oriented by orienting \(L{\cal
  B}^\Lambda\). Since the deformation operators for \(\cm^\Lambda\)
are stabilizations of those for \(\cm\), this also orients the fiber
moduli spaces \(\cm_\lambda\) for \(\lambda\in \Lambda\).

Notice that the above definition does not 
require nondegeneracy of the moduli
spaces \(\cm\), and hence we make no such assumptions in this 
subsection.
Nevertheless, when \(\cm\) is nondegenerate, the determinant line
for the relevant deformation operator \(\det \mathfrak{D}_u=\det T_u\cm\)
at any \(u\in\cm\). In this case this definition agrees
with the usual definition of the orientation of a manifold.

We do, however, assume nondegeneracy of the spaces of critical points.
Namely, we assume (FS1) for a Floer theory \((\lc,  \mathfrak{H}, \ind;
{\cal Y}_\chi, {\cal V}_\chi)\), and assume
(RHFS1*) for a CHFS throughout this subsection.

We begin by some general discussion on abstract Floer theories 
in \S7.2.1--4.

\subsubsection{General strategy for orientability.}
Below we roughly outline a scheme to establish orientability of
\(L{\cal B}\), which is particularly useful for symplectic Floer theories, 
when the configuration spaces have complicated topology. 
To begin, construct a map \[m: {\cal B}\to \Sigma/G,\] 
where \(\Sigma \) is a contractible space parameterizing certain
operators, and \(G\) is a suitable automorphism group.
The map \(m\) is typically defined by identifying the deformation
operator at \(u\in {\cal B}\) to an operator in \(\Sigma\), after
certain trivialization is chosen. \(G\) is usually the group of
automorphisms relating different possible trivializations. 

The space \(\Sigma\) parameterizes a trivial determinant line bundle
\(L\Sigma\), over which the action by \(G\) extends. Moreover,
\((L\Sigma)/G=L(\Sigma/G)\) and \(L{\cal B}=m^*L(\Sigma/G)\). One next
shows that \(G\) induces trivial actions 
on the determinant lines. Thus \(L(\Sigma/G)\), and hence also
\(L{\cal B}\), are trivial.

In family settings, \({\cal B}\) and \(\Sigma\) above are both
replaced by bundles \({\cal B}^\Lambda\), \(\Sigma^\Lambda\) over the
parameter space \(\Lambda\), and \(m\) above will be a bundle map.

In the case \({\cal B}={\cal B}_P(x, y)\), in order for the
deformation operator to be Fredholm, \(x\), \(y\in {\cal P}\) have to be
nondegenerate. More generally, one may consider exponentially weighted
versions of deformation operators \(E_u^{(\sigma_1, \sigma_2)}\)
(cf. \S I.3.2.3) instead of \(E_u\). When \(x\), \(y\) are
respectively \(\sigma_1\)-weighted nondegenerate and
\(\sigma_2\)-weighted nondegenerate, this defines another determinant
line bundle over \({\cal B}_P(x, y)\), which we denote by
\(L^{(\sigma_1, \sigma_2)}{\cal B}_P(x, y)\). Under this
weighted-nondegeneracy condition on \(x\), \(y\), the determinant line
bundle \(L^{(\sigma_1, \sigma_2)}{\cal B}_P(x, y)\) is independent of
small perturbations to the weights \(\sigma_1, \sigma_2\).

These weighted versions are useful for dealing with the case when one
of \(x\), \(y\) is minimally degenerate: in this case, the deformation
operator \(E_u\) is defined as a map between complicated polynomially
weighted Sobolev spaces (cf. \S I.5). However, we showed in \S I.5.2.5
that \(E_u\) has identical kernel and cokernel as \(E_u^{(-\sigma,
  \sigma)}\) for any small positive \(\sigma\). As we are only
concerned with the {\em linear} aspect (Kuranishi structure) of the
Floer theory, there is thus no harm in replacing \(E_u\) by the
simpler \(E_u^{(-\sigma, \sigma)}\): the orientation of \(\cm_P(x,
y)\) in this case will be given by an orientation of \(L^{(-\sigma,
  \sigma)}{\cal B}_P(x, y)\).

Turning now to the family situation of a CHFS \(\{{\cal
  V}_\lambda\}_{\lambda\in\Lambda}\) satisfying (RHFS1*), given \({\bf
  x}, {\bf y}\in \aleph_\Lambda\) and an interval \(S\subset
\Lambda\), let \({\cal B}_P^S(\underline{\bf x}, \underline{\bf y})=\bigcup_{\lambda\in S\cap
  \bar{\Lambda}_{\bf x}\cap \bar{\Lambda}_{\bf y}} {\cal
  B}_P(x_\lambda, y_\lambda)\). Under the assumption (RHFS1*), one may
choose a set of intervals \(\{S_i\}\) covering \(\Lambda\), such that
each \(S_i\) contains at most one death-birth, and all overlaps
\(S_i\cap S_j\) for different \(i, j\) contains no death-birth. Over
each \(S_i\), one may choose appropriate weights \(\sigma_{x,i}\),
\(\sigma_{y,i}\) with small absolute value, such that \(x_\lambda\) is
\(\sigma_{x, i}\)-weighted nondegenerate and \(y_\lambda\) is
\(\sigma_{y, i}\)-weighted nondegenerate for all \(\lambda\in S_i\cap
  \bar{\Lambda}_{\bf x}\cap \bar{\Lambda}_{\bf y}\). An orientation of
  \(L^{(\sigma_{x, i}, \sigma_{y_i})}{\cal B}_P^{S_i}(\underline{\bf x},
  \underline{\bf y})\) determines an orientation of
  \(\cm_P^{S_i}(\underline{\bf x}, \underline{\bf y})\) as well as
  ones for its fibers \(\cm_{P, \lambda}(\underline{\bf x}, \underline{\bf y})\),
  which agree with our previous discussion on orienting \(\cm_P(x,
  y)\) for non-degenerate or minimally degenerate \(x\), \(y\).
Note again that the precise values of the weights \(\sigma_{x, i}\),
\(\sigma_{y, i}\) are immaterial; in particular,
  when \(S_i\) contains no death-birth, they can be chosen to be
  0. Otherwise, only the signs of these weights matter.  

Lastly, we may patch up the determinant line bundles 
\(L^{(\sigma_{x, i}, \sigma_{y_i})}{\cal B}_P^{S_i}(\underline{\bf x},
\underline{\bf y})\) for all intervals \(S_i\) to define a determinant
line bundle \(L{\cal B}_P^{\Lambda}(\underline{\bf x},
\underline{\bf y})\) over \({\cal B}_P^{\Lambda}(\underline{\bf x},
\underline{\bf y})\), by observing that, since for all \(\lambda \in \bigcup_{i,
  j}\overline{S_i\cap S_j}\), \(x_\lambda, y_\lambda\) are
nondegenerate, determinant line bundles with different weights over
\({\cal B}_P^{\overline{S_i\cap S_j}}(\underline{\bf x},
\underline{\bf y})\) can be identified. 

More concretely, in \S7.2.5 we shall apply the above general scheme to the
specific Floer theory described in \S I.3. See also \cite{L2} for its
application in other versions of symplectic Floer theories. In gauge
theoretic settings, the configuration space \({\cal B}\) itself has
the structure of 
\({\cal A}/{\cal G}\), where \({\cal A}\) is an affine space, and
\({\cal G}\) is the gauge group, which is often connected under the
assumption of simple-connectivity of the underlying manifold.
Thus, much of the above scheme also carry over to this context.

\subsubsection{Coherent Orientations for \(L{\cal B}_P\). }
Assuming the orientability of \(L{\cal B}_P^k(x, y)\) and \(L{\cal
  B}^{\Lambda, k+1}_P({\bf x}, {\bf y})\) for any pair of \(x, y\in {\cal
  P}\) or \({\bf x}, {\bf y}\in \aleph_\Lambda\) and any \(k\in \Z\),
we explain in \S7.2.2--3 how
the orientations of all these should be related, so as to endow the
moduli spaces of broken trajectories with a correct oriented
manifold-with-corners structure. 
\begin{notation*}
Given a determinant line bundle \(LQ\) over a parameter space \(Q\),
\(LQ\backslash \text{zero section}\) contracts to a \(\Z/2\Z\)-bundle
over \(Q\), which we denote by \(\op{\frak{Or}} (LQ)\). \(LQ\) is
orientable if \(\op{\frak{Or}} (LQ)\) is trivial; in this case, the
space of sections of \(\op{\frak{Or}} (LQ)\) is denoted
\(\op{\frak{or}} (LQ)\): this \(\Z/2\Z\)-torsor is the space of
possible orientations for \(LQ\).
\end{notation*}
Recall the continuity of the gluing homomorphism
\(\op{\frak{g}}_R(\mathfrak{D}_1, \mathfrak{D}_2)\) in
\(\mathfrak{D}_1\), \(\mathfrak{D}_2\), and \(R\). Thus, it defines a 
gluing homomorphism 
\begin{equation*}\label{eq:glue-or}
\op{\frak{g}}: \op{\frak{or}}(L{\cal B}_P^{k_1}(z_1, z_2))\times\op{\frak{or}}(L{\cal
  B}_P^{k_2}(z_2, z_3))  \to \op{\frak{or}}(L{\cal B}_P^{k_1+k_2}(z_1, z_3))
\end{equation*}
for any \(z_1, z_2, z_3\in {\cal P}\) and \(k_1, k_2\in \Z\).
We write \(\op{\frak{g}}(o_1, o_2)=o_1\# o_2\).
\begin{definition*}Let \((\lc, {\frak H}, \ind; {\cal Y}_\chi, {\cal V}_\chi)\) be a
Floer theory satisfying (FS1); in particular, \({\cal P}\) consists of
finitely many nondegenerate elements.
A {\em coherent orientation} of \[L{\cal B}_P=\coprod_{k\in \Z}\coprod_{x, y\in {\cal
    P}}L{\cal B}_P^k(x, y)\] is a section, \({\frak s}\), of 
\(\op{\frak{Or}}(L{\cal B}_P)\), so that for all \(k_1, k_2\in \Z\), \(z_1, z_2, z_3\in {\cal P}\),
\begin{equation}\label{cond:coh}
\frak{s}|_{{\cal B}_P^{k_1}(z_1, z_2)}\#{\frak s}|_{{\cal B}_P^{k_2}(z_2, z_3)}={\frak s}|_{{\cal B}_P^{k_1+k_2}(z_1, z_3)}.\end{equation}
\end{definition*}
A moment's thought (or cf. \cite{FH}) reveals that coherent
orientations always exist. 
Fixing an \(x\in {\cal P}\), a coherent orientation for \(L{\cal
  B}_P\) is determined by choosing an element in
\(\op{\frak{or}}(L{\cal B}_P^{k_y}(x, y))\) for each \(y\neq
x\) and an integer \(k=\op{gr}(x, y) \mod 2\mathbb{N}_\psi\), and in
the case when \(\mathbb{N}_\psi\neq0\), an additional 
element of \(\op{\frak{or}}(L{\cal B}_P^{2\mathbb{N}_\psi}(x,
x))\).
The cases of 
\(\op{card}({\cal P})=0\), or \(\op{card}({\cal P})=1\) and
\(\mathbb{N}_\psi=0\) are excluded: in the first case, \({\cal B}_P\) is
empty, while in the second case, there is no nonconstant connecting
flow line. Thus, there is nothing to orient in these cases.

The following fact is immediate from the definition of coherent
orientation, but shall be important later.

\begin{lemma*}Let \({\frak s}\) be an arbitrary coherent
orientation of \(L{\cal B}_P\). Then for any \(x\in {\cal P}\), 
\({\frak s}|_{{\cal B}_P^0(x, x)}\)
is the canonical orientation of \(L{\cal B}_P^0(x, x)\). 
\end{lemma*}

In the above, the {\em canonical
orientation} of \(L{\cal B}_P^0(x, x)\) is that determined by the  
canonical orientation of \(\det E_{\bar{x}}\), the latter being due to
the identification of the kernel and cokernel of 
\(E_{\bar{x}}=d/ds+A_{x}\) via  the facts that \(\ker E_{\bar{x}}=\ker
A_x\),  \(\cok E_{\bar{x}}=\cok A_x\), and that \(A_x\) is self-adjoint. 
\smallbreak
\begin{proof}
The coherence condition (\ref{cond:coh}) requires the gluing maps
\begin{gather*}
\op{\frak g} ({\frak s}|_{{\cal B}_P^0(x, x)}, -): \op{\frak{or}}(L{\cal
  B}_P^{k}(x, y))  \to \op{\frak{or}}(L{\cal B}_P^{k}(x, y)); \\
\op{\frak g} (-, {\frak s}|_{{\cal B}_P^0(x, x)}): \op{\frak{or}}(L{\cal
  B}_P^{k'}(z, x))  \to \op{\frak{or}}(L{\cal B}_P^{k'}(z, x))
\end{gather*}
to be the identity map. 
\end{proof}

\subsubsection{Coherent Orientations for \(L{\cal B}_P^\Lambda\). }
Given a CHFS \(\{(\lc, \mathfrak{H}, \ind; {\cal Y}_\lambda, {\cal V}_\lambda )\}_{\lambda\in \Lambda}\) satisfying
(RHFS1*), we aim to orient \(L{\cal B}_P^\Lambda\), where
\[{\cal B}_P^\Lambda=\coprod_{k\in
  \Z}\coprod_{{\bf x}, {\bf y}\in \aleph_\Lambda}{\cal
  B}_P^{\Lambda, k}(\underline{\bf x}, \underline{\bf y}).\]
There is a natural projection map \(\Pi_\Lambda: {\cal B}_P^\Lambda\to
\Lambda\), whose fiber over \(\lambda\in \Lambda\) is:
\begin{itemize}
\item\({\cal B}_{P, \lambda}=\coprod_{x_\lambda, y_\lambda\in {\cal
    P}_\lambda} {\cal B}_P(x_\lambda, y_\lambda)\), when \(\lambda\)
is not a death-birth,
\item \({\cal B}_{P, \lambda}=\coprod_{x_\lambda, y_\lambda\in {\cal
    P}_\lambda\backslash \{z_\lambda\}\cup \{z_{\lambda+},
  z_{\lambda-}\}} {\cal B}_P(x_\lambda, y_\lambda)\), when \({\cal
  P}_\lambda\) contains a unique minimally degenerate critical point \(z_\lambda\).
\end{itemize}
The elements \(z_{\lambda+}, z_{\lambda-}\) should be regarded as the
end points of the two path components \({\bf z}_+, {\bf z}_-\) of
\({\cal P}^\Lambda\backslash {\cal P}^{\Lambda, deg}\) connected at
\(z_\lambda\), with \(\ind({\bf z}_+)=\ind_+(z)\), \(\ind({\bf
  z}_-)=\ind_-(z)\) respectively. 
We write 
\[
L{\cal B}_P(x_\lambda, z_{\lambda\pm})=\rho_\lambda^*L{\cal
  B}_P^\Lambda(\underline{\bf x}, \underline{\bf z}_\pm)=L^{(\sigma_x,
  \mp\sigma)}{\cal B}_P(x_\lambda, z_\lambda), \quad \text{for \(0<\sigma\ll1\),}
\]
where \(\rho_\lambda: {\cal B}_{P, \lambda}\hookrightarrow {\cal
  B}_P^\Lambda\) is the inclusion. 

First, observe that it is also useful to identify \(\det
E_u^{(\sigma_1, \sigma_2)}\) with \(\det E_u^{[\sigma_1, \sigma_2]}\), where
\[
E_u^{[\sigma_1, \sigma_2]}:=\varsigma^{\sigma_1, \sigma_2}E_u(\varsigma^{\sigma_1, \sigma_2})^{-1}=E_u+(\sigma_2s\beta(s)+\sigma_1s\beta(-s))'.
\]
is a map between ordinary Sobolev spaces. 
With this identification, one may extend the gluing homomorphism to
the weighted case:
\[
\op{\frak{g}}: \op{\frak{or}}(L^{(\sigma_1, \sigma_2)}{\cal
  B}_P^{k_1}(z_1, z_2))\times\op{\frak{or}}(L^{(\sigma_2, \sigma_3)}{\cal
  B}_P^{k_2}(z_2, z_3))  \to \op{\frak{or}}(L^{(\sigma_1, \sigma_3)}{\cal B}_P^{k_1+k_2}(z_1, z_3)).
\]
For any triple \({\bf z}_1, {\bf z}_2, {\bf z}_3\in \aleph_\Lambda\)
with \(\bar{\Lambda}_{{\bf z}_1}\cap \bar{\Lambda}_{{\bf z}_2} \cap \bar{\Lambda}_{{\bf
    z}_3}\neq \emptyset\), and any pair of integers \(k_1, k_2\),
one has also the parameterized version of gluing homomorphism
\[
\op{\frak{g}}^\Lambda: \op{\frak{or}}(L{\cal B}_P^{\Lambda, k_1}(\underline{\bf z}_1, \underline{\bf
  z}_2))\times\op{\frak{or}}(L{\cal
  B}_P^{\Lambda, k_2}(\underline{\bf z}_2, \underline{\bf z}_3))  \to \op{\frak{or}}(L{\cal
  B}_P^{\Lambda, k_1+k_2}(\underline{\bf z}_1, \underline{\bf z}_3))
\]
extending the gluing homomorphism \(\mathfrak{g}\) over the fibers
\(L{\cal B}_{P, \lambda}\).
\begin{definition*}Given a CHFS satisfying (RHFS1*) as above, a 
{\em coherent orientation} of  \(L{\cal B}_P^\Lambda\)
is a section, \({\frak s}\), of \(\op{\frak{Or}}(L{\cal B}_P^\Lambda)\), so that:
\begin{enumerate} 
\item For any triple \({\bf z}_1, {\bf z}_2, {\bf z}_3\in \aleph_\Lambda\)
with \(\bar{\Lambda}_{{\bf z}_1}\cap \bar{\Lambda}_{{\bf z}_2} \cap \bar{\Lambda}_{{\bf
    z}_3}\neq \emptyset\), and any pair of integers \(k_1, k_2\),
\[\frak{s}|_{{\cal
  B}_P^{\Lambda, k_1}(\underline{\bf z}_1, \underline{\bf z}_2)}\#{\frak s}|_{{\cal
  B}_P^{\Lambda, k_2}(\underline{\bf z}_2, \underline{\bf z}_3)}={\frak s}|_{
{\cal B}_P^{\Lambda, k_1+k_2}(\underline{\bf z}_1, \underline{\bf z}_3)};\]
\item For all \(y_\lambda\in {\cal P}^{\Lambda, deg}\), \({\frak
    s}|_{{\cal B}_P^1(y+, y-))}\) is the {\em standard orientation} of
  \(L^{(-\sigma, \sigma)}{\cal B}_P^1(y_\lambda, y_\lambda)\), namely,
  the orientation given by the oriented K-model
  \([E_{\bar{y}_\lambda}^{(-\sigma, \sigma)}: \R{\bf e}_{y_\lambda}\to *]_B\).
\end{enumerate}
\end{definition*}

Again, it is easy to see that such coherent orientation always exists.
When \({\mathbb N}_\psi\neq0\), there are \(\op{card}(\aleph _\Lambda)\)
possible coherent orientations. When \({\mathbb N}_\psi=0\), there are 
\(\op{card}(\aleph _\Lambda)-1\) of them.
Condition 2 in the above definition is imposed so that the short flow
line between the two new critical points \(y_{\lambda\pm}\) described
in \S5.3 has positive sign.

\subsubsection{Grading-compatible orientation for \(L{\cal B}_O^1\).}
The definition of our invariant \(I_F\) involves both \(\cm_O^1\) and
\(\cm_P^1\), which are related by gluing elements in \(\cm_P^0(x,x)\)
during a CHFS. It is thus crucial to orient
\(\cm_P^0(x, x)\) and \(\cm_O^1\) consistently. The notion of
``grading compatible orientation'' describes such a suitable
compatibility relation. More generally, one may consider compatibility
conditions relating the orientations of 
higher dimensional moduli spaces \(\cm_O^{2k{\mathbb N}_\psi+1}\), 
and \(\cm_P(x, x)^{2k{\mathbb N}_\psi}\), but this does not concern
us, since our invariant involves only low dimensional moduli spaces.

Let \({\cal B}_O^{k}\subset {\cal B}_O\) be the subset consisting of
elements \((T, u)\) with \(\op{gr}(u)=k-1\), and \(L{\cal B}_O^k\) be
the determinant line bundle of the family of deformation operators
\(\tilde{D}_{(T, u)}\). Assume that \(L{\cal B}_O^1\) is orientable.

Recall that the relative grading \(\op{gr}\) in a Floer theory is
typically defined via spectral flow by identifying deformation 
operators \(A_x\) with elements in a space of self-adjoint operators
\(\Sigma_C\) (cf. \S I.3.1.4 for the version relevant to this
article). On the other hand, the orientation of \(L{\cal B}_O^1\) is
defined by a map \(m_O: {\cal B}_O^1\to \Sigma_O^1/G_O\), where
\(\Sigma_O^1\) is a space of Fredholm operators of index 1, which
includes rank-1 stablizations of the operator
\[D_{\mathbb{A}; T}:=\partial_s+\mathbb{A}, \quad s\in S^1_T\]
for any surjective \(\mathbb{A}\in \Sigma_C\) and \(T\in \R^+\).
\begin{definition*}
For a Floer theory \((\lc, \mathfrak{H}, \ind; {\cal Y}_\chi, {\cal V}_\chi)\), 
the {\em grading-compatible} orientation of \(L{\cal B}_O^1\), or more
generally \(L\Sigma_O^1/G_O\) (also called the orientation {\em compatible} with
the absolute \(\Z/2\Z\)-grading \(\ind\)), is the orientation given by
the canonical orientation of \(\det \tilde{D}_{\mathbb{A}, T}\), where
\( \tilde{D}_{\mathbb{A}, T}\in \Sigma_O\) is the rank 1 stabilization
of \(D_{\mathbb{A}, T}\) by the
zero map, and \(\mathbb{A}\) is a surjective operator in
\(\Sigma_C\) of even index. 
\end{definition*}
In the above, the canonical orientation of \(\tilde{D}_{\mathbb{A},
  T}\) is the stabilization of the canonical orientation of
\(D_{\mathbb{A}, T}\), which in turn is defined in the same way as the
canonical orientation of \(E_{\bar{x}}\) (cf. \S7.2.2). 
Note that the choice of \(\mathbb{A}\) and \(T\) do not matter in the
above definition: as one varies \(T\), \(D_{\mathbb{A}, T}\)
remains surjective; one the other hand, the independence of the choice
of \(\mathbb{A}\) is a consequence of the following basic Lemma. 
\begin{lemma*}
For any two surjective operators \(\mathbb{A}, \mathbb{A}'\in \Sigma_C\), 
the canonical orientations of \(\det D_{\mathbb{A}, T}\)
and \(\det D_{\mathbb{A}', T}\) are of relative sign
\((-1)^{\op{gr}(\mathbb{A}, \mathbb{A}')}\) with respective to the family
\(\{D_{\mathbb{A}, T} | \mathbb{A}\in \Sigma_C\}\), where
\(\op{gr}(\mathbb{A}, \mathbb{A}')\)
denotes the relative index between \(\mathbb{A}\) and \(\mathbb{A}'\).
\end{lemma*}
An immediate corollary is:
\begin{corollary*}
Suppose that \(L\Sigma_O^1/G_O\) is orientable.
Then for any surjective \(\mathbb{A}\in \Sigma_C\) and \(T\in \R^+\),
the relative sign between the grading-compatible orientation of
\(L\Sigma_O^1/G_O\) and
the canonical orientation of \(\det D_{\mathbb{A}, T}\)
is \((-1)^{\ind A}\).  
\end{corollary*}
\subsubsection{Orientability in symplectic Floer theory.}
We now apply the general strategy described in \S7.2.1 to establish
the orientability of moduli spaces for the specific version of Floer
theory considered in this article.
\medbreak
\noindent{\bf (1)} {\sc Orienting \(L{\cal B}_P\).} This follows from
  \cite{FH}, which we now review in our terminology: 
Let \(J\in {\cal
    J}_K^{reg}\), \(X\) be \(J\)-nondegnerate (cf. \S I.3.2.1), and 
\(\Sigma_C\) be as in \S I.3.1.4.
Given two self-adjoint, surjective operators \(A_-, A_+\in
\Sigma_C\), let \(\Sigma_P(A_-, A_+)\) be the
space of operators of the form:
\footnote{\(\Sigma_P(A_-, A_+)\) is basically the space \(\Theta\) in Proposition 7 of
\cite{FH}.} 
\[
\partial_s+J(s, t)\partial_t+\nu(s,t) \, : \, L^p_1(\R\times S^1; \R^{2n})\to
L^p(\R\times S^1; \R^{2n}),
\]
where \(J\) is a smooth complex structure on the trivial
\(\R^{2n}\)-bundle over \(\R\times S^1\), compatible with the standard
symplecic structure on \(\R^{2n}\). \(\nu\) is a smooth matrix-valued
function, and both \(J\) and \(\nu\) extend smoothly over the
cylinder \([-\infty, +\infty]\times S^1\) that compactifies \(\R\times
S^1\). Furthermore, over the two circles at infinity, 
\[
J(-\infty, t)\partial_t+\nu(-\infty,t)=A_-; \quad 
J(\infty, t)\partial_t+\nu(\infty,t)=A_+.
\]

The contractibility of \(\Sigma_P(A_-, A_+)\) follows from well-known contractibility of the
space of complex structures, and the fact that \(\nu\) lies in a
vector space. 

Next, denote by \(\mathfrak{T}_x\) the space of unitary trivializations of
\(x^*K\) for \(x\in {\cal P}\) and let \(\mathfrak{T}=\coprod_{x\in
  {\cal P}}\mathfrak{T}_x\). This is a \(C^\infty(S^1, U(n))\)-bundle
  over \({\cal P}\). Fix a \(g\in C^\infty(S^1, U(n))\) and a section 
\(\Phi: {\cal P}\to \mathfrak{T}\), such that the inclusion
\(\mathfrak{S}:=\{g^k\Phi (x)\,
|\, k\in \Z, x\in {\cal P}\}\hookrightarrow \mathfrak{T}\) induces an isomorphism
\(i_\pi: \mathfrak{S}=\pi_0\mathfrak{S}\to
\pi_0\mathfrak{T}\). 

Recall from \S I.3.1.4 that we have a bundle map (over \({\cal P}\)) from
\(\tilde{\cal P}\) to \(\pi_0\mathfrak{T}\): from a fixed unitary
trivialization of \(\gamma_0^*K\) and a path \(w\subset \lc\) from \(\gamma_0\)
to \(x\in {\cal P}\), we extend the trivialization over \(w^*K\) to obtain a homotopy
class of trivializations of \(x^*K\). If \(w'\) is another path in the same
equivalence class, i.e. \(\op{im}[w-w']=0\), then \((w-w')^*K\) is
trivial, since \(c_1^f(\op{im}[w-w'])=0\). Hence \(w, w'\) induce the
same homotopy class of trivializations of \(x^*K\).
Composing this map with \(i_\pi^{-1}\), we have a map assigning to each
\((x, [w])\in \tilde{\cal P}\) a trivialization \(\Phi_{x, [w]}\in \mathfrak{S}_x\).
Let 
\[
{\mathbb A}_{(x, [w])}:=\Phi_{x, [w]}A_x\Phi_{x, [w]}^{-1}.
\]
We have a map \(m_P: {\cal B}_P^k((x, [w]), (y, [v]))\to
\Sigma_P({\mathbb A}_{(x, [w])}, {\mathbb A}_{(y, [v])})/G_P\),
defined as follows. 

A \(u\in {\cal B}_P((x, [w]), (y, [v]))\), together with a
trivialization \(\Phi_u\) of the
symplectic vector bundle \(u^*K\) that restricts respectively to \(\Phi_{(x, [w])}\)
and \(\Phi_{(y,[v])}\) at the circles at \(-\infty\) and \(\infty\),
assigns an element in \(\Sigma_P\).  Namely, \(\mathbb{E}_u:=\Phi_{u*}E_u\Phi_{u*}^{-1}\).

The space of such trivializations \(\Phi_u\) is an affine space under 
\[
G_P=\Big\{\Psi \, \Big|\, \Psi\in C^\infty([-\infty, \infty]\times S^1,
\op{Sp}_n), \Psi |_{\{\pm\infty\}\times S^1}=1\Big\}.
\]
Let \(m_P(u)\) be the \(G_P\)-orbit of \(\mathbb{E}_u\). 
It is shown in Lemma 13 of \cite{FH} that any orbit of \(G_P\) in
\(\op{\frak{Or}}(L\Sigma_P({\mathbb A}_{(x,[w])}, {\mathbb A}_{(y, [v])}))\) is contained in a single path component
of \(\op{\frak{Or}}(L\Sigma_P({\mathbb A}_{(x,[w])}, {\mathbb A}_{(y,
  [v])}))\); thus, \(L(\Sigma_P({\mathbb A}_{(x,[w])}, {\mathbb A}_{(y,
  [v])})/G_P)\) is trivial; hence so is \(L{\cal B}_P((x, [w]), (y,
[v]))\).

Notice that by definition, for any \((z_1, [v_1]), (z_2, [v_2]), (z_2,
[v_3])\in \tilde{\cal P}\), \((\mathbb{E}_1,
\mathbb{E}_2)\in\Sigma_P({\mathbb A}_{(z_1, [v_1])}, {\mathbb
  A}_{(z_2, [v_2])}))\times \Sigma_P({\mathbb A}_{(z_2, [v_2])},
{\mathbb A}_{(z_3, [v_3])})) \) is glueable. This gives rise to a
gluing homomorphism from \(\op{\frak or}(L{\cal B}_P((z_1,
[v_1]), (z_2, [v_2]))\times \op{\frak or}(L{\cal B}_P((z_2,
[v_2]), (z_3, [v_3]))\) to \(\op{\frak or}(L{\cal B}_P((z_1,
[v_1]), (z_3, [v_3]))\).

Lastly, observe that if \(\op{gr} ((x, [w), (y,[v]))=\op{gr} ((x, [w'],
(y,[v']))=k\), the spaces \(\Sigma_P({\mathbb A}_{(x,[w])}, {\mathbb A}_{(y,
  [v])}))\) and \(\Sigma_P({\mathbb A}_{(x,[w'])}, {\mathbb A}_{(y,
  [v'])}))\) may be identified via conjugation by \(\bar{g}^i\) for
some \(i\in \Z\) and \(\bar{g}\in C^\infty([-\infty, \infty]\times S^1,
\op{Sp}_n)\), \(\bar{g}(s,t):=g(t)\).
Thus, the above discussion in fact verifies the orientability of
\(L{\cal B}_P^k(x, y)\) for any \(x, y\in {\cal P}\), \(k\in \Z\), and
the gluing homomorphism above gives the gluing homomorphism
\(\mathfrak{g}\) described in \S7.2.3.

\medbreak
\noindent{\bf (2)} {\sc Orienting \(L{\cal B}_P^{\Lambda}.\)}
Suppose \((J^\Lambda, X^\Lambda)\) generates a CHFS satisfying
the properties (RHFS1*). 
Let \(({\bf x, [w]}), ({\bf y, [v]})\) be two path components of \(\tilde{\cal
  P}^\Lambda/\tilde{\cal P}^{\Lambda, deg}\). 

The deformation operator of \(\cm_P^\Lambda\)
at \(u_\lambda\),
 \(\hat{E}_{u_\lambda}\), may be regarded as a stabilization of
 \(E_{u_\lambda}\). Because of the stabilization isomorphism for
 families, to orient the determinant line bundle 
 \(\det \{\hat{E}_{u_\lambda}\}_{u_\lambda\in {\cal B}_P^{\Lambda
    }(({\bf x, [w]}), ({\bf y, [v]}))}\), it is equivalent to orient the
 determinant line bundle 
\[
L{\cal B}_P^{\Lambda}(({\bf x, [w]}), ({\bf y, [v]})):= \det\{E_{u_\lambda}\}_{u_\lambda\in {\cal B}_P^{\Lambda}(({\bf x, [w]}), ({\bf y, [v]}))}.
\]
This can be oriented by repeating part (1) above, replacing
\(\Sigma_P\) by the parameterized version:
\[
\Sigma_P({\bf A}_-, {\bf A}_+):=\bigcup_{\lambda\in
\Lambda_{\bf x}\cap \Lambda_{\bf y}}\Sigma_P({\mathbb A}_{(x_\lambda,
[w_\lambda])}, {\mathbb A}_{(y_\lambda, [v_\lambda])}),
\]
which is a \(\Sigma_P\)-bundle over \(\Lambda_{\bf x}\cap \Lambda_{\bf
  y}\). In the above, \(\mathbb{A}_{(x_\lambda, [w_\lambda])}\) \(\forall
(x_\lambda, [w_\lambda])\in \tilde{\cal P}\) is defined via
a smooth section \(\Phi^\Lambda: {\cal P}^\Lambda\to
\mathfrak{T}^\Lambda\),
\(\mathfrak{T}^\Lambda:=\bigcup_{x_\lambda\in{\cal
    P}^\Lambda}\mathfrak{T}_{x_\lambda}\).
This is again a contractible space, since it is a bundle with
contractible fibers and base. 

As in (1), this in turn demonstrates the orientability of \(L{\cal
  B}^{\Lambda, k+1}({\bf x}, {\bf y})\), and defines the parameterized
version of gluing homomorphism \(\mathfrak{g}^\Lambda\) described in
\S7.2.3. Now one may follow the arguments in \S7.2.3 to define a
coherent orientation of \(L{\cal B}_P^\Lambda\).

\medbreak
\noindent{\bf (3)} {\sc Orienting $L{\cal B}_O^{1}$ and \(L{\cal
    B}_O^{\Lambda, 2}\).}
Since \(\tilde{D}_{(T, u)}\) is a rank 1 stabilization of \(D_{(T, u)}\), it is
equivalent to orient \(\underline{L}{\cal B}_O^1:=\det \{D_{(T,
  u)}\}_{(T, u)\in{\cal B}_O^1}\).

Similarly to parts (1), (2) above, we
introduce a map \(m_O: {\cal B}_O^1\to \Sigma_O^1/G_O\), where
\(\Sigma_O^1\) is the space of rank-1 stabilizations of operators of the form:
\[\bar{\partial}_{J, \nu; T}:=
\partial_s+J(s, t)\partial_t+\nu(s,t) \, : \, L^p_1(S^1_T\times S^1; \R^{2n})\to
L^p(S^1_T\times S^1; \R^{2n})\quad \text{for some \(T\in \R^+\),}
\]
with \(J, \nu\) defined similarly to part (1). The determinant line bundle
\(L\Sigma_O^1\) is canonically oriented as follows: 
Note that \(\Sigma_O^1\) contracts to the
subspace consisting of complex linear \(\bar{\partial}_{J, \nu; T}\), 
which we denote by
\(\Sigma'_O\). However, \(L\Sigma'_O\) is canonically oriented by the complex
linearity of kernels and cokernels.
Next, note that \(u^*K\) is trivial for any \((T, u)\in {\cal
  B}_O^1\). Given a \((u, T)\in {\cal B}_O^1\) and a
trivialization \(\Phi_u\) of \(u^*K\), one has
\[\tilde{\mathbb{D}}_{(T, u)}:=\Phi_{u*}\tilde{D}_{(T,
    u)}\tilde{\Phi}_{u*}^{-1}\in \Sigma_O, 
\] 
where \(\tilde{\Phi}_{u*}^{-1}:=1\oplus \Phi_{u*}^{-1}\) as an
endomorphism of \(\R\oplus L^p_1(S^1_T\times S^1; \R^{2n})\).
This defines \(m_O\).

It is not hard to see that \(G_O\) acts trivially on the \(\Z/2\Z\)-bundle 
\(\op{\frak{Or}}(L\Sigma_O^1)\) by conjugation:  by continuation
(cf. the commutative diagram in p.28 of \cite{FH}), it suffices to check this
for \(\op{\frak{Or}}(L\Sigma'_O)\). However, for 
\(\op{\frak{Or}}(L\Sigma'_O)\) this is obvious, again by the complex
linearity of elements in \(\Sigma'_O\).

The orientability of \(L{\cal B}_O^{\Lambda, 2}\) follows immediately
from that of \(L{\cal B}_O^{1}\), 
since \({\cal B}_O^{\Lambda, 2}={\cal B}_O^1\times \Lambda\) by definition.

Finally, note that for this version of symplectic Floer theory, the
canonical orientation of \(L\Sigma_O^1/G_O\) is compatible with the
mod 2 Conley-Zehnder index \(\ind\): the former is given by the
canonical orientation of \(\det \tilde{D}_{\mathbb{A}_0, T}\), where
\(\mathbb{A}_0\) is such that \(D_{\mathbb{A}_0, T}\) is complex
linear. By the definition of Conley-Zehnder index (cf. \S I.3.1.4),
\(\op{CZ}(\mathbb{A}_0)\) is even.

\subsection{The Signs.}
It was shown in \cite{FH} that with a coherent orientation for
\(\cm_P\), the Floer complex indeed satisfies \(\tilde{\partial}_F^2=0\). 
In this subsection, we generalize this result to the setting of CHFSs and
verify the second statement of (RHFS4). Namely, we
show that with \(\cm_P^\Lambda\) endowed with
coherent orientations and \(\cm_O^{\Lambda,2}\) endowed with
the grading-compatible orientation, the various
0-dimensional strata \(\mathbb{J}_P\), \(\mathbb{T}_{P, hs-s}\), \(\mathbb{T}_{P, db}\)
in \(\hat{\cm}_P^{\Lambda, 1, +}\) and their analogs for broken orbits
are expressed in terms of products of 0-dimensional 
moduli spaces, with relative signs
given by the formulae (I.28--33). 

As the signs for \(\mathbb{J}_P\), \(\mathbb{J}_O\) given in (I.28,
31) follow immediately from the definition of coherent orientation, 
we shall concentrate on the signs for
\(\mathbb{T}_{P, hs-s}\), \(\mathbb{T}_{O,
  hs-s}\), \(\mathbb{T}_{P, db}\), and \(\mathbb{T}_{O, db}\): the
formulae (I.30, 33, 29, 32) are respectively 
rephrased in terms of the gluing theorems 
Propositions \ref{8.1}, \ref{handleslide} in Lemmas  
7.3.2--7.3.6 below.

We assume throughout this subsection that \(L\Sigma_P/G_P,
L\Sigma_O/G_O\) and their parameterized versions are endowed with
coherent orientations/ grading compatible orientations, and all the
oriented K-models are compatible with these orientations, unless
otherwise specified. 

The results and arguments in this subsection
apply to general Floer theory, in which the relevant moduli spaces are
oriented according to the scheme in \S7.2.1--4 above.

\subsubsection{Preparations.}

\noindent {\bf (1)} {\sc Signs of flowlines.}
The sign of a flow in a 0-dimensional reduced moduli space, \(\hat{u}\in \cm^1/\R\),
in general means the relative sign \([u']/\ker \mathfrak{D}_u\) for any
representative \(u\in \cm^1\) in the un-reduced moduli space, where
\(\mathfrak{D}_u\) is the deformation operator of \(\cm^1\). It will
be denotedy by \(\sign (u)\).

\medbreak
\noindent{\bf (2)} {\sc Trivializations of deformation operators.}
Instead of working with the deformation operators \(E_u,
\tilde{D}_{(T,u)}\) and their parameterized versions, it is often more
convenient to work with their corresponding operators in \(\Sigma_P\)
or \(\Sigma_O\) via lifts of the maps \(m_P, m_O\). These will be
denoted by boldface letters such as \({\mathbb E}_u, \tilde{\mathbb
  D}_{(T, u)}\). When \(L\Sigma_P/G_P\), \(L\Sigma_O/G_O\) are
orientable, the choice of liftings does not matter.
We shall also omit specifying
the class \([v]\) in \({\mathbb A}_{(y, [v])}=:{\mathbb
  A}_y\), when the precise choice is immaterial.

For the symplectic Floer theory discussed in this article, this means
replacing the deformation operators by their conjugates by
trivializations of \(u^*K\), namely \(\Phi_{u}\) (cf. the definition
of \({\mathbb E}_u\), \({\mathbb D}_u\) in \S7.1.5). 
We write \((f)_\Phi:=\Phi_{u*}f\), e.g. \(f=u'\) for \(u\in {\cal B}_P\) or \({\cal B}_O\).

The families of operators considered in the rest of this subsection
will always be subfamilies of various versions of \(\Sigma_P\),
\(\Sigma_O\). Thus, we shall refer to the correlation and relative
signs of orientations of determinant lines, or mutual co-orientation
of K-models without specifying the family. 

The following consistency conditions on the choice of
liftings will be assumed in the following discussion: 
\begin{description}\itemsep -2pt
\item [{\rm (a)}] For a subfamily \(U\subset {\cal B}\), the lifting
  \(\tilde{m}: U\to \Sigma\) is continuous;
\item [{\rm (b)}] the liftings are ``coherent'' in the sense that they
  are consistent with gluing.
\end{description}

\subsubsection{Signs for \(\mathbb{T}_{P, hs-s}\).}
To verify the sign in (I.30), we need to examine oriented K-models for
the gluing theorem, Proposition \ref{handleslide} (a). 
Let \((J^\Lambda, X^\Lambda)\) be an admissible \((J,
X)\)-homotopy, for any \({\bf x}_1, {\bf x}_2\in \aleph_\Lambda\) and \(R_0\gg1\), the
(omitted) proof of Proposition \ref{handleslide} (a) defines a gluing
map
\[
\op{Gl}_{P, hs}(\mathbf{x}_1,\mathbf{x}_2;\Re): 
\mathbb{T}_{P, hs-s}(\mathbf{x_1},\mathbf{x_2};\Re)\times
(R_0, \infty)\to \cm_P^{\Lambda,2} (\mathbf{x}_1,\mathbf{x}_2;
\op{wt}_{-{\cal Y}, e_{\cal P}}\leq \Re).
\]
Let \((\lambda_0, \hat{u})\in \hat{\cm}_P^{\Lambda, 0} (\mathbf{x}, \mathbf{y}; J^\Lambda, X^\Lambda)\) be a handleslide.
Without loss of generality, assume \(\lambda_0=0\).
Let $\mathbf{q},\mathbf{z} \in \aleph_\Lambda$ be of
indices $\ind \mathbf{y}+1$ and $\ind \mathbf{y}-1$ respectively.
Let \[\hat{v}_-\in \hat{\cm}_P^0 (q_0,x_0; J_0, X_0), \quad  
\hat{v}_+\in \hat{\cm}_P^0 (y_0,z_0;J_0,X_0)\]
and \(v_-, v_+, u\) be centered representatives of \(\hat{v}_-,
\hat{v}_+, \hat{u}\) respectively. Let
\[w_{\# -}(R)=v_-\#_Ru,\quad w_{\#, +}(R)=u\#_Rv_+\]
be the pregluings defined in \S1.2.2, and let
\[\begin{split}
(\lambda_-(R), w_-(R)):= &\op{Gl}_{P,
  hs}(\mathbf{q},\mathbf{y};\Re)(\{\hat{v}_-, \hat{u}\}, R), \\
(\lambda_+(R),w_+(R)):=&\op{Gl}_{P,
  hs}(\mathbf{x},\mathbf{z};\Re)(\{\hat{u}, \hat{v}_+\},
R).\end{split}\]
be the images of the gluing map obtained by further perturbing \(w_{\#
  -}(R)\) and \(w_{\# +}(R)\) respectively.
To simplify notation, we shall omit \(R\) when there is
no danger of confusion.
\begin{lemma*}\label{lemma:handleslide-sign}
Let \(u, v_\pm, (\lambda_\pm, w_\pm)\) be as above. Then
\begin{equation}\begin{split}{\rm (1)} &\, - \sign (\lambda_-)\sign(w_-)=\sign(v_-)\sign(u)\\
\label{sign:2}
{\rm (2)} &\, \sign(\lambda_+)\sign(w_+)=\sign(u)\sign(v_+).
\end{split}\end{equation}
\end{lemma*}
\begin{proof}
We shall focus on case (1) below, since case (2) is entirely parallel.
According to \S7.1.3, 7.1.5, and the choice of coherent orientation,
we have the oriented K-models:
\begin{description}
\item[{\rm (i\(\pm\))}] \([\, \hat{\mathbb{E}}_{(0, w_{\#\pm})}: \hat{K}_{\#\pm}\to C_{\#\pm}]\),
where
\begin{gather*}
\hat{K}_{\#-}=-\R\oplus (\ker \mathbb{E}_{v_{- }}\#_{R}\ker
\mathbb{E}_{u}), \quad C_{\#-}=*\#_R\cok \mathbb{E}_{u};\\
\hat{K}_{\#+}=\R\oplus (\ker \mathbb{E}_{u}\#_R\ker
\mathbb{E}_{v_+}), \quad C_{\#+}=\cok \mathbb{E}_{u}\#_R*.
\end{gather*}
\item[{\rm (ii)}] \([\, \hat{\mathbb{E}}_{(0, u)}:  \R\oplus (\ker \mathbb{E}_u)
 \to \cok \mathbb{E}_u]\).
\end{description}
Since \(u\) is by assumption a nondegenerate element of
\(\cm_P^{\Lambda, 1}\), the standard oriented K-model for
\(\hat{\mathbb{E}}_{(0, u)}\) may be viewed as a reduction of the
oriented K-model (ii) by \(-\R\), taking \[\begin{split}
\cok \mathbb{E}_u &=-\R (Y_{(0, u)})_\Phi,\quad \text{and}\\
\ker \mathbb{E}_u &=\ker \hat{\mathbb{E}}_{(0, u)}=\sign
(u) \R (u')_\Phi\end{split}\] 
as oriented spaces. (Recall that \(Y_{(0, u)}\) is a cutoff version of
\(\partial_\lambda {\cal V}_\lambda\) appearing in the definition of
\(\hat{\mathbb{E}}_{(0, u)}\), cf. I.6.1.5.) 

Next, decompose \(\hat{K}_{\#-}\) into the direct sum of the ordered
triple of oriented subspaces 
\[*\oplus (\ker \mathbb{E}_{v_{-}}\#_{R}*), \quad  \R\oplus *, \quad
\text{and} \,\, *\oplus (*\#_{R}\ker
\mathbb{E}_{u}), .\]
By Lemma \ref{glue-K} (2), for large \(R\), 
the restriction of \(\Pi_{C_{\#-}}\hat{\mathbb{E}}_{(0,w_{\#-})}\) 
to the first and last subspaces are small, while its restriction to the second
subspace approximates the multiplication by \(\Pi_{\cok
  \mathbb{E}_u}\tilde{Y}\) (under the natural identification of 
the domain and range spaces), where \(\tilde{Y}\) is another cutoff
version of \(\partial_\lambda {\cal V}_\lambda\) which agrees with
\((Y_{(0, u)})_\Phi\) except in the region where \(s\ll-1\). Let \(Y_\nu:=\nu
(Y_{(0, u)})_\Phi+(1-\nu)\tilde{Y}\) for \(\nu\in [0,1]\), and \(\hat{\mathbb{E}}_\nu\) be the rank 1
stabilization of \(\mathbb{E}_u\) by multiplication by \(Y_\nu\). By
the surjectivity of \(\partial_s+\mathbb{A}_y\) and
\(\hat{\mathbb{E}}_0=\mathbb{E}_{(0, u)}\), and an excision argument
(as outlined in \S1.2.5), \(\hat{\mathbb{E}}_\nu\) has uniformly
bounded right inveres. 
Thus, we may conclude that \(\Pi_{\cok \mathbb{E}_u}\tilde{Y}\), and hence also
\(\Pi_{C_{\#-}}\hat{\mathbb{E}}_{(0,w_{\#-})}|_{\R\oplus *}\), are
positive of \(O(1)\). This implies that the reduction of the
oriented K-model (i\(-\)) by \(-\R\oplus *\) is equivalent to the
standard oriented K-model of \(\hat{\mathbb{E}}_{(0,w_{\#-})}\), which
is in turn equivalent to the standard oriented K-model of 
\(\hat{\mathbb{E}}_{(\lambda_-,w_-)}\), due to the the proximity
between \(w_-\) and \(w_{\#-}\). In other words, the projection 
\begin{equation*}
\Pi_{K_{\#-}}: \ker \hat{\mathbb{E}}_{(\lambda_-,w_-)}\to \ker
\mathbb{E}_{v_{-}}\#_R\ker\hat{\mathbb{E}}_{(0,u)}=:K_{\#-}.\label{isom:w-}
\end{equation*}
is an orientation-preserving isomorphism.
We have the following ordered bases compatible with the former and
latter oriented spaces above: 
\[
\Big\{\sign (w_-)(0, w_-'), \, (\p_R\lambda_-)^{-1}(\p_R\lambda_-,\p_R w_-)\Big\}, \quad \Big\{\sign
(v_-)(0, v'_-), \, \sign (u) (0,u')\Big\}.
\] 
Observing that \(\sign (\p_R\lambda_-)=-\sign (\lambda_-)\), and recalling the
description of \(\Pi_{\iota_\#^jK_J}\) in \S1.2.4, one finds that
with respect to these bases, \(\Pi_{K_{\#-}}\) is a matrix of the form
\[\left(\begin{array}{cc}
C_1\sign (w_-)\sign(v_-) &C_1'\sign (\lambda_-)\sign(v_-) \\
C_2 \sign (w_-)\sign(u) & -C_2'\sign (\lambda_-)\sign(u) 
\end{array}\right)
\quad \text{for \(C_1, C_1', C_2, C_2'\in \R^+\).}
\]
Equation (\ref{sign:2}.1) follows from the requirement that this
matrix has positive determinant. 
\end{proof}

\subsubsection{Signs for \(\mathbb{T}_{O, hs-s}\).}
To verify the sign in (I.33), we examine the oriented K-model for the
gluing theorem Proposition \ref{handleslide} (b).
Let $y, u$ be as in \S\ref{lemma:handleslide-sign}, but now assume
that the handleslide \(u\) is of Type II, namely, \(x=y\). Let
\(w_\#(R)=u_{\#R}\) be the glued orbit introduced in \S1.2.2, and 
\((\lambda (R), (T(R), w(R))):=\op{Gl}_{O, hs}(\hat{u}, R)\) be the
image of the gluing map obtained by perturbing \(w_\#(R)\),  
\(\op{Gl}_{O, hs}: \mathbb{T}_{O, hs-s}(\Re)\times (R_0, \infty)\to
\hat{\cm}_O^{\Lambda, 1}(\op{wt}_{-{\cal Y}, e_{\cal P}}\leq \Re)\) being the
gluing map in the (omitted) proof of Proposition \ref{handleslide}
(b).
\begin{lemma*}
In the above notation, 
\(\label{handleslide-o-sign}
\sign (w)=(-1)^{\ind y} \sign (\lambda)\sign(u).
\)
\end{lemma*}
\begin{proof}
According to Corollary 7.2.4, \(\sign (w)=(-1)^{\ind y}
[(w')_\Phi]/\ker^{\mathfrak{o}_y} \tilde{\mathbb{D}}_{(T,w)}\), where
\(\ker^{\mathfrak{o}_y} \tilde{\mathbb{D}}_{(T,w)}\) denotes \(\ker
\tilde{\mathbb{D}}_{(T,w)}\) endowed with the orientation correlated
to the canonical orientation of \(\tilde{\mathbb{D}}_{\mathbb{A}_y,
  T}\).
We compute \([(w')_\Phi]/\ker^{\mathfrak{o}_y}
\tilde{\mathbb{D}}_{(T,w)}\) in two steps.

\medbreak

\noindent{\sc Step 1. The relative sign \([(w')_\Phi]/\ker^{\mathfrak{o}_y}
  \acute{\mathbb{D}}_{(\lambda,w)}\).} Let \(\acute{D}_{(\lambda, w)}\) be
the rank 1 stabilization of \(D_w\) defined by 
\[
\acute{D}_{(\lambda, w)}(\alpha, \xi)=
\alpha\p_\lambda {\cal V}_\lambda(w)+D_{w}\xi.\]
Perform cyclic gluing to the oriented K-model \([
\hat{\mathbb{E}}_{(\lambda, u)}: \R\oplus \ker\mathbb{E}_u\to \cok
\mathbb{E}_u]\), we obtain an oriented K-model for
\(\acute{\mathbb{D}}_{w_\#}\), a rank 1 stabilization of
\(\mathbb{D}_{w_\#}\) by multiplication with a cutoff version of
\(\partial_\lambda{\cal V}_\lambda (w_\#)\). The argument in \S7.3.2
shows that a reduction of this oriented K-model by \(\R\) is
equivalent to a standard K-model for \(\acute{\mathbb{D}}_{w_\#}\),
which is in turn equivalent to a standard K-model for 
\(\acute{\mathbb{D}}_{(\lambda, w)}\). Moreover, 
according to the continuity of gluing homomorphisms and Lemma 7.2.2, 
the orientation of this standard K-model is correlated to the
canonical orientation of \(\tilde{\mathbb{D}}_{\mathbb{A}_y, T}\). In
other words, 
\[[(w')_\Phi]/\ker^{\mathfrak{o}_y} \acute{\mathbb{D}}_{(\lambda,w)}=\sign (u).  \]

\noindent{\sc Step 2. The relative sign \(\ker^{\mathfrak{o}_y} \acute{D}_{(\lambda,w)}/\ker^{\mathfrak{o}_y} \tilde{D}_{(T,w)}\).} Notice that the operators 
\(\acute{D}_{(\lambda, w)}\), \(\tilde{D}_{(T,w)}\) have a common
stabilization, \(\hat{D}_{(\lambda,(T,w))}=(\partial_\lambda{\cal
  V}_\lambda, (-w'/T, D_w))\). 
The two dimensional space \(\ker\hat{D}_{(\lambda,(T,w))}\) is spanned
by \(\{\partial_R(\lambda , (T, w)), \, (0, (0, w'))\}\). This means 
\(\Pi_{\cok D_w}\Big(\partial_R\lambda\partial_\lambda{\cal V}_\lambda+\partial_RT(-w'/T)\Big)=0\), and hence the relative
sign is computed by
\[
\ker^{\mathfrak{o}_y} \acute{D}_{(\lambda,w)}/\ker^{\mathfrak{o}_y} \tilde{D}_{(T,w)}
=-\sign(\partial_R\lambda)/\sign(\partial_R T)=\sign (\lambda).
\]

Finally, the Lemma is obtained taking the product of the relative
signs obtained in Steps 1 and 2 above with \((-1)^{\ind y}\).
\end{proof}

\subsubsection{Signs for \(\mathbb{T}_{P, db}\).}
To verify the signs in (I.29), we need to analyze the orientation of
the K-model for the gluing theorem Proposition
\ref{8.1} (a). Let \(\lambda, u_0, \ldots, u_{k+1}\) be as in \S2.2,
and let \((\lambda,w)\) be the image of \((\{u_0, u_1, \ldots,
u_{k+1}\}, \lambda)\) under the gluing map defined in \S\ref{gluing-map:2.1}.  
\begin{lemma*} 
Under the assumptions in \S2.1 and in the above notation, 
\[\label{sign:deg-gluing-a}
\sign(w)=(-1)^{k+1}
\prod_{i=0}^{k+1}\sign (u_i).
\]
\end{lemma*}
\begin{proof}
As explained earlier in this section, since we work with the ordinary
Sobolev norms instead of the complicated polynomially-weighted ones, 
it is convenient to replace the delicate pregluing \(w_\chi\) defined
in \S2.2 by ordinary glued trajectories or orbits: 
Choose \(R_i\), \(R'_i\), \(i\in\{1, \ldots, k+1\}\) and \(L\)
appropriately so that:
\[
w_\#:=\tau_L \Big( u_0\#_{R_1}\bar{y}\#_{R_1'}u_1\#_{R_2}\cdots\#_{R_{k+1}'}u_{k+1}\Big)
\]
is pointwise close to \(w\) and \(w_\chi\): more precisely, \(w_\#(s)\),
\(w(s)\), \(w_\chi(s)\) are \(C_\epsilon\)-close to each other
\(\forall s\), and
\[
w_\chi(\gamma_{u_i}^{-1}(0))=w_{\#}(\gamma_{u_i}^{-1}(0)).
\] 
As explained in \S7.2.1 and \S7.2.3, the deformation
operator in \(\Sigma_P\) corresponding to \(w_\#\) is:
\[
\mathbb{E}_{w_\#}^{\sigma}:=\mathbb{E}_{u_0}^{[0,
  \sigma]}\#_{R_1}\mathbb{E}_{\bar{y}}^{[\sigma,
  -\sigma]}\#_{R_1'}\mathbb{E}_{u_1}^{[-\sigma,
  \sigma]}\#_{R_2}\mathbb{E}_{\bar{y}}^{[\sigma,
  -\sigma]}\#_{R_2'}\cdots\#_{R_k'}\mathbb{E}_{u_{k+1}}^{[-\sigma,
  0]}\quad \text{for a small \(\sigma>0\).}
\]
Let \(\bar{\bf e}_y^\sigma(s):= \varsigma^{-\sigma, \sigma}(s){\bf
  e}_y\), and recall that 
\begin{equation*}\label{K-C-y-sigma}
\ker \mathbb{E}_{\bar{y}}^{[-\sigma, \sigma]}=\cok \mathbb{E}_{\bar{y}}^{[\sigma,
  -\sigma]}=\R \bar{\bf e}_y^\sigma; \quad 
\cok \mathbb{E}_{\bar{y}}^{[-\sigma, \sigma]}=\ker \mathbb{E}_{\bar{y}}^{[\sigma,
  -\sigma]}=*.
\end{equation*}
Then by Lemma \ref{glue-K}, we have
the following oriented K-model for \(\mathbb{E}_{w_\#}\) compatible with
the coherent orientation:
\[\begin{split}
[\, \mathbb{E}_{w_\#}^{\sigma}: K_\#\to C_\#], \quad \text{where
  \(C_\#=*\#_{R_1}\R\bar{\bf e}_y^\sigma\#_{R_1'}* \cdots \#_{R_k'}*\)}, \qquad \\
K_\#=(-1)^{k+1}\ker \mathbb{E}_{u_0}^{[0, \sigma]}\#_{R_1}*\#_{R_1'}\ker \mathbb{E}_{u_1}^{[-\sigma,
  \sigma]} \#_{R_2}*\#_{R_2'}\cdots\#_{R_k'}\ker \mathbb{E}_{u_{k+1}}^{[-\sigma,
  0]} 
\end{split}\]
On the other hand, in section 3, we constructed
the following K-model: 
\[
[\mathbb{E}_{w_\chi}: K_\chi\to C_\chi]=\Big[\R{\mathfrak
  e}_{u_0}\oplus\cdots\oplus\R{\mathfrak
  e}_{u_{k+1}} \to \R{\mathfrak f}_1\oplus\cdots\oplus\R{\mathfrak
  f}_{k+1}\Big].
\] 
(\(\mathbb{E}_{w_\chi}\) is now considered as an operator between
ordinary Sobolev spaces.
As remarked before, the polynomially weighted spaces
are commensurate with the ordinary Sobolev spaces, and we do not need
uniform boundedness of right inverses in this section. We have also
suppressed the subscript \(\Phi\) and written \({\mathfrak
  e}_{u_i}=({\mathfrak
  e}_{u_i})_\Phi, {\mathfrak f}_j=({\mathfrak f}_j)_\Phi\) above for
simplicity.) 

Using the descriptions of \(\Pi_{K_\chi}\) and \(\Pi_{C_\chi}\) given
in \S\ref{glue-K}  and Proposition \ref{prop:P-j} and the proximity between
\(\mathbb{E}_{w_\#}^{\sigma}\) and \(\mathbb{E}_{w_\chi}\), 
one may easily check that the oriented K-model 
\([K_\#\to C_\#]\) is equivalent to 
\[
\Big[\, (-1)^{k+1}\prod_{i=0}^{k+1}\sign (u_i) K_\chi\to C_\chi\, \Big],
\]
implying that the latter is also compatible with the coherent orientation. 

Next, observe that \((w')_\Phi\) projects positively to all
\({\mathfrak e}_{u_i}\). This, together with the form of
\(\Pi_{C_\chi}\mathbb{E}_\chi|_{K_\chi}\) given in Lemma \ref{est:cok} (b),
implies that the reduction of the above oriented K-model, 
\[\Big[(-1)^{k+1}\prod_{i=0}^{k+1}\sign (u_i)\, \R (w')_\Phi\to\, 
*\, \Big], \] is equivalent to the standard oriented K-model for 
\(\mathbb{E}_{w_\chi}\), which is in turn equivalent to the standard
oriented K-model for \(\mathbb{E}_w\), due to the proximity between
\(w\) and \(w_\chi\). These observations immediately imply the Lemma.
\end{proof}

\subsubsection{Signs for \(\mathbb{T}_{O, db}\).}
To verify the sign in (I.32), we examine the orientation of the
K-model in the proof of Proposition \ref{8.1} (b). Let \(\{\hat{u}_1,
\ldots, \hat{u}_k\}\) be a broken orbit, \(u_i\) be the centered
representative of \(\hat{u}_i\), and \((\lambda,(T,w))\) be the image of \((\{\hat{u}_1, \ldots,
\hat{u}_{k}\}, \lambda')\) under the gluing map defined in \S\ref{gluing-map:2.1b}. 
\begin{lemma*}\label{sign:deg-gluing-b}
Under the assumptions in \S2.1 and in the above notation, 
\[
\sign(w)=(-1)^{\ind_- y+k}
\prod_{i=1}^{k}\sign (u_i).
\]
\end{lemma*}
\begin{proof}
As argued in the proof of Lemma 7.3.3, \begin{equation}\label{sign4}
\sign (w)=(-1)^{\ind_- y}\sign
(\lambda)(w')_\Phi/\ker^{\mathfrak{o}_{y-}}\acute{\mathbb{D}}_{(\lambda,w)},\end{equation}
where the superscript \(\mathfrak{o}_{y-}\) indicates the orientation
correlated with the canonical orientation of
\(\tilde{\mathbb{D}}_{\mathbb{A}_y+\sigma, T}\). (Recall the definition
\(\ind_-y=\ind (\mathbb{A}_y+\sigma)\).) According to the assumption
of \S2.1, \(\sign(\lambda)=1\).

Instead of working with the standard K-model for
\(\tilde{\mathbb{D}}_{\mathbb{A}_y+\sigma, T}\), we find it easier to work
with the following mutually co-oriented K-model: 
\([\acute{\mathbb{D}}_y: -({\bf e}_y)_\Phi\to *]\), where
\(\acute{\mathbb{D}}_y\) is the stabilization of
\(\mathbb{D}_{\mathbb{A}_y, T}\) by multiplication with \(({\bf e}_y)_\Phi\).
To see that they are indeed co-oriented, observe that the
interpolation between them, 
\(
 \acute{\mathbb D}_{\nu}=
((1-\nu) ({\bf e}_{y})_\Phi, 
 {\mathbb D}_{\mathbb{A}_y+\nu \sigma, T})
\) are surjective $\forall \nu\in [0,1]$,  and has the following
continuous basis for the kernel: \(\xi_\nu:=(\nu\sigma,
-(1-\nu)({\bf e}_{y})_\Phi).\)

We now consider mutually co-oriented K-models for two operators 
approximating \(\acute{\mathbb{D}}_{(\lambda, w)}\) and
\(\acute{\mathbb{D}}_y\) respectively. The proximity of the
operators implies that these K-models also form
mutually co-oriented K-models for to \(\acute{\mathbb{D}}_{(\lambda, w)}\) and
\(\acute{\mathbb{D}}_y\) respectively. Choose an glued orbit
\[w_{\#}=\tau_L(\bar{y}\#_{R_1}u_1\#_{R_1'}\bar{y}\#_{R_2}\cdots\#_{R_k}u_k\#_{R_k'})\]
appproximating \(w\) and \(w_\chi\) pointwise in the sense of
\S7.3.4, and let 
\begin{gather*}\label{gl-Dw}
\begin{split}\acute{\mathbb{D}}_{w\#}^\sigma:=\Big(\p_\lambda{\cal V}_\lambda(w_{\#}), \, \mathbb{E}_{\bar{y}}^{[\sigma,-\sigma]}\#_{R_1} \mathbb{E}_{u_{1}}^{[-\sigma, \sigma]}\#_{R_1'}
\mathbb{E}_{\bar{y}}^{[\sigma,-\sigma]}\#_{R_2} \mathbb{E}_{u_{2}}^{[-\sigma, \sigma]}\cdots
\mathbb{E}_{u_{k}}^{[-\sigma, \sigma]}\#_{R_k'}\Big)
\end{split}\\\label{gl-Dy}
\begin{split}\acute{\mathbb{D}}_{y\#}^\sigma:=\Big(-({\bf e}_y)_\Phi,
  \, 
\mathbb{E}_{\bar{y}}^{[\sigma, -\sigma]}\#_{R_1} \mathbb{E}_{\bar{y}}^{[-\sigma, \sigma]}\#_{R_1'}\mathbb{E}_{\bar{y}}^{[\sigma,
  -\sigma]}\#_{R_2} \mathbb{E}_{\bar{y}}^{[-\sigma, \sigma]}\cdots
\mathbb{E}_{\bar{y}}^{[-\sigma, \sigma]}\#_{R_k'}\Big). \\
\end{split}
\end{gather*}
Since \([\, \mathbb{E}_{u_i}^{[-\sigma, \sigma]}: \ker \mathbb
\mathbb{E}_{u_i}^{[-\sigma, \sigma]}\to  *]\) and
\([\mathbb{E}_{\bar{y}}^{[-\sigma, \sigma]}: \R \bar{\bf
  e}_y^\sigma\to  *]\) are mutually co-oriented K-models (by coherent
orientation), the continuity of gluing homorphisms and stabilization
imply that we have the mutually co-oriented K-models
\([\acute{\mathbb{D}}_{w\#}^\sigma: \hat{K}_{w\#}\to C_{\#}]^{\mathfrak{gl}}\), 
\([\acute{\mathbb{D}}_{y\#}^\sigma: \hat{K}_{y\#}\to C_{\#}]^{\mathfrak{gl}}\), where
\begin{gather*}\begin{split}
\hat{K}_{w\#}= &
(-1)^{k+1}\R\oplus (*\#_{R_1}\ker \mathbb{E}^{[-\sigma, \sigma]}_{u_{1}}\#_{R_1'}*\#_{R_2}\ker \mathbb{E}^{[-\sigma, \sigma]}_{u_{2}}\cdots\ker
\mathbb{E}^{[-\sigma, \sigma]}_{u_{k}}\#_{R_k'})
\\
C_{\#}= & (\R \hat{\bf
  e}_y^\sigma)\#_{R_1}*\#_{R_1'}(\R \hat{\bf e}_y^\sigma)\cdots*\#_{R_k'}, \\
\hat{K}_{y\#}= &
(-1)^{k+1}\R\oplus (*\#_{R_1}(\R \hat{\bf e}_y^\sigma)\#_{R_1'} *\#_{R_2}(\R \hat{\bf
  e}_y^\sigma) \cdots(\R \hat{\bf e}_y^\sigma)\#_{R_k'}).\\
\end{split}\end{gather*}
Note that the orientation of these K-models is different from the
grading-compatible orientation, or the \(\mathfrak{o}_{y-}\)
orientation. We call it the ``glued orientation'', indicated by the
superscript \({\mathfrak{gl}}\) above. 

We now compute the sign of \((w')_\Phi\) relative to the
\(\mathfrak{gl}\)-orientation above. 
As in \S7.3.4, this is done by comparing the glued K-model above to
\(
[\acute{\mathbb{D}}_{(\lambda, w_\chi)}: \hat{K}_\chi\to
C_\chi]:=[\R_\alpha  \oplus \R{\mathfrak e}_{u_1}\oplus\cdots\oplus\R{\mathfrak
  e}_{u_{k}}\to  \, \R{\mathfrak f}_1\oplus\cdots\oplus\R{\mathfrak f}_{k}]
\), constructed previously in \S4.3. 
In this case, we find \([\acute{\mathbb{D}}_{w\#}^\sigma: \hat{K}_{w\#}\to
C_{w\#}]^{\mathfrak{gl}}\) equivalent to
\[[\acute{\mathbb{D}}_{(\lambda,w_\chi)}: (-1)^{k+1}\Pi_i(\sign (u_i))\, \hat{K}_\chi\to C_\chi].\]
On the other hand, in Lemma \ref{est:cok} (b),  
$\Pi _{C_\chi}\acute{\mathbb{D}}_{w_\chi}|_{\hat{K}_\chi}$ is computed in
the bases \(\{1, {\frak e}_{u_1}, \ldots{\frak e}_{u_k}\}\),
\(\{{\frak f}_j\}\) to be of the form:
\[
\left(\begin{array}{cccccc}
 +&- &0&\cdots& \cdots&+\\
+&+ &-& 0&\cdots &0\\
 +&0 &+ &\ddots&\cdots & 0 \\
\vdots &0 &0 &\ddots & \ddots&\vdots \\
 +&0 &\cdots &\cdots& + &-
\end{array}\right)\quad \text{(\(+/-\) denote positive/negative numbers.)} 
\] 
modulo ignorable terms. 
Combining this with the fact that, in terms of the same basis, 
\[\Pi_{\hat{K}_\chi} (w')_\Phi=(0, +, +, \cdots, +),\]
we see that
\begin{equation}\label{sgn-Dw}
[(w')_\Phi]/\ker^\mathfrak{gl} \acute{\mathbb{D}}_{(\lambda, w)}=-\prod_{i=1}^k\op{sign}(u_i).
\end{equation}

Next, we need to find the relative sign between the \(\mathfrak{gl}\)
and \(\mathfrak{o}_y\) orientations. For this purpose, 
we compute explicitly the form of the operator 
$\Pi_{C_{\#}}\acute{\mathbb{D}}_{y\#}|_{\hat{K}_{y\#}}$.
In terms of the bases \(\{1, e_1, \ldots, e_k, \}\) and \(\{f_j\}\), where
\[
e_i:= *\#_{R_1}*\cdots*\#_{R_i}\bar{\bf
  e}_y^\sigma\#_{R_i'}*\cdots*\#_{R_k'}; \quad
f_j:=*\#_{R_1}*\cdots*\#_{R_{j-1}'}\bar{\bf
  e}_y^\sigma\#_{R_j}*\cdots*\#_{R_k'}, 
\] 
it has the following form:
\[
\left(\begin{array}{cccccc}
 +&+&0&\cdots& \cdots&-\\
+&- &+& 0&\cdots &0\\
+&  0 &- &\ddots&\cdots & 0\\
\vdots& \vdots & \ddots &\ddots &\ddots & \vdots\\
+&0 &\cdots & 0& - &+
\end{array}\right).
\] 
Combining this with the facts that, in the same basis,
\[\Pi_{\hat{K}_{y\#}}(-({\bf e}_y)_\Phi)=(0, -, -, \cdots, -),\] 
we have by the proximity between $\acute{\mathbb{D}}_{y\#}$ and
$\acute{\mathbb{D}}_{y}$ that
\begin{equation*}\label{sgn-Dy}
\sign (\mathfrak{o}_{y-}/\mathfrak{gl})=-[({\bf e}_y)_\Phi]/\ker^\mathfrak{gl} 
\acute{\mathbb{D}}_{y}= (-1)^{k+1}.
\end{equation*}
The Lemma now follows by combining this with (\ref{sgn-Dw}) and (\ref{sign4}).
\end{proof}

\end{document}